\documentclass[12pt, a4paper]{amsart}

\usepackage{pdfsync}
\usepackage{amsmath,amsthm,amsfonts,amscd,amssymb,eucal,latexsym,mathrsfs}
\usepackage[all]{xy}
\usepackage{latexsym, amssymb, amscd, mathrsfs, graphics, fullpage}
\usepackage[backref=page,colorlinks,bookmarksopen=true]{hyperref}

\def\gG{\mathfrak{G}}
\def\gH{\mathfrak{H}}

\def\gM{\mathfrak{M}}
\def\gN{\mathfrak{N}}
\def\gT{\mathfrak{T}}
\def\ga{\mathfrak{a}}
\def\gb{\mathfrak {b}}

\def\gS{\mathfrak{S}}

\begin{document}

\title[Morita equivalence of measured quantum groupoids] 
{Morita equivalence of measured quantum groupoids\\
Application to deformation of measured quantum groupoids by 2-cocycles}
\author{Michel Enock}
\address{Institut de Math\'ematiques de Jussieu, Unit\'{e} Mixte Paris 6 / Paris 7 /
CNRS de Recherche 7586 \\Case 247, 4 place Jussieu, 75252 Paris Cedex}
 \email{enock@math.jussieu.fr}
\date{september 12}

\begin{abstract}
In a recent article of Kenny De Commer, was investigated a Morita equivalence between locally compact quantum groups, in which a measured quantum groupoid, of basis $\mathbb{C}^2$, was constructed as a linking object. Here, we generalize all these constructions and concepts to the level of measured quantum groupoids. As for locally compact quantum groups, we apply this construction to the deformation of a measured quantum groupoid by a 2-cocycle. 
 \end{abstract}

\maketitle
\newpage
\section{Introduction}
\label{intro}
\subsection{}
 In two articles (\cite{Val1}, \cite{Val2}), J.-M. Vallin has introduced two notions (pseudo-multiplicative
unitary, Hopf-bimodule), in order to generalize, up to the groupoid
case, the classical notions of multiplicative unitary \cite{BS} and of Hopf-von Neumann algebras \cite{ES}
which were introduced to describe and explain duality of groups, and leaded to appropriate notions
of quantum groups (\cite{ES}, \cite{W1}, \cite{W2}, \cite{BS}, \cite{MN}, \cite{W3}, \cite{KV1}, \cite{KV2}, \cite{MNW}). 
\\ In another article \cite{EVal}, J.-M. Vallin and the author have constructed, from a depth 2 inclusion of
von Neumann algebras $M_0\subset M_1$, with an operator-valued weight $T_1$ verifying a regularity
condition, a pseudo-multiplicative unitary, which leaded to two structures of Hopf bimodules, dual
to each other. Moreover, we have then obtained an action of one of these structures on the algebra $M_1$ such that $M_0$
is the fixed point subalgebra, the algebra $M_2$ given by the basic construction being then
isomorphic to the crossed-product. There is on $M_2$ an action of the other structure, which
can be considered as the dual action.
\\  If the inclusion
$M_0\subset M_1$ is irreducible, we recovered quantum groups, as proved and studied in former papers
(\cite{EN}, \cite{E2}).
\\ Therefore, this construction leads to a notion of "quantum groupoid", and a
duality within "quantum groupoids". 
\subsection{}
In a finite-dimensional setting, this construction can be
mostly simplified, and is studied in \cite{NV2}, \cite{BSz1},
\cite{BSz2}, \cite{Sz},\cite{Val3}, \cite{Val4}, \cite{Val5}, and examples are described. In \cite{NV2}, the link between these "finite quantum
groupoids" and depth 2 inclusions of
$II_1$ factors is given, and in \cite{D} had been proved that any finite-dimensional connected $\mathbb{C}^*$-quantum groupoid can act outerly on the hyperfinite $II_1$ factor. 
\subsection{}
In \cite{E3}, the author studied, in whole generality, the notion of pseudo-multiplicative unitary introduced par J.-M. Vallin in \cite{Val2}; following the strategy given by \cite{BS}, with the help of suitable fixed vectors, he introduced a notion of "measured quantum groupoid of compact type". Then F. Lesieur in \cite{L}, starting from a Hopf-bimodule (as introduced in \cite{Val1}), when there exist a left-invariant operator-valued weight and a right-invariant operator-valued weight, mimicking in this wider setting the technics of Kustermans and Vaes (\cite{KV1}, \cite{KV2}), obtained a pseudo-multiplicative unitary, which, as in the quantum group case, "contains" all the informations about the object (the von Neuman algebra, the coproduct, the antipod, the co-inverse). Lesieur gave the name of "measured quantum groupoids" to these objects. A new set of axioms for these had been given in an appendix of \cite{E5}. Moreover, in \cite{E4} had been shown that, with suitable conditions, the objects constructed in \cite{EVal} from depth 2 inclusions, are "measured quantum groupoids" in the sense of Lesieur.   
\subsection{}
In \cite{E5} have been developped the notions of action (already introduced in \cite{EVal}), crossed-product, etc, following what had been done for locally compact quantum groups in (\cite{E1}, \cite{ES1}, \cite{V2}); a biduality theorem for actions had been obtained in (\cite{E5}, 11.6). Moreover, we proved in (\cite{E5} 13.9) that, for any action of a measured quantum groupoid, the inclusion of the initial algebra (on which the measured quantum groupoid is acting) into the crossed-product is depth 2, which leads, thanks to \cite{E4}, to the construction of another measured quantum groupoid (\cite{E5} 14.2).  In \cite{E6} was proved a generalization of Vaes' theorem (\cite{V2}, 4.4) on the standard implementation of an action of a locally compact quantum group; namely, we had obtained such a result when there exists a normal semi-finite faithful operator-valued weight from the von Neumann algebra on which the measured quantum groupoid is acting, onto the copy of the basis of this measured quantum groupoid which is put inside this algebra by the action. 
\subsection{}
In \cite{E7} was studied outer actions of measured quantum groupoids. This notion was used to prove that any measured quantum groupoid can be constructed from a depth 2 inclusion. 
\subsection{}
In \cite{DC1}, Kenny De Commer introduced a notion of monoidal equivalence between two locally compact quantum groups, and constructed, in that situation, a measured quantum groupoid of basis $\mathbb{C}^2$ as a linking object between these two locally compact quantum groups. More precisely, from a locally compact quantum group $\bf{G}_1$ having a specific action $\ga_1$, called a Galois action, on a von Neumann algebra $A$, he was able to construct an important bunch of structures on $A$, and, by a reflexion technic, inspired by the work of P. Shauenburg in an algebraic context (\cite{Sc}), a second locally compact quantum group $\bf{G}_2$, and, more precisely, a measured quantum groupoid linking $\bf{G}_1$ and $\bf{G}_2$. This leads to an equivalence relation between locally compact quantum groups. 
\subsection{}
In that article, we generalize De Commer's construction to measured quantum groupoids. We call Morita equivalence this equivalence relation; two measured quantum groupoids $\gG_1$ and $\gG_2$ are Morita equivalent if there exists a von Neumann algebra on which $\gG_1$ acts on the right, $\gG_2$ acts on the left, and the two actions commute and being Galois, roughly speaking in a similar sense as de Commer's. This von Neumann algebra is then called an imprimitivity bi-comodule for these two measured quantum groupoids. This definition is similar to Renault's equivalence of locally compact groupoids, as defined in \cite{R1}, and developped in \cite{R2}, in which he proved that the $\bf{C}^*$-algebras of these two locally compact groupoids are then Morita equivalent. This is why we had chosen this terminology of "Morita equivalence". In \cite{DC2}, De Commer uses also this terminology, but two quantum groups are Morita equivalent in his sense if and only if their duals are Morita equivalent in ours. 
\subsection{}
In fact, De Commer's technics remain unchanged in the measured quantum groupoid context, if we start from a measured quantum groupoid $\gG$, and a Galois action $\ga$ of $\gG$ on a von Neumann algebra $A$, such that the invariant subalgebra $A^\ga$ is a finite sum of factors. This was remarked also in \cite{DC4}. In the general context, some extra hypothesis is needed, and we had to introduce what we called a "Galois system", which is, roughly speaking, a Galois action, equipped with an invariant weight. 
\subsection{}
De Commer used his construction to solve the problem of deforming a locally compact quantum group by a 2-cocycle. Namely, if $\bf{G}$ is a locally compact quantum group, and $\Omega$ a 2-cocycle, it had been observed since years that it is possible to deform the coproduct by using $\Omega$. Is the deformation still a locally compact quantum group ? or, equivalently, is there, in that case, an existence theorem for a left (resp. right) Haar weight ? This problem was solved in several particular cases and examples (\cite {EV}, \cite{V}, \cite{FV}) and De Commer answered positively to this question in whole generality. Of course, the same problem holds for measured quantum groupoids, and the answer is still positive when the basis of the measured quantum groupoid is a finite sum of factors. In the general case, we were able to give different sufficient conditions on the 2-cocycle, and give some examples, based on the construction of matched pairs of groupoids (\cite{Val6}). 

\subsection{}
This article is organized as follows : 
\newline
In chapter \ref{pre}, we recall as quickly as possible all the notations and results needed in that article; we emphazise that this article should be understood as the continuation of \cite{E5} and \cite{E6}, and that reading this article needs having \cite {E5} in hand. 
\newline
In chapter \ref{integrable}, inspired by \cite{V2} and \cite{DC1}, we prove specific results on integrable actions of a measured quantum groupoid $\gG$ and define Galois actions of $\gG$ and Galois systems for $\gG$ .
\newline
In the chapter \ref{fromGtoG}, inspired by \cite{DC1}, we associate to a Galois action of $\gG$ several data which will be usefull in the sequel. In particular we discuss how it is possible to construct a Galois system from a Galois action. 
\newline
In the chapter \ref{TLG}, we use the reflexion technique introduced in \cite{DC1}, in order to construct, "through the Galois system", another measured quantum groupoid $\gG_1$, and, more precisely, a measured quantum groupoid linking $\gG$ and $\gG_1$.  
\newline
The chapter \ref{Morita} is devoted to several equivalent definitions of Morita equivalence of measured quantum groupoids. We finish that chapter by giving some examples and constructions of locally compact quantum groups being Morita equivalent to measured quantum groupoids (\ref{exLCQG}). 
\newline
In the chapter \ref{cocycles}, following K. De Commer, we tried to use Morita equivalence to solve the problem of deforming a measured quantum groupoid by a 2-cocycle. This problem is here solved if the basis of the measured quantum groupoid is a finite sum of factors. In the general case, we obtain sufficient conditions, which will help, in chapter \ref{atleast}, to give a new example of construction of measured quantum groupoids, using J.-M. Vallin's construction of matched pairs of groupoids (\cite{Val6}).

 \section{Preliminaries}
 \label{pre}
This article is the continuation of \cite{E5}; preliminaries are to be found in \cite{E5}, and we just recall herafter the following definitions and notations :

\subsection{Spatial theory; relative tensor products of Hilbert spaces and fiber products of von Neumann algebras (\cite{C1}, \cite{S}, \cite{T}, \cite{EVal})}
\label{spatial}
 Let $N$ be a von Neumann algebra, $\psi$ a normal semi-finite faithful weight on $N$; we shall denote by $H_\psi$, $\gN_\psi$, ... the canonical objects of the Tomita-Takesaki theory associated to the weight $\psi$; let $\alpha$ be a non degenerate faithful representation of $N$ on a Hilbert space $\mathcal H$; the set of $\psi$-bounded elements of the left-module $_\alpha\mathcal H$ is :
\[D(_\alpha\mathcal{H}, \psi)= \lbrace \xi \in \mathcal{H};\exists C < \infty ,\| \alpha (y) \xi\|
\leq C \| \Lambda_{\psi}(y)\|,\forall y\in \gN_{\psi}\rbrace\]
Then, for any $\xi$ in $D(_\alpha\mathcal{H}, \psi)$, there exists a bounded operator
$R^{\alpha,\psi}(\xi)$ from $H_\psi$ to $\mathcal{H}$,  defined, for all $y$ in $\gN_\psi$ by :
\[R^{\alpha,\psi}(\xi)\Lambda_\psi (y) = \alpha (y)\xi\]
which intertwines the actions of $N$. 
\newline
If $\xi$, $\eta$ are bounded vectors, we define the operator product 
\[<\xi,\eta>_{\alpha,\psi} = R^{\alpha,\psi}(\eta)^* R^{\alpha,\psi}(\xi)\]
belongs to $\pi_{\psi}(N)'$, which, thanks to Tomita-Takesaki theory, will be identified to the opposite von Neumann algebra $N^o$. 
\newline
If now $\beta$ is a non degenerate faithful antirepresentation of $N$ on a Hilbert space $\mathcal K$, the relative tensor product $\mathcal K\underset{\psi}{_\beta\otimes_\alpha}\mathcal H$ is the completion of the algebraic tensor product $K\odot D(_\alpha\mathcal{H}, \psi)$ by the scalar product defined,  if $\xi_1$, $\xi_2$ are in $\mathcal{K}$, $\eta_1$, $\eta_2$ are in $D(_\alpha\mathcal{H},\psi)$, by the following formula :
\[(\xi_1\odot\eta_1 |\xi_2\odot\eta_2 )= (\beta(<\eta_1, \eta_2>_{\alpha,\psi})\xi_1 |\xi_2)\]
If $\xi\in \mathcal{K}$, $\eta\in D(_\alpha\mathcal{H},\psi)$, we shall denote $\xi\underset{\psi}{_\beta\otimes_\alpha}\eta$ the image of $\xi\odot\eta$ into $\mathcal K\underset{\psi}{_\beta\otimes_\alpha}\mathcal H$, and, writing $\rho^{\beta, \alpha}_\eta(\xi)=\xi\underset{\psi}{_\beta\otimes_\alpha}\eta$, we get a bounded linear operator from $\mathcal H$ into $\mathcal K\underset{\nu}{_\beta\otimes_\alpha}\mathcal H$, which is equal to $1_\mathcal K\otimes_\psi R^{\alpha, \psi}(\eta)$. 
\newline
Changing the weight $\psi$ will give an isomorphic Hilbert space, but the isomorphism will not exchange elementary tensors !
\newline

We shall denote $\sigma_\psi$ the relative flip, which is a unitary sending $\mathcal{K}\underset{\psi}{_\beta\otimes_\alpha}\mathcal{H}$ onto $\mathcal{H}\underset{\psi^o}{_\alpha\otimes _\beta}\mathcal{K}$, defined, for any $\xi$ in $D(\mathcal {K}_\beta ,\psi^o )$, $\eta$ in $D(_\alpha \mathcal {H},\psi)$, by :
\[\sigma_\psi (\xi\underset{\psi}{_\beta\otimes_\alpha}\eta)=\eta\underset{\psi^o}{_\alpha\otimes_\beta}\xi\]
In $x\in \beta(N)'$, $y\in \alpha(N)'$, it is possible to define an operator $x\underset{\psi}{_\beta\otimes_\alpha}y$ on $\mathcal K\underset{\psi}{_\beta\otimes_\alpha}\mathcal H$, with natural values on the elementary tensors. As this operator does not depend upon the weight $\psi$, it will be denoted $x\underset{N}{_\beta\otimes_\alpha}y$. 
\newline
We define a relative flip $\varsigma_N$ from $\mathcal L(\mathcal K)\underset{N}{_\beta*_\alpha}\mathcal L(\mathcal H)$ onto $\mathcal L(\mathcal H)\underset{N^o}{_\alpha*_\beta}\mathcal L(\mathcal K)$ by $\varsigma_N(X)=\sigma_\psi X(\sigma_{\psi})^*$, for any $X\in \mathcal L(\mathcal K)\underset{N}{_\beta*_\alpha}\mathcal L(\mathcal H)$ and any normal semi-finite faithful weight $\psi$ on $N$. 
\newline
If $P$ is a von Neumann algebra on $\mathcal H$, with $\alpha(N)\subset P$, and $Q$ a von Neumann algebra on $\mathcal K$, with $\beta(N)\subset Q$, then we define the fiber product $Q\underset{N}{_\beta*_\alpha}P$ as $\{x\underset{N}{_\beta\otimes_\alpha}y, x\in Q', y\in P'\}'$. 
\newline
Moreover, this von Neumann algebra can be defined independantly of the Hilbert spaces on which $P$ and $Q$ are represented; if $(i=1,2)$, $\alpha_i$ is a faithful non degenerate homomorphism from $N$ into $P_i$, $\beta_i$ is a faithful non degenerate antihomomorphism from $N$ into $Q_i$, and $\Phi$ (resp. $\Psi$) an homomorphism from $P_1$ to $P_2$ (resp. from $Q_1$ to $Q_2$) such that $\Phi\circ\alpha_1=\alpha_2$ (resp. $\Psi\circ\beta_1=\beta_2$), then, it is possible to define an homomorphism $\Psi\underset{N}{_{\beta_1}*_{\alpha_1}}\Phi$ from $Q_1\underset{N}{_{\beta_1}*_{\alpha_1}}P_1$ into $Q_2\underset{N}{_{\beta_2}*_{\alpha_2}}P_2$. 
\newline
The operators $\theta^{\alpha, \psi}(\xi, \eta)=R^{\alpha, \psi}(\xi)R^{\alpha, \psi}(\eta)^*$, for all $\xi$, $\eta$ in $D(_\alpha\mathcal H, \psi)$, generates a weakly dense ideal in $\alpha(N)'$. Moreover, there exists a family $(e_i)_{i\in I}$ of vectors in $D(_\alpha\mathcal H, \psi)$ such that the operators $\theta^{\alpha, \psi}(e_i, e_i)$ are 2 by 2 orthogonal projections ($\theta^{\alpha, \psi}(e_i, e_i)$ being then the projection on the closure of $\alpha(N)e_i$). Such a family is called an orthogonal $(\alpha, \psi)$-basis of $\mathcal H$.

\subsection{Measured quantum groupoids (\cite{L}, \cite{E5})}
\label{MQG}
 A quintuplet $(N, M, \alpha, \beta, \Gamma)$ will be called a Hopf-bimodule, following (\cite{Val2}, \cite{EVal} 6.5), if
$N$,
$M$ are von Neumann algebras, $\alpha$ a faithful non-degenerate representation of $N$ into $M$, $\beta$ a
faithful non-degenerate anti-representation of
$N$ into $M$, with commuting ranges, and $\Gamma$ an injective involutive homomorphism from $M$
into
$M\underset{N}{_\beta *_\alpha}M$ such that, for all $X$ in $N$ :
\newline
(i) $\Gamma (\beta(X))=1\underset{N}{_\beta\otimes_\alpha}\beta(X)$
\newline
(ii) $\Gamma (\alpha(X))=\alpha(X)\underset{N}{_\beta\otimes_\alpha}1$ 
\newline
(iii) $\Gamma$ satisfies the co-associativity relation :
\[(\Gamma \underset{N}{_\beta *_\alpha}id)\Gamma =(id \underset{N}{_\beta *_\alpha}\Gamma)\Gamma\]
This last formula makes sense, thanks to the two preceeding ones and
\ref{spatial}. The von Neumann algebra $N$ will be called the basis of $(N, M, \alpha, \beta, \Gamma)$\vspace{5mm}.\newline
If $(N, M, \alpha, \beta, \Gamma)$ is a Hopf-bimodule, it is clear that
$(N^o, M, \beta, \alpha,
\varsigma_N\circ\Gamma)$ is another Hopf-bimodule, we shall call the symmetrized of the first
one. (Recall that $\varsigma_N\circ\Gamma$ is a homomorphism from $M$ to
$M\underset{N^o}{_\alpha*_\beta}M$).
\newline
If $N$ is abelian, $\alpha=\beta$, $\Gamma=\varsigma_N\circ\Gamma$, then the quadruplet $(N, M, \alpha, \alpha,
\Gamma)$ is equal to its symmetrized Hopf-bimodule, and we shall say that it is a symmetric
Hopf-bimodule\vspace{5mm}.\newline
A measured quantum groupoid is an octuplet $\mathfrak {G}=(N, M, \alpha, \beta, \Gamma, T, T', \nu)$ such that (\cite{E5}, 3.8) :
\newline
(i) $(N, M, \alpha, \beta, \Gamma)$ is a Hopf-bimodule, 
\newline
(ii) $T$ is a left-invariant normal, semi-finite, faithful operator valued weight $T$ from $M$ to $\alpha (N)$ (to be more precise, from $M^{+}$ to the extended positive elements of $\alpha(N)$ (cf. \cite{T} IX.4.12)), which means that, for any $x\in\gM_T^+$, we have $(id\underset{\nu}{_\beta*_\alpha}T)\Gamma(x)=T(x)\underset{N}{_\beta\otimes_\alpha}1$. 
\newline
(iii) $T'$ is a right-invariant normal, semi-finite, faithful operator-valued weight $T'$ from $M$ to $\beta (N)$, which means that, for any $x\in\gM_{T'}^+$, we have $(T'\underset{\nu}{_\beta*_\alpha}id)\Gamma(x)=1\underset{N}{_\beta\otimes_\alpha}T'(x)$. 
\newline
(iv) $\nu$ is normal semi-finite faitfull weight on $N$, which is relatively invariant with respect to $T$ and $T'$, which means that the modular automorphisms groups of the weights $\Phi=\nu\circ\alpha^{-1}\circ T$ and $\Psi=\nu^o\circ\beta^{-1}\circ T'$ commute. 
\newline
We shall write $H=H_\Phi$, $J=J_\Phi$, and, for all $n\in N$, $\hat{\beta}(n)=J\alpha(n^*)J$, $\hat{\alpha}(n)=J\beta(n^*)J$.  The weight $\Phi$ will be called the left-invariant weight on $M$. 
\newline
Examples are described and explained in \ref{exMQG}. 
\newline
Then, $\mathfrak {G}$ can be equipped with a pseudo-multiplicative unitary $W$ which is a unitary from $H\underset{\nu}{_\beta\otimes_\alpha}H$ onto $H\underset{\nu^o}{_\alpha\otimes_{\hat{\beta}}}H$ (\cite{E5}, 3.6), which intertwines $\alpha$, $\hat{\beta}$, $\beta$  in the following way : for all $X\in N$, we have :
\[W(\alpha
(X)\underset{N}{_\beta\otimes_\alpha}1)=
(1\underset{N^o}{_\alpha\otimes_{\hat{\beta}}}\alpha(X))W\]
\[W(1\underset{N}{_\beta\otimes_\alpha}\beta
(X))=(1\underset{N^o}{_\alpha\otimes_{\hat{\beta}}}\beta (X))W\]
\[W(\hat{\beta}(X) \underset{N}{_\beta\otimes_\alpha}1)=
(\hat{\beta}(X)\underset{N^o}{_\alpha\otimes_{\hat{\beta}}}1)W\]
\[W(1\underset{N}{_\beta\otimes_\alpha}\hat{\beta}(X))=
(\beta(X)\underset{N^o}{_\alpha\otimes_{\hat{\beta}}}1)W\]
and the operator $W$ satisfies :
\[(1\underset{N^o}{_\alpha\otimes_{\hat{\beta}}}W)
(W\underset{N}{_\beta\otimes_\alpha}1_{\gH})
=(W\underset{N^o}{_\alpha\otimes_{\hat{\beta}}}1)
\sigma^{2,3}_{\alpha, \beta}(W\underset{N}{_{\hat{\beta}}\otimes_\alpha}1)
(1\underset{N}{_\beta\otimes_\alpha}\sigma_{\nu^o})
(1\underset{N}{_\beta\otimes_\alpha}W)\]
Here, $\sigma^{2,3}_{\alpha, \beta}$
goes from $(H\underset{\nu^o}{_\alpha\otimes_{\hat{\beta}}}H)\underset{\nu}{_\beta\otimes_\alpha}H$ to $(H\underset{\nu}{_\beta\otimes_\alpha}H)\underset{\nu^o}{_\alpha\otimes_{\hat{\beta}}}H$, 
and $1\underset{N}{_\beta\otimes_\alpha}\sigma_{\nu^o}$ goes from $H\underset{\nu}{_\beta\otimes_\alpha}(H\underset{\nu^o}{_\alpha\otimes_{\hat{\beta}}}H)$ to $H\underset{\nu}{_\beta\otimes_\alpha}H\underset{\nu}{_{\hat{\beta}}\otimes_\alpha}H$. 
\newline
All the intertwining properties properties allow us to write such a formula, which will be called the
"pentagonal relation". Moreover, $W$, $M$ and $\Gamma$ are related by the following results :
\newline
(i) $M$ is the weakly closed linear space generated by all operators of the form $(id*\omega_{\xi, \eta})(W)$, where $\xi\in D(_\alpha H, \nu)$, and $\eta\in D(H_{\hat{\beta}}, \nu^o)$ (\cite{E5}, 3.8(vii)).
\newline
(ii) for any $x\in M$, we have $\Gamma(x)=W^*(1\underset{N^o}{_\alpha\otimes_{\hat{\beta}}}x)W$ (\cite{E5}, 3.6). 
\subsubsection{Lemma}
\label{Gamma}
{\it Let $\gG$ be a measured quantum groupoid, $W$ its pseudo-multiplicative unitary, $\xi\in D(_\alpha H, \nu)$ and $\eta\in D(H_{\hat{\beta}}, \nu^o)$; then :}
\[\Gamma((id*\omega_{\xi, \eta})(W))=(id\underset{N}{_\beta*_\alpha}id*\omega_{\xi, \eta})(\sigma^{2,3}_{\alpha, \beta}(W\underset{N}{_{\hat{\beta}}\otimes_\alpha}1)
(1\underset{N}{_\beta\otimes_\alpha}\sigma_{\nu^o})
(1\underset{N}{_\beta\otimes_\alpha}W))\]
\begin{proof}
This is clear, using the pentagonal relation, and the formula linking $\Gamma$ and $W$. \end{proof}

Moreover, it is also possible to construct many other data, namely a co-inverse $R$, a scaling group $\tau_t$, an antipod $S$, a modulus $\delta$, a scaling operator $\lambda$, a managing operator $P$, and a canonical one-parameter group $\gamma_t$ of automorphisms on the basis $N$ (\cite{E5}, 3.8). Instead of $\mathfrak {G}$, we shall mostly use $(N, M, \alpha, \beta, \Gamma, T, RTR, \nu)$ which is another measured quantum groupoid, denoted $\underline{\mathfrak {G}}$, which is equipped with the same data ($W$, $R$, ...) as $\gG$. 
\newline
A dual measured quantum group $\widehat{\mathfrak{G}}$, which is denoted $(N, \widehat{M}, \alpha, \hat{\beta}, \widehat{\Gamma}, \widehat{T}, \widehat{R}\widehat{T}\widehat{R}, \nu)$, can be constructed, and we have $\widehat{\widehat{\mathfrak {G}}}=\underline{\mathfrak {G}}$. 
\newline
In particular, from the fact that $\nu$ is relatively invariant with respect to $T$ and $R\circ T\circ R$, is obtained the definition of the modulus and the scaling operator by the formula :
\[(D\Phi\circ R:D\Phi)_t=\lambda^{it^2/2}\delta^{it}\]
Then, thanks to \cite{V1}, we obtain that, if $a\in M$ is such that the operator $a\delta^{1/2}$ is bounded and its closure $\overline{a\delta^{1/2}}$ belongs to $\gN_\Phi$, then $a$ belongs to $\gN_{\Phi\circ R}$, and that we can identify $H_{\Phi\circ R}$ with $H$ by writing then $\Lambda_{\Phi\circ R}(a)=\Lambda_{\Phi}(\overline{a\delta^{1/2}})$. 
\newline
Canonically associated to $\mathfrak {G}$, can be defined also the opposite measured quantum groupoid is $\mathfrak{G}^o=(N^o, M, \beta, \alpha, \varsigma_N\Gamma, RTR, T, \nu^o)$ and the commutant measured quantum groupoid $\mathfrak{G}^c=(N^o, M', \hat{\beta}, \hat{\alpha}, \Gamma^c, T^c, R^cT^cR^c, \nu^o)$; we have $(\mathfrak{G}^o)^o=(\mathfrak{G}^c)^c=\underline{\mathfrak{G}}$, $\widehat{\mathfrak{G}^o}=(\widehat{\mathfrak {G}})^c$, $\widehat{\mathfrak {G}^c}=(\widehat{\mathfrak {G}})^o$, and $\mathfrak{G}^{oc}=\mathfrak {G}^{co}$ is canonically isomorphic to $\underline{\mathfrak {G}}$ (\cite{E5}, 3.12). 
\newline
The pseudo-multiplicative unitary of $\widehat{\mathfrak{G}}$ (resp. $\mathfrak{G}^o$, $\mathfrak{G}^c$) will be denoted $\widehat{W}$ (resp. $W^o$, $W^c$). The left-invariant weight on $\widehat{\mathfrak{G}}$ (resp. $\mathfrak{G}^o$, $\mathfrak{G}^c$) will be denoted $\widehat{\Phi}$ (resp. $\Phi^o$, $\Phi^c$). For simplification, we shall write $\hat{J}$ for $J_{\widehat{\Phi}}$. 
\newline
We have $\widehat{W}=\sigma_{\nu^o}W^*\sigma_{\nu}$ which is a unitary from $H\underset{\nu}{_{\hat{\beta}}\otimes_\alpha}H$ onto $H\underset{\nu^o}{_\alpha\otimes_\beta}H$. The algebra $\widehat{M}$ is generated by the operators $(\omega_{\xi, \eta}*id)(W)$, where $\xi$ belongs to $D(H_{\hat{\beta}}, \nu^o)$ and $\eta$ belongs to $D(_\alpha H, \nu)$. In (\cite{E5}4.8) was proved that such an element belongs to $\gN_{\widehat{\Phi}}$ if and only if $\xi$ belongs to $\mathcal D(\pi'(\eta)^*)$, where $\pi'(\eta)^*$ is the adjoint of the (densely defined) operator $\pi'(\eta)$ defined on $\Lambda_\Phi(\gN_\Phi)$ by $\pi'(\eta)\Lambda_\Phi(x)=x\eta$, and we have then :
\[\widehat{\Phi}[(\omega_{\xi, \eta}*id)(W)^*(\omega_{\xi, \eta}*id)(W)]=\|\pi'(\eta)^*\xi\|^2\]
which allows us to identify $H_{\widehat{\Phi}}$ with $H$ by writing $\Lambda_{\widehat{\Phi}}((\omega_{\xi, \eta}*id)(W))=\pi'(\eta)^*\xi$, or, for any $x\in \gN_\Phi$ :
\[(\Lambda_\Phi(x)|\Lambda_{\widehat{\Phi}}((\omega_{\xi, \eta}*id)(W)))=(x\eta|\xi)\]
The pseudo-multiplicative unitary $W^o$ is equal to $(\hat{J}\underset{N^o}{_\alpha\otimes_{\hat{\beta}}}\hat{J})W(\hat{J}\underset{N^o}{_\alpha\otimes_\beta}\hat{J})$ (\cite{E5}, 3.12(v)), which is a unitary from $H\underset{\nu^o}{_\alpha\otimes_\beta}H$ onto $H\underset{\nu}{_\beta\otimes_{\hat{\alpha}}}H$, where, for all $n\in N$, $\hat{\alpha}(n)=J\beta(n^*)J$. Therefore, applying this result about $\gN_{\widehat{\Phi}}$ to the duality between $\gG^o$ and $\widehat{\gG}^c$, we obtain that the operator $(\omega_{\xi, \eta}*id)(W^o)=(\omega_{\xi, \eta}*id)[(\hat{J}\underset{N^o}{_\alpha\otimes_{\hat{\beta}}}\hat{J})W(\hat{J}\underset{N^o}{_\alpha\otimes_\beta}\hat{J})]$ belongs to $\gN_{\widehat{\Phi}^c}$ if and only if $\hat{J}\xi$ belongs to $\mathcal D(\pi'(\hat{J}\eta)^*)$, and then, we get :
\[\Lambda_{\widehat{\Phi}^c}(\omega_{\xi, \eta}*id)(W^o)=\Lambda_{\widehat{\Phi}^c}[(\omega_{\xi, \eta}*id)[(\hat{J}\underset{N^o}{_\alpha\otimes_{\hat{\beta}}}\hat{J})W(\hat{J}\underset{N^o}{_\alpha\otimes_\beta}\hat{J})]]=\hat{J}\pi'(\hat{J}\eta)^*\hat{J}\xi\]
and, if, moreover, $\eta$ belongs to $\mathcal D(\delta^{-1/2})$, we get, for any $x\in\gN_\Phi$ :
\[(\Lambda_\Phi(x)|\Lambda_{\widehat{\Phi}^c}(\omega_{\xi, \eta}*id)(W^o))=(\Lambda_\Phi(x)|\Lambda_{\widehat{\Phi}^c}[(\omega_{\xi, \eta}*id)[(\hat{J}\underset{N^o}{_\alpha\otimes_{\hat{\beta}}}\hat{J})W(\hat{J}\underset{N^o}{_\alpha\otimes_\beta}\hat{J})]])=(x\delta^{-1/2}\eta|\xi)\]
Let $_a\gH_b$ be a $N-N$-bimodule, i.e. an Hilbert space $\gH$ equipped with a normal faithful non degenerate representation $a$ of $N$ on $\gH$ and a normal faithful non degenerate anti-representation $b$ on $\gH$, such that $b(N)\subset a(N)'$. A corepresentation of $\gG$ on $_a\gH_b$ is a unitary $V$ from $\gH\underset{\nu^o}{_a\otimes_\beta}H$ onto 
 $\gH\underset{\nu}{_b\otimes_\alpha}H$, satisfying, for all $n\in N$ :
 \[V(b(n)\underset{N^o}{_a\otimes_\beta}1)=(1\underset{N}{_b\otimes_\alpha}\beta(n))V\]
 \[V(1\underset{N^o}{_a\otimes_\beta}\alpha(n))=(a(n)\underset{N}{_b\otimes_\alpha}1)V\]
 \[V(1\underset{N^o}{_a\otimes_\beta}\hat{\beta}(n))=(1\underset{N}{_b\otimes_\alpha}\hat{\beta}(n))V\]
such that, for any $\xi\in D(_a\gH, \nu)$ and $\eta\in D(\gH_b, \nu^o)$, the operator $(\omega_{\xi, \eta}*id)(V)$ belongs to $M$ (then, it is possible to define $(id*\theta)(V)$, for any $\theta$ in $M_*^{\alpha, \beta}$ which is the linear set generated by the $\omega_\xi$, with $\xi\in D(_\alpha H, \nu)\cap D(H_\beta, \nu^o)$), and such that the application $\theta\rightarrow (id*\theta)(V)$ from $M_*^{\alpha, \beta}$ into $\mathcal L(\gH)$ is multiplicative (\cite{E5} 5.1, 5.5).

\subsubsection{Lemma}
\label{lemW}
{\it Let $\gG$ be a measured quantum groupoid; we have, for any $\xi\in D(H_{\beta}, \nu^o)$ and $\eta\in D(_\alpha H, \nu)$, $t\in\mathbb{R}$ :}
\[\sigma_t^{\widehat{\Phi}}[(\omega_{\xi, \eta}*id)(W)]=(\omega_{P^{it}\xi, \delta^{-it}P^{it}\eta}*id)(W)\]

\begin{proof}
Let $\zeta_1\in D(_\alpha H, \nu)$, and $\zeta_2\in D(H_{\hat{\beta}}, \nu^o)$; we have, using successively \cite{E5}, 3.10(vii), 3.8(vii), 3.11(iii), 3.8(vi) and again 3.11(iii):
\begin{eqnarray*}
\sigma_t^{\widehat{\Phi}}((\omega_{\xi, \eta}*id)(W)\zeta_1|\zeta_2)
&=&
(W(\xi\underset{\nu}{_\beta\otimes_\alpha}\Delta_{\widehat{\Phi}}^{-it}\zeta_1)|\eta\underset{\nu^o}{_\alpha\otimes_{\hat{\beta}}}\Delta_{\widehat{\Phi}}^{-it}\zeta_2)\\
&=&
(W(\xi\underset{\nu}{_\beta\otimes_\alpha}P^{-it}J\delta^{it}J\zeta_1)|\eta\underset{\nu^o}{_\alpha\otimes_{\hat{\beta}}}P^{-it}J\delta^{it}J\zeta_2)\\
&=&
(W(P^{it}\xi\underset{\nu}{_\beta\otimes_\alpha}J\delta^{it}J\zeta_1)|P^{it}\eta\underset{\nu^o}{_\alpha\otimes_{\hat{\beta}}}J\delta^{it}J\zeta_2)\\
&=&
(\hat{J}P^{it}\eta\underset{\nu}{_\beta\otimes_\alpha}\delta^{it}J\zeta_2|W^*(\hat{J}P^{it}\xi\underset{\nu^o}{_\alpha\otimes_{\hat{\beta}}}\delta^{it}J\zeta_1))\\
&=&
(\hat{J}P^{it}\eta\underset{\nu}{_\beta\otimes_\alpha}\delta^{it}J\zeta_2|(\delta^{it}\underset{N}{_\beta\otimes_\alpha}\delta^{it})W^*(\hat{J}P^{it}\xi\underset{\nu^o}{_\alpha\otimes_{\hat{\beta}}}J\zeta_1))\\
&=&
(\hat{J}\delta^{-it}P^{it}\eta\underset{\nu}{_\beta\otimes_\alpha}J\zeta_2|W^*(\hat{J}P^{it}\xi\underset{\nu^o}{_\alpha\otimes_{\hat{\beta}}}J\zeta_1))\\
&=&
(W(P^{it}\xi\underset{\nu}{_\beta\otimes_\alpha}\zeta_1)|\delta^{-it}P^{it}\eta\underset{\nu^o}{_\alpha\otimes_{\hat{\beta}}}\zeta_2)
\end{eqnarray*}
from which we get the result. \end{proof}

\subsubsection{Lemma}
\label{lemW2}
{\it Let $\gG$ be a measured quantum groupoid, and $m\in\widehat{M}'$; then, we have :}
\[\widehat{W}(1\underset{N}{_{\hat{\beta}}\otimes_\alpha}m)\widehat{W}^*=W^{o*}(\hat{J}\hat{R}^c(m^*)\hat{J}\underset{N}{_\beta\otimes_{\hat{\alpha}}}1)\widehat{W}^o\]

\begin{proof}
By definition, we have :
\begin{eqnarray*}
\widehat{W}(1\underset{N}{_{\hat{\beta}}\otimes_\alpha}m)\widehat{W}^*
&=&
\sigma W^*\sigma ((1\underset{N}{_{\hat{\beta}}\otimes_\alpha}m)\sigma W\sigma\\
&=&
\sigma W^*(m\underset{N^o}{_\alpha\otimes_{\hat{\beta}}}1) W\sigma
\end{eqnarray*}
and, using (\cite{E5}, 3.11(iii), 3.10 (iii), 3.12(v) and 3.11(iii) again), we get it is equal to :
\begin{eqnarray*}
\sigma (\hat{J}\underset{N}{_\beta\otimes_\alpha}J)W(\hat{J}m\hat{J}\underset{N}{_\beta\otimes_\alpha}1)W^*(\hat{J}\underset{N}{_\beta\otimes_{\hat{\alpha}}}J)\sigma
&=&
\sigma  (\hat{J}\underset{N}{_\beta\otimes_\alpha}J)\varsigma \widehat{\Gamma}(\hat{J}m\hat{J})(\hat{J}\underset{N}{_\beta\otimes_{\hat{\alpha}}}J)\sigma\\
&=&
(J\underset{N^o}{_\alpha\otimes_\beta}\hat{J})\widehat{\Gamma}(\hat{J}m\hat{J})(J\underset{N^o}{_{hat{\alpha}}\otimes_\beta}\hat{J})\\
&=&
W^{o*}(\hat{J}\hat{R}^c(m^*)\hat{J}\underset{N}{_\beta\otimes_{\hat{\alpha}}}1)\widehat{W}^o
\end{eqnarray*}
\end{proof}

\subsection{Examples of measured quantum groupoids}
\label{exMQG}
Examples of measured quantum groupoids are the following :
\newline
(i) locally compact quantum groups, as defined and studied by J. Kustermans and S. Vaes (\cite{KV2}, \cite {KV2}, \cite{V2}); these are, trivially, the measured quantum groupoids with the basis $N=\mathbb{C}$. 
\newline
(ii) measured groupoids, equipped with a left Haar system and a quasi-invariant measure on the set of units, as studied mostly by T. Yamanouchi (\cite{Y1}, \cite{Y2}, \cite{Y3}, \cite{Y4}); it was proved in \cite{E8} that these measured quantum groupoids are exactly those whose underlying von Neumann algebra is abelian. This example had been presented in full details in (\cite{E5}, 3.4 and 3.13). 
\newline
(iii) the finite dimensional case had been studied by D. Nikshych and L. Vainermann (\cite{NV2}, \cite{NV2}) and J.-M. Vallin (\cite{Val3}, \cite{Val4}); in that case, non trivial examples are given. 
\newline
(iv) continuous fields of ($\bf{C}^*$-version of) locally compact quantum groups, as studied by E. Blanchard in (\cite{Bl1}, \cite{Bl2}); it was proved in \cite{E8} that these measured quantum groupoids are exactly those whose basis is central in the underlying von Neumann algebras of both the measured quantum groupoid and its dual. 
\newline
(v) in (\cite{L}, 17.1), be given a family $\gG_i=(N_i, M_i, \alpha_i, \beta_i, \Gamma_i, T_i, T'_i, \nu_i)$ a measured quantum groupoids, Lesieur showed that it is possible to construct another measured quantum groupoid $\gG=\oplus_{i\in I}\gG_i=(\oplus_{i\in I}N_i, \oplus_{i\in I}M_i, \oplus_{i\in I}\alpha_i, \oplus_{i\in I}\beta_i, \oplus_{i\in I}\Gamma_i, \oplus_{i\in I}T_i, \oplus_{i\in I}T'_i, \oplus_{i\in I}\nu_i)$. 
\newline
(vi) in \cite{DC1}, K. De Commer proved that, in the case of a monoidal equivalence between two locally compact quantum groups (which means that each of these locally compact quantum group has an ergodic and integrable action on the other one), it is possible to construct a measured quantum groupoid of basis $\mathbb{C}^2$ which contains all the data. Moreover, he proved that such measured quantum groupoids are exactly those whose basis $\mathbb{C}^2$ is central in the underlying von Neumann algebra of the measured quantum groupoid, but not in the underlying von Neumann algebra of the dual measured quantum groupoid. 
\newline
(vii) in \cite{E5} was described how, from an action $(b, \ga)$ of a measured quantum groupoid $\gG$, it is possible to construct another measured quantum groupoid $\gG(\ga)$; as a particular case, this allows to canonically associate to any action $\ga$ of a locally compact quantum group $\bf{G}$ on a von Neumann algebra $A$, a measured quantum groupoid $\gG(\ga)$. 
\newline
(viii) in \cite{VV} was given a specific procedure to construct locally compact quantum groups, starting from a locally compact group $G$, whose almost all elements belong to the product $G_1G_2$ (where $G_1$ and $G_2$ are closed subgroups of $G$ such that $G_1\cap G_2=\{e\}$, where $e$ is the neutral element of $G$); such $(G_1, G_2)$ is called a "matched pair" of locally compact groups. Then, $G_1$ acts naturally on $L^\infty(G_2)$ (and vice versa), and the two crossed-products obtained bear the structure of two locally compact quantum groups in duality. In \cite{Val5}, J.-M. Vallin generalizes this constructions up to groupoids, and, then, obtains examples of measured quantum groupoids; more specific examples are then given by the action of a matched pair of groups on a locally compact space, and also more exotic examples. 
\newline
(ix) in \cite{L}, 9.5.5, was given the following exemple, called "quantum space quantum groupoid"; let $N$ be a von Neumann algebra; let us consider $M=N^o\otimes_{Z(N)}N$, the representation $\alpha$ of $N$ into $M$ given by ($n\in N$) $\alpha(n)=1\otimes_{Z(N)}n$, and the anti-representation $\beta$ given by $\beta(n)=n^o\otimes_{Z(N)}1$. Then if $\tau$ is a normal semi-finite faithful trace on $Z(N)$, $\nu$ a normal faithful semi-finite weight on $N$, let $T_\nu$ be the normal faithful semi-finite operator-valued weight from $N$ onto $Z(N)$ such that $\nu=\tau\circ T_\nu$, we can easily get that the relative tensor product $(H_\nu\otimes_\tau H_\nu)\underset{\nu}{_\beta\otimes_\alpha}(H_\nu\otimes_\tau H_\nu)$ is canonically isomorphic to $H_\nu\otimes_\tau H_\nu\otimes_\tau H_\nu$ and, this isomorphism sends $M\underset{N}{_\beta*_\alpha}M$ onto $N^o\otimes_\tau Z(N)\otimes_\tau N$; we can therefore identify $M\underset{N}{_\beta*_\alpha}M$ with $M$, and verify that $(N, M, \alpha, \beta, id)$ is a Hopf-bimodule. 
\newline
Moreover, we can get that $\gG(N)=(N, M, \alpha, \beta, id, T_\nu^o\otimes_{Z(N)}id, id\otimes_{Z(N)}T_\nu, \nu)$ is a measured quantum groupoid. We shall call it the $N$-measured quantum groupoid. 
\newline
The dual measured quantum groupoid is $\widehat{\gG(N)}=(N, Z(N)', \alpha, \hat{\beta}, id, (T_\nu^o)^{-1}, T_\nu^{-1}, \nu)$, where $\hat{\beta}(n)=J_\nu n^*J_\nu$, $T_\nu^{-1}$ is the canonical operator-valued weight from $Z(N)'$ to $N'$ given from $T_\nu$, and $(T_\nu^o)^{-1}$ is the canonical operator-valued weight from $Z(N)'$ to $N$ given from $T_\nu^o$. This measured quantum groupoid will be called the dual $N$-measured quantum groupoid. 
\newline
(x) If $\gG_1=(N_1, M_1, \alpha_1, \beta_1, \Gamma_1, T_1, T'_1, \nu_1)$ and $\gG_2=(N_2, M_2, \alpha_2, \beta_2, \Gamma_2, T_2, T'_2, \nu_2)$ are two measured quantum groupoids, then we can define another measured quantum groupoid :
\[\gG_1\otimes\gG_2=(N_1\otimes N_2, M_1\otimes M_2, \alpha_1\otimes\alpha_2, \beta_1\otimes\beta_2, (id\otimes\varsigma\otimes id)(\Gamma_1\otimes\Gamma_2), T_1\otimes T_2, T'_1\otimes T'_2, \nu_1\otimes\nu_2)\]
 Moreover, it easy to get that $\widehat{\gG_1\otimes\gG_2}=\widehat{\gG_1}\otimes\widehat{\gG_2}$.
\newline
(xi) The $SU(2)$ dynamical quantum group, as studied in particular by E. Koelink and H. Rosengren (\cite{KR}) can be lifted, thanks to \cite{Ti}, to the level of operator algebras, and give another example of a measured quantum groupoid. 
\newline
(xii) last, but not least, De Commer studied Morita equivalence between the quantum group $SU_q(2)$, and various quantum groups (\cite{DC2}, \cite{DC3}). In a new work (\cite{DC4}), he obtains an integrable Galois action of $SU_q(2)$ which is not ergodic. Therefore, this leads to a measured quantum groupoid (\ref{KDC}).

\subsection{Action of a measured quantum groupoid (\cite{E5})}
\label{action}

An action (\cite{E5}, 6.1) of $\mathfrak{G}$ on a von Neumann algebra $A$ is a couple $(b, \mathfrak a)$, where :
\newline
(i) $b$ is an injective $*$-antihomomorphism from $N$ into $A$; 
\newline
(ii) $\mathfrak a$ is an injective $*$-homomorphism from $A$ into $A\underset{N}{_b*_\alpha}M$; 
\newline
(iii) $b$ and $\mathfrak a$ are such that, for all $n$ in $N$:
\[\mathfrak a (b(n))=1\underset{N}{_b\otimes_\alpha}\beta(n)\]
(which allow us to define $\mathfrak a\underset{N}{_b*_\alpha}id$ from $A\underset{N}{_b*_\alpha}M$ into $A\underset{N}{_b*_\alpha}M\underset{N}{_\beta*_\alpha}M$)
and such that :
\[(\mathfrak a\underset{N}{_b*_\alpha}id)\mathfrak a=(id\underset{N}{_b*_\alpha}\Gamma)\mathfrak a\]
If we start from a measured groupoid, we get the usual notion of action of a groupoid (\cite{E5}, 6.3).
\newline
The invariant subalgebra $A^\ga$ is defined by :
\[A^\ga=\{x\in A\cap b(N)'; \ga(x)=x\underset{N}{_b\otimes_\alpha}1\}\]
As $A^\ga\subset b(N)'$, $A$ (and $L^2(A)$) is a $A^\ga-N^o$-bimodule.
\newline
Let us write, for any $x\in A^+$, $T_\ga(x)=(id\underset{\nu}{_b*_\alpha}\Phi)\ga(x)$; this formula defines a normal faithful operator-valued weight from $A$ onto $A^\ga$; the action $\ga$ will be said integrable if $T_\ga$ is semi-finite (\cite{E5}, 6.11, 12, 13 and 14). 
\newline
If the von Neumann algebra $A$ acts on a Hilbert space $\gH$, and if there exists a representation $a$ of $N$ on $\gH$ such that $b(N)\subset A\subset a(N)'$, a corepresentation $V$ of $\gG$ on the bimodule $_a\gH_b$ will be called an implementation of $\ga$ if we have $\ga(x)=V(x\underset{N^o}{_a\otimes_b}1)V^*$ , for all $x\in A$ (\cite{E5}, 6.6); moreover, if $\psi$ is a normal semi-finite faithful weight on $A$, we shall define a representation $a$ of $N$ on $H_\psi$ by $a(n)=J_\psi b(n^*)J_\psi$, for all $n\in N$, and we shall look after an implementation $V$ of $\ga$ on $_a(H_\psi)_b$ such that (\cite{E5}, 6.9):
\[V^*=(J_\psi\underset{\nu^o}{_\alpha\otimes_\beta}J_{\widehat{\Phi}})V(J_\psi\underset{\nu}{_b\otimes_\alpha} J_{\widehat{\Phi}})\]
If the weight $\psi$ is $\delta$-invariant, which means that, for all $\eta\in D(_\alpha H, \nu)\cap\mathcal D(\delta^{1/2})$ such that $\delta^{1/2}\eta$ belongs to $D(H_\beta, \nu^o)$, and $x\in\gN_\psi$, we have :
\[\psi[(id\underset{N}{_b*_\alpha}\omega_\eta)\ga(x^*x)]=\|\Lambda_{\psi}(x)\underset{\nu^o}{_a\otimes_\beta}\delta^{1/2}\eta\|^2\]
and if, moreover, $\psi$ bears the density property, which means that $D((H_\psi)_b, \nu^o)\cap D(_a H_\psi, \nu)$ is dense in $H_\psi$, then such an implementation $V_\psi$ was constructed in \cite{E5}, 8.8); more precisely (\cite{E5}, 8.4), if $x\in\gN_\psi$, $\xi\in D(_\alpha H, \nu)$ and $\eta$ is as above, we get that $(id\underset{N}{_b*_\alpha}\omega_{\eta, \xi})\ga(x)$ belongs to $\gN_\psi$ and that :
\[\Lambda_\psi[(id\underset{N}{_b*_\alpha}\omega_{\eta, \xi})\ga(x)]=(id*\omega_{\delta^{1/2}\eta, \xi})(V_\psi)\Lambda_\psi(x)\]
In (\cite{E6}, 7.6) was introduced the notion of invariant weight by an action; a normal faithful semi-finite weight $\phi$ on $A$ will be called invariant by $\ga$ if, for all $\eta\in D(_\alpha H, \nu)\cap D(H_\beta, \nu^o)$, and $x\in\gN_\phi$, we have :
\[\phi[(id\underset{N}{_b*_\alpha}\omega_\eta)\ga(x^*x)=\|\Lambda_{\phi}(x)\underset{\nu^o}{_a\otimes_\beta}\eta\|^2\]
If, moreover, $\phi$ bear the density property, a similar implementation $V'_\phi$ was constructed also in (\cite{E6}, 7.7). Moreover, with these hypothesis, it is possible to prove that there exists a normal semi-finite operator-valued weight $\gT$ from $A$ onto $b(N)$ (we shall say that the action is "weighted"), such that $\phi=\nu^o\circ b^{-1}\circ\gT$. This operator-valued weight $\gT$ satisfies, for all positive $x$ in $A$ :
\[(\gT\underset{N}{_b*_\alpha}id)\ga(x)=\ga(\gT (x))=1\underset{N}{_b\otimes_\alpha}\beta\circ b^{-1}\gT(x)\]
\newline
Let us remark that, if we define, for $n\in N$, $b^o(n)=b(n)^o$, we obtain a $*$-homomorphism from $N$ into $A^o$; moreover, for $x\in A$, let us write $\ga^o(x^o)=(.^o\underset{N}{_b*_\alpha}R)\circ a(x)$; it is straightforward to get that $(b^o, \ga^o)$ is an action of $\gG^o$ on $A^o$. 
\newline
Of course, one should write in this paragraph "right action" instead of simply action. At some stage of this paper, we shall need left-actions. A left-action of $\gG$ on a von Neumann algebra $A$ is a couple $(a, \gb)$, where :
\newline
(i) $a$ is an injective *-homomorphism from $N$ into $A$;
\newline
(ii) $\gb$ is an injective $*$-homomorphism from $A$ into $M\underset{N}{_\beta*_a}A$; 
\newline
(iii) $a$ and $\gb$ are such that, for all $n$ in $N$ :
\[\gb(a(n))=\alpha(n)\underset{N}{_\beta\otimes_a}1\]
\[(id\underset{N}{_\beta*_a}\gb)\gb=(\Gamma\underset{N}{_\beta*_a}id)\gb\]
Then, it is clear that $(a, \varsigma_N\gb)$ is an action (a right action) of $\gG^o$ on $A$, and $(a^o, (\sigma_{N}\gb)^o)$ is an action (a right action) of $\gG$ on $A^o$. Conversely, if $(b, \ga)$ is an action of $\gG$ on $A$, then, $(b^o, \sigma_{N^o}\ga^o)$ is a left-action of $\gG$ on $A^o$. 
\newline
The invariant subalgebra $A^\gb$ is defined by :
\[A^\gb=\{x\in A\cap a(N)'; \gb(x)=1\underset{N}{_\beta\otimes_a}x\}\]
and $T_\gb=(\Phi\underset{\nu}{_\beta*_a}id)\gb$ is a normal faithful operator-valued weight from $A$ onto $A^\gb$; the action $\gb$ will be said integrable if $T_\gb$ is semi-finite. It is clear that $\gb$ is integrable if and only if $(\sigma_N\gb)^o$ is integrable. 
\newline
If $(b, \ga)$ is an action of $\gG_1=(N_1, M_1, \alpha_1, \beta_1, \Gamma_1, T_1, T'_1, \nu_1)$ on a von Neumann algebra $A$, and $(a, \gb)$ a left-action of $\gG_2=(N_2, M_2, \alpha_2, \beta_2, \Gamma_2, T_2, T'_2, \nu_2)$ on $A$, such that $a(N_2)\subset b(N_1)'$, then, we shall say that the actions $\ga$ and $\gb$ commute if we have :
\[b(N_1)\subset A^\gb\]
\[a(N_2)\subset A^\ga\]
\[(\gb\underset{N_1}{_b*_{\alpha_1}}id)\ga=(id\underset{N_2}{_{\beta_2}*_a}\ga)\gb\]
Let us remark that the two first properties allow us to write the fiber products $\gb\underset{N_1}{_b*_{\alpha_1}}id$ and $id\underset{N_2}{_{\beta_2}*_a}\ga$. 

\subsection{Crossed-product (\cite{E5})}
\label{crossed}
The crossed-product of $A$ by $\mathfrak {G}$ via the action $\mathfrak a$ is the von Neumann algebra generated by $\mathfrak a(A)$ and $1\underset{N}{_b\otimes_\alpha}\widehat{M}'$ (\cite{E5}, 9.1) and is denoted $A\rtimes_\mathfrak a\mathfrak {G}$; then there exists (\cite{E5}, 9.3) an integrable action $(1\underset{N}{_b\otimes_\alpha}\hat{\alpha}, \tilde{\mathfrak a})$ of $(\widehat{\mathfrak {G}})^c$ on $A\rtimes_\mathfrak a\mathfrak {G}$. 
\newline
The biduality theorem (\cite{E5}, 11.6) says that the bicrossed-product $(A\rtimes_\mathfrak a\mathfrak {G})\rtimes_{\tilde{\mathfrak a}}\widehat{\mathfrak {G}}^o$ is canonically isomorphic to $A\underset{N}{_b*_\alpha}\mathcal L(H)$; more precisely, this isomorphism is given by :
\[\Theta (\ga\underset{N}{_b*_\alpha}id)(A\underset{N}{_b*_\alpha}\mathcal L(H))=(A\rtimes_\mathfrak a\mathfrak {G})\rtimes_{\tilde{\mathfrak a}}\widehat{\mathfrak {G}}^o\]
where $\Theta$ is the spatial isomorphism between $\mathcal L(\gH\underset{\nu}{_b\otimes_\alpha}H\underset{\nu}{_\beta\otimes_\alpha}H)$ and $\mathcal L(\gH\underset{\nu}{_b\otimes_\alpha}H\underset{\nu^o}{_{\hat{\alpha}}\otimes_\beta}H)$ implemented by $1_\gH\underset{\nu}{_b\otimes_\alpha}\sigma_\nu W^o\sigma_\nu$; the biduality theorem says also that this isomorphism sends  the action $(1\underset{N}{_b\otimes_\alpha}\hat{\beta}, \underline{\mathfrak a})$ of $\gG$ on $A\underset{N}{_b*_\alpha}\mathcal L(H)$, defined, for any $X\in A\underset{N}{_b*_\alpha}\mathcal L(H)$, by :
\[\underline{\mathfrak a}(X)=(1\underset{N}{_b\otimes_\alpha}\sigma_{\nu^o}W\sigma_{\nu^o})(id\underset{N}{_b*_\alpha}\varsigma_N)(\mathfrak a\underset{N}{_b*_\alpha}id)(X)(1\underset{N}{_b\otimes_\alpha}\sigma_{\nu^o}W\sigma_{\nu^o})^*\]
on the bidual action (of $\mathfrak{G}^{co}$) on $(A\rtimes_\mathfrak a\mathfrak {G})\rtimes_{\tilde{\mathfrak a}}\widehat{\mathfrak {G}}^o$. 
\newline
We have $(A\rtimes_\ga\gG)^{\tilde{\ga}}=\ga(A)$ (\cite{E5} 11.5), and, therefore, the normal faithful semi-finite operator-valued weight $T_{\tilde{\ga}}$ sends $A\rtimes_\mathfrak a\mathfrak {G}$ onto $\ga(A)$; therefore, starting with a  normal semi-finite weight $\psi$ on $A$, we can construct a dual weight $\tilde{\psi}$ on $A\rtimes_\mathfrak a\mathfrak {G}$ by the formula $\tilde{\psi}=\psi\circ\ga^{-1}\circ T_{\tilde{\ga}}$ (\cite{E5} 13.2). 
\newline
Moreover (\cite{E5} 13.3), the linear set generated by all the elements $(1\underset{N}{_b\otimes_\alpha}a)\mathfrak a(x)$, for all $x\in\gN_\psi$, $a\in\gN_{\widehat{\Phi}^c}\cap\gN_{\hat{T}^c}$, is a core for $\Lambda_{\tilde{\psi}}$, and it is possible to identify the GNS representation of $A\rtimes_\mathfrak a\gG$ associated to the weight $\tilde{\psi}$ with the natural representation on $H_\psi\underset{\nu}{_b\otimes_\alpha}H$ by writing :
\[\Lambda_\psi(x)\underset{\nu}{_b\otimes_\alpha}\Lambda_{\widehat{\Phi}^c}(a)=\Lambda_{\tilde{\psi}}[(1\underset{N}{_b\otimes_\alpha}a)\mathfrak a(x)]\]
which leads to the identification of $H_{\tilde{\psi}}$ with $H_\psi\underset{\nu}{_b\otimes_\alpha}H$. 
\newline
If the weight $\psi$ is $\delta$-invariant (resp. invariant) and bears the density property, then the implementation $V_\psi$ (resp. $V'_\psi$) recalled in \ref{action} is equal to $J_{\tilde{\psi}}(J_\psi\underset{N^o}{_a\otimes_\beta}J_{\widehat{\Phi}})$ (\cite{E6}, 3.2). More generally, if we write $V=J_{\tilde{\psi}}(J_\psi\underset{N^o}{_a\otimes_\beta}J_{\widehat{\Phi}})$, we have :
\[V^*=(J_\psi\underset{\nu^o}{_\alpha\otimes_\beta}J_{\widehat{\Phi}})V(J_\psi\underset{\nu}{_b\otimes_\alpha} J_{\widehat{\Phi}})\]
and, if it is an implementation of $\ga$, we shall call it a standard implementation of $\ga$. It had been proved that it is the case, for any normal semi-finite faithful weight $\psi$ on $A$, whenever the action is weighted (i.e. if there exists a normal semi-finite faithful operator-valued weight from $A$ onto $b(N)$). 
\newline
If $(a, \gb)$ is a left action of $\gG$ on $A$, we shall define the crossed product $\gG\ltimes_{\gb}A$ as the von Neumann algebra generated by $\widehat{M}\underset{N}{_\beta\otimes_a} 1$ and $b(A)$; therefore, it is 
the image under $\sigma_N$ of the crossed product $A\rtimes_{\sigma_N\gb}\gG^o$. 

\subsection{Basic Construction}
\label{basic}
Let $M_0\subset M_1$ be an inclusion of $\sigma$-finite von Neumann algebras, equipped with a normal faithful semi-finite operator-valued weight $T_1$ from $M_1$ to $M_0$. Let $\psi_0$ be a normal faithful semi-finite weight on $M_0$, and $\psi_1=\psi_0\circ T_1$.
 \newline
 Following (\cite{J}, 3.1.5(i)), the von Neumann algebra $M_2=J_{\psi_1}M'_0J_{\psi_1}$ defined on the Hilbert space $H_{\psi_1}$ will be called the basic construction made from the inclusion $M_0\subset M_1$. We have $M_1\subset M_2$, and we shall say that the inclusion $M_0\subset M_1\subset M_2$ is standard. 
 \newline
  Let us write $r$ for the inclusion of $M_0$ into $M_1$ (or the representation of $M_0$ on $H_{\psi_1}$ given by the restriction of $\pi_{\psi_1}$ to $M_0$), and let us define $s$, for any $x\in M_0$, by $s(x)=J_{\psi_1}r(x)^*J_{\psi_1}$; $s$ is a normal faithful anti-representation of $M_0$ on $H_{\psi_1}$, and $M_2=s(M_0)'$. Therefore (\ref{spatial}), the operators $\theta^{s, \psi_0^o}(\xi, \eta)$, for all $\xi$, $\eta$ in $D((H_{\psi_1})_s, \psi_0^o)$ generate a dense ideal in $M_2$. 
 \newline
 Following (\cite{EN} 10.6), for $x$ in $\gN_{T_1}$, we shall define $\Lambda_{T_1}(x)$ by the following formula, for all $z$ in $\gN_{\psi_{0}}$ :
\[\Lambda_{T_1}(x)\Lambda_{\psi_{0}}(z)=\Lambda_{\psi_1}(xz)\]
This operator belongs to $Hom_{M_{0}^o}(H_{\psi_0}, H_{\psi_1})$; if $x$, $y$ belong to $\gN_{T_1}$, then $\Lambda_{T_1}(x)\Lambda_{T_1}(y)^*$ belongs to $M_{2}$, and
$\Lambda_{T_1}(x)^*\Lambda_{T_1}(y)=T_1(x^*y)\in M_0$. 
\newline
 Using then Haagerup's construction (\cite{T}, IX.4.24), it is possible to construct a normal semi-finite faithful operator-valued weight $T_2$ from $M_2$ to $M_1$ (\cite{EN}, 10.7), which will be called the basic construction made from $T_1$. If $x$, $y$ belong to $\gN_{T_1}$, then the operators $\Lambda_{T_1}(x)\Lambda_{T_1}(y)^*$ form a dense sub-$*$algebra of $M_2$, included into
$\gM_{T_{2}}$, and we have $T_{2}(\Lambda_{T_1}(x)\Lambda_{T_1}(y)^*)=xy^*$.  The operator-valued weight $T_2$ is characterized by the equality (\cite{EN}, 10.3) :
\[\frac{d\psi_1\circ T_2}{d\psi_0^o}=\frac {d\psi_1}{d(\psi_0\circ T_1)^o}=\Delta_{\psi_1}\]
from which, writing $\psi_2=\psi_1\circ T_2$, we get that :
\[\sigma_t^{\psi_2}(\Lambda_{T_1}(x)\Lambda_{T_1}(y)^*)=\Lambda_{T_1}(\sigma_t^{\psi_1}(x))\Lambda_{T_1}(\sigma_t^{\psi_1}(y^*))^*\] 
The operator-valued weight $T_2$ from $M_2$ to $M_1$ will be called the basic construction made from the operator-valued weight $T_1$ from $M_1$ to $M_0$. Using (\cite{EN}, 3.7 and 10.6 (v)), we easily get that, for any $x$, $y$ in $\gN_{T_1}\cap\gN_{\psi_1}\cap\gN_{T_1}^*\cap\gN_{\psi_1}^*$, we have $T_2(\Lambda_{T_1}(x)\Lambda_{T_1}(y^*)^*)=xy$, and :
\[\|\Lambda_{\psi_2}(\Lambda_{T_1}(x)\Lambda_{T_1}(y^*)^*)\|=\|\Lambda_{\psi_1}(x)\underset{\psi_0}{_s\otimes_r}\Lambda_{\psi_1}(y)\|\]
where $r$ is the inclusion of $M_0$ into $M_1$, and, for $a\in M_0$, $s(a)=J_{\psi_1}a^*J_{\psi_1}$; so, we can identify $H_{\psi_2}$ with $H_{\psi_1}\underset{\psi_0}{_s\otimes_r}H_{\psi_1}$ by writing $\Lambda_{\psi_2}(\Lambda_{T_1}(x)\Lambda_{T_1}(y^*)^*)=\Lambda_{\psi_1}(x)\underset{\psi_0}{_s\otimes_r}\Lambda_{\psi_1}(y)$; then, we identify $\Delta_{\psi_2}^{it}$ with $\Delta_{\psi_1}^{it}\underset{\psi_0}{_s\otimes_r}\Delta_{\psi_1}^{it}$ (here, this relative tensor product of operators means that there exists a bounded operator with natural values on elementary tensors) and $J_{\psi_2}$ with $\sigma_{M_0^o}(J_{\psi_1}\underset{M_0}{_s\otimes_r}J_{\psi_1})$.
\newline
Then, for any $\xi\in D((H_{\psi_1})_s, \psi_0^o)$ and $\eta\in D((H_{\psi_1})_s, \psi_0^o)\cap \mathcal D(\Delta_{\psi_1}^{1/2})$ such that $\Delta_{\psi_1}^{-1/2}\eta$ belongs to $D((H_{\psi_1})_s, \psi_0^o)$, we have $\Lambda_{\psi_2}(\theta^{s, \psi_0^o}(\xi, \eta))=\xi\underset{\psi_0}{_s\otimes_r}J_{\psi_1}\Delta_{\psi_1}^{1/2}\eta$.
\newline
Using similar arguments as in (\cite{E6}, 4.7(ii)), we can prove that there exists a family $(e_i)_{i\in I}$, which is an orthogonal $(s, \psi_0^o)$-basis of $H_{\psi_1}$ , such that each vector $e_i$ belongs to $\mathcal D(\Delta_{\psi_1}^{1/2})$; we can prove then, as in (\cite{E6}, 4.7(iii)), that $\psi_2=\sum_i\omega_{\Delta_{\psi_1}^{1/2}e_i}$. 
\newline
Let $\mathcal T_{\psi_1, T_1}$ be the Tomita algebra associated to the operator-valued weight $T_1$ and the weight $\psi_1$ (\cite{EN}, 10.12, and \cite{E5}, 2.2.1), which is made of elements $x$ in $\gN_{\psi_1}\cap\gN_{\psi_1}^*\cap\gN_{T_1}\cap\gN_{T_1}^*$, which are analytic with respect to $\sigma_t^{\psi_1}$, and such that, for any $z\in\mathbb{C}$, $\sigma_z^{\psi_1}(x)$ belongs to $\gN_{\psi_1}\cap\gN_{\psi_1}^*\cap\gN_{T_1}\cap\gN_{T_1}^*$; such elements are a dense $*$ subalgebra of $M_1$. Moreover, it is possible to prove (\cite{DC1}, 1.4) that an element $X\in M_2$ belongs to $\gN_{\psi_2}$ if and only if there exists $\Xi\in H_{\psi_2}$ such that, for any $x$, $y$ in $\mathcal T_{\psi_1, T_1}$, we have :
\[(\Lambda_{\psi_1}(x)\underset{\psi_0}{_s\otimes_r}\Lambda_{\psi_1}(y)|\Xi)=
(\Lambda_{\psi_1}(x)|X\Lambda_{\psi_1}(\sigma_{-i}^{\psi_1}(y^*)))\]
and, then, we have $\Xi=\Lambda_{\psi_2}(X)$.

\section{Integrable actions of a measured quantum groupoid}
\label{integrable}
In that chapter are generalized, up to measured quantum groupoids, results about integrable actions ((\cite{V2}, 5.3, \cite{DC1}, 2.1); namely, if $(b,\ga)$ is an integrable action of $\gG$ on a von Neumann algebra $A$ (the definition had been given in \ref{action}), we construct then a representation $\pi_\ga$ of the crossed product on the Hilbert space $L^2(A)$(\ref{pia}), whose image is the von Neumann algebra $s(A^\ga)'$ given by the standard construction made from the inclusion $A^\ga\subset A$; moreover is constructed an isometry $G$ from  $L^2(s(A^\ga)')$ into $L^2(A\rtimes_\ga\gG)$ (\ref{G}), which is a unitary if and only if the representation $\pi_\ga$ is faithful; following (\cite{DC1}, 2.7), we say that the integrable action $(b, \ga)$ is then Galois (\ref{defGalois}).

\subsection{Lemma}
\label{lemintegrable}
{\it Let $(b, \ga)$ be an integrable action of a measured quantum groupoid on a von Neumann algebra $A$; let $\psi_0$ be a normal faithful semi-finite weight on $A^\ga$, and $\psi_1=\psi_0\circ T_\ga$ be the lifted normal semi-finite faithful weight on $A$; let $(\xi_i)_{i\in I}$ be a family of vectors in $ D((H_{\psi_1})_b, \nu^o)$ such that $\psi_0(x)=\sum_i\omega_{\xi_i}(x)$, for all positive $x\in A^\ga$; then, there exists an isometry $\mathcal V$ from $H_{\psi_1}$ into $\oplus_i (H_{\psi_1}\underset{\nu}{_b\otimes_\alpha}H)=(H_{\psi_1}\underset{\nu}{_b\otimes_\alpha}H)\otimes l^2(I)$ such that :
\newline
(i) for all $y\in A$, $(\ga(y)\otimes 1_{l^2(I)})\mathcal V=\mathcal V y$; 
\newline
(ii) for all $n\in N$, $(1\underset{N}{_b\otimes_\alpha}\hat{\alpha}(n)\otimes 1_{l^2(I)})\mathcal V=\mathcal V a(n)$. }
\begin{proof}
Let $(\eta_j)_{j\in J}$ be an orthogonal $(b, \nu^o)$ basis of $H_{\psi_1}$; we have, for any $i$, $j$ and $x\in\gN_{\psi_1}$:
\[[(\omega_{\xi_i, \eta_j}\underset{N}{_b*_\alpha}id)\ga(x)]^*[(\omega_{\xi_i, \eta_j}\underset{N}{_b*_\alpha}id)\ga(x)]=(\omega_{\xi_i}\underset{N}{_b*_\alpha}id)[\ga(x^*)(\theta^{b, \nu^o}(\eta_j, \eta_j)\underset{N}{_b\otimes_\alpha}1)\ga(x)]\]
and, therefore :
\[\Phi([(\omega_{\xi_i, \eta_j}\underset{N}{_b*_\alpha}id)\ga(x)]^*[(\omega_{\xi_i, \eta_j}\underset{N}{_b*_\alpha}id)\ga(x)])\leq \Phi[(\omega_{\xi_i}\underset{N}{_b*_\alpha}id)\ga(x^*x)]=\omega_{\xi_i}\circ T_\ga(x^*x)\leq \psi_1(x^*x)\]
So, for any $i$, $j$ and $x\in\gN_{\psi_1}$, $(\omega_{\xi_i, \eta_j}\underset{N}{_b*_\alpha}id)\ga(x)$ belongs to $\gN_\Phi$; moreover, we have :
\[\Phi([\sum_j(\omega_{\xi_i, \eta_j}\underset{N}{_b*_\alpha}id)\ga(x)]^*[\sum_j(\omega_{\xi_i, \eta_j}\underset{N}{_b*_\alpha}id)\ga(x)])
=\omega_{\xi_i}\circ T_\ga(x^*x)\]
and, therefore :
\[\sum_i \Phi([\sum_j(\omega_{\xi_i, \eta_j}\underset{N}{_b*_\alpha}id)\ga(x)]^*[\sum_j(\omega_{\xi_i, \eta_j}\underset{N}{_b*_\alpha}id)\ga(x)])=\psi_1(x^*x)\]
which proves that we can define now $\mathcal V$, for all $x\in\gN_{\psi_1}$ by :
\[\mathcal V\Lambda_{\psi_1}(x)=\oplus_i\sum_j \eta_j\underset{\nu}{_b\otimes_\alpha}\Lambda_\Phi((\omega_{\xi_i, \eta_j}\underset{N}{_b*_\alpha}id)\ga(x))\]
As, for $x\in\gN_{\psi_1}$, we have $\|\mathcal V\Lambda_{\psi_1}(x)\|^2=\psi_1(x^*x)$, we can extend $\mathcal V$ to an isometry from $H_{\psi_1}$ into $\oplus_{i\in I}(H_{\psi_1}\underset{\nu}{_b\otimes_\alpha}H)=(H_{\psi_1}\underset{\nu}{_b\otimes_\alpha}H)\otimes l^2(I)$. 
Let now $y$ be in $A$; we have :
\begin{eqnarray*}
(\ga(y)\otimes 1_{l^2(I)})\mathcal V\Lambda_{\psi_1}(x)
&=&
\oplus_i\sum_j \ga(y)(\eta_j\underset{\nu}{_b\otimes_\alpha}\Lambda_\Phi((\omega_{\xi_i, \eta_j}\underset{N}{_b*_\alpha}id)\ga(x)))\\
&=&
\oplus_i\sum_j\sum_k \eta_k\underset{N}{_b\otimes_\alpha}(\omega_{\eta_k, \eta_j}\underset{N}{_b*_\alpha}id)\ga(y)\Lambda_\Phi((\omega_{\xi_i, \eta_j}\underset{N}{_b*_\alpha}id)\ga(x))\\
&=&
\oplus_i\sum_j\sum_k \eta_k\underset{N}{_b\otimes_\alpha}\Lambda_\Phi((\omega_{\eta_k, \eta_j}\underset{N}{_b*_\alpha}id)\ga(y)(\omega_{\xi_i, \eta_j}\underset{N}{_b*_\alpha}id)\ga(x))\\
&=&
\oplus_i\sum_k\eta_k\underset{N}{_b\otimes_\alpha}\Lambda_\Phi(\omega_{\xi_i, \eta_k}\underset{N}{_b*_\alpha}id)\ga(yx))\\
&=&
\mathcal V \Lambda_{\psi_1}(yx)
\end{eqnarray*}
and, therefore, $(\ga(y)\otimes 1_{l^2(I)})\mathcal V=\mathcal V y$. 
\newline
Let $n$ be in the Tomita algebra of the weight $\nu$; we have :
\begin{eqnarray*}
(1\underset{N}{_b\otimes_\alpha}\hat{\alpha}(n))\otimes  1_{l^2(I)})\mathcal V\Lambda_{\psi_1}(x)
&=&
\oplus_i\sum_j \eta_j\underset{\nu}{_b\otimes_\alpha}\hat{\alpha}(n)\Lambda_\Phi((\omega_{\xi_i, \eta_j}\underset{N}{_b*_\alpha}id)\ga(x))\\
&=&
\oplus_i\sum_j \eta_j\underset{\nu}{_b\otimes_\alpha}\Lambda_\Phi((\omega_{\xi_i, \eta_j}\underset{N}{_b*_\alpha}id)\ga(x)\beta(\sigma_{-i/2}^\nu(n))\\
&=&
\oplus_i\sum_j \eta_j\underset{\nu}{_b\otimes_\alpha}\Lambda_\Phi((\omega_{\xi_i, \eta_j}\underset{N}{_b*_\alpha}id)\ga(xb(\sigma_{-i/2}^\nu(n))))\\
&=&
\mathcal V \Lambda_{\psi_1}(xb(\sigma_{-i/2}^\nu(n))\\
&=&
\mathcal V a(n)\Lambda_{\psi_1}(x)
\end{eqnarray*}
which, by contnuity, remains true for all $n\in N$. 
\end{proof}

\subsection{Theorem}
\label{thintegrable}
{\it Let $(b, \ga)$ be an integrable action of a measured quantum groupoid on a von Neumann algebra $A$; let $\psi_0$ be a normal faithful semi-finite weight on $A^\ga$, and $\psi_1=\psi_0\circ T_\ga$ be the lifted normal semi-finite faithful weight on $A$; then, the weight $\psi_1$ is $\delta$-invariant, and bears the density property, in the sense of \ref{action}}

\begin{proof}
Let's use the notations of \ref{lemintegrable}; let $x$ be in $\gN_{\psi_1}$, and $\eta\in D(_\alpha H, \nu)\cap\mathcal D(\delta^{1/2})$, such that $\delta^{1/2}\eta$ belongs to $D((H)_\beta, \nu^o)$; we have, using the isometry $\mathcal V$ and (\cite{E5}, 8.2) :
\begin{eqnarray*}
\|\Lambda_{\psi_1}(x)\underset{\nu}{_b\otimes_\alpha}\delta^{1/2}\eta\|^2
&=&
\oplus_i\|\sum_j\eta_j\underset{\nu}{_b\otimes_\alpha}\Lambda_{\Phi}((\omega_{\xi_i, \eta_j}\underset{N}{_b*_\alpha}id)\ga(x))\underset{\nu^o}{_{\hat{\alpha}}\otimes_\beta}\delta^{1/2}\eta\|^2\\
&=&
\oplus\|\Lambda_\Phi(\alpha(<\eta_j, \eta_j>_{b, \nu^o})(\omega_{\xi_i, \eta_j}\underset{N}{_b*_\alpha}id)\ga(x))\underset{\nu^o}{_{\hat{\alpha}}\otimes_\beta}\delta^{1/2}\eta\|^2\\
&=&
\oplus\|\Lambda_\Phi((\omega_{\xi_i, \eta_j}\underset{N}{_b*_\alpha}id)\ga(x))\underset{\nu^o}{_{\hat{\alpha}}\otimes_\beta}\delta^{1/2}\eta\|^2\\
&=&
\sum_i\Phi (\sum_j(id\underset{N}{_b*_\alpha}\omega_\eta)\Gamma[(\omega_{\xi_i, \eta_j}\underset{N}{_b*_\alpha}id)\ga(x)^*(\omega_{\xi_i, \eta_j}\underset{N}{_b*_\alpha}id)\ga(x)]\\
&=&
\sum_i\Phi[(id\underset{N}{_b*_\alpha}\omega_\eta)\Gamma(\omega_{\xi_i}\underset{N}{_b*_\alpha}*id)\ga(x^*]\\
&=&
\sum_i \omega_{\xi_i}\circ T_\ga(id\underset{N}{_b*_\alpha}\omega_\eta)\ga(x^*x))\\
&=&
\psi_1[(id\underset{N}{_b*_\alpha}\omega_\eta)\ga(x^*x)]
\end{eqnarray*}
which proves that $\psi_1$ is $\delta$-invariant; moreover, if we take the Tomita algebra relatively to the weight $\psi_1$ and the operator-valued weight $T_\ga$, we get that the weight $\psi_1$ bears the density property. 
\end{proof}

\subsection{Proposition}
\label{biaction}
{\it Let $\gG_i=(N_i, M_i, \alpha_i, \beta_i, \Gamma_i, T_i, T'_i, \nu_i)$ ($ i=1,2$) be two measured quantum groupoids, $(b, \ga)$ an action of $\gG_1$ on a von Neumann algebra $A$, and $(a, \gb)$ a left-action of $\gG_2$ on $A$; let us suppose that the actions $\ga$ and $\gb$ commute; then :
\newline
(i) the operator-valued weight $T_\gb$ from $A$ onto $A^\gb$ satisfies :
\[(T_\gb\underset{\nu_1}{_b*_{\alpha_1}}id)\ga=\ga\circ T_\gb\]
(ii) if $\gb$ is integrable and if $A^\gb=b(N_1)$, the weight $\phi_1=\nu_1\circ b^{-1}\circ T_\gb$ is a normal semi-finite faithful weight on $A$, invariant under the action $\ga$, $\delta_{\gG_2}$-invariant under the action $\gb$, and bears the density property. }

\begin{proof}
Result (i) is straightforward, using the definition of commuting actions. With the hypothesis of (ii), we get that $T_\gb$ is a normal semi-finite faithful operator-valued weight from $A$ onto $b(N_1)$, that $\phi_1$ is a normal semi-finite faithful weight on $A$ which satisfies, for all $x\in A^+$ :
\[(T_b\underset{\nu_1}{_b*_{\alpha_1}}id)\ga(x)=\ga\circ T_\gb(x)=1\underset{N_1}{_b\otimes_{\alpha_1}}\beta_1\circ b^{-1}\circ T_\gb(x)\]
\[(\phi_1\underset{\nu_1}{_b*_{\alpha_1}}id)\ga(x)=\beta_1\circ b^{-1}\circ T_\gb(x)\]
 from which we get that $\phi_1$ is invariant by $\ga$. On the other hand, $\phi_1$ is $\delta_{\gG_2}$-invariant under the action $\gb$, and bears the density property by \ref{thintegrable}. \end{proof}

\subsection{Lemma}
\label{lemV}
{\it Let $\gG$ be a measured quantum groupoid; let $V$ be a corepresentation of $\gG$ on a $N-N$ bimodule $_a\gH_b$ (\cite{E5}, 5.1), and let $(b, \ga)$ the canonical action implemented by $V$ on $a(N)'$ by $\ga(x)=V(x\underset{N^o}{_a\otimes_\beta}1)V^*$ (\cite{E5}, 6.6); then we have :}
\[(a(N)')^\ga=a(N)'\cap b(N)'\cap\{(id*\omega_{\xi, \eta})(V), \xi\in D(_\alpha H, \nu), \eta\in D((H_\phi)_\beta, \nu^o)\}'\]

\begin{proof}
Clear\end{proof}

\subsection{Lemma}
\label{A2} 
{\it Let $(b, \ga)$ be an integrable action of a measured quantum groupoid $\gG$ on a von Neumann algebra $A$, and let $\psi_0$ be a normal semi-finite faithful weight on $A^\ga$, and $\psi_1=\psi_0\circ T_\ga$ be the normal semi-finite faithful lifted weight on $A$; let $V_{\psi_1}$ be the standard implementation of $\ga$ defined in \ref{action}:  let $s(A^\ga)'=J_{\psi_1} (A^\ga)'J_{\psi_1}$ the basic construction made from the inclusion $A^\ga\subset A$ (cf. \ref{basic}); then, we have :}
\[s(A^\ga)'=(A\cup a(N)\cup \{(id*\omega_{\eta, \xi})(V_{\psi_1}^*), \xi\in D(_\alpha H, \nu), \eta\in D((H)_\beta, \nu^o)\})"\]
\begin{proof}
Using \ref{lemV}, we get $A^\ga=A\cap b(N)'\cap\{(id*\omega_{\xi, \eta})(V_{\psi_1}), \xi\in D(_\alpha H, \nu), \eta\in D((H_\phi)_\beta, \nu^o)\}'$, and, therefore :
\[J_{\psi_1}A^\ga J_{\psi_1}=A'\cap a(N)'\cap J_{\psi_1}\{(id*\omega_{\xi, \eta})(V_{\psi_1}), \xi\in D(_\alpha H, \nu), \eta\in D((H_\phi)_\beta, \nu^o)\}'J_{\psi_1}\]
As $V_{\psi_1}(J_{\psi_1}\underset{N}{_b\otimes_\alpha}\hat{J})=(J_{\psi_1}\underset{N}{_b\otimes_\alpha}\hat{J})V_{\psi_1}^*$, we have :
\[J_{\psi_1}(id*\omega_{\xi, \eta})(V_{\psi_1})J_{\psi_1}=(id*\omega_{\hat{J}\xi, \hat{J}\eta})(V_{\psi_1}^*)\]
and we get :
\[J_{\psi_1}A^\ga J_{\psi_1}=A'\cap a(N)'\cap \{(id*\omega_{\eta, \xi})(V_{\psi_1}^*), \xi\in D(_\alpha H, \nu), \eta\in D((H)_\beta, \nu^o)\}'\]
from which we get the result. 
\end{proof}

\subsection{Theorem}
\label{pia}
{\it Let $(b, \ga)$ be an integrable action of a measured quantum groupoid $\gG$ on a von Neumann algebra $A$, let $V_{\psi_1}$ be the standard implementation of $\ga$, as defined in \ref{crossed}; let us denote by $r$ the injection of $A^\ga$ into $A$, and let us write $s(x)=J_{\psi_1}r(x)^*J_{\psi_1}$ for any $x\in A^\ga$; then $s(A^\ga)'$ is the basic construction made from the inclusion $A^\ga\subset A$ (cf. \ref{basic}); then, there exists a normal surjective $*$-homomorphism $\pi_\ga$ from the crossed-product $A\rtimes_\ga\gG$ onto $s(A^\ga)'$, called the Galois homomorphism associated to the integrable action $(b, \ga)$, such that, for all $x\in A$, $n\in N$, $\xi\in D(_\alpha H, \nu)$, $\eta\in D((H)_\beta, \nu^o)$:
\[\pi_\ga(\ga(x))=x\]
\[\pi_\ga(1\underset{N}{_b\otimes_\alpha}\hat{\alpha}(n))=a(n)\]
\[\pi_\ga(1\underset{N}{_b\otimes_\alpha}(\omega_{\eta, \xi}*id)[(W^o)^*])=\pi_\ga(1\underset{N}{_b\otimes_\alpha}(\omega_{\eta, \xi}*id)[(\hat{J}\underset{N}{_\beta\otimes_\alpha}\hat{J})W^*(\hat{J}\underset{N}{_\beta\otimes_{\hat{\alpha}}}\hat{J})])=(id*\omega_{\eta, \xi})(V_{\psi_1})\]
For simplification, we shall write $\mu(m)=\pi_\ga(1\underset{N}{_b\otimes_\alpha}m)$, for any $m\in\widehat{M}'$, and we obtain this way a representation of $\widehat{M}'$ on $\mathcal L(H_{\psi_1})$.}
\begin{proof}
Let us use the notations of \ref{lemintegrable} and \ref{thintegrable}; let's suppose that $\eta$ belongs also to $\mathcal D(\delta^{-1/2})$ and that $\delta^{-1/2}\eta$ belongs to $D(_\alpha H, \nu)$; then, we have :
\begin{multline*}
(1\underset{N}{_b\otimes_\alpha}(\omega_{\eta, \xi}*id)[(\hat{J}\underset{N}{_\beta\otimes_\alpha}\hat{J})W^*(\hat{J}\underset{N}{_\beta\otimes_{\hat{\alpha}}}\hat{J})])\otimes 1_{l^2(I)})\mathcal V\Lambda_{\psi_1}(x)=\\
\oplus_i\sum_j \eta_j\underset{\nu}{_b\otimes_\alpha}(\omega_{\eta, \xi}*id)[(\hat{J}\underset{N}{_\beta\otimes_\alpha}\hat{J})W^*(\hat{J}\underset{N}{_\beta\otimes_{\hat{\alpha}}}\hat{J})])\Lambda_\Phi((\omega_{\xi_i, \eta_j}\underset{N}{_b*_\alpha}id)\ga(x))
\end{multline*}
which, using (\cite{E5}, 3.10(ii) applied to $\gG^o$, 3.8(vi)), and the identification of $H_{\Phi\circ R}$ with $H$ made in \ref{MQG}) is equal to :
\[\oplus_i\sum_j \eta_j\underset{\nu}{_b\otimes_\alpha}\Lambda_\Phi((id\underset{N}{_\beta*_\alpha}\omega_{\delta^{-1/2}\eta, \xi})\Gamma[(\omega_{\xi_i, \eta_j}\underset{N}{_b*_\alpha}id)\ga(x)])\]
or, to :
\[\oplus_i\sum_j \eta_j\underset{\nu}{_b\otimes_\alpha}\Lambda_\Phi((\omega_{\xi_i, \eta_j}\underset{N}{_b*_\alpha}id)\ga[(id\underset{N}{_\beta*_\alpha}\omega_{\delta^{-1/2}\eta, \xi})\ga(x)])\]
which is, using (\cite{E5}, 8.4), equal to :
\[\mathcal V\Lambda_{\psi_1}[(id\underset{N}{_\beta*_\alpha}\omega_{\delta^{-1/2}\eta, \xi})\ga(x)])=
\mathcal V (id*\omega_{\eta, \xi})(V_{\psi_1})\Lambda_{\psi_1}(x)\]
from which, by density, we get that :
\[(1\underset{N}{_b\otimes_\alpha}(\omega_{\eta, \xi}*id)[(\hat{J}\underset{N}{_\beta\otimes_\alpha}\hat{J})W^*(\hat{J}\underset{N}{_\beta\otimes_{\hat{\alpha}}}\hat{J})])\otimes 1_{l^2(I)})\mathcal V=\mathcal V (id*\omega_{\eta, \xi})(V_{\psi_1})\]
which, by density and continuity, remains true for any $\eta$ in $D(H)_\beta, \nu^o)$. Using now \ref{MQG}, we get that the weak closure of the linear span of all operators of the form 
\[(\omega_{\eta, \xi}*id)[(\hat{J}\underset{N}{_\beta\otimes_\alpha}\hat{J})W^*(\hat{J}\underset{N}{_\beta\otimes_{\hat{\alpha}}}\hat{J})]\]
for all $\xi\in D(_\alpha H, \nu)$, $\eta\in D((H)_\beta, \nu^o)$, is equal to the von Neumann algebra $\widehat{M}'$; therefore, we get that, for any $y\in \widehat{M}'$,  the image of $(1\underset{N}{_b\otimes_\alpha}y\otimes 1_{l^2(I)})\mathcal V$ is included into the image of $\mathcal V$, which means that $\mathcal V\mathcal V^*(1\underset{N}{_b\otimes_\alpha}y\otimes 1_{l^2(I)})\mathcal V=(1\underset{N}{_b\otimes_\alpha}y\otimes 1_{l^2(I)})\mathcal V$; therefore, we have, for any $y\in\widehat{M}'$, $\mathcal V\mathcal V^*(1\underset{N}{_b\otimes_\alpha}y\otimes 1_{l^2(I)})\mathcal V\mathcal V^*=(1\underset{N}{_b\otimes_\alpha}y\otimes 1_{l^2(I)})\mathcal V\mathcal V^*$, which proves that $\mathcal V\mathcal V^*$ commutes with $1\underset{N}{_b\otimes_\alpha}\widehat{M}'\otimes 1_{l^2(I)}$. Using \ref{lemintegrable}, we easily get that $\mathcal V\mathcal V^*$ commutes also with $\ga(A)\otimes 1_{l^2(I)}$, and, therefore, that it commutes with $A\rtimes_\ga\gG\otimes 1_{l^2(I)}$. Let us write now, for any $z\in A\rtimes_\ga\gG$ :
\[\pi_\ga(z)=\mathcal V^*(z\otimes 1_{l^2(I)})\mathcal V\]
Thanks to this commutation property, $\pi_\ga$ is a $*$-homomorphism from $A\rtimes_\ga\gG$ into $\mathcal L(H_{\psi_1})$. Using now \ref{A2}, we get that the image of $\pi_\ga$ is $s(A^\ga)'$. 
\end{proof}

\subsection{Lemma}
\label{lemR}
{\it With the notations of \ref{pia}, we have, for all $m\in\widehat{M}'$ :}
\[\pi_\ga(1\underset{N}{_b\otimes_\alpha}\hat{R}^c(m))=J_{\psi_1}\pi_\ga(1\underset{N}{_b\otimes_\alpha}m^*)J_{\psi_1}\]

\begin{proof}
Let $\xi\in D(_\alpha H, \nu)$, $\eta\in D(H_\beta, \nu^o)$; using \ref{pia}, we get that :
\[J_{\psi_1}\pi_\ga(1\underset{N}{_b\otimes_\alpha}(\omega_{\eta, \xi}*id)[(\hat{J}\underset{N}{_\beta\otimes_\alpha}\hat{J})W^*(\hat{J}\underset{N}{_\beta\otimes_{\hat{\alpha}}}\hat{J})])^*J_{\psi_1}\]
is equal to $J_{\psi_1}(id*\omega_{\eta, \xi})(V_{\psi_1})^*J_{\psi_1}$, which, using \ref{action}, is equal to $(i*\omega_{\hat{J}\xi, \hat{J}\eta})(V_{\psi_1})$, and, using \ref{pia} again, is equal to :
\[\pi_\ga(1\underset{N}{_b\otimes_\alpha}(\omega_{\hat{J}\xi, \hat{J}\eta}*id)[(\hat{J}\underset{N}{_\beta\otimes_\alpha}\hat{J})W^*(\hat{J}\underset{N}{_\beta\otimes_{\hat{\alpha}}}\hat{J})])\]
which is $\pi_\ga((1\underset{N}{_b\otimes_\alpha}\hat{J}(\omega_{\xi, \eta}*id)(W^*)\hat{J})$, and, using (\cite{E5}, 3.11(iii)), is equal to :
\[\pi_\ga(1\underset{N}{_b\otimes_\alpha}J(\omega_{\xi, \eta}*id)[(\hat{J}\underset{N^o}{_\alpha\otimes_{\hat{\beta}}}\hat{J})W(\hat{J}\underset{N^o}{_\alpha\otimes_{\hat{\beta}}}\hat{J})]J)\]
which is :
\[\pi_\ga(1\underset{N}{_b\otimes_\alpha}J(\omega_{\eta, \xi}*id)[(\hat{J}\underset{N}{_\beta\otimes_\alpha}\hat{J})W^*(\hat{J}\underset{N}{_\beta\otimes_{\hat{\alpha}}}\hat{J})])^*J)\]
which is $\pi_\ga(1\underset{N}{_b\otimes_\alpha}\hat{R}^c[(\omega_{\eta, \xi}*id)[(\hat{J}\underset{N}{_\beta\otimes_\alpha}\hat{J})W^*(\hat{J}\underset{N}{_\beta\otimes_{\hat{\alpha}}}\hat{J})]])$; we get then the result by density. \end{proof}
\subsection{Theorem}
\label{G}
{\it Let $(b, \ga)$ an integrable action of $\gG$ on a von Neumann algebra $A$, $\pi_\ga$ the Galois homomorphism associated by \ref{pia}; let $\psi_0$ be a normal semi-finite faithful weight on $A^\ga$, and $\psi_1=\psi_0\circ T_\ga$; let $a$ be the representation of $N$ on $H_{\psi_1}$ defined, for $n\in N$, by  :
\[a(n)=J_{\psi_1}b(n^*)J_{\psi_1}\]
Let us write $r$ for the injection of $A^\ga$ into $A$, and $s$ for the antirepresentation of $A^\ga$ on $H_{\psi_1}$ given, for $a\in A^\ga$, by $s(a)=J_{\psi_1}r(a^*)J_{\psi_1}$. Then :
\newline
(i) there exists an isometry $G$ from $H_{\psi_1}\underset{\psi_0}{_s\otimes_r}H_{\psi_1}$ into $H_{\psi_1}\underset{\nu}{_b\otimes_\alpha}H$, such that :
\[G(\Lambda_{\psi_1}(x)\underset{\psi_0}{_s\otimes_r}\zeta)=\sum_i e_i\underset{\nu}{_b\otimes_\alpha}\Lambda_\Phi[(\omega_{\zeta, e_i}\underset{N}{_b*_\alpha}id)\ga(x)]\]
for all $x$ in $\gN_{T_\ga}\cap\gN_{\psi_1}$, $\zeta\in D((H_{\psi_1})_b, \nu^o)$, and for all $(b, \nu^o)$-orthogonal basis $(e_i)_{i\in I}$ of $H_{\psi_1}$. Moreover, for any $n\in N$, $a\in A^\ga$, we have :
\[G(b(n)\underset{A^\ga}{_s\otimes_r}1)=(1\underset{N}{_b\otimes_\alpha}\beta(n))G\]
\[G(1\underset{A^\ga}{_s\otimes_r}b(n))=(1\underset{N}{_b\otimes_\alpha}\hat{\beta}(n))G\]
\[G(1\underset{A^\ga}{_s\otimes_r}a(n))=(a(n)\underset{A^\ga}{_b\otimes_\alpha}1)G\]
\[G(r(a)\underset{A^\ga}{_s\otimes_r}1)=(r(a)\underset{A^\ga}{_b\otimes_\alpha}1)G\]
\[G(1\underset{A^\ga}{_s\otimes_r}s(a))=(s(a)\underset{A^\ga}{_b\otimes_\alpha}1)G\]
(ii) for any $e\in \gN_\Phi$, we have :
\[(1\underset{N}{_b\otimes_\alpha}J_\Phi eJ_\Phi)G(\Lambda_{\psi_1}(x)\underset{\psi_0}{_s\otimes_r}\zeta)
=\ga(x)(\zeta\underset{\nu}{_b\otimes_\alpha}J_\phi\Lambda_\Phi(e))\]
(iii) for all $\zeta'$ in $D((H_{\psi_1})_b, \nu^o)$, $(\omega_{\zeta, \zeta'}\underset{N}{_b*_\alpha}id)\ga(x)$ belongs to $\gN_\Phi$, and we have :
\[\Lambda_\Phi[(\omega_{\zeta, \zeta'}\underset{N}{_b*_\alpha}id)\ga(x)]=(\omega_{\Lambda_{\psi_1}(x), \zeta'}*id)(G)\zeta\]
(iv) for any $a'\in A$, $Y\in\widehat{M}'$, we have :
\[\ga(a')G=G(a'\underset{A^\ga}{_s\otimes_r}1)\]
\[(1\underset{N}{_b\otimes_\alpha}Y)G=G(\pi_\ga(1\underset{N}{_b\otimes_\alpha}Y){_s\otimes_r}1)\]
(v) the projection $GG^*$ commutes with $A\rtimes_\ga\gG$, and, for any $X\in A\rtimes_\ga\gG$, we have :
\[\pi_\ga(X)\underset{A^\ga}{_s\otimes_r}1_{H_{\psi_1}}=G^*XG\]
(vi) for any $t\in\mathbb{R}$, we have $\Delta_{\tilde{\psi_1}}^{it}G=G\Delta_{\psi_2}^{it}$.
\newline
(vii) for all $t\in\mathbb{R}$, we have :
\[\tilde{G}(\Delta_{\psi_1}^{it}\underset{\psi_0}{_s\otimes_r}\Delta_{\psi_1}^{it})=
((\delta\Delta_{\widehat{\Phi}})^{-it}\underset{\nu^o}{_\alpha\otimes_b}\Delta_{\psi_1}^{it})\tilde{G}\]}

\begin{proof}
As :
\[[(\omega_{\zeta, e_i}\underset{N}{_b*_\alpha}id)\ga(x)]^*[(\omega_{\zeta, e_i}\underset{N}{_b*_\alpha}id)\ga(x)]\leq (\omega_{\zeta}\underset{N}{_b*_\alpha}id)\ga(x^*x)\]
we get that :
\[\Phi( [(\omega_{\zeta, e_i}\underset{N}{_b*_\alpha}id)\ga(x)]^*[(\omega_{\zeta, e_i}\underset{N}{_b*_\alpha}id)\ga(x)])
\leq \Phi[(\omega_{\zeta}\underset{N}{_b*_\alpha}id)\ga(x^*x)]=
\omega_{\zeta}\circ T_\ga (x^*x)\]
and we get that $[(\omega_{\zeta, e_i}\underset{N}{_b*_\alpha}id)\ga(x)]$ belongs to $\gN_\Phi$; defining $G$ by the formula given in (i), we obtain, for $x$, $x'$, in $\gN_{T_\ga}\cap\gN_{\psi_1}$, $\zeta$, $\zeta'$ in $D(H_{\psi_1})_b, \nu^o)$, that :
\[(G(\Lambda_{\psi_1}(x)\underset{\psi_0}{_s\otimes_r}\zeta)|G(\Lambda_{\psi_1}(x')\underset{\psi_0}{_s\otimes_r}\zeta'))\]
is equal to :
\[\sum_i(\Lambda_\Phi[(\omega_{\zeta, e_i}\underset{N}{_b*_\alpha}id)\ga(x)]|\Lambda_\Phi[(\omega_{\zeta', e_i}\underset{N}{_b*_\alpha}id)\ga(x')])
=
(T_\ga(x'^*x)\zeta|\zeta')\]
or, to $(\Lambda_{\psi_1}(x)\underset{\psi_0}{_s\otimes_r}\zeta|\Lambda_{\psi_1}(x')\underset{\psi_0}{_s\otimes_r}\zeta')$
which implies that this formula defines an isometry which can be extended by continuity to $H_{\psi_1}\underset{\psi_0}{_b\otimes_a}H_{\psi_1}$ and does not depend upon the choice of the basis, which is the first result of (i). If $n$ is a unitary in $N$, $(a(n)e_i)_{i\in I}$ is another orthogonal $(b, \nu^o)$-basis of $H_{\psi_1}$, and the independance of $G$ from the basis gives the second and the third formula of (i); let us remark that, for all $n\in N$, $b(n)$ belongs to $A$ and, therefore, commutes with $s$, and that it commutes with $A^\ga$, and, therefore to $r$; moreover, as $\ga(b(n))=1\underset{N}{_b\otimes_\alpha}\beta(n)$, we easily get the first formula linking $G$ with $b(n)$. If we suppose now that $n$ is analytic with respect to $\nu$, we obtain :
\begin{eqnarray*}
G(1\underset{A^\ga}{_s\otimes_r}b(n))(\Lambda_{\psi_1}(x)\underset{\psi_0}{_s\otimes_r}\zeta)
&=&
\sum_i e_i\underset{\nu}{_b\otimes_\alpha}\Lambda_\Phi[(\omega_{b(n)\zeta, e_i}\underset{N}{_b*_\alpha}id)\ga(x)]\\
&=&
\sum_i e_i\underset{\nu}{_b\otimes_\alpha}\Lambda_\Phi[(\omega_{\zeta, e_i}\underset{N}{_b*_\alpha}id)\ga(x)\alpha(\sigma_{-i/2}^\nu(n))]\\
&=&
\sum_i e_i\underset{\nu}{_b\otimes_\alpha}J_\Phi\alpha(n^*)J_\Phi\Lambda_\Phi[(\omega_{\zeta, e_i}\underset{N}{_b*_\alpha}id)\ga(x)]\\
&=&
(1\underset{N}{_b\otimes_\alpha}\hat{\beta}(n))G(\Lambda_{\psi_1}(x)\underset{\psi_0}{_s\otimes_r}\zeta)
\end{eqnarray*}
which, by continuity, finishes the proof of (i). We then obtain :
\begin{eqnarray*}
(1\underset{N}{_b\otimes_\alpha}J_\Phi eJ_\Phi)G(\Lambda_{\psi_1}(x)\underset{\psi_0}{_s\otimes_r}\zeta)
&=&
\sum_i e_i\underset{\nu}{_b\otimes_\alpha}J_\Phi eJ_\Phi \Lambda_\Phi((\omega_{\zeta, e_i}\underset{N}{_b*_\alpha}\ga(x))\\
&=&
\sum_i e_i\underset{\nu}{_b\otimes_\alpha}(\omega_{\zeta, e_i}\underset{N}{_b*_\alpha}id)\ga(x)J_\Phi\Lambda_{\Phi}(e)\\
&=&
\ga(x)(\zeta\underset{\nu}{_b\otimes_\alpha}J_\Phi\Lambda_\Phi(e))
\end{eqnarray*}
which is (ii). We then get :
\[(\omega_{\zeta, \zeta'}\underset{N}{_b*_\alpha}id)\ga(x)J_\Phi\Lambda_\Phi(e)=J_\Phi eJ_\Phi (id*\omega_{\Lambda_{\psi_1}(x), \zeta'})(G)\zeta\]
from which we deduce (iii).
\newline
On the other hand, using again (ii), we get :
\begin{eqnarray*}
(1\underset{N}{_b\otimes_\alpha}J_\Phi eJ_\Phi)\ga(a')G(\Lambda_{\psi_1}(x)\underset{\psi_0}{_s\otimes_r}\zeta)
&=&
\ga(a')(1\underset{N}{_b\otimes_\alpha}J_\Phi eJ_\Phi)G(\Lambda_{\psi_1}(x)\underset{\psi_0}{_s\otimes_r}\zeta)\\
&=&
\ga(a')\ga(x)(\zeta\underset{\nu}{_b\otimes_\alpha}J_\phi\Lambda_\Phi(e))\\
&=&
\ga(a'x)(\zeta\underset{\nu}{_b\otimes_\alpha}J_\phi\Lambda_\Phi(e))\\
&=&
(1\underset{N}{_b\otimes_\alpha}J_\Phi eJ_\Phi)G(\Lambda_{\psi_1}(a'x)\underset{\psi_0}{_s\otimes_r}\zeta)\\
&=&
(1\underset{N}{_b\otimes_\alpha}J_\Phi eJ_\Phi)G(a'\underset{A^\ga}{_s\otimes_r}1)(\Lambda_{\psi_1}(x)\underset{\psi_0}{_s\otimes_r}\zeta)
\end{eqnarray*}
from which we get, by continuity :
\[(1\underset{N}{_b\otimes_\alpha}J_\Phi eJ_\Phi)\ga(a')(G)=(1\underset{N}{_b\otimes_\alpha}J_\Phi eJ_\Phi)G(a'\underset{A^\ga}{_s\otimes_r}1)\]
and, making $e$ going weakly to $1$, we get the first result of (iv). 
\newline
Let $\xi\in D(_\alpha H, \nu)$, $\eta\in D((H)_\beta, \nu^o)\cap \mathcal D(\delta^{-1/2})$, such that $\delta^{-1/2}\eta$ belongs to $D(_\alpha H, \nu)$; using (iii), we get that :
\[(\omega_{\eta, \xi}*id)[(\hat{J}\underset{N}{_\beta\otimes_\alpha}\hat{J})W^*(\hat{J}\underset{N}{_\beta\otimes_{\hat{\alpha}}}\hat{J})](\omega_{\Lambda_{\psi_1}(x), \zeta'}*id)(G)\zeta\]
is equal to :
\[(\omega_{\eta, \xi}*id)[(\hat{J}\underset{N}{_\beta\otimes_\alpha}\hat{J})W^*(\hat{J}\underset{N}{_\beta\otimes_{\hat{\alpha}}}\hat{J})]\Lambda_\Phi[(\omega_{\zeta, \zeta'}\underset{N}{_b*_\alpha}id)\ga(x)]\]
which, using (\cite{E5}, 4.3), is equal to :
\[\Lambda_\Phi[(id\underset{N}{_\beta*_\alpha}\omega_{\delta^{-1/2}\eta, \xi})\Gamma((\omega_{\zeta, \zeta'}\underset{N}{_b*_\alpha}id)\ga(x))]=\Lambda_\Phi[(\omega_{\zeta, \zeta'}\underset{N}{_b*_\alpha}id)\ga((id\underset{N}{_\beta*_\alpha}\omega_{\delta^{-1/2}\eta, \xi})\ga(x)]\]
which, using (iii) again, and (\cite{E5}, 8.4), is equal to :
\[(\omega_{\Lambda_{\psi_1}[(id\underset{N}{_\beta*_\alpha}\omega_{\delta^{-1/2}\eta, \xi})\ga(x)], \zeta'}*id)(G)\zeta
=(\omega_{(id*\omega_{\eta, \xi})(V_\psi)\Lambda_{\psi_1}(x), \zeta'}*id)(G)\zeta\]
from which we get, by continuity :
\[(\omega_{\eta, \xi}*id)[(\hat{J}\underset{N}{_\beta\otimes_\alpha}\hat{J})W^*(\hat{J}\underset{N}{_\beta\otimes_{\hat{\alpha}}}\hat{J})](\omega_{\Lambda_{\psi_1}(x), \zeta'}*id)(G)=(\omega_{(id*\omega_{\eta, \xi})(V_\psi)\Lambda_{\psi_1}(y), \zeta'}*id)(G)\]
which, by continuity and density, remains true for any $\eta\in D(H_\beta, \nu^o)$. 
\newline
Using \ref{pia}, we get, by continuity and density :
\begin{multline*}
(1\underset{N}{_b\otimes_\alpha}(\omega_{\eta, \xi}*id)[(\hat{J}\underset{N}{_\beta\otimes_\alpha}\hat{J})W^*(\hat{J}\underset{N}{_\beta\otimes_{\hat{\alpha}}}\hat{J})])G\\
=G(\pi_\ga(1\underset{N}{_b\otimes_\alpha}(\omega_{\eta, \xi}*id)[(\hat{J}\underset{N}{_\beta\otimes_\alpha}\hat{J})W^*(\hat{J}\underset{N}{_\beta\otimes_{\hat{\alpha}}}\hat{J})])\underset{A^\ga}{_s\otimes_r}1)
\end{multline*}
Using now \ref{MQG}, we get that the weak closure of the linear span of all operators of the form $(\omega_{\eta, \xi}*id)[(\hat{J}\underset{N}{_\beta\otimes_\alpha}\hat{J})W^*(\hat{J}\underset{N}{_\beta\otimes_{\hat{\alpha}}}\hat{J})]$, for all $\xi\in D(_\alpha H, \nu)$, $\eta\in D((H)_\beta, \nu^o)$, is equal to the von Neumann algebra $\widehat{M}'$; therefore, we get, for all $Y\in \widehat{M}'$, that :
\[(1\underset{N}{_b\otimes_\alpha}Y)G=G(\pi_\ga(1\underset{N}{_b\otimes_\alpha}Y)\underset{A^\ga}{_s\otimes_r}1)\]
which finishes the proof of (iv). 
\newline
From (iv), we get that :
\[(1\underset{N}{_b\otimes_\alpha}Y)GG^*=G(\pi_\ga(1\underset{N}{_b\otimes_\alpha}Y)\underset{A^\ga}{_s\otimes_r}1)G^*\]
and that the projection $GG^*$ commutes with $1\underset{N}{_b\otimes_\alpha}\widehat{M}'$; using same arguments, we get that $GG^*$ commutes with $\ga(A)$, and, therefore, it commutes with $A\rtimes_\ga \gG$. So, we get that the application which sends $Z\in A\rtimes_\ga\gG$ on $G^*ZG$ is a $*$-homomorphism, which is equal to $\pi_\ga(Z)$ for any $Z=1\underset{N}{_b\otimes_\alpha}Y$, with $Y\in \widehat{M}'$; using \ref{pia}, we get that the same property holds if $Z=\ga(a')$; therefore, it is true for any $Z\in A\rtimes_\ga\gG$, which is (v). 
\newline
Let us first remark that, because $\psi_1$ is $\delta$-invariant (\ref{thintegrable}), we have (\cite{E6}, 3.2(ii)) : 
\[\Delta_{\tilde{\psi_1}}^{it}=\Delta_{\psi_1}^{it}\underset{N}{_b\otimes_\alpha}(\delta \Delta_{\widehat{\Phi}})^{-it}\]
where this relative tensor product of operators means that it is possible to define a bounded operator with natural values on elementary tensors. With the same definition of relative tensors of operators, we have (\ref{basic}) $\Delta_{\psi_2}^{it}=\Delta_{\psi_1}^{it}\underset{A^\ga}{_s\otimes_r}\Delta_{\psi_1}^{it}$. Using these remaks, we get, for any $x$ in $\gN_{\psi_1}\cap \gN_{T_\ga}$ and $\zeta$, $\zeta'$ in $D((H_{\psi_1})_b, \nu^o)$, that :
\[(\omega_{\Lambda_{\psi_1}(x), \zeta'}*id)(\Delta_{\tilde{\psi_1}}^{it}G\Delta_{\psi_2}^{-it})\zeta\]
is equal to, using (iii) :
\[(\delta\Delta_{\widehat{\Phi}})^{-it}(\omega_{\Delta_{\psi_1}^{-it}\Lambda_{\psi_1}(x), \Delta_{\psi_1}^{-it}\zeta'}*id)(G)\Delta_{\psi_1}^{-it}\zeta
=
(\delta\Delta_{\widehat{\Phi}})^{-it}\Lambda_\Phi((\omega_{\Delta_{\psi_1}^{-it}\zeta, \Delta_{\psi_1}^{-it}\zeta'}\underset{N}{_b*_\alpha}id)\ga(\sigma_{-t}^{\psi_1}(x)))\]
As $\psi_1$ is $\delta$-invariant, we have (\cite{E5}, 88(iii)), for all $t\in\mathbb{R}$ and $x\in A$ :
\[\ga(\sigma_t^{\psi_1}(x))=(\sigma_t^{\psi_1}\underset{N}{_b*_\alpha}\tau_{-t}\sigma_{-t}^{\Phi\circ R}\sigma_t^{\Phi})\ga(x)\]
and, therefore, we have, using (iii) :
\[(\omega_{\Lambda_{\psi_1}(x), \zeta'}*id)(\Delta_{\tilde{\psi_1}}^{it}G\Delta_{\psi_2}^{-it})\zeta
=
(\delta\Delta_{\widehat{\Phi}})^{-it}\Lambda_\Phi(\tau_t\sigma_t^{\Phi\circ R}\sigma_{-t}^{\Phi}[(\omega_{\zeta, \zeta'}\underset{N}{_b*_\alpha}id)\ga(x)])\]
which, using using (\cite{E5}3.8(vii) and (vi)) is equal to :
\begin{multline*}
(\delta\Delta_{\widehat{\Phi}})^{-it}\lambda^{-t/2}P^{it}\Lambda_\Phi(\sigma_t^{\Phi\circ R}\sigma_{-t}^{\Phi}[(\omega_{\zeta, \zeta'}\underset{N}{_b*_\alpha}id)\ga(x)])\\
=
(\delta\Delta_{\widehat{\Phi}})^{-it}\lambda^{-t/2}P^{it}\lambda^{t/2}\delta^{it}J_\Phi\delta^{it}J_\Phi\Lambda_\Phi[(\omega_{\zeta, \zeta'}\underset{N}{_b*_\alpha}id)\ga(x)])
\end{multline*}
which, using  (\cite{E5} 3.10 (vii)) and again (iii), is equal to :
\[\Lambda_\Phi[(\omega_{\zeta, \zeta'}\underset{N}{_b*_\alpha}id)\ga(x)])
=
(\omega_{\Lambda_{\psi_1}(x), \zeta'}*id)(G)\zeta\]
which gives (vi). As (vii) had been proved as well, this finishes the proof. 
\end{proof}

\subsection{Theorem}
\label{G*}
{\it Let $(b, \ga)$ an integrable action of $\gG$ on a von Neumann algebra $A$, $\pi_\ga$ the Galois homomorphism associated by \ref{pia} from the crossed-product $A\rtimes_\ga\gG$ onto the von Neumann algebra $s(A^\ga)'$ obtained by the basic construction made from the inclusion $A^\ga\subset A$; let $\psi_0$ be a normal semi-finite faithful weight on $A^\ga$, and $\psi_1=\psi_0\circ T_\ga$; let us define the representation $a$ of $N$ on $H_{\psi_1}$ by, for $n\in N$ :
\[a(n)=J_{\psi_1}b(n^*)J_{\psi_1}\]
Let us write $r$ for the injection of $A^\ga$ into $A$, and $s$ for the antirepresentation of $A^\ga$ on $H_{\psi_1}$ given, for $a\in A^\ga$, by $s(a)=J_{\psi_1}r(a^*)J_{\psi_1}$, and let $G$ be the isometry from $H_{\psi_1}\underset{\psi_0}{_s\otimes_r}H_{\psi_1}$ into $H_{\psi_1}\underset{\nu}{_b\otimes_\alpha}H$ constructed in \ref{G}; let $\psi_2$ be the weight $\psi_1\circ T_2$ where $T_2$ is the operator-valued weight from $s(A^\ga)'$ onto $A$ obtained by the basic construction (\ref{basic}). Then :
\newline
(i) for any $\xi\in D(_\alpha H, \nu)$, $\eta\in D(H_\beta, \nu^o)$ such that $(\omega_{\xi, \eta}*id)(W^o)$ belongs to $\gN_{\widehat{\Phi}^c}$, and for any $z$ in $\gN_{\psi_1}$, we have :
\[G^*(\Lambda_{\psi_1}(z)\underset{N}{_b\otimes_\alpha}\Lambda_{\widehat{\Phi}^c}[(\omega_{\xi, \eta}*id)(W^o)]
=\Lambda_{\psi_2}[\pi_\ga((1\underset{N}{_b\otimes_\alpha}(\omega_{\xi, \eta}*id)(W^o))\ga(z)]\]
(ii) for any $X\in \gN_{\tilde{\psi_1}}$, $\pi_\ga (X)$ belongs to $\gN_{\psi_2}$, and :
\[G^*\Lambda_{\tilde{\psi_1}}(X)=\Lambda_{\psi_2}(\pi_\ga(X))\]
(iii) $G^*J_{\tilde{\psi_1}}=J_{\psi_2}G^*$
\newline
(iv) the projection $GG^*$ is equal to the support $p$ of $\pi_\ga$; let us consider $\pi_\ga$ as an isomorphism between $(A\rtimes_\ga\gG)_p$ and $s(A^\ga)'$; then, this isomorphism sends the weight $\tilde{\psi_1}_p$ on $\psi_2$.}

\begin{proof}
Let $x$, $y$ be in the Tomita algebra $\mathcal T_{\psi_1, T_1}$; we get that the scalar product 
\[(G^*(\Lambda_{\psi_1}(z)\underset{N}{_b\otimes_\alpha}\Lambda_{\widehat{\Phi}^c}[(\omega_{\xi, \eta}*id)[(\hat{J}\underset{N^o}{_\alpha\otimes_{\hat{\beta}}}\hat{J})W(\hat{J}\underset{N^o}{_\alpha\otimes_\beta}\hat{J})]])|\Lambda_{\psi_1}(x)\underset{\psi_0}{_s\otimes_r}\Lambda_{\psi_1}(y))\]
is equal to :
\[(\Lambda_{\psi_1}(z)\underset{N}{_b\otimes_\alpha}\Lambda_{\widehat{\Phi}^c}[(\omega_{\xi, \eta}*id)[(\hat{J}\underset{N^o}{_\alpha\otimes_{\hat{\beta}}}\hat{J})W(\hat{J}\underset{N^o}{_\alpha\otimes_\beta}\hat{J})]]|G(\Lambda_{\psi_1}(x)\underset{\psi_0}{_s\otimes_r}\Lambda_{\psi_1}(y)))\]
or, to :
\[(\Lambda_{\widehat{\Phi}^c}[(\omega_{\xi, \eta}*id)[(\hat{J}\underset{N^o}{_\alpha\otimes_{\hat{\beta}}}\hat{J})W(\hat{J}\underset{N^o}{_\alpha\otimes_\beta}\hat{J})]]|(\omega_{\Lambda_{\psi_1}(x), \Lambda_{\psi_1}(z)}*id)(G)\Lambda_{\psi_1}(y))\]
which, using \ref{G}(iii), is equal to :
\[(\Lambda_{\widehat{\Phi}^c}[(\omega_{\xi, \eta}*id)[(\hat{J}\underset{N^o}{_\alpha\otimes_{\hat{\beta}}}\hat{J})W(\hat{J}\underset{N^o}{_\alpha\otimes_\beta}\hat{J})]]|\Lambda_\Phi((\omega_{\Lambda_{\psi_1}(y), \Lambda_{\psi_1}(z)}\underset{N}{_b*_\alpha}id)\ga(x))\]
If, moreover, $\eta$ belongs to $\mathcal D(\delta^{-1/2})$, we get, using \ref{MQG}, that it is equal to :
\[(\xi|(\omega_{\Lambda_{\psi_1}(y), \Lambda_{\psi_1}(z)}\underset{N}{_b*_\alpha}id)\ga(x)\delta^{-1/2}\eta)\]
and, if, moreover, $\delta^{-1/2}\eta$ belongs to $D(_\alpha H, \nu)$, this is equal to :
\begin{eqnarray*}
(\Lambda_{\psi_1}(z)|(id\underset{N}{_b*_\alpha}\omega_{\delta^{-1/2}\eta, \xi})\ga(x)\Lambda_{\psi_1}(y)
&=&
\psi_1(y^*(id\underset{N}{_b*_\alpha}\omega_{\xi, \delta^{-1/2}\eta})\ga(x^*)z)\\
&=&
\psi_1((id\underset{N}{_b*_\alpha}\omega_{\xi, \delta^{-1/2}\eta})\ga(x^*)z\sigma_{-i}^{\psi_1}(y^*))\\
&=&
(z\Lambda_{\psi_1}(\sigma_{-i}^{\psi_1}(y^*))|\Lambda_{\psi_1}((id\underset{N}{_b*_\alpha}\omega_{\delta^{-1/2}\eta, \xi})\ga(x))
\end{eqnarray*}
which, using \ref{action} and the standard implementation associated to the weight $\psi_1$, thanks to \ref{thintegrable}, is equal to :
\[(z\Lambda_{\psi_1}(\sigma_{-i}^{\psi_1}(y^*))|(id*\omega_{\eta, \xi})(V_{\psi_1})\Lambda_{\psi_1}(x))\]
and, by continuity, we get that the equality :
\begin{multline*}
(G^*(\Lambda_{\psi_1}(z)\underset{N}{_b\otimes_\alpha}\Lambda_{\widehat{\Phi}^c}[(\omega_{\xi, \eta}*id)[(\hat{J}\underset{N^o}{_\alpha\otimes_{\hat{\beta}}}\hat{J})W(\hat{J}\underset{N^o}{_\alpha\otimes_\beta}\hat{J})]])|\Lambda_{\psi_1}(x)\underset{\psi_0}{_s\otimes_r}\Lambda_{\psi_1}(y))\\
=(z\Lambda_{\psi_1}(\sigma_{-i}^{\psi_1}(y^*))|(id*\omega_{\eta, \xi})(V_{\psi_1})\Lambda_{\psi_1}(x))\\
=((id*\omega_{\xi, \eta})(V_{\psi_1}^*)z\Lambda_{\psi_1}(\sigma_{-i}^{\psi_1}(y^*))|\Lambda_{\psi_1}(x))
\end{multline*}
remains true for the initial hypothesis on $\xi$ and $\eta$. Therefore, we get, using \ref{pia}, that this scalar product is equal to :
\[(\pi_\ga[(1\underset{N}{_b\otimes_\alpha}[(\omega_{\xi, \eta}*id)[(\hat{J}\underset{N^o}{_\alpha\otimes_{\hat{\beta}}}\hat{J})W(\hat{J}\underset{N^o}{_\alpha\otimes_\beta}\hat{J})]]\ga(z)]\Lambda_{\psi_1}(\sigma_{-i}^{\psi_1}(y^*))|\Lambda_{\psi_1}(x))\]
which, using \ref{basic}, is equal to :
\[\Lambda_{\psi_2}(\pi_\ga[(1\underset{N}{_b\otimes_\alpha}[(\omega_{\xi, \eta}*id)[(\hat{J}\underset{N^o}{_\alpha\otimes_{\hat{\beta}}}\hat{J})W(\hat{J}\underset{N^o}{_\alpha\otimes_\beta}\hat{J})]]\ga(z)])|\Lambda_{\psi_1}(x)\underset{\psi_0}{_s\otimes_r}\Lambda_{\psi_1}(y))\]
which, by continuity and density, gives (i). By density, we get, using \ref{crossed}, that, for any $z\in\gN_{\psi_1}$ and $a\in\gN_{\widehat{\Phi}^c}\cap\gN_{\hat{T}^c}$, we have :
\[G^*\Lambda_{\tilde{\psi_1}}((1\underset{N}{_b\otimes_\alpha}a)\ga(z))=G^*(\Lambda_{\psi_1}(z)\underset{N}{_b\otimes_\alpha}\Lambda_{\widehat{\Phi}^c}(a))=
\Lambda_{\psi_2}(\pi_\ga[(1\underset{N}{_b\otimes_\alpha}a)\ga(z)])\]
The linear set generated by elements of the form $(1\underset{N}{_b\otimes_\alpha}a)\ga(z)$, with $a\in\gN_{\widehat{\Phi}^c}\cap\gN_{\hat{T}^c}$ and $z\in\gN_{\psi_1}$, is a core for $\Lambda_{\tilde{\psi_1}}$ (\cite{E5}, 10.8(ii)). So, if $X\in\gN_{\tilde{\psi_1}}$, there exists elements $a_i\in\gN_{\widehat{\Phi}^c}\cap\gN_{\hat{T}^c}$ and $z_i\in\gN_{\psi_1}$ such that the finite sums $\sum_i(1\underset{N}{_b\otimes_\alpha}a_i)\ga(z_i)$ are weakly converging to $X$, and $\sum_i\Lambda_{\tilde{\psi_1}}[(1\underset{N}{_b\otimes_\alpha}a_i)\ga(z_i)]$ is converging to $\Lambda_{\tilde{\psi_1}}(X)$. But then, on one hand, $\pi_\ga(\sum_i(1\underset{N}{_b\otimes_\alpha}a_i)\ga(z_i))$ is converging to $\pi_\ga(X)$, and, on the other hand,  the finite sums $\sum_i\Lambda_{\psi_2}(\pi_\ga[(1\underset{N}{_b\otimes_\alpha}a_i)\ga(z_i)])=\sum_i G^*(\Lambda_{\tilde{\psi_1}}[(1\underset{N}{_b\otimes_\alpha}a_i)\ga(z_i)]$ are converging. So, applying the closed graph theorem to the closed application $\Lambda_{\psi_2}$, we get (ii). 
If now $X\in \gN_{\tilde{\psi_1}}\cap\gN_{\tilde{\psi_1}}^*$, we get that $\pi_\ga(X)$ belongs to $\gN_{\psi_2}\cap\gN_{\psi_2}^*$, and that $G^*S_{\tilde{\psi_1}}\Lambda_{\tilde{\psi_1}}(X)=S_{\psi_2}G^*\Lambda_{\tilde{\psi_1}}(X)$. So, we have $G^*S_{\tilde{\psi_1}}\subset S_{\psi_2}G^*$; using now \ref{G}(vii), we get (iii). 
\newline
Using (iii), we get that $J_{\tilde{\psi_1}}GG^*J_{\tilde{\psi_1}}=GG^*$; as, in \ref{G}(v), we had obtained that $GG^*$ belongs to $(A\rtimes_\ga\gG)'$, we get that $GG^*\in Z(A\rtimes_\ga\gG)$. Using then \ref{G}(v) again, we see that, for any $X\in A\rtimes_\ga\gG$, we have $XGG^*=G\pi_\ga(X)G^*$, and, therefore, that $\pi_\ga(X)=0$ if and only if $XGG^*=0$, from which we get that $GG^*$ is equal to the support of $\pi_\ga$. Using then (ii), we finish the proof. 
\end{proof}

\subsection{Lemma}
\label{lemG*}
{\it With the hypothesis and notations of \ref{G*}, we get, for any $m\in\widehat{M}'$ :}
\[1\underset{A^\ga}{_s\otimes_r}\pi_\ga(1\underset{N}{_b\otimes_\alpha}m)=G^*V_{\psi_1}[1\underset{N}{_a\otimes_\beta}\hat{J}\hat{R}^c(m^*)\hat{J}]V_{\psi_1}^*G\]

\begin{proof}
Using \ref{lemR}, \ref{G*}(iii), \ref{G}(iv) and \ref{crossed}, we get :
\begin{eqnarray*}
G[1\underset{A^\ga}{_s\otimes_r}\pi_\ga(1\underset{N}{_b\otimes_\alpha}m)]
&=&
G[1\underset{A^\ga}{_s\otimes_r}J_{\psi_1}\pi_\ga(1\underset{N}{_b\otimes_\alpha}\hat{R}^c(m^*))J_{\psi_1}]\\
&=&
GJ_{\psi_2}[\pi_\ga(1\underset{N}{_b\otimes_\alpha}\hat{R}^c(m^*))\underset{N}{_s\otimes_r}1]J_{\psi_2}\\
&=&
J_{\tilde{\psi_1}}G[\pi_\ga(1\underset{N}{_b\otimes_\alpha}\hat{R}^c(m^*))\underset{N}{_s\otimes_r}1]J_{\psi_2}\\
&=&
J_{\tilde{\psi_1}}[1\underset{N}{_b\otimes_\alpha}\hat{R}^c(m^*)]GJ_{\psi_2}\\
&=&
J_{\tilde{\psi_1}}[1\underset{N}{_b\otimes_\alpha}\hat{R}^c(m^*)]J_{\tilde{\psi_1}}G\\
&=&
V_{\psi_1}[1\underset{N^o}{_a\otimes_\beta}\hat{J}\hat{R}^c(m^*)\hat{J}]V_{\psi_1}^*G
\end{eqnarray*}
from which, $G$ being an isometry, we get the result. \end{proof}

\subsection{Definitions}
\label{defGalois}
Let $(b, \ga)$ be an integrable action of a measured quantum groupoid $\gG$ on a von Neumann algebra $A$, $A\rtimes_\ga\gG$ be the crossed product, $\pi_\ga$ be the Galois homomorphism from $A\rtimes_\ga\gG$ onto the algebra $s(A^\ga)'$ obtained by the standard construction made from the inclusion $A^\ga\subset A$ (\ref{pia}), and $G$ be the isometry constructed in \ref{G}; then, using \ref{G*}(iv), we get that the following properties are equivalent :
\newline
(i) $\pi_\ga$ is an isomorphism between $A\rtimes_\ga\gG$ and $s(A^\ga)'$;
\newline
(ii) the isometry $G$ is a unitary; 
\newline
(iii) the inclusion $A^\ga\underset{N}{_b\otimes_\alpha}1_H\subset \ga(A)\subset A\rtimes_\ga\gG$ is standard, and the operator-valued weight $T_{\tilde{\ga}}$ is obtained from $T_\ga$ by this standard construction. 
\newline
In that situation, following \cite{DC1}, we shall say that $(b, \ga)$ is a Galois action of $\gG$, and that the $A^\ga-N^o$-bimodule $A$ will be called a Galois bimodule for $\gG$, and the unitary $\tilde{G}=\sigma_\nu G$ from $H_{\psi_1}\underset{\psi_0}{_s\otimes_r}H_{\psi_1}$ onto  $H\underset{\nu^o}{_\alpha\otimes_b}H_{\psi_1}$ will be called its Galois unitary. Then, it is clear that the representation $\mu$ of $\widehat{M}'$ on $H_{\psi_1}$, defined in \ref{pia},  is faithful. 
\newline
Moreover, a normal semi-finite faithful weight $\psi_0$ on $A^\ga$ will be said $\ga$-relatively invariant, if there exists a normal semi-finite faithful weight $\phi$ on $A$, invariant by $\ga$, and bearing the density property, such that the two automorphism groups $\sigma^\phi$ and $\sigma^{\psi_1}$ on $A$ commute (where $\psi_1=\psi_0\circ T_\ga$). In that situation, we shall say that the 5-uple $(A, b, \ga, \phi, \psi_0)$ is a Galois system for $\gG$. Then, thanks to \cite{V1}, we know that there exists a positive operator $\delta_A$ affiliated to $A$, and a positive operator $\lambda_A$ affiliated to $Z(A)$ such that :
\[(D\phi :D\psi_1)_t=\lambda_A^{it^2/2}\delta_A^{it}\]
We shall call $\delta_A$ the modulus of the action $(b,\ga)$, and $\lambda_A$ the scaling operator of the action. 
\newline
Starting from a left-action $(a, \gb)$, we get the notion of left Galois system and that $(A, a, \gb, \phi, \psi_0)$ is left Galois if and only if $(A^o, a^o, (\sigma_N\gb)^o, \phi^o, \psi_0^o)$ is Galois, or, conversely, that $(A, b, \ga, \phi, \psi_0)$ is Galois if and only if $(A^o, b^o, \sigma_N\ga^o, \phi^o, \psi_0^o)$ is left-Galois.

\subsection{Examples}
\label{exGalois}
(i) Let $(b, \ga)$ be any action of $\gG$ on a von Neumann algebra $A$; then (\ref{crossed}), there exists an action $(1\underset{N}{_b\otimes_\alpha}\hat{\alpha}, \tilde{\ga})$ of $\widehat{\gG}^c$ on the crossed product $A\rtimes_\ga\gG$. This action is integrable (\cite{E5}, 9.8); we have $(A\rtimes_\ga\gG)^{\tilde{\ga}}=\ga(A)$ (\cite{E5}, 11.5), and, as the inclusion $\ga(A)\subset A\rtimes_\ga\gG\subset A\underset{N}{_b*_\alpha}\mathcal L(H)$ is depth 2 (\cite{E5}, 13.8), we obtain, by (\cite{E5}, 13.9) that the dual action $(1\underset{N}{_b\otimes_\alpha}\hat{\alpha}, \tilde{\ga})$ is a Galois action of $\widehat{\gG}^c$, with $\alpha(A)\subset A\rtimes_\ga\gG$ and $1\underset{N}{_b\otimes_\alpha}\hat{\alpha}(N)\subset A\rtimes_\ga\gG$ as Galois bi-module. 
\newline
(ii) in particular (\cite{E5}, 9.5) we get that $(\beta, \Gamma)$ is a Galois action of $\gG$, with $\alpha(N)\subset M$ and $\beta(N)\subset M$ as Galois bi-module. Moreover, we get that $(M, \beta, \Gamma, \Phi\circ R, \nu)$ is a Galois system for $\gG$. Then, we can easily check that $M^\Gamma=\alpha(N)$, that the operator-valued weight $T_\Gamma$ is equal to the left-invariant weight $T_L$, and, therefore, that $\psi_1=\Phi$, $a=\hat{\alpha}$, $r=\alpha$, $s=\hat{\beta}$, $V_{\psi_1}=\sigma W^{o*}\sigma$, $\pi_\Gamma=id$, and $\tilde{G}=\widehat{W}$. 
\newline
(iii) if $(b, \ga)$ is an integrable outer action of $\gG$ on $A$, then, $(b, \ga)$ is Galois : let $G$ be the isometry constructed in \ref{G}; as, by definition (\cite{E7}), we have $A\rtimes_\ga\gG\cap\ga(A)'=1\underset{N}{_b\otimes_\alpha}\hat{\alpha}(N)$, we get, by \ref{G*}, that there exists a projection $p\in Z(N)$ such that $GG^*=1\underset{N}{_b\otimes_\alpha}\hat{\alpha}(p)$ is the support of $\pi_\ga$; using \ref{G}, we get that $p=1$, which gives the injectivity of $\pi_\ga$. 

\subsection{Lemma}
\label{lemG}
{\it (i) Let $(b, \ga)$ be a Galois action of the measured quantum groupoid $\gG$ on the von Neumann algebra $A$; let $\psi_1=\psi_0\circ T_\ga$, $V_{\psi_1}$ the standard implementation of $\ga$ associated to $\psi_1$(\ref{crossed}), and $\tilde{G}$ the Galois unitary of the Galois system. Then, we have :
\[\tilde{G}\sigma_{\psi_0^o}(J_{\psi_1}\underset{A^\ga}{_s\otimes_r}J_{\psi_1})=\sigma_\nu V_{\psi_1}\sigma_{\nu^o}(\hat{J}\underset{N^o}{_\alpha\otimes_b}J_{\psi_1})\tilde{G}\]
(ii) Let $W$ be the pseudo-multiplicative unitary of $\gG$, $W^o$ be the pseudo-multiplicative unitary of $\gG^o$, $\widehat{W}$ be the pseudo-multiplicative unitary of $\widehat{\gG}$. We have :}
\[\widehat{W}\sigma_{\nu^o}(J\underset{N}{_{\hat{\beta}}\otimes_\alpha}J)=W^{o *}(\hat{J}\underset{N^o}{_\alpha\otimes_\beta}J)\widehat{W}\]
\begin{proof}
Using \ref{G*}(iii) and \ref{crossed}, we have :
\begin{eqnarray*}
\sigma_\nu \tilde{G}\sigma_{\psi_0^o}(J_{\psi_1}\underset{A^\ga}{_s\otimes_r}J_{\psi_1})
&=&
GJ_{\psi_2}\\
&=&
J_{\tilde{\psi_1}}G\\
&=&
V_{\psi_1}(J_{\psi_1}\underset{N}{_b\otimes_\alpha}\hat{J})G\\
&=&
V_{\psi_1}\sigma_{\nu^o}(\hat{J}\underset{N^o}{_\alpha\otimes_b}J_{\psi_1})\tilde{G}
\end{eqnarray*}
from which we get (i). Using \ref{exGalois}(ii), we obtain (ii). \end{proof}

\section{From Galois actions to Galois systems and back}
\label{fromGtoG}
In this chapter, we suppose that we have a Galois action $(b, \ga)$ of a measured quantum groupoid $\gG$ on a von Neumann algebra $A$, and a normal semi-finite faithful weight $\psi_0$ on $A^\ga$, such that the subspace $D((H_{\psi_1})_b, \nu^o)\cap D(_rH_{\psi_1}, \psi_0)$ is dense in $H_{\psi_1}$, where $\psi_1=\psi_0\circ T_\ga$. We then prove that right leg of the Galois unitary introduced in \ref{defGalois} generates $A$, and that this unitary satisfies a pentagonal relation (\ref{tildeG}). This allows us to prove, in some particular cases (\ref{inner}, \ref{corsum}) that there exits then a normal semi-finite faithful weight $\phi$ on $A$ such that $(A, b, \ga, \phi, \psi_0)$ is a Galois system for $\gG$. Conversely, if there exists a Galois system  $(A, b, \ga, \phi, \psi_0)$ for $\gG$, then the weight $\psi_0$ satisfies this density property (\ref{dens}). 

\subsection{Definition}
\label{defdens}
Let $(b, \ga)$ be a Galois action of the quantum groupoid $\gG$ on a von Neumann algebra $A$; let $\psi_0$ a normal semi-finite faithful weight on $A^\ga$; let us write $\psi_1=\psi_0\circ T_\ga$, $r$ for the injection of $A^\ga$ into $A$. We shall say that the weight $\psi_0$ bear the Galois density property if the subspace $D((H_{\psi_1})_b, \nu^o)\cap D(_rH_{\psi_1}, \psi_0)$ is dense in $H_{\psi_1}$. 

\subsection{Theorem}
\label{tildeG}
{\it Let $(b, \ga)$ be a Galois action of the measured quantum groupoid $\gG$ on a von Neumann algebra $A$; let $\psi_0$ a normal semi-finite faithful weight on $A^\ga$, bearing the Galois density property, in the sense of \ref{defdens}; ,let $\tilde{G}$ be the Galois unitary of $(b, \ga)$, from $H_{\psi_1}\underset{\psi_0}{_s\otimes_r}H_{\psi_1}$ onto  $H\underset{\nu^o}{_\alpha\otimes_b}H_{\psi_1}$, as defined in \ref{defGalois}. We have :
\newline
(i) for any $x$ in $\gN_{\psi_1}\cap\gN_{T_\ga}$, $\zeta\in D((H_{\psi_1})_b, \nu^o)\cap D(_r H_{\psi_1}, \psi_0)$ , $\zeta'$ in $D((H_{\psi_1})_b, \nu^o)$, $\eta\in H$, we have :
\[(\tilde{G}(\Lambda_{\psi_1}(x)\underset{\psi_0}{_s\otimes_r}\zeta)|\eta\underset{\nu^o}{_\alpha\otimes_b}\zeta')=(\Lambda_\Phi[(\omega_{\zeta, \zeta'}\underset{N}{_b*_\alpha}id)\ga(x)]|\eta)\]
and, therefore :
\[(id*\omega_{\zeta, \zeta'})(\tilde{G})\Lambda_{\psi_1}(x)=\Lambda_\Phi[(\omega_{\zeta, \zeta'}\underset{N}{_b*_\alpha}id)\ga(x)]\]
(ii) for any $x\in\gN_{\psi_1}\cap\gN_{T_\ga}$, $y\in\gN_\Phi\cap \gN_T$, $\xi\in D(_\alpha H, \nu)$, we have :
\[(\omega_{\Lambda_{\psi_1}(x), J_\Phi y^*J_\Phi\xi}*id)(\tilde{G})=(id\underset{N}{_b*_\alpha}\omega_{J_\Phi\Lambda_\Phi(y), \xi})\ga(x)\]
(iii), for any $x\in\gN_{\psi_1}\cap\gN_{\psi_1}^*$, $y$, $z$ in $\gN_\Phi\cap\gN_T$, we have :
\[(\omega_{\Lambda_{\psi_1}(x), J_\Phi\Lambda_\Phi(y^*z)}*id)(\tilde{G})^*=(\omega_{\Lambda_{\psi_1}(x^*), J_\Phi\Lambda_\Phi(z^*y)}*id)(\tilde{G})\]
(iv) the two unitaries $(1\underset{N^o}{_\alpha\otimes_b}\tilde{G})(\tilde{G}\underset{A^\ga}{_s\otimes_r}1)$ and $(\widehat{W}\underset{N^o}{_\alpha\otimes_b}1)\sigma^{2,3}_{\alpha, \hat{\beta}}(\tilde{G}\underset{N}{_b\otimes_\alpha}1)(1\underset{A^\ga}{_s\otimes_r}\sigma_{\nu^o})(1\underset{A^\ga}{_s\otimes_r}\tilde{G})$, from $H_{\psi_1}\underset{\psi_0}{_s\otimes_r}H_{\psi_1}\underset{\psi_0}{_s\otimes_r}H_{\psi_1}$ to $H\underset{\nu^o}{_\alpha\otimes_\beta}H\underset{\nu^o}{_\alpha\otimes_b}H_{\psi_1}$, are equal.}

\begin{proof}
Using \ref{G}(iii) and the definition of $\tilde{G}$, we get (i) by a direct calculation. 
Using (i), we then get :
\begin{eqnarray*}
((\omega_{\Lambda_{\psi_1}(x), J_\Phi y^*J_\Phi\xi}*id)(\tilde{G})\zeta|\zeta')
&=&
(\tilde{G}(\Lambda_{\psi_1}(x)\underset{\psi_0}{_s\otimes_r}\zeta)|J_\Phi y^*J_\Phi\xi\underset{\nu^o}{_\alpha\otimes_b}\zeta')\\
&=&
(\Lambda_\Phi(\omega_{\zeta, \zeta'}\underset{N}{_b*_\alpha}id)\ga(x))|J_\Phi y^*J_\Phi\xi)\\
&=&
(J_\Phi yJ_\Phi\Lambda_\Phi(\omega_{\zeta, \zeta'}\underset{N}{_b*_\alpha}id)\ga(x))|\xi)\\
&=&
((\omega_{\zeta, \zeta'}\underset{N}{_b*_\alpha}id)\ga(x)J_\Phi\Lambda_\Phi(y)|\xi)\\
&=&
((id\underset{N}{_b*_\alpha}\omega_{J_\Phi\Lambda_\Phi(y), \xi})\ga(x)\zeta|\zeta')
\end{eqnarray*}
from which we get (ii). Using (ii), we get :
\begin{eqnarray*}
(\omega_{\Lambda_{\psi_1}(x), J_\Phi\Lambda_\Phi(y^*z)}*id)(\tilde{G})^*
&=&
(id\underset{N}{_b*_\alpha}\omega_{J_\Phi\Lambda_\Phi(y), J_\Phi\Lambda_{\Phi}(z)})\ga(x)^*\\
&=&
(id\underset{N}{_b*_\alpha}\omega_{J_\Phi\Lambda_\Phi(z), J_\Phi\Lambda_{\Phi}(y)})\ga(x^*)\\
&=&
(\omega_{\Lambda_{\psi_1}(x^*), J_\Phi\Lambda_\Phi(z^*y)}*id)(\tilde{G})
\end{eqnarray*}
from which we get (iii). 
\newline
Let $v\in D(H_\beta, \nu^o)$, $w\in D(_\alpha H,\nu)\cap D(H_\beta, \nu^o)$, $\zeta\in D((H_{\psi_1})_b, \nu^o)\cap D(_rH_{\psi_1}, \psi_0)$, $x\in \gN_{\psi_1}\cap \gN_{T_{\ga}}$; we have, using (i) and (\cite{E5}, 3.10 (i)) :
\begin{eqnarray*}
(i*\omega_{v,w})(\widehat{W})(id*\omega_{\zeta, \zeta'})(\tilde{G})\Lambda_{\psi_1}(x)
&=&
(i*\omega_{v,w})(\widehat{W})\Lambda_\Phi[(\omega_{\zeta, \zeta'}\underset{N}{_b*_\alpha}id)\ga(x)]\\
&=&
\Lambda_\Phi[(\omega_{v, w}\underset{N}{_\beta*_\alpha}id)\Gamma[(\omega_{\zeta, \zeta'}\underset{N}{_b*_\alpha}id)\ga(x)]]\\
&=&
\Lambda_\Phi[(\omega_{\zeta, \zeta'}\underset{N}{_b*_\alpha}\omega_{v, w}\underset{N}{_\beta*_\alpha}id)(\ga\underset{N}{_b*_\alpha}id)\ga(x)]\\
&=&
\Lambda_\Phi[((\omega_{v, w}\underset{N^o}{_\alpha*_b}\omega_{\zeta, \zeta'})\circ \varsigma_N\ga)\ga(x)]
\end{eqnarray*}
Using \ref{G}(iv), we get that, for all $y\in A$ :
\[(\omega_{v, w}\underset{N^o}{_\alpha*_b}\omega_{\zeta, \zeta'})\circ \varsigma_N\ga)(y)=
(\tilde{G}(y\underset{A^\ga}{_s\otimes_r}1)\tilde{G}^*(v\underset{\nu^o}{_\alpha\otimes_b}\zeta)|w\underset{\nu^o}{_\alpha\otimes_b}\zeta')\]
Let $(e_i)_{i\in I}$ be an orthogonal $(r, \psi_0)$-basis for $H_{\psi_1}$; there exists $(v_i)_{i\in I}$ and $(w_i)_{i\in I}$ in $H_{\psi_1}$, such that :
\[\tilde{G}^*(v\underset{\nu^o}{_\alpha\otimes_b}\zeta)=\sum_i v_i\underset{\psi_0}{_s\otimes_r}e_i\]
\[\tilde{G}^*(w\underset{\nu^o}{_\alpha\otimes_b}\zeta')=\sum_i w_i\underset{\psi_0}{_s\otimes_r}e_i\]
Using the intertwining properties given in \ref{G}(i), we get, for all $n\in N$ :
\begin{multline*}
\sum_i\|b(n)v_i\|^2=\|\sum_i b(n)v_i\underset{\psi_0}{_s\otimes_r}e_i\|^2=\|(b(n)\underset{A^\ga}{_s\otimes_r}1)\tilde{G}^*(v\underset{\nu^o}{_\alpha\otimes_b}\zeta)\|^2\\
=\|\tilde{G}^*(\beta(n)v\underset{\nu^o}{_\alpha\otimes_b}\zeta)\|^2=\|\beta(n)v\underset{\nu^o}{_\alpha\otimes_b}\zeta\|^2
\end{multline*}
and, therefore, as $\zeta$ is in $D((H_{\psi_1})_b, \nu^o)$ and $v$ is in $D(H_\beta, \nu^o)$, we get that each $v_i$ is in $D((H_{\psi_1})_b, \nu^o)$, and the same result holds for the $w_i$s.
\newline
On the other hand, for any $z\in A^\ga$, we get :
\begin{multline*}
\sum_i\|r(z)v_i\|^2=\|\sum_i r(z)v_i\underset{\psi_0}{_s\otimes_r}e_i\|^2=
\|(r(z)\underset{A^\ga}{_s\otimes_r}1)\tilde{G}^*(v\underset{\nu^o}{_\alpha\otimes_b}\zeta)\|^2\\
=\|\tilde{G}^*(v\underset{\nu^o}{_\alpha\otimes_b}r(z)\zeta)\|^2=\|v\underset{\nu^o}{_\alpha\otimes_b}r(z)\zeta\|^2
\end{multline*}
and, therefore, as $v$ is in $D(_\alpha H, \nu)$ and $\zeta$ is in $D(_rH_{\psi_1}, \psi_0)$, we get that each $v_i$ is in $D(_rH_{\psi_1}, \psi_0)$. 
\newline
So, we get that there exists $v_i$ in $D((H_{\psi_1})_b, \nu^o)\cap D(_rH_{\psi_1}, \psi_0)$, and $w_i$ in $D((H_{\psi_1})_b, \nu^o)$, such that, for all $y\in A$, we have :
\[(\omega_{v, w}\underset{N^o}{_\alpha*_b}\omega_{\zeta, \zeta'})\circ \varsigma_N\ga(y)=
\sum_i(yv_i|w_i)\]
and, therefore $(\omega_{v, w}\underset{N^o}{_\alpha*_b}\omega_{\zeta, \zeta'})\circ \varsigma_N\ga=\sum_i\omega_{v_i, w_i}$. So, we get :
\begin{eqnarray*}
(i*\omega_{v,w})(\widehat{W})(id*\omega_{\zeta, \zeta'})(\tilde{G})\Lambda_{\psi_1}(x)
&=&
\Lambda_\Phi[((\omega_{v, w}\underset{N^o}{_\alpha*_b}\omega_{\zeta, \zeta'})\circ \varsigma_N\ga)\ga(x)]\\
&=&
\sum_i\Lambda_\Phi [(\omega_{v_i, w_i}\underset{N}{_b*_\alpha}id)\ga(x)]\\
&=&
\sum_i (i*\omega_{v_i,w_i})(\tilde{G})\Lambda_{\psi_1}(x)
\end{eqnarray*}
and, therefore, for all $\eta\in H$, we have :
\begin{eqnarray*}
((i*\omega_{v,w})(\widehat{W})(id*\omega_{\zeta, \zeta'})(\tilde{G})\Lambda_{\psi_1}(x)|\eta)
&=&
\sum_i ((i*\omega_{v_i,w_i})(\tilde{G})\Lambda_{\psi_1}(x)|\eta)\\
&=&
\sum_i(\tilde{G}(\Lambda_{\psi_1}(x)\underset{\psi_0}{_s\otimes_r}v_i)|\eta\underset{\nu^o}{_\alpha\otimes_b}w_i)
\end{eqnarray*}
which is equal to :
\begin{multline*}
\sum_i((\tilde{G}\underset{A^\ga}{_s\otimes_r}1)(\Lambda_{\psi_1}(x)\underset{\psi_0}{_s\otimes_r}v_i\underset{\psi_o}{_s\otimes_r}e_i)|\eta\underset{\nu^o}{_\alpha\otimes_b}w_i\underset{\psi_0}{_s\otimes_r}e_i)\\
=
((\tilde{G}\underset{A^\ga}{_s\otimes_r}1)(1\underset{A^\ga}{_s\otimes_r}\tilde{G}^*)(\Lambda_{\psi_1}(x)\underset{\psi_0}{_s\otimes_r}(v\underset{\nu^o}{_\alpha\otimes_b}\zeta))|(1\underset{N^o}{_\alpha\otimes_b}\tilde{G}^*)(\eta\underset{\nu^o}{_\alpha\otimes_{\hat{\beta}}}w\underset{\nu^o}{_\alpha\otimes_b}\zeta'))
\end{multline*}
On the other hand, we get that :
\begin{eqnarray*}
((i*\omega_{v,w})(\widehat{W})(id*\omega_{\zeta, \zeta'})(\tilde{G})\Lambda_{\psi_1}(x)|\eta)
&=&
((id*\omega_{\zeta, \zeta'})(\tilde{G})\Lambda_{\psi_1}(x)|(i*\omega_{v,w})(\widehat{W})^*\eta)\\
&=&
(\tilde{G}(\Lambda_{\psi_1}(x)\underset{\psi_0}{_s\otimes_r}\zeta)|(i*\omega_{v,w})(\widehat{W})^*\eta\underset{\nu^o}{_\alpha\otimes_b}\zeta')
\end{eqnarray*}
is, using again \ref{G}(i), equal to :
\[(\sigma^{2,3}_{\alpha, \hat{\beta}}(\tilde{G}\underset{N}{_b\otimes_\alpha}1)(1\underset{A^\ga}{_s\otimes_r}\sigma_{\nu^o})[\Lambda_{\psi_1}(x)\underset{\psi_0}{_s\otimes_r}(v\underset{\nu^o}{_\alpha\otimes_b}\zeta)]|(\widehat{W}^*\underset{N^o}{_\alpha\otimes_b}1)(\eta\underset{\nu^o}{_\alpha\otimes_{\hat{\beta}}}w\underset{\nu^o}{_\alpha\otimes_b}\zeta')\]
where $\sigma^{2,3}_{\alpha, \hat{\beta}}$ is the flip from $(H\underset{N^o}{_\alpha\otimes_b}H_{\psi_1})\underset{N}{_{\hat{\beta}}\otimes_\alpha}H$ onto $(H\underset{N}{_{\hat{\beta}}\otimes_\alpha}H)\underset{N^o}{_\alpha\otimes_b}H_{\psi_1}$, exchanging the second and the third leg, and $\alpha$ with $\hat{\beta}$. 
\end{proof}
From which we get that :
\[((1\underset{N^o}{_\alpha\otimes_b}\tilde{G})(\tilde{G}\underset{A^\ga}{_s\otimes_r}1)(1\underset{A^\ga}{_s\otimes_r}\tilde{G}^*)([\Lambda_{\psi_1}(x)\underset{\psi_0}{_s\otimes_r}(v\underset{\nu^o}{_\alpha\otimes_b}\zeta)]|\eta\underset{\nu^o}{_\alpha\otimes_{\hat{\beta}}}w\underset{\nu^o}{_\alpha\otimes_b}\zeta')\]
is equal to :
\[(\widehat{W}\underset{N^o}{_\alpha\otimes_b}1)\sigma^{2,3}_{\alpha, \hat{\beta}}(\tilde{G}\underset{N}{_b\otimes_\alpha}1)(1\underset{A^\ga}{_s\otimes_r}\sigma_{\nu^o})[\Lambda_{\psi_1}(x)\underset{\psi_0}{_s\otimes_r}(v\underset{\nu^o}{_\alpha\otimes_b}\zeta)]|\eta\underset{\nu^o}{_\alpha\otimes_{\hat{\beta}}}w\underset{\nu^o}{_\alpha\otimes_b}\zeta')\]
and, therefore, that :
\[(\widehat{W}\underset{N^o}{_\alpha\otimes_b}1)\sigma^{2,3}_{\alpha, \hat{\beta}}(\tilde{G}\underset{N}{_b\otimes_\alpha}1)(1\underset{A^\ga}{_s\otimes_r}\sigma_{\nu^o})(1\underset{A^\ga}{_s\otimes_r}\tilde{G})=(1\underset{N^o}{_\alpha\otimes_b}\tilde{G})(\tilde{G}\underset{A^\ga}{_s\otimes_r}1)\]

\subsection{Corollary}
\label{cortildeG}
{\it Let $(b, \ga)$ be a Galois action of the measured quantum groupoid $\gG$ on a von Neumann algebra $A$; let $\psi_0$ a normal semi-finite faithful weight on $A^\ga$, bearing the Galois density property defined in \ref{defdens}; let us write $\psi_1=\psi_0\circ T_\ga$, and let us write $r$ for the injection of $A^\ga$ in $A$, and $s$ for the anti-representation $s(x)=J_{\psi_1}r(x^*)J_{\psi_1}$ of $A^\ga$ on $H_{\psi_1}$; let $\tilde{G}$ be the Galois unitary of $(b, \ga)$, from $H_{\psi_1}\underset{\psi_0}{_s\otimes_r}H_{\psi_1}$ onto  $H\underset{\nu^o}{_\alpha\otimes_b}H_{\psi_1}$, as defined in \ref{defGalois}. Then, the linear space generated by the elements of the form $(\omega_{\zeta, \zeta'}*id)(\tilde{G})$, for all $\zeta$ in $D((H_{\psi_1})_s, \psi_0^o)$, and $\zeta'\in D(_\alpha H, \nu)$ is weakly dense in $A$. }

\begin{proof}
Let us first look at the product of two elements of that form. Let $\zeta_1\in D((H_{\psi_1})_s, \psi_0^o)$, and $\zeta'_1\in D(_\alpha H, \nu)$; let $\xi$, $\eta$ in $H_{\psi_1}$. Then, we have :
\begin{multline*}
((\omega_{\zeta, \zeta'}*id)(\tilde{G})(\omega_{\zeta_1, \zeta_1'}*id)(\tilde{G})\xi|\eta)\\
=([\sigma^{2,3}_{\alpha, \hat{\beta}}(\tilde{G}\underset{N}{_b\otimes_\alpha}1)(1\underset{A^\ga}{_s\otimes_r}\sigma_{\nu^o})(1\underset{A^\ga}{_s\otimes_r}\tilde{G})](\zeta\underset{\psi_0}{_s\otimes_r}\zeta_1\underset{\psi_0}{_s\otimes_r}\xi)|(\zeta'\underset{\nu}{_{\hat{\beta}}\otimes_\alpha}\zeta_1')\underset{\psi_0^o}{_\alpha\otimes_b}\eta)
\end{multline*}
which, using \ref{tildeG}(iv), is equal to :
\[((1\underset{N^o}{_\alpha\otimes_b}\tilde{G})(\tilde{G}\underset{A^\ga}{_s\otimes_r}1)(\zeta\underset{\psi_0}{_s\otimes_r}\zeta_1\underset{\psi_0}{_s\otimes_r}\xi)|(\widehat{W}\underset{(A^\ga)^o}{_\alpha\otimes_b}1)[(\zeta'\underset{\nu}{_{\hat{\beta}}\otimes_\alpha}\zeta_1')\underset{\psi_0^o}{_\alpha\otimes_b}\eta])\]
Let $(e_i)_{i\in I}$ be an orthogonal $(\alpha, \nu)$-basis. As in (\cite{E3}, 3.4), we can prove that there exist $(\zeta_i)_{i\in I}\in D((H_{\psi_1})_s, \psi_0^o)$ and $(\zeta'_i)_{i\in I}\in D(_\alpha H, \nu)$ such that :
\[\tilde{G}(\zeta\underset{\psi_0}{_s\otimes_r}\zeta_1)=\sum_i e_i\underset{\nu^o}{_\alpha\otimes_b}\zeta_1\]
\[\widehat{W}(\zeta'\underset{\nu}{_{\hat{\beta}}\otimes_\alpha}\zeta'_1)=\sum_i e_i\underset{\nu^o}{_\alpha\otimes_\beta}\zeta'_i\]
and, therefore, we get that :
\[((\omega_{\zeta, \zeta'}*id)(\tilde{G})(\omega_{\zeta_1, \zeta_1'}*id)(\tilde{G})\xi|\eta)
=\sum_i(\omega_{\zeta_i, \zeta'_i}*id)(\tilde{G})\xi|\eta)\]
which proves that the product $(\omega_{\zeta, \zeta'}*id)(\tilde{G})(\omega_{\zeta_1, \zeta_1'}*id)(\tilde{G})$ is the weak limit of the finite sums $(\omega_{\zeta_i, \zeta'_i}*id)(\tilde{G})$. So, by continuity, we get that the weak closure of the linear space generated by the elements of the form $(\omega_{\zeta, \zeta'}*id)(\tilde{G})$, for all $\zeta$ in $D((H_{\psi_1})_s, \psi_0^o)$, and $\zeta'\in D(_\alpha H, \nu)$ is an algebra. 
\newline
Using \ref{tildeG}(ii), we get, on one hand, that all the operators of the form $(\omega_{\zeta_2, \zeta'}*id)(\tilde{G})$ (with $\zeta_2\in D((H_{\psi_1})_s, \psi_0^o)$) belong to $A$, and, on the other hand, that the closure of the linear set generated by these operators is the closure of the set of all operators $(id\underset{N}{_b*_\alpha}\omega_{\zeta'', \zeta'})\ga(x)$, for all $\zeta', \zeta''$ in $D(_\alpha H, \nu)$, and $x\in A$, and is therefore invariant by taking the adjoint, and that it contains all operators $b(<\zeta'', \zeta''>_{b, \nu^o})=(id\underset{N}{_b*_\alpha}\omega_{\zeta'', \zeta'})\ga(1)$. Therefore, it is a sub-von Neumann algebra $B$ of $A$ which contains $b(N)$. If now $X\in B'$, we get that $X\underset{N}{_b\otimes_\alpha}1$ belongs to $\ga(A)'\cap b(N)'\underset{N}{_b\otimes_\alpha}1=\ga(A)'\cap \mathcal L(H_{\psi_1})\underset{N}{_b*_\alpha}\alpha(N)$, which is equal to the commutant $(\ga(A)\cup 1_{H_{\psi_1}}\underset{N}{_b\otimes_\alpha}\alpha(N)')'$. Thanks to (\cite{E5}, 11.5(ii)), we have :
\[(\ga(A)\cup 1_{H_{\psi_1}}\underset{N}{_b\otimes_\alpha}\alpha(N)')''=A\underset{N}{_b*_\alpha}\mathcal L(H)\]
and, therefore, we get that $(\ga(A)\cup 1_{H_{\psi_1}}\underset{N}{_b\otimes_\alpha}\alpha(N))'=A'\underset{N}{_b\otimes_\alpha}1$. So, any $X\in B'$ belongs to $A'$, and we finally get that $B=A$, which finishes the proof. 
\end{proof}

\subsection{Proposition}
\label{propK}
{\it With the assumptions of \ref{cortildeG}, let us write $K^{it}=\tilde{G}^*(J_\Phi\delta^{it}J_\Phi\underset{N^o}{_\alpha\otimes_b}1)\tilde{G}$. Let $T_2$ be the canonical operator-valued weight from $s(A^\ga)'$ onto $A$ obtained by the basic construction from $T_\ga$, and $\psi_2=\psi_1\circ T_2$. Then :
\newline
(i) the one-parameter group of unitaries $K^{it}$ on $H_{\psi_1}\underset{\psi_0}{_s\otimes_r}H_{\psi_1}$ belongs to $A'\underset{A^\ga}{_s*_r}A$. 
\newline
(ii) there exists a one-parameter group of automorphisms $\rho_t$ of $s(A^\ga)'$ such that, for all $y\in s(A^\ga)'$, we have :
\[K^{it}(y\underset{A^\ga}{_s\otimes_r}1)K^{-it}=\rho_t(y)\underset{A^\ga}{_s\otimes_r}1\]
with $\rho_t(x)=x$, for all $x\in A$. 
\newline
(iii) for any $y\in s(A^\ga)'^+$, we have $\psi_2\circ\rho_t(y)=\psi_2(b(q)^{-t}y)$, where $q$ belongs to $Z(N)$ and is such that the scaling operator of $\gG$ is $\lambda=\alpha(q)=\beta(q)$, and $b(q)\in Z(A)$, and $T_2(\rho_t(y))=b(q)^{-t}T_2(y)$. 
\newline
(iv) Let us identify $H_{\psi_2}$ and $H_{\psi_1}\underset{\psi_0}{_s\otimes_r}H_{\psi_1}$ (\ref{basic}); then $K^{it}$ is the standard implementation of $\rho_t$, and, therefore, we have :
\[K^{it}\Lambda_{\psi_2}(X)=\Lambda_{\psi_2}(b(q)^{t/2}\rho_t(X))\]
for any $X\in\gN_{\psi_2}$, and :
\[\sigma_{\psi_0^o}(J_{\psi_1}\underset{A^\ga}{_s\otimes_r}J_{\psi_1})K^{it}(J_{\psi_1}\underset{(A^\ga)^o}{_r\otimes_s}J_{\psi_1})\sigma_{\psi_0}=K^{it}\]
(v) It is possible to define a one-parameter group of unitaries $K^{it}\underset{\nu}{_{1\underset{A^\ga}{_s\otimes_r}b}\otimes_\alpha}\delta^{it}$ on $H_{\psi_1}\underset{\psi_0}{_s\otimes_r}H_{\psi_1}\underset{\nu}{_b\otimes_\alpha}H$, with natural values on elementary tensors; moreover, we have :
\[(id\underset{A^\ga}{_s*_r}\ga)(K^{it})=K^{it}\underset{\nu}{_{1\underset{A^\ga}{_s\otimes_r}b}\otimes_\alpha}\delta^{it}\]
(vi) for any $s$,$t$ in $\mathbb{R}$, we have $(\Delta_{\psi_1}^{it}\underset{\psi_0}{_s\otimes_r}\Delta_{\psi_1}^{it})(K^{is})(\Delta_{\psi_1}^{-it}\underset{\psi_0}{_s\otimes_r}\Delta_{\psi_1}^{-it})=K^{is}$. }

\begin{proof}
By a straightforward application of \ref{tildeG}(iv), we get that $K^{it}$ belongs to $\mathcal L(H_{\psi_1})\underset{A^\ga}{_s*_r}A$; moreover, using \ref{G}(v), we get, for any $X\in A\rtimes_\ga\gG$, that :
\[K^{it}(\pi_\ga(X)\underset{A^\ga}{_s\otimes_r}1)K^{-it}=K^{it}G^*XGK^{-it}\]
and, therefore, if $X=\ga(x)$, with $x\in A$, we have, using \ref{G}(iv) :
\[K^{it}(x\underset{A^\ga}{_s\otimes_r}1)K^{-it}=G^*(1\underset{N}{_b\otimes_\alpha}J_\Phi\delta^{it}J_\Phi)\ga(x)(1\underset{N}{_b\otimes_\alpha}J_\Phi\delta^{-it}J_\Phi)G=G^*\ga(x)G=x\underset{A^\ga}{_s\otimes_r}1\]
from which we finish the proof of (i). 
\newline
Using (\cite{E5}, 3.11(ii)), we get that, for any $a\in M$, that $\hat{\delta}^{it}a\hat{\delta}^{-it}=\tau_{-t}\sigma_{-t}^{\Phi\circ R}(a)$, and, applying this result to $\widehat{\gG}$, we get that $\delta^{it}b\delta^{-it}=\hat{\tau}_{-t}\sigma_{-t}^{\widehat{\Phi}\circ \hat{R}}(b)$, for any $b$ in $\widehat{M}$; moreover, using now (\cite{E5}, 3.10(iv)), we get that :
\[J_\Phi\delta^{it}J_\Phi bJ_\Phi\delta^{-it}J_\Phi=\hat{\tau}_{-t}\sigma_t^{\widehat{\Phi}}(b)\]
and that, for any $c$ in $\widehat{M}'$,  $J_\Phi\delta^{it}J_\Phi cJ_\Phi\delta^{-it}J_\Phi$ belongs to $\widehat{M}'$, and, more precisely, that :
 \[J_\Phi\delta^{it}J_\Phi cJ_\Phi\delta^{-it}J_\Phi=\hat{\tau}_{-t}^c\sigma_t^{\widehat{\Phi}^c}(c)\]
from which we infer that the one-parameter group of unitaries $1\underset{N}{_b\otimes_\alpha}J_\Phi\delta^{it}J_\Phi$ implements a one-parameter group of automorphisms of $A\rtimes_\ga\gG$; which gives (ii), thanks to \ref{G}(iv). 
\newline
Then, using (\cite{E5}, 13.4, and 3.8(vii) applied to $\widehat{\gG}^c$), we get that, for any $x\in\gN_{\psi_1}$ and $c\in \gN_{\widehat{\Phi}^c}$, we have :
\begin{eqnarray*}
\widetilde{\psi_1}[(1\underset{N}{_b\otimes_\alpha}J_\Phi\delta^{it}J_\Phi)\ga(x^*)(1\underset{N}{_b\otimes_\alpha}c^*c)\ga(x)(1\underset{N}{_b\otimes_\alpha}J_\Phi\delta^{-it}J_\Phi)]
&=&
\|\Lambda_{\psi_1}(x)\underset{\nu}{_b\otimes_\alpha}\Lambda_{\widehat{\Phi}^c}(\hat{\tau}_{-t}^c\sigma_t^{\widehat{\Phi}^c}(c))\|^2\\
&=&
\|\Lambda_{\psi_1}(x)\underset{\nu}{_b\otimes_\alpha}\Lambda_{\widehat{\Phi}^c}(\lambda^{-t/2}c)\|^2\\
&=&
\widetilde{\psi_1}[\ga(x^*)(1\underset{N}{_b\otimes_\alpha}\lambda^{-t}c^*c)\ga(x)]
\end{eqnarray*}
from which we get, using again (\cite{E5}, 13.4), that :
\[\widetilde{\psi_1}[(1\underset{N}{_b\otimes_\alpha}J_\Phi\delta^{it}J_\Phi)X(1\underset{N}{_b\otimes_\alpha}J_\Phi\delta^{-it}J_\Phi)]=\widetilde{\psi_1}((1\underset{N}{_b\otimes_\alpha}\lambda^{-t})X)\]
for any $X\in\gM_{\tilde{\psi_1}}^+$. 
\newline
As $\lambda=\alpha(q)=\beta(q)$ is affiliated to $Z(M)$ (\cite{E5}, 3.8(vi)), we get that $1\underset{N}{_b\otimes_\alpha}\lambda=\ga(b(q))$, and that $b(q)$ is affiliated to $Z(A)$. Using now \ref{G*}, we get the result. 
\newline
 Using now (\cite{E5}, 3.10(vii) and again 13.4), we get that $(1\underset{N}{_b\otimes_\alpha}J_\Phi\delta^{it}J_\Phi)$ is the standard implementation of $Ad((1\underset{N}{_b\otimes_\alpha}J_\Phi\delta^{it}J_\Phi)_{|A\rtimes_\ga\gG}$ on $H_{\psi_1}\underset{\nu}{_b\otimes_\alpha}H$ (which is identified with $H_{\tilde{\psi_1}}$ by (\cite{E5}, 13.4). Therefore, using again \ref{G*}, we get that $K^{it}$ is the standard implementation of $\rho_t$; thanks to (iii), we finish the proof of (iv). 
 \newline
 Similarly, using again (\cite{E5}, 3.8 (i) and (v)), we get $\alpha(n)\delta^{it}=\delta^{it}\alpha(\gamma_t\sigma_t^\nu(n))$; as $\nu$ is invariant under $\gamma_t$ (\cite{E5}, 3.8 (v)), there exists a one-parameter group of unitaries $h^{it}$ on $H_\nu$, such that, for all $n\in\gN_\nu$, we gave $\Lambda_\nu(\gamma_t\sigma_t^\nu(n))=h^{it}\Lambda_\nu(n)$, and $h^{it}mh^{-it}=\gamma_t\sigma_t^\nu(m)$, for all $m\in N$; therefore, if $\eta$ is in $D(_\alpha H, \nu)$, it is straightforward to get that $\delta^{it}\eta$ belongs also to $D(_\alpha H, \nu)$; more precisely, we have then :
\[R^{\alpha, \nu}(\delta^{it}\eta)\Lambda_\nu(n)=\alpha(n)\delta^{it}\eta=\delta^{it}\alpha(\gamma_t\sigma_t^\nu(n))\eta=\delta^{it}R^{\alpha, \nu}(\eta)h^{it}\Lambda_\nu(n)\]
from which we infer that $R^{\alpha, \nu}(\delta^{it}\eta)=\delta^{it}R^{\alpha, \nu}(\eta)h^{it}$, and, if $\eta$, $\eta'$ belong to $D(_\alpha H, \nu)$, we get that $<\delta^{it}\eta, \delta^{it}\eta'>_{\alpha, \nu}^o=\gamma_{-t}\sigma_{-t}^\nu(<\eta, \eta'>_{\alpha, \nu}^o)$. From which we get, for all $\xi$, $\xi'$ in $H_{\psi_1}\underset{\psi_0}{_s\otimes_r}H_{\psi_1}$, using (ii), that :
\begin{eqnarray*}
(K^{it}\xi\underset{\nu}{_{1\underset{A^\ga}{_s\otimes_r}b}\otimes_\alpha}\delta^{it}\eta|K^{it}\xi'\underset{\nu}{_{1\underset{A^\ga}{_s\otimes_r}b}\otimes_\alpha}\delta^{it}\eta')
&=&
((1\underset{A^\ga}{_s\otimes_r}b(\gamma_{-t}\sigma_{-t}^\nu(<\eta, \eta'>_{\alpha, \nu}^o)))K^{it}\xi|\xi')\\
&=&
(K^{it}(1\underset{A^\ga}{_s\otimes_r}b(<\eta, \eta'>_{\alpha, \nu}^o))\xi|K^{it}\xi')\\
&=&
(\xi\underset{\nu}{_{1\underset{A^\ga}{_s\otimes_r}b}\otimes_\alpha}\eta|\xi'\underset{\nu}{_{1\underset{A^\ga}{_s\otimes_r}b}\otimes_\alpha}\eta')
\end{eqnarray*}
from which we get the first result of (v). 
\newline 
Using now \ref{G}(iv) and \ref{tildeG}(iv), we get that $(id\underset{A^\ga}{_s*_r}\varsigma_N\ga)(K^{it})$ is equal to :
\[[(\widehat{W}\underset{N^o}{_\alpha\otimes_b}1)\sigma^{2,3}_{\alpha, \hat{\beta}}(\tilde{G}\underset{N}{_b\otimes_\alpha}1)(1\underset{A^\ga}{_s\otimes_r}\sigma_{\nu^o})]^*(1\underset{N}{_\alpha\otimes_b}\tilde{G})(J_\Phi\delta^{it}J_\Phi\underset{N}{_\alpha\otimes_b}1\underset{A^\ga}{_s\otimes_r}1)(1\underset{N}{_\alpha\otimes_b}\tilde{G})^*...\]
\[...(\widehat{W}\underset{N^o}{_\alpha\otimes_b}1)\sigma^{2,3}_{\alpha, \hat{\beta}}(\tilde{G}\underset{N}{_b\otimes_ \alpha}1)(1\underset{A^\ga}{_s\otimes_r}\sigma_{\nu^o})\]
and is therefore equal to 
\[[\sigma^{2,3}_{\alpha, \hat{\beta}}(\tilde{G}\underset{N}{_b\otimes_\alpha}1)(1\underset{A^\ga}{_s\otimes_r}\sigma_{\nu^o})]^*(\widehat{W}^*(J_\Phi\delta^{it}J_\Phi\underset{N^o}{_\alpha\otimes_\beta}1)\widehat{W}\underset{N^o}{_\alpha\otimes_b}1)[\sigma^{2,3}_{\alpha, \hat{\beta}}(\tilde{G}\underset{N}{_b\otimes_\alpha}1)(1\underset{A^\ga}{_s\otimes_r}\sigma_{\nu^o})]\]
Using (\cite{E5}, successively 3.10 (vii), 3.11 (iii), 3.6 and 3.8 (vi)) we get :
\begin{eqnarray*}
\widehat{W}^*(J_\Phi\delta^{it}J_\Phi\underset{N^o}{_\alpha\otimes_\beta}1)\widehat{W}
&=&
\sigma_{\nu^o} W(1\underset{N}{_\beta\otimes_\alpha}J_\Phi\delta^{it}J_\Phi)W^*\sigma_{\nu^o}\\
&=&
\sigma_{\nu^o} (\hat{J}\underset{\nu}{_\beta\otimes_\alpha} J_\Phi)W^*(1\underset{N}{_\beta\otimes _\alpha}\delta^{it}) W(\hat{J}\underset{\nu^o}{_\alpha\otimes_{\hat{\beta}}} J_\Phi)\sigma_{\nu^o}\\
&=&
\sigma_{\nu^o} (\hat{J}\underset{\nu}{_\beta\otimes_\alpha} J_\Phi)\Gamma(\delta^{it})(\hat{J}\underset{\nu^o}{_\alpha\otimes_{\hat{\beta}}} J_\Phi)\sigma_{\nu^o}\\
&=&
\sigma_{\nu^o} (\hat{J}\underset{\nu}{_\beta\otimes_\alpha} J_\Phi)(\delta^{it}\underset{\nu}{_\beta\otimes_\alpha}\delta^{it})(\hat{J}\underset{\nu^o}{_\alpha\otimes_{\hat{\beta}}} J_\Phi)\sigma_{\nu^o}\\
&=&
J_\Phi\delta^{it}J_\Phi\underset{N}{_{\hat{\beta}}\otimes_\alpha}\hat{J}\delta^{it}\hat{J}\\
&=&
J_\Phi\delta^{it}J_\Phi\underset{N}{_{\hat{\beta}}\otimes_\alpha}\delta^{it}
\end{eqnarray*}
from which we get the second result of (v). 
\newline
Using (iv) and (\cite{E5}, 3.11(ii)), we get :
\begin{eqnarray*}
(\Delta_{\psi_1}^{it}\underset{\psi_0}{_s\otimes_r}\Delta_{\psi_1}^{it})(K^{is})(\Delta_{\psi_1}^{-it}\underset{\psi_0}{_s\otimes_r}\Delta_{\psi_1}^{-it})
&=&
(\Delta_{\psi_1}^{it}\underset{\psi_0}{_s\otimes_r}\Delta_{\psi_1}^{it})\tilde{G}^*(J_\Phi\delta^{is}J_\Phi\underset{\nu^o}{_\alpha\otimes_b}1)\tilde{G}(\Delta_{\psi_1}^{-it}\underset{\psi_0}{_s\otimes_r}\Delta_{\psi_1}^{-it})\\
&=&
\tilde{G}^*((\delta\Delta_{\widehat{\Phi}})^{-it}J_\Phi\delta^{is}J_\Phi(\delta\Delta_{\widehat{\Phi}})^{it}\underset{\nu^o}{_\alpha\otimes_b}1)\tilde{G}\\
&=&
\tilde{G}^*(J_\Phi\delta^{is}J_\Phi\underset{\nu^o}{_\alpha\otimes_b}1)\tilde{G}
\end{eqnarray*}
which is (vi). 
\end{proof}

\subsection{Theorem}
\label{thK}
{\it Let's suppose again the assumptions of  \ref{cortildeG} and \ref{propK}; let $K^{it}$ be the one-parameter group of unitaries on $H_{\psi_1}\underset{\psi_0}{_s\otimes_r}H_{\psi_1}$ defined in \ref{propK}; let us suppose that there exists a positive non singular operator $\delta_A$ affiliated to $A$ such that we have, for all $t\in\mathbb{R}$ :
\[K^{it}=J_{\psi_1}\delta_A^{it}J_{\psi_1}\underset{A^\ga}{_s\otimes_r}\delta_A^{it}\]
Then :
\newline
(i) it is possible to define a one parameter group of unitaries $\delta_A^{it}\underset{\nu}{_b\otimes_\alpha}\delta^{it}$ on $H_{\psi_1}\underset{\nu}{_b\otimes_\alpha}H$, with natural values on elementary tensors; moreover, we have :
\[\ga(\delta_A^{it})=\delta_A^{it}\underset{\nu}{_b\otimes_\alpha}\delta^{it}\]
(ii) for all $s$, $t$ in $\mathbb{R}$, we have $\sigma_s^{\psi_1}(\delta_A^{it})=b(q)^{ist}\delta_A^{it}$.
\newline
(iii) there exists a normal semi-finite faithful weight $\phi$ on $A$ such that $(A, b, \ga, \phi, \psi_0)$ is a Galois system. Moreover, the modulus of this Galois action is operator $\delta_A$, and the scaling operator is equal to $b(q)$, where $q\in Z(N)$ is such that $\alpha(q)=\beta(q)=\lambda$, the scaling operator of $\gG$.}

\begin{proof}
Using \ref{propK}(v), we easily get (i). Using now \cite{E5}, 8.8(iii), we get that $\ga(\sigma_s^{\psi_1}(\delta_A)^{it})=\sigma_s^{\psi_1}(\delta_A)^{it}\underset{\nu}{_b\otimes_\alpha}\delta^{it}$, and, therefore, that $\sigma_s^{\psi_1}(\delta_A)^{it})\delta_A^{-it}$ belongs to $r(A^\ga)$. 
\newline
So, there exists $k\eta r(A^\alpha)$ such that $\sigma_s^{\psi_1}(\delta_A^{it})=k^{ist}\delta_A^{it}=\delta_A^{it}k^{ist}$. 
\newline
Let us write $k=\int_0^\infty \lambda de_\lambda$, and let us put $f_n=\int_{1/n}^n de_\lambda$; then using (\cite{E5}, 2.2.2), we get that, for any $x\in\gN_{\psi_1}\cap\gN_{T_\ga}$, $xf_nk^{-t/2}$ is bounded and belongs to $\gN_{T_\ga}\cap \gN_{\psi_1}$, and, with same arguments, we get that $xf_nk^{-t/2}\delta_A^{-it}$ belongs also to $\gN_{T\ga}\cap\gN_{\psi_1}$. We then get that :
\[J_{\psi_1}\delta_A^{it}J_{\psi_1}\Lambda_{T_\ga}(xf_n)=\Lambda_{T_\ga}(xf_nk^{-t/2}\delta_A^{-it})\]
and, therefore, with the notations of \ref{propK}(ii) :
\[\rho_t(\Lambda_{T_\ga}(xf_n)\Lambda_{T_\ga}(xf_n)^*)=\Lambda_{T_\ga}(xf_nk^{-t/2}\delta_A^{-it})\Lambda_{T_\ga}(xf_nk^{-t/2}\delta_A^{-it})^*\]
from which we get that :
\[T_2\rho_t(\Lambda_{T_\ga}(xf_n)\Lambda_{T_\ga}(xf_n)^*)=xf_nk^{-t}x^*\]
and, on the other hand, using \ref{propK}(iii), we have, using the fact that $b(q)$ is affiliated to $Z(A)$ :
\[T_2\rho_t(\Lambda_{T_\ga}(xf_n)\Lambda_{T_\ga}(xf_n)^*)=b(q)^{-t}xf_nx^*=xf_nb(q)^{-t}x^*\]
from which we easily deduce that $k=b(q)$, which finishes the proof of (ii). 
\newline
Using \cite{V1}, we get that there is a normal semi-finite faithful weight $\phi$ on $A$, such that $(D\phi : D\psi_1)_t=b(q)^{it^2/2}\delta_A^{it}$, and that the modular groups $\sigma^\phi$ and $\sigma^{\psi_1}$ commute. If $x\in N_\phi$ is such that $x\delta_A^{1/2}$ is bounded, then this last operator belongs to $\gN_{\psi_1}$ and we can identify $\Lambda_\phi(x)$ with $\Lambda_{\psi_1}(x\delta_A^{1/2})$ and $J_\phi$ with $b(q)^{i/4}J_{\psi_1}$; we shall denote, for $n\in N$, $a(n)=J_{\psi_1}b(n^*)J_{\psi_1}=J_\phi b(n^*)J_\phi$. 
\newline
For $x\in\gN_\phi$ and $\eta\in D(_\alpha H, \nu)\cap D(H_\beta, \nu^o)\cap \mathcal D(\delta^{1/2})$, such that $\delta^{-1/2}\eta$ belongs to $D(_\alpha H, \nu)$, we have, using these remarks, and the fact that $\psi_1$ is $\delta$-invariant, \ref{inner}(iii) :
\begin{eqnarray*}
\|\Lambda_\phi(x)\underset{\nu^o}{_a\otimes_\beta}\eta\|^2
&=&
\|\Lambda_{\psi_1}(x\delta_A^{1/2})\underset{\nu^o}{_a\otimes_\beta}\eta\|^2\\
&=&
(\psi_1\underset{N}{_b*_\alpha}\omega_{\delta^{-1/2}\eta})(\ga(\delta_A^{1/2}x^*x\delta_A^{1/2})]\\
&=&
(\psi_1\underset{N}{_b*_\alpha}\omega_{\delta^{-1/2}\eta})[(\delta_A^{1/2}\underset{\nu}{_b\otimes_\alpha}\delta^{1/2})\ga(x^*x)(\delta_A^{1/2}\underset{\nu}{_b\otimes_\alpha}\delta^{1/2})]\\
&=&
(\phi\underset{N}{_b*_\alpha}\omega_\eta)[\ga(x^*x)]
\end{eqnarray*}
which remains true for any $\eta\in D(_\alpha H, \nu)\cap D(H_\beta, \nu^o)$, and gives then, using (\cite{E6}7.6) that the weight $\phi$ is invariant under $\ga$, which finishes the proof. 

\end{proof}

\subsection{Corollary}
\label{inner}
{\it Let $(b, \ga)$ be a Galois action of the measured quantum groupoid $\gG$ on a von Neumann algebra $A$; let $\psi_0$ be a normal semi-finite faithful weight on $A^\ga$, bearing the Galois density property defined in \ref{defdens}; let $\rho_t$ be the one-parameter group of automorphisms of $s(A^\ga)'$ defined in \ref{propK}. Let us suppose that this one-parameter group is inner; then :
\newline
(i) there exists a non singular positive operator $\delta_A$ affiliated to $A\cap r(A^\ga)'$ such that :
\[K^{it}=J_{\psi_1}\delta_A^{it}J_{\psi_1}\underset{A^\ga}{_s\otimes_r}\delta_A^{it}\]
\newline
(ii) there exists a normal semi-finite faithful weight $\phi$ on $A$ such that $(A, b, \ga, \phi, \psi_0)$ is a Galois system.}
\begin{proof}
As $\rho_t(x)=x$ for all $x\in A$, we get that there exists a positive non-singular operator $\delta_A$ affiliated to $A\cap r(A^\ga)'$ such that, for all $x\in s(A^\ga)'$, we have :
\[\rho_t(x)=J_{\psi_1}\delta_A^{it}J_{\psi_1}xJ_{\psi_1}\delta_A^{-it}J_{\psi_1}\]
and then, using \ref{propK}(iv), we have (i). Result (ii) is then a direct corollary of \ref{thK}(iii). 
\end{proof}

\subsection{Corollary}
\label{corsum}
{\it Let $(b, \ga)$ be a Galois action of the measured quantum groupoid $\gG$ on a von Neumann algebra $A$; let us suppose that the invariant subalgebra $A^\ga$ is a finite sum of factors (in particular, if $A^\ga$ is finite dimensional); let $\psi_0$ be a normal semi-finite faithful weight on $A^\ga$, bearing the Galois density property defined in \ref{defdens}; then, there exists a normal semi-finte faithful weight $\phi$ on $A$ such that $(A, b, \ga, \phi, \psi_0)$ is a Galois system. }
\begin{proof}
The center $Z(s(A^\ga)')$ is equal to $s(Z(A^\ga))$; if $A^\ga=\oplus_{i\in I}F_i$ (with a finite set $I$), we get that $Z(A^\ga)=\oplus_{i\in I}\mathbb{C}$, and that any automorphism of $Z(s(A^\ga)')$ gives a permutation in the set $I$; therefore, the restriction of the one parameter group $\rho_t$ defined in \ref{propK} to this center gives a continuous function $\rho$ from $\mathbb{R}$ to the set $\gS(I)$ (with the pointwise topology); as $\rho^{-1}({id})$ is open, closed, non empty, we get that $\rho_t$ acts identically on the center $Z(s(A^\ga)')$; therefore, $\rho_t$ is inner, by (\cite{StZ} 8.11), and we get the result by \ref{inner}. 
\end{proof}

\subsection{Theorem}
\label{tau}
{\it Let $(b, \ga)$ be a Galois action of the measured quantum groupoid $\gG$ on a von Neumann algebra $A$; let $\psi_0$ be a normal semi-finite faithful weight on $A^\ga$, bearing the Galois density property defined in \ref{defdens};  let $\tilde{G}$ be the Galois unitary of $(b, \ga)$, as defined in \ref{defGalois}; let $\psi_1=\psi_0\circ T_\ga$, and $\tilde{\psi_1}$ its dual weight on the crossed-product $A\rtimes_\ga\gG$. Then :
\newline
(i) for all $s$, $t$ in $\mathbb{R}$, we have :
\[\sigma_t^{\tilde{\psi_1}}(1\underset{N}{_b\otimes_\alpha}\hat{J}\hat{\delta}^{is}\hat{J})=1\underset{N}{_b\otimes_\alpha}\hat{J}\hat{\delta}^{is}\hat{J}\]
(ii) there exists a one-parametrer group of unitaries $P_A^{it}=\Delta_{\psi_1}^{it}\pi_\ga((1\underset{N}{_b\otimes_\alpha}\hat{J}\hat{\delta}^{-it}\hat{J})$ on $H_{\psi_1}$, which defines a one-parameter group $\tau^A_t$ of automorphism of $A$ defined, for all $X\in A$, by $\tau^A_t(X)=P_A^{it}AP_A^{-it}$. 
For any $x\in A^\ga$, we have $\tau^A_t(x)=\sigma_t^{\psi_0}(x)$, and, for all $n\in N$, we have $\tau^A_t(b(n))=b(\sigma^\nu_{t}(n))$. 
\newline
(iii) for all $t\in\mathbb{R}$, we have :
\[\ga(\tau^A_t(X))=(\sigma_t^{\psi_1}\underset{N}{_b*_\alpha}\sigma_t^{\Phi\circ R})\ga(X)=(\tau^A_t\underset{N}{_b*_\alpha}\tau_t)\ga(X)\]
\[\ga(\sigma_t^{\psi_1}(X))=(\tau^A_t\underset{N}{_b*_\alpha}\sigma_t^\Phi)\ga(X)\]
(iv) we have, for any positive $X\in A$ :
\[\psi_1\circ\tau^A_t(X)=\psi_1(b(q)^{-t}X)\]
where $q$ is the positive non singular operator affiliated to $N$ such that the scaling operator $\lambda$ of $\gG$ satisfies $\lambda=\alpha(q)=\beta(q)$ (\cite{E6}, 3.8(vi)).
\newline
(v) for any $X\in\gN_{\psi_1}$, we have $P_A^{it}\Lambda_{\psi_1}(X)=b(q)^{t/2}\Lambda_{\psi_1}(\tau^A_t(X))$. So, $P_A^{it}$ is the standard implementation of $\tau_t^A$, and $J_{\psi_1}P_A^{it}=P_A^{it}J_{\psi_1}$.
\newline
(vi) there exists a one parameter group of unitaries $P_A^{it}\underset{A^\ga}{_s\otimes_r}P_A^{it}$ on $H_{\psi_1}\underset{\psi_0}{_s\otimes_r}H_{\psi_1}$, with natural values on elementary tensors, and a one parametrer group of unitaries $P^{it}\underset{N^o}{_\alpha\otimes_b}P_A^{it}$ on $H\underset{\nu^o}{_\alpha\otimes_b}H_{\psi_1}$, with natural values on elementary tensors, and we have, for all $t\in\mathbb{R}$ :
\[\tilde{G}(P_A^{it}\underset{A^\ga}{_s\otimes_r}P_A^{it})=(P^{it}\underset{N^o}{_\alpha\otimes_b}P_A^{it})\tilde{G}\]
}
\begin{proof}
We know (\cite{E6}, 3.2) that $\Delta_{\tilde{\psi_1}}^{it}=\Delta_{\psi_1}^{it}\underset{N}{_b\otimes_\alpha}(\delta\Delta_{\widehat{\Phi}}^{-it})$; from which, using (\cite{E5}, 3.11 (ii), (vii) and (iv)) we get that :
\begin{eqnarray*}
\sigma_t^{\tilde{\psi_1}}(1\underset{N}{_b\otimes_\alpha}\hat{J}\hat{\delta}^{is}\hat{J})
&=&
(\Delta_{\psi_1}^{it}\underset{N}{_b\otimes_\alpha}(\delta\Delta_{\widehat{\Phi}}^{-it}))(1\underset{N}{_b\otimes_\alpha}\hat{J}\hat{\delta}^{is}\hat{J})(\Delta_{\psi_1}^{-it}\underset{N}{_b\otimes_\alpha}(\delta\Delta_{\widehat{\Phi}}^{it}))\\
&=&
1\underset{N}{_b\otimes_\alpha}(\hat{\delta}\Delta_{\Phi})^{it}(\hat{J}\hat{\delta}^{is}\hat{J})(\hat{\delta}\Delta_{\Phi})^{-it}\\
&=&
1\underset{N}{_b\otimes_\alpha}\Delta_{\Phi}^{it}(\hat{J}\hat{\delta}^{is}\hat{J})\Delta_{\Phi}^{-it}\\
&=&
1\underset{N}{_b\otimes_\alpha}P^{it}(\hat{J}\hat{\delta}^{it}\hat{J})(\hat{J}\hat{\delta}^{is}\hat{J})(\hat{J}\hat{\delta}^{-it}\hat{J})P^{-it}\\
&=&
1\underset{N}{_b\otimes_\alpha}P^{it}(\hat{J}\hat{\delta}^{is}\hat{J})P^{-it}\\
&=&
1\underset{N}{_b\otimes_\alpha}\hat{J}\hat{\tau}_t(\hat{\delta}^{is})\hat{J}\\
&=&
1\underset{N}{_b\otimes_\alpha}\hat{J}\hat{\delta}^{is}\hat{J}
\end{eqnarray*}
which gives (i). From (i), and \ref{G*}(iv) and \ref{defGalois}, we get that :
\[\sigma_t^{\psi_2}(\pi_\ga(1\underset{N}{_b\otimes_\alpha}\hat{J}\hat{\delta}^{is}\hat{J}))=\pi_\ga(1\underset{N}{_b\otimes_\alpha}\hat{J}\hat{\delta}^{is}\hat{J})\]
from which we get :
\[\Delta_{\psi_1}^{it}(\pi_\ga(1\underset{N}{_b\otimes_\alpha}\hat{J}\hat{\delta}^{is}\hat{J})\Delta_{\psi_1}^{-it}=\pi_\ga(1\underset{N}{_b\otimes_\alpha}\hat{J}\hat{\delta}^{is}\hat{J})\]
which gives the commutation of the two one-parameter groups of unitaries $\Delta_{\psi_1}^{it}$ and 
$\pi_\ga(1\underset{N}{_b\otimes_\alpha}\hat{J}\hat{\delta}^{is}\hat{J})$, and the existence of the one-parameter group of unitaries $P_A^{it}$. 
\newline
We easily get  that $\tilde{\ga}[(1\underset{N}{_b\otimes_\alpha}\hat{J}\hat{\delta}^{-it}\hat{J})\ga(X)(1\underset{N}{_b\otimes_\alpha}\hat{J}\hat{\delta}^{it}\hat{J})]$ is equal to :
\begin{multline*}
(1\underset{N}{_b\otimes_\alpha}\hat{J}\hat{\delta}^{-it}\hat{J}\underset{N^o}{_{\hat{\alpha}}\otimes_\beta}\hat{J}\hat{\delta}^{-it}\hat{J})(\ga(X){_{\hat{\alpha}}\otimes_\beta}1)(1\underset{N}{_b\otimes_\alpha}\hat{J}\hat{\delta}^{it}\hat{J}\underset{N^o}{_{\hat{\alpha}}\otimes_\beta}\hat{J}\hat{\delta}^{it}\hat{J})\\
=
[(1\underset{N}{_b\otimes_\alpha}\hat{J}\hat{\delta}^{-it}\hat{J})\ga(X)(1\underset{N}{_b\otimes_\alpha}\hat{J}\hat{\delta}^{it}\hat{J})]\underset{N^o}{_{\hat{\alpha}}\otimes_\beta}1
\end{multline*}
from which we get, using (\cite{E5}, 10.12), that $(1\underset{N}{_b\otimes_\alpha}\hat{J}\hat{\delta}^{-it}\hat{J})\ga(X)(1\underset{N}{_b\otimes_\alpha}\hat{J}\hat{\delta}^{it}\hat{J})$ belongs to $\ga(A)$, and, therefore, that $\pi_\ga(1\underset{N}{_b\otimes_\alpha}\hat{J}\hat{\delta}^{-it}\hat{J})X\pi_\ga(1\underset{N}{_b\otimes_\alpha}\hat{J}\hat{\delta}^{it}\hat{J})$ belongs to $A$, from which it is straightforward to get that $P_A^{it}XP_A^{-it}$ belongs to $A$, and gives the existence of $\tau^A$. 
\newline
We have :
\begin{eqnarray*}
\ga(\tau^A_t(x))
&=&
\ga[\pi_\ga(1\underset{N}{_b\otimes_\alpha}\hat{J}\hat{\delta}^{-it}\hat{J})\sigma_t^{\psi_0}(x)\pi_\ga(1\underset{N}{_b\otimes_\alpha}\hat{J}\hat{\delta}^{it}\hat{J})]\\
&=&(1\underset{N}{_b\otimes_\alpha}\hat{J}\hat{\delta}^{-it}\hat{J})(\sigma_t^{\psi_0}(x)\underset{N}{_b\otimes_\alpha}1)(1\underset{N}{_b\otimes_\alpha}\hat{J}\hat{\delta}^{it}\hat{J})\\
&=&
\sigma_t^{\psi_0}(x)\underset{N}{_b\otimes_\alpha}1\\
&=&
\ga(\sigma_t^{\psi_0}(x))
\end{eqnarray*}
from which we get that $\tau^A_t(x)=\sigma_t^{\psi_0}(x)$.
\newline
We have, for all $n\in N$ :
\begin{eqnarray*}
\ga(\tau^A_t(b(n)))
&=&
\ga[\pi_\ga(1\underset{N}{_b\otimes_\alpha}\hat{J}\hat{\delta}^{-it}\hat{J})\sigma_t^{\psi_1}(b(n))\pi_\ga(1\underset{N}{_b\otimes_\alpha}\hat{J}\hat{\delta}^{it}\hat{J})]\\
&=&
\ga[\pi_\ga(1\underset{N}{_b\otimes_\alpha}\hat{J}\hat{\delta}^{-it}\hat{J})\delta_A^{-it}\sigma_t^{\phi}(b(n))\delta_A^{it}\pi_\ga(1\underset{N}{_b\otimes_\alpha}\hat{J}\hat{\delta}^{it}\hat{J})] \\
&=&
\ga[\pi_\ga(1\underset{N}{_b\otimes_\alpha}\hat{J}\hat{\delta}^{-it}\hat{J})\delta_A^{-it}b(\sigma_{-t}^\nu(n))\delta_A^{it}\pi_\ga(1\underset{N}{_b\otimes_\alpha}\hat{J}\hat{\delta}^{it}\hat{J})]\\
&=&
(1\underset{N}{_b\otimes_\alpha}\hat{J}\hat{\delta}^{-it}\hat{J})(\delta_A^{-it}\underset{N}{_b\otimes_\alpha}\delta^{-it})(1\underset{N}{_b\otimes_\alpha}\beta((\sigma_{-t}^\nu(n)))(\delta_A^{it}\underset{N}{_b\otimes_\alpha}\delta^{it})(1\underset{N}{_b\otimes_\alpha}\hat{J}\hat{\delta}^{it}\hat{J})\\
&=&
1\underset{N}{_b\otimes_\alpha}\hat{J}\hat{\delta}^{-it}\hat{J}\delta^{-it}\beta((\sigma_{-t}^\nu(n)))\delta^{it}\hat{J}\hat{\delta}^{it}\hat{J}\\
&=&
1\underset{N}{_b\otimes_\alpha}\hat{J}\hat{\delta}^{-it}\hat{J}\sigma_{-t}^{\Phi\circ R}\sigma_{t}^\Phi(\beta((\sigma_{-t}^\nu(n))))\hat{J}\hat{\delta}^{it}\hat{J}\\
&=&
1\underset{N}{_b\otimes_\alpha}\hat{J}\hat{\delta}^{-it}\hat{J}\beta(\gamma_t(n))\hat{J}\hat{\delta}^{it}\hat{J}
\end{eqnarray*}
and, therefore, we get :
\begin{multline*}
\ga(\tau^A_t(b(n)))=
1\underset{N}{_b\otimes_\alpha}\hat{J}\hat{\delta}^{-it}\alpha(\gamma_t(n^*))\hat{\delta}^{it}\hat{J}
=
1\underset{N}{_b\otimes_\alpha}\hat{J}\sigma_{-t}^{\widehat{\Phi}\circ \hat{R}}\sigma_t^{\widehat{\Phi}}(\alpha(\gamma_t(n^*)))\hat{J}\\
=
1\underset{N}{_b\otimes_\alpha}\hat{J}\alpha(\sigma_t^\nu(n^*))\hat{J}
=
1\underset{N}{_b\otimes_\alpha}\beta(\sigma_t^\nu(n))
=
\ga(b(\sigma_t^\nu(n)))
\end{multline*}
 which finishes the proof of (ii). 
\newline
As, using (\cite{E5}, 3.10(vi)), we get that $\hat{J}\hat{\delta}^{-it}\hat{J}$ is equal to $P^{it}\Delta_\Phi^{-it}$, and, therefore, implements $\tau_t\sigma_{-t}^\Phi$, we get, using (\cite{E5}, 8.8), that :
\[\ga(\tau^A_t(X))=(\sigma_t^{\psi_1}\underset{N}{_b*_\alpha}\sigma_{-t}^{\Phi\circ R})\ga(X)\]
Using (\cite{E5}, 3.8 (i) and (ii)), we get that $\Gamma\circ\sigma_{-t}^{\Phi\circ R}=(\sigma_{-t}^{\Phi\circ R}\underset{N}{_b*_\alpha}\tau_t)\Gamma$, from which, we infer :
\begin{eqnarray*}
(\ga\underset{N}{_b*_\alpha}id)\ga(\tau_t^A(x))
&=&
(id\underset{N}{_b*_\alpha}\Gamma)\ga(\tau^A_t(X))\\
&=&
(id\underset{N}{_b*_\alpha}\Gamma)(\sigma_t^{\psi_1}\underset{N}{_b*_\alpha}\sigma_{-t}^{\Phi\circ R})\ga(X)\\
&=&
(\sigma_t^{\psi_1}\underset{N}{_b*_\alpha}\sigma_{-t}^{\Phi\circ R}\underset{N}{_b*_\alpha}\tau_t)(id\underset{N}{_b*_\alpha}\Gamma)\ga(X)\\
&=&
(\sigma_t^{\psi_1}\underset{N}{_b*_\alpha}\sigma_{-t}^{\Phi\circ R}\underset{N}{_b*_\alpha}\tau_t)(\ga\underset{N}{_b*_\alpha}id)\ga(X)\\
&=&
(\ga\underset{N}{_b*_\alpha}id)(\tau^A_t\underset{N}{_b*_\alpha}\tau_t)\ga(X)
\end{eqnarray*}
from which we get that $\ga(\tau^A_t(X))=(\tau^A_t\underset{N}{_b*_\alpha}\tau_t)\ga(X)$. 
\newline
Finally, we have :
\begin{eqnarray*}
\ga(\sigma_t^{\psi_1}(X))
&=&
(1\underset{N}{_b\otimes_\alpha}\hat{J}\hat{\delta}^{it}\hat{J})\ga(\tau^A_t(X))(1\underset{N}{_b\otimes_\alpha}\hat{J}\hat{\delta}^{-it}\hat{J})\\
&=&
(id\underset{N}{_b*_\alpha}\tau_{-t}\sigma_t^\Phi)(\tau^A_t\underset{N}{_b*_\alpha}\tau_t)\ga(X)\\
&=&
(\tau^A_t\underset{N}{_b*_\alpha}\sigma^\Phi_t)\ga(X)
\end{eqnarray*}
which finishes the proof of (iii). 
\newline
We have, for any positive $X\in A$ :
\begin{eqnarray*}
T_\ga(\tau^A_t(X))
&=&
(id\underset{N}{_b*_\alpha}\Phi)\ga(\tau^A_t(X))\\
&=&
(id\underset{N}{_b*_\alpha}\Phi)(\tau^A_t\underset{N}{_b*_\alpha}\tau_t)\ga(X)\\
&=&
\tau^A_t(id\underset{N}{_b*_\alpha}\Phi\circ \tau_t)\ga(X)\\
\end{eqnarray*}
As $\Phi\circ\tau_t(Y)=\Phi(\lambda^{-t}Y)$ for any positive $Y\in M$ (\cite{E5}, 3.8 (vii)), we get that it is equal to $\tau_t^A[T_\ga(b(q)^{-t}X)]$, and, using (ii), to $\sigma_t^{\psi_0}[T_\ga(b(q)^{-t}X)]$; from which we get (iv).
\newline
If $\zeta$ is in $D(_rH_{\psi_1}, \psi_0)$, we have, using (\ref{G}(i) and (iv)), and then (iii) :
\begin{eqnarray*}
G(\pi_\ga(1\underset{N}{_b\otimes_\alpha}\hat{J}\hat{\delta}^{-it}\hat{J})(\Lambda_{\psi_1}(X)\underset{\psi_0}{_s\otimes_r}\zeta)
&=&
(1\underset{N}{_b\otimes_\alpha}\hat{J}\hat{\delta}^{-it}\hat{J})\sum_i e_i\underset{N}{_b\otimes_\alpha}\Lambda_\Phi((\omega_{\zeta, e_i}\underset{N}{_b*_\alpha}id)\ga(X))\\
&=&
\sum_i e_i\underset{N}{_b\otimes_\alpha}P^{it}\Delta_\Phi^{-it}\Lambda_\Phi((\omega_{\zeta, e_i}\underset{N}{_b*_\alpha}id)\ga(X))\\
&=&
\sum_i e_i\underset{N}{_b\otimes_\alpha}\lambda^{-t/2}\Lambda_\Phi[\tau_t\sigma_{-t}^\Phi(\omega_{\zeta, e_i}\underset{N}{_b*_\alpha}id)\ga(X)]\\
&=&
\sum_i e_i\underset{N}{_b\otimes_\alpha}\Lambda_\Phi((\omega_{\zeta, e_i}\underset{N}{_b*_\alpha}id)\ga(\tau^A_t\sigma_t^{\psi_1}X))\\
&=&
(1\underset{N}{_b\otimes_\alpha}\lambda^{-t/2})G(\Lambda_{\psi_1}(\tau^A_t\sigma_{-t}^{\psi_1}(X))\underset{\psi_0}{_s\otimes_r}\zeta)
\end{eqnarray*}
from which we get that :
\[\pi_\ga(1\underset{N}{_b\otimes_\alpha}\hat{J}\hat{\delta}^{-it}\hat{J})\Lambda_{\psi_1}(X)=
b(q)^{-t/2}\Lambda_{\psi_1}(\tau^A_t\sigma_{-t}^{\psi_1}(X))\]
which gives (v). 
\newline
Thanks to (ii), one can easily get that if $\zeta$ belongs to $D(_rH_{\psi_1}, \psi_0)$, so does, for all $t\in\mathbb{R}$, $P_A^{it}\zeta$, and that $R^{r, \psi_0}(P_A^{it}\zeta)=P_A^{it}R^{r,\psi_0}(\zeta)\Delta_{\psi_0}^{-it}$. Using then (v), we obtain easily the existence of the first one parameter group of unitaries. Using again (ii), we get that, if $\zeta'$ belongs to $D((H_{\psi_1})_b, \nu^o)$, so does $P_A^{it}\zeta'$, and that $R^{b, \nu^o}(P_A^{it}\zeta')=P_A^{it}R^{b, \nu^o}(\zeta')\Delta_\nu^{-it}$, from which one gets the existence of the second one parameter group of unitaries. Moreover, using successively (v), \ref{tildeG}(i), (iii), \cite{E6}, 3.8(vii) and (vi), and again \ref{tildeG}(i), we get, for all $\zeta\in D(_rH_{\psi_1}, \psi_0)\cap D((H_{\psi_1})_b, \nu^o)$, $\zeta'\in D((H_{\psi_1})_b, \nu^o)$, $x\in\gN_{T_{\ga}}\cap\gN_{\psi_1}$ :
\begin{eqnarray*}
(id*\omega_{P_A^{it}\zeta, P_A^{it}\zeta'}(\tilde{G})P_A^{it}\Lambda_{\psi_1}(x)=
&=&
(id*\omega_{P_A^{it}\zeta, P_A^{it}\zeta'}(\tilde{G})b(q)^{t/2}\Lambda_{\psi_1}(\tau_t^A(x))\\
&=&
\Lambda_\Phi((\omega_{P_A^{it}\zeta, P_A^{it}\zeta'}\underset{N}{_b*_\alpha}id)\ga(b(q)^{t/2}\tau_t^A(x))\\
&=&
\Lambda_\Phi(\beta(q)^{t/2}\tau_t^A(\omega_{\zeta, \zeta'}\underset{N}{_b*_\alpha}id)\ga(x))\\
&=&
P^{it}\Lambda_\Phi (\omega_{\zeta, \zeta'}\underset{N}{_b*_\alpha}id)\ga(x))\\
&=&
P^{it}(id*\omega_{\zeta, \zeta'})(\tilde{G})\Lambda_{\psi_1}(x)
\end{eqnarray*}
from which we get the formula we were looking for, and which finishes the proof. 
\end{proof}

\subsection{Theorem}
\label{thgalois}
{\it Let $(b, \ga)$ be a Galois action of the measured quantum groupoid $\gG$ on a von Neumann algebra $A$; let $\psi_0$ be a normal semi-finite faithful weight on $A^\ga$, bearing the Galois density property defined in \ref{defdens};  let $\psi_1=\psi_0\circ T_\ga$, and $\tilde{G}$ be its Galois unitary,. Let us suppose that there exists two strongly commuting positive non-singular operator $\delta_A$ and $\lambda_A$, affiliated to $A$, such that the normal semi-finite faithful weight $\phi$ on $A$ defined by $(D\phi:D\psi_1)_t=\lambda_A^{it^2/2}\delta_A^{it}$ (by \cite{V1} 5.1) is invariant under $\ga$; then :
\newline
(i) there exists a one-parametrer group of unitaries $\delta_A^{it}\underset{N}{_b\otimes_\alpha}\delta^{it}$ on $H_{\psi_1}\underset{\nu}{_b\otimes_\alpha}H$, having natural values on elementary tensors, such that, for all $t\in\mathbb{R}$ :
\[\ga(\delta_A^{it})=\delta_A^{it}\underset{N}{_b\otimes_\alpha}\delta^{it}\]
(ii) there exists a one-parametrer group of unitaries $J_{\psi_1}\delta_A^{it}J_{\psi_1}\underset{A^\ga}{_s\otimes_r}\delta_A^{it}$ on $H_{\psi_1}\underset{\psi_0}{_s\otimes_r}H_{\psi_1}$, having natural values on elementary tensors, such that, for all $t\in\mathbb{R}$ :
\[J_{\psi_1}\delta_A^{it}J_{\psi_1}\underset{A^\ga}{_s\otimes_r}\delta_A^{it}=K^{it}=\tilde{G}^*(J\delta^{it}J\underset{N^o}{_\alpha\otimes_b}1)\tilde{G}\]
(iii) we have $\lambda_A=b(q)$, where $q$ is the positive non-singular operator affiliated to $Z(N)$, such that $\lambda=\alpha(q)=\beta(q)$ (\cite{E5}3.8(vi)); the operator $\lambda_A$ is affiliated to $Z(A)$, and $(A,b,\ga,\phi, \psi_0)$ is a Galois system for $\gG$. 
\newline
(iv) we have $\tau^A_t(\delta_A^{is})=\delta_A^{is}$.}

\begin{proof}
By definition, $(D\phi:D\psi_1)_t=\lambda_A^{it^2/2}\delta_A^{it}$, and, therefore, $\ga(\lambda_A^{it^2/2})\ga(\delta_A^{it})=(D\overline{\phi}:D\overline{\psi_1})_t$ where $\overline{\phi}$ (resp. $\overline{\psi_1}$) is the weight on $A\underset{N}{_b*_\alpha}\mathcal L(H)$ given by the bidual weight on the bicrossed-product (which is isomorphic to $A\underset{N}{_b*_\alpha}\mathcal L(H)$) (\cite{E5}, 11.6). As the weight $\phi$ is invariant with respect to $\ga$, the weight $\overline{\phi}$ is equal to another weight $\underline{\phi}$ (\cite{E6},7.7(x)), which is defined by the formula $\frac{d\underline{\phi}}{d\phi^o}=\Delta_\phi^{1/2}\underset{N}{_b\otimes_\alpha}\Delta_{\widehat{\Phi}}^{-1/2}$ (\cite{E6},4.4). On the other hand, using (\cite{E5}, 13.7), and (\cite{E6}, 3.2), we get that $\frac{d\overline{\psi_1}}{d\psi_1^o}=\Delta_{\psi_1}^{1/2}\underset{N}{_b\otimes_\alpha}(\delta\Delta_{\widehat{\Phi}})^{-1/2}$. 
\newline
Finally, we get :
\begin{eqnarray*}
\ga(\lambda_A^{it^2/2})\ga(\delta_A^{it})
&=&
(\Delta_\phi^{it}\underset{N}{_b\otimes_\alpha}\Delta_{\widehat{\Phi}}^{-it})[(D\phi^o:D\psi_1^o)_t\underset{N}{_b\otimes_\alpha}1](\Delta_{\psi_1}^{-it}\underset{N}{_b\otimes_\alpha}(\delta\Delta_{\widehat{\Phi}})^{it})\\
&=&
(D\phi:D\psi_1)_t\underset{N}{_b\otimes_\alpha}\delta^{it}\\
&=&
\lambda_A^{it^2/2}\delta_A^{it}\underset{N}{_b\otimes_\alpha}\delta^{it}
\end{eqnarray*}
from which we get (i), and that $\lambda_A$ is affiliated to $A^\ga$. 
\newline
It is straightforward to get that there exists on $H_{\psi_1}\underset{\psi_0}{_s\otimes_r}H_{\psi_1}$ a one-parametrer group of unitaries $J_{\psi_1}\delta_A^{it}J_{\psi_1}\underset{A^\ga}{_s\otimes_r}\delta_A^{it}$, having natural values on elementary tensors; 
using \ref{G}(i), we get, for any $x\in\gN_{T_\ga}\cap\gN_{\psi_1}$, $\zeta\in D((H_{\psi_1})_b, \nu^o)$ and $(e_i)_{i\in I}$ an orthogonal $(b, \nu^o)$-basis of $H_{\psi_1}$ :
\begin{eqnarray*}
\tilde{G}(J_{\psi_1}\delta_A^{it}J_{\psi_1}\underset{A^\ga}{_s\otimes_r}\delta_A^{it})(\Lambda_{\psi_1}(x)\underset{\psi_0}{_s\otimes_r}\zeta)
&=&
\tilde{G}(\Lambda_{\psi_1}(x\lambda_A^{-t/2}\delta_A^{-it})\underset{\psi_0}{_s\otimes_r}\delta_A^{it}\zeta)\\
&=&
\sum_i\Lambda_\Phi(\omega_{\delta_A^{it}\zeta, e_i}\underset{N}{_b*_\alpha}id)\ga(x\lambda_A^{-t/2}\delta_A^{-it})\underset{\nu^o}{_\alpha\otimes_b}e_i
\end{eqnarray*}
which, using (i) and the fact that $\lambda_A$ is affiliated to $A^\ga$, is equal to :
\begin{eqnarray*}
\sum_i\Lambda_\Phi((\omega_{\lambda_A^{-t/2}\zeta, e_i}\underset{N}{_b*_\alpha}id)\ga(x)\delta^{-it})\underset{\nu^o}{_\alpha\otimes_b}e_i
&=&
\sum_iJ\delta^{it}J\lambda^{t/2}\Lambda_\Phi((\omega_{\lambda_A^{-t/2}\zeta, e_i}\underset{N}{_b*_\alpha}\ga(x))\underset{\nu^o}{_\alpha\otimes_b}e_i\\
&=&
(J\delta^{it}J\lambda^{t/2}\underset{N^o}{_\alpha\otimes_b}1)\tilde{G}(\Lambda_{\psi_1}(x)\underset{\psi_0}{_s\otimes_r}\lambda_A^{-t/2}\zeta)\\
&=&
(J\delta^{it}J\underset{N^o}{_\alpha\otimes_b}1)\tilde{G}(\pi_\ga(1\underset{N}{_b\otimes_\alpha}\lambda^{t/2})\Lambda_{\psi_1}(x)\underset{\psi_0}{_s\otimes_r}\lambda_A^{-t/2}\zeta)
\end{eqnarray*}
from which we get :
\[(J_{\psi_1}\delta_A^{it}J_{\psi_1}\underset{A^\ga}{_s\otimes_r}\delta_A^{it})(\Lambda_{\psi_1}(x)\underset{\psi_0}{_s\otimes_r}\zeta)=\tilde{G}^*(J\delta^{it}J\underset{N^o}{_\alpha\otimes_b}1)\tilde{G}(\pi_\ga(1\underset{N}{_b\otimes_\alpha}\lambda^{t/2})\Lambda_{\psi_1}(x)\underset{\psi_0}{_s\otimes_r}\lambda_A^{-t/2}\zeta)\]
Therefore, the application $Q$ which sends $\Lambda_{\psi_1}(x)\underset{\psi_0}{_s\otimes_r}\zeta$ on $\pi_\ga(1\underset{N}{_b\otimes_\alpha}\lambda^{t/2})\Lambda_{\psi_1}(x)\underset{\psi_0}{_s\otimes_r}\lambda_A^{-t/2}\zeta$ is bounded; as it is clearly positive, by the unicity of polar decomposition, we get (ii), and the fact that $Q=1$; from which one gets that $\lambda_A$ is affiliated to $Z(A^\ga)$, and that $\lambda_A^{t/2}=\pi_\ga(1\underset{N}{_b\otimes_\alpha}\lambda^{t/2})$, and, therefore, that :
\[\lambda_A^{t/2}\underset{N}{_b\otimes_\alpha}1=\ga(\lambda_A^{t/2})=1\underset{N}{_b\otimes_\alpha}\lambda^{t/2}=1\underset{N}{_b\otimes_\alpha}\beta(q^{t/2})=\ga(b(q^{t/2}))\]
from which we get that $\lambda_A=b(q)$. Then, we get, for all $x\in A$ :
\[\ga(\lambda_A^{it}x\lambda_A^{-it})=(1\underset{N}{_b\otimes_\alpha}\lambda^{it})\ga(x)(1\underset{N}{_b\otimes_\alpha}\lambda^{-it})=\ga(x)\]
because $\lambda$ is affiliated to $Z(M)$; so we get that $\lambda_A$ is affiliated to $Z(A)$, and, by \cite{V1} 5.2, that the modular groups of $\sigma^\phi$ and $\sigma^{\psi_1}$ commute; which, thanks to \ref{thK}, gives (iii). 
\newline
We have :
\begin{eqnarray*}
\tau_t^A(\delta_A^{is})
&=&
\pi_\ga(1\underset{N}{_b\otimes_\alpha}\hat{J}\hat{\delta}^{-it}\hat{J})\sigma_t^{\psi_1}(\delta_A^{is})\pi_\ga(1\underset{N}{_b\otimes_\alpha}\hat{J}\hat{\delta}^{-it}\hat{J})\\
&=&
\pi_\ga(1\underset{N}{_b\otimes_\alpha}\hat{J}\hat{\delta}^{-it}\hat{J})\lambda_A^{ist}\delta_A^{is}\pi_\ga(1\underset{N}{_b\otimes_\alpha}\hat{J}\hat{\delta}^{-it}\hat{J})
\end{eqnarray*}
and, therefore, using \ref{thgalois} (i) and (iii) :
\begin{eqnarray*}
\ga(\tau_t^A(\delta_A^{is}))
&=&
(1\underset{N}{_b\otimes_\alpha}\hat{J}\hat{\delta}^{-it}\hat{J})(1\underset{N}{_b\otimes_\alpha}\lambda^{ist})(\delta_A^{is}\underset{N}{_b\otimes_\alpha}\delta^{is})(1\underset{N}{_b\otimes_\alpha}\hat{J}\hat{\delta}^{-it}\hat{J})\\
&=&
\delta_A^{is}\underset{N}{_b\otimes_\alpha}\lambda^{ist}\tau_t\sigma_{-t}^\Phi(\delta^{is})\\
&=&
\delta_A^{is}\underset{N}{_b\otimes_\alpha}\lambda^{ist}\lambda^{-ist}\delta^{is}\\
&=&
\delta_A^{is}\underset{N}{_b\otimes_\alpha}\delta^{is}\\
&=&
\ga(\delta_A^{is})
\end{eqnarray*}
from which we get (iv) and finish the proof. 
\end{proof}
\subsection{Corollary}
\label{corgalois}
{\it Let $(b, \ga)$ be a Galois action of the measured quantum groupoid $\gG$ on a von Neumann algebra $A$; let $\psi_0$ be a normal semi-finite faithful weight on $A^\ga$, bearing the Galois density property defined in \ref{defdens};  let $\psi_1=\psi_0\circ T_\ga$; then, are equivalent :
\newline
(i) there exists a positive non singular operator $\delta_A$ affiliated to $A$ such that, for all $t\in\mathbb{R}$, we have $K^{it}=J_{\psi_1}\delta_A^{it}J_{\psi_1}\underset{A^\ga}{_s\otimes_r}\delta_A^{it}$. 
\newline
(ii) there exists two strongly commuting positive non-singular operator $\delta_A$ and $\lambda_A$, affiliated to $A$, such that the normal semi-finite faithful weight $\phi$ on $A$ defined by $(D\phi:D\psi_1)_t=\lambda_A^{it^2/2}\delta_A^{it}$ (by \cite{V1} 5.1) is invariant under $\ga$.
\newline
(iii) there exists a normal semi-finite faithful weight $\phi$ on $A$, such that $(A, b, \ga, \phi, \psi_0)$ is a Galois system. 
\newline
Then, $\delta_A$ is the modulus of the action $(b, \ga)$, and $\lambda_A=b(q)$, where $q\eta Z(N)$ is such that $\lambda=\alpha(q)=\beta(q)$. }
\begin{proof}
We had obtained in \ref{thK} that (i) implies (iii); in \ref{thgalois}, we had obtained that (ii) implies (i) and (iii); and applying \ref{thgalois} to (iii), we obtain (ii). 
\end{proof}
\subsection{Proposition}
\label{dens}
{\it Let $(A, b, \ga, \phi, \psi_0)$ be a Galois system for the measured quantum groupoid $\gG$; let $\psi_1=\psi_0\circ T_\ga$, $\gT$ be the normal semi-finite faithful weight from $A$ onto $b(N)$ such that $\phi=\nu^o\circ b^{-1}\circ \gT$, and $r$ the canonical injection of $A^\ga$ into $A$. Then :
\newline
(i) the left ideal $\gN_{\psi_1}\cap\gN_{T_\ga}\cap\gN_\phi\cap\gN_{\gT}$ is dense in $A$. 
\newline
(ii) the subspace $\Lambda_{\phi}(\gN_{\psi_1}\cap\gN_{T_\ga}\cap\gN_\phi\cap\gN_{\gT})$ is dense in $H_{\phi}$. 
\newline
(iii) the subspace $D((H_{\psi_1})_b, \nu^o)\cap D(_rH_{\psi_1}, \psi_0)$ is dense in $H_{\psi_1}$, i.e. $\psi_0$ satisfies the Galois density property defined in \ref{defdens}. }

\begin{proof}
Using \cite{V1}, we know that, if $x$ in $A$ is such that $x\delta_A^{1/2}$ is bounded and its closure $\overline{x\delta_A^{1/2}}$ belongs to $\gN_\phi$, then $x$ belongs to $\gN_{\psi_1}$; we can then (and we shall) identify $\Lambda_{\psi_1}(x)$ with $\Lambda_\phi(\overline{x\delta_A^{1/2}})$ and $J_{\psi_1}$ with $\lambda_A^{i/4}J_\phi$. In particular, using the selfadjoint elements of $A$ given by the formula :
\[e_n=\frac{2n^2}{\Gamma(1/2)\Gamma(1/4)}\int_{\mathbb{R}^2}e^{-n^2x^2-n^4y^4}\lambda_A^{ix}\delta_A^{iy}dxdy\]
which are analytic with respect to $\sigma^\phi$ and such that, for any $z\in \mathbb{C}$,  the sequence $\sigma_z^\phi(e_n)$ is bounded and strongly converges to $1$, we get that for any $x\in\gN_\phi$, $x(e_n\delta_A^{1/2})$ belongs to $\gN_{\psi_1}$. 
\newline
Let $\gT$ be the normal faithful semi-finite operator-valued weight from $A$ onto $b(N)$ such that $\phi=\nu^o\circ b^{-1}\circ \gT$. Let us suppose that $x$ is positive in the Tomita algebra $\mathcal T_{\phi, \gT}$ (\cite{E5}, 2.2.1) associated to $\phi$ and $\gT$ (i.e. $x$ belongs to $\gN_\phi\cap\gN_{\gT}$, is analytical with respect to $\sigma^\phi$, and, for all $z\in\mathbb{C}$, $\sigma_z^\phi(x)$ belongs to $\gN_\phi\cap\gN_\phi^*\cap\gN_{\gT}\cap\gN_{\gT}^*$. As in (\cite{L} 5,17), let us define :
\[x_{p,q}=f_p\sqrt\frac{d}{\pi}\int_{-\infty}^{+\infty}e^{-qt^2}\sigma_t^{\psi_1}(x)dt\]
with $f_p=\int_{1/p}^pde_t$, where $\lambda_A=\int_0^\infty tde_t$.
and we get that $x_{p,q}$ belongs to $T_{\phi, \gT}$, is analytical with respect to $\sigma^{\psi_1}$,and that, for all $z\in\mathbb{C}$, $\sigma_z^{\psi_1}(x_{p,q})$ belongs to $T_{\phi, \gT}$. As $\sigma_t^{\psi_1}\sigma_{-t}^\phi=Ad\delta_A^{it}$, we get that, for all $z$ in $\mathbb{C}$, $\delta_A^{iz}x_{p,q}\delta_A^{-iz}$ belongs to $T_{\phi, \gT}$; in particular, $\delta_A^{-1/2}x_{p,q}\delta_A^{1/2}$ belongs to $T_{\phi, \gT}$ and $e_nx_{p,q}\delta_A^{1/2}=(\delta_A^{1/2}e_n)\delta_A^{-1/2}x_{p,q}\delta_A^{1/2}$ belongs to $\gN_\phi\cap\gN_\gT$. We prove this way that the set $T^{\psi_1}_{\phi, \gT}$ of elements $x$ in $\gN\phi\cap\gN_\phi^*\cap\gN_\gT\cap\gN_\gT^*$ which are analytic with respect, both, of $\sigma^\phi$ and $\sigma^{\psi_1}$, and such that $x\delta_A^{1/2}$ is bounded and belongs to $\gN_\phi\cap\gN_\gT$ is weakly dense in $A$, and its image under $\Lambda_\phi$ is a dense subspace of $H_\phi$.
\newline
Let us take $x\in T^{\psi_1}_{\phi, \gT}$; using the fact that $\psi_1$ is a $\delta$-invariant weight with respect to $\ga$ (\ref{crossed}), we get (where $a$ is the representation of $N$ on $H_{\psi_1}$ given by $a(n)=J_{\psi_1}b(n^*)J_{\psi_1}$) that :
\begin{eqnarray*}
(T_\ga\underset{N}{_b*_\alpha}id)\ga(x^*x)
&=&
\delta^{1/2}\beta(<\Lambda_{\psi_1}(x), \Lambda_{\psi_1}(x)>_{a, \nu})\delta^{1/2}\\
&=&
\delta^{1/2}\beta(<J_{\psi_1}\Lambda_{\psi_1}(x), J_{\psi_1}\Lambda_{\psi_1}(x)>_{b, \nu^o})\delta^{1/2}\\
&=&
\delta^{1/2}\beta(<\lambda_A^{i/4}J_{\phi}\Lambda_{\phi}(x\delta_A^{1/2}), \lambda_A^{i/4}J_{\phi}\Lambda_{\phi}(x\delta_A^{1/2})>_{b, \nu^o})\delta^{1/2}
\end{eqnarray*}
and, therefore, if $\eta$ belongs to $D(_\alpha H, \nu)\cap \mathcal D(\delta^{1/2})$, we get that $(id\underset{N}{_b*_\alpha}\omega_\eta)\ga(x^*x)$ belongs to $\gM_{T_\ga}^+$ (and to $\gM_{\psi_1}^+$ by similar arguments).
\newline
Using now the fact that $\phi$ is invariant with respect to $\ga$, we get (where $a$ means here the representation of $N$ on $H_\phi$ given by $a(n)=J_\phi b(n^*)J_\phi$) that :
\begin{eqnarray*}
(\gT\underset{N}{_b*_\alpha}id)\ga(x^*x)
&=&
\beta(<\Lambda_\phi(x), \Lambda_\phi(x)>_{a, \nu})\\
&=&
\beta(<J_\phi \Lambda_\phi(x), J_\phi \Lambda_\phi(x)>_{b, \nu^o})
\end{eqnarray*}
and we get that $(id\underset{N}{_b*_\alpha}\omega_\eta)\ga(x^*x)$ belongs also to $\gM_{\gT}^+\cap\gM_\phi^+$. So, we get that, for any $x\in T^{\psi_1}_{\phi, \gT}$, $\xi\in D(_\alpha H, \nu)$, $\eta \in D(_\alpha H, \nu)\cap \mathcal D(\delta^{1/2})$, the operator $(id\underset{N}{_b*_\alpha}\omega_{\eta, \xi})\ga(x)$ belongs to $\gN_{T_\ga}\cap\gN_{\psi_1}\cap\gN_{\gT}\cap\gN_\phi$. 
\newline
So, we get that the weak closure of $\gN_{T_\ga}\cap\gN_{\psi_1}\cap\gN_{\gT}\cap\gN_\phi$ contains all elements of the form $(id\underset{N}{_b*_\alpha}\omega_{\eta, \xi}\ga(x)$ for any $\xi$, $\eta$ in $D(_\alpha H, \nu)$ and $x\in A$; using now \cite{E5}, 11.5 (ii), we get (i). 
\newline
Let us suppose now that $\zeta\in H_\phi$ is orthogonal to $\Lambda_\phi(\gN_{T_\ga}\cap\gN_{\psi_1}\cap\gN_{\gT}\cap\gN_\phi)$. Using \ref{action} and (\cite{E6} 7.7), we get that :
\[(V'_\phi(\Lambda_\phi(x)\underset{N^o}{_a\otimes_\beta}\eta)|\zeta\underset{N}{_b\otimes_\alpha}\xi)=0\]
for all $x\in T^{\psi_1}_{\phi, \gT}$, $\eta\in D{_\alpha H, \nu})\cap \mathcal D(\delta^{1/2})$, $\xi\in D(_\alpha, \nu)$. As $\Lambda_\phi(T^{\psi_1}_{\phi, \gT})$ is dense in $H_\phi$ and $V'_\phi$ is a unitary, we get that $\zeta\underset{N}{_b\otimes_\alpha}\xi)=0$ for all $\xi\in D(_\alpha H, \nu)$, and, therefore, that $\zeta=0$; which is (ii). 
\newline
We know that $J_\phi\Lambda_\phi(\gN_\phi\cap \gN_\gT)\subset D((H_\phi)_b, \nu^o)$ and that $J_{\psi_1}\Lambda_{\psi_1}(\gN_{\psi_1}\cap\gN_{T_\ga})\subset D(_r H_{\psi_1}, \psi^o)$. As the canonical isomorphism between $H_\phi$ and $H_{\psi_1}$ exchanges the representations of $A$ (and, therefore, of $b(N)$), and sends $J_\phi$ on $J_{\psi_1}\lambda_A^{i/4}$, we get that :
\[J_{\psi_1}\Lambda_{\psi_1}(\gN_{\psi_1}\cap\gN_{T_\ga}\cap\gN_\phi\cap\gN_{\gT})\subset D((H_{\psi_1})_b, \nu^o)\cap D(_rH_{\psi_1}, \psi_0)\]
from which we get (iii). 
\end{proof}

\section{Through the looking-glass}
\label{TLG}
In this chapter, we use the reflection technic introduced by De Commer in \cite{DC1}; if we start from a Galois action $(b, \ga)$ of a measured quantum groupoid $\gG$ on a von Neumann algebra $A$, we obtain a co-involutive Hopf-bimodule which has $A^\ga$ as basis (\ref{coproduct}). If we start from a Galois system $(A, b, \ga, \phi, \psi_0)$, we then construct a left-invariant operator-valued weight on this co-involutive Hopf-bimodule, and obtain this way, "throught the Galois system", another measured quantum groupoid. More precisely, we get in fact two measured quantum groupoids, one with basis $A^\ga$, called the reflected measured quantum groupoid of $\gG$, throught the Galois system, whose underlying von Neumann algebra acts on $H_{\psi_1}$ (\ref{thP2}), and another one which will be a von Neumann algebra acting on $H_{\psi_1}\oplus H$, with the basis $A^\ga\oplus N$, and will be called the linking measured quantum groupoid, between the two preceeding ones (\ref{dual}).

\subsection{Notations}
\label{notTLG}
Let $(b, \ga)$ be a Galois action of a measured quantum groupoid $\gG$ on a von Neumann algebra $A$; let $\psi_0$ be a normal semi-finite faithful weight on $A^\ga$ bearing the density condition. Let us now consider the von Neumann algebra $\widetilde{N}=A^\ga\oplus N$, equipped with a normal faithful semi-finite weight $\psi_0\oplus\nu$, its representation $\widetilde{\alpha}=r\oplus \alpha$, and its anti-representation $\widetilde{\hat{\beta}}=s\oplus \hat{\beta}$ on the Hilbert space $H_{\psi_1}\oplus H$. For any $m'\in\widehat{M}'$, let us write $\mu(m')=\pi_\ga(1\underset{N}{_b\otimes_\alpha}m')$, and consider the operator $\varpi(m')=\mu(m')\oplus m'$ on $H_{\psi_1}\oplus H$; we define this way a normal faithful representation $\varpi$ of $\widehat{M}'$ on $H_{\psi_1}\oplus H$, and a faithful normal anti-representation $\varpi^o$ of $\widehat{M}$ given, for any $m\in\widehat{M}$ by :
\[\varpi^o(m)=\mu(\hat{J}m^*\hat{J})\oplus\hat{J}m^*\hat{J}\]
We shall denote by $\widehat{Q}$ the commutant $\varpi^o(\widehat{M})'$. We shall use matrix notation for elements in $\widehat{Q}$, or, more generally, in $\mathcal L(H_{\psi_1}\oplus H)$. In particular, we shall write 
\[\widehat{Q}=\left(\begin{matrix}{\widehat{P}}&{\widehat{I}}\cr{\widehat{I}^*}&{\widehat{M}}\end{matrix}\right)\]
where $\widehat{P}=\pi_\ga(1\underset{N}{_b\otimes_\alpha}\widehat{M}')'$, and $\widehat{I}$ is the following closed linear set of intertwinners :
\[\widehat{I}=\{X\in \mathcal L(H, H_{\psi_1}), Xm=\pi_\ga(1\underset{N}{_b\otimes_\alpha} m)X, \forall m\in\widehat{M}'\}\]
We see that $r(A^\ga)\subset \widehat{P}$ and $s(A^\ga)\subset \widehat{P}$, and, therefore, $\widetilde{\alpha}(\tilde{N})\subset \widehat{Q}$ and $\widetilde{\hat{\beta}}(\tilde{N})\subset \widehat{Q}$. 
\newline
Let us remark that, for any $\xi\in D((H_{\psi_1})_{\mu^o}, \widehat{\Phi}^o)$, the operator $R^{\mu^o, \widehat{\Phi}^o}(\xi)$ belongs to $\widehat{I}$ (which implies that $\widehat{I}$ is not reduced to $\{0\}$). 
Using \ref{lemR}, we get that, if $X\in\widehat{I}$, we have, for any $m\in\widehat{M}'$ :
\[J_{\psi_1}XJm=J_{\psi_1}X\hat{R}^c(m^*)J=J_{\psi_1}\pi_\ga(1\underset{N}{_b\otimes_\alpha}\hat{R}^c(m^*))XJ=\pi_\ga(1\underset{N}{_b\otimes_\alpha}m)J_{\psi_1}XJ\]
from which we get that $J_{\psi_1}XJ$ belongs to $\widehat{I}$. 
\newline
In particular, for any $n\in N$, we get that $X\in\widehat{I}$ satisfies $b(n)X=X\beta(n)$, and we can define $1_H\underset{N^o}{_\alpha\otimes_\beta}X$ from $H\underset{N}{_\alpha\otimes_\beta}H$ to $H\underset{N}{_\alpha\otimes_b}H_{\psi_1}$. Applying this result to $J_{\psi_1}XJ$, we get that :
\[X\hat{\alpha}(n)=XJ\beta(n^*)J=J_{\psi_1}b(n^*)J_{\psi_1}X=a(n)X\]
and we can define $X\underset{N}{_{\hat{\alpha}}\otimes_\beta}1_H$ from $H\underset{\nu^o}{_{\hat{\alpha}}\otimes_\beta}H$ to $H_{\psi_1}\underset{\nu^o}{_a\otimes_\beta}H$, and $X\underset{N}{_\beta\otimes_\alpha}1_H$ from $H\underset{\nu}{_\beta\otimes_\alpha}H$ to $H_{\psi_1}\underset{\nu}{_b\otimes_\alpha}H$. 
\newline
Using then \ref{pia}, we get that :
\[V_{\psi_1}^*(X\underset{N}{_\beta\otimes_\alpha}1_H)=(X\underset{N}{_{\hat{\alpha}}\otimes_\beta}1_H)(\sigma(\hat{J}\otimes\hat{J})W(\hat{J}\otimes\hat{J})\sigma)\]
\[V_{\psi_1}(X\underset{N}{_{\hat{\alpha}}\otimes_\beta}1_H)=(X\underset{N}{_\beta\otimes_\alpha}1_H)(\sigma(\hat{J}\otimes\hat{J})W^*(\hat{J}\otimes\hat{J})\sigma)\]
Let us denote $e_1=1_{A^\ga}\in \widetilde{N}$, and $e_2=1_N\in\widetilde{N}$; we get that $\widetilde{\alpha}(e_1)=\widetilde{\hat{\beta}}(e_1)=P_{H_{\psi_1}}\in \widehat{P}$, and that $\widetilde{\alpha}(e_2)=\widetilde{\hat{\beta}}(e_2)=P_H\in \widehat{M}$, and $\widehat{P}=\widehat{Q}_{\widetilde{\alpha}(e_1)}$, $\widehat{M}=\widehat{Q}_{\widetilde{\alpha}(e_2)}$. We can verify that $s(A^\ga)\subset\widehat{P}$ and $r(A^\ga)\subset \widehat{P}$. 
\newline
Let us describe now the fiber product $\widehat{Q}\underset{\widetilde{N}}{_{\widetilde{\hat{\beta}}}*_{\widetilde{\alpha}}}\widehat{Q}$. This von Neumann algebra is defined on the Hilbert space :
\[(H_{\psi_1}\oplus H)\underset{\psi_0\oplus\nu}{_{\widetilde{\hat{\beta}}}\otimes_{\widetilde{\alpha}}}(H_{\psi_1}\oplus H)=(H_{\psi_1}\underset{\psi_0}{_s\otimes_r}H_{\psi_1})\oplus (H\underset{\nu}{_{\hat{\beta}}\otimes_\alpha}H)\]
where this direct sum decomposition can be seen with the projections :
\[P_{H_{\psi_1}\underset{\psi_0}{_s\otimes_r}H_{\psi_1}}=\widetilde{\alpha}(e_1)\underset{\psi_0\oplus\nu}{_{\widetilde{\hat{\beta}}}\otimes_{\widetilde{\alpha}}}1=1\underset{\psi_0\oplus\nu}{_{\widetilde{\hat{\beta}}}\otimes_{\widetilde{\alpha}}}\widetilde{\alpha}(e_1)=\widetilde{\alpha}(e_1)\underset{\psi_0\oplus\nu}{_{\widetilde{\hat{\beta}}}\otimes_{\widetilde{\alpha}}}\widetilde{\alpha}(e_1)\]
\[P_{H\underset{\nu}{_{\hat{\beta}}\otimes_\alpha}H}=\widetilde{\alpha}(e_2)\underset{\psi_0\oplus\nu}{_{\widetilde{\hat{\beta}}}\otimes_{\widetilde{\alpha}}}1=1\underset{\psi_0\oplus\nu}{_{\widetilde{\hat{\beta}}}\otimes_{\widetilde{\alpha}}}\widetilde{\alpha}(e_2)=\widetilde{\alpha}(e_2)\underset{\psi_0\oplus\nu}{_{\widetilde{\hat{\beta}}}\otimes_{\widetilde{\alpha}}}\widetilde{\alpha}(e_2)\]
So, we can also use matrix notations for elements in $\widehat{Q}\underset{\widetilde{N}}{_{\widetilde{\hat{\beta}}}*_{\widetilde{\alpha}}}\widehat{Q}$, or, more generally, in $\mathcal L((H_{\psi_1}\oplus H)\underset{\psi_0\oplus\nu}{_{\widetilde{\hat{\beta}}}\otimes_{\widetilde{\alpha}}}(H_{\psi_1}\oplus H))$. In particular, we shall get :
\[\widehat{Q}\underset{\widetilde{N}}{_{\widetilde{\hat{\beta}}}*_{\widetilde{\alpha}}}\widehat{Q}=
\left(\begin{matrix}{\widehat{P}\underset{A^\ga}{_s*_r}\widehat{P}}&{\widehat{I}\underset{N}{_{\hat{\beta}}*_\alpha}\widehat{I}}\cr{(\widehat{I}\underset{N}{_{\hat{\beta}}*_\alpha}\widehat{I})^*}&{\widehat{M}\underset{N}{_{\hat{\beta}}*_\alpha}\widehat{M}}\end{matrix}\right)\]
where $\widehat{I}\underset{N}{_{\hat{\beta}}*_\alpha}\widehat{I}$ is the closed set of intertwinners :
\begin{multline*}
\widehat{I}\underset{N}{_{\hat{\beta}}*_\alpha}\widehat{I}=\{Y\in\mathcal L(H\underset{\nu}{_{\hat{\beta}}\otimes_\alpha}H, H_{\psi_1}\underset{\psi_0}{_s\otimes_r}H_{\psi_1}), Y(m_1\underset{N}{_{\hat{\beta}}\otimes_\alpha}m_2)\\=(\mu(m_1)\underset{A^\ga}{_s\otimes_r}\mu(m_2))Y, \forall m_1, m_2\in \widehat{M}'\}
\end{multline*}

\subsection{Lemma}
\label{lemTLG}
{\it Let's use the notations of \ref{notTLG}, then, we have, for all $X\in\widehat{I}$ and $m\in\widehat{M}'$ :
\[\tilde{G}^*(1\underset{N^o}{_\alpha\otimes_\beta}X)\widehat{W}(m\underset{N}{_{\hat{\beta}}\otimes_\alpha}1)=[\pi_\ga(1\underset{N}{_b\otimes_\alpha}m]\underset{A^\ga}{_s\otimes_r}1)\tilde{G}^*(1\underset{N^o}{_\alpha\otimes_\beta}X)\widehat{W}\]
\[\tilde{G}^*(1\underset{N^o}{_\alpha\otimes_\beta}X)\widehat{W}(1\underset{N}{_{\hat{\beta}}\otimes_\alpha}m)=[1\underset{A^\ga}{_s\otimes_r}\pi_\ga(1\underset{N}{_b\otimes_\alpha}m)]\tilde{G}^*(1\underset{N^o}{_\alpha\otimes_\beta}X)\widehat{W}\]
and, therefore, $\tilde{G}^*(1\underset{N^o}{_\alpha\otimes_\beta}X)\widehat{W}$ belongs to $\widehat{I}\underset{N}{_{\hat{\beta}}*_\alpha}\widehat{I}$. We obtain also that, for any $v\in D((H_{\psi_1})_s, \psi_0^o)$ and $\xi\in D(H_{\hat{\beta}}, \nu^o)$, the operator $(\omega_{v, \xi}*id)[\tilde{G}^*(1\underset{N^o}{_\alpha\otimes_\beta}X)\widehat{W}]$ belongs to $\widehat{I}$. }
\begin{proof}
Using (\cite{E5}, 3.6(ii)) applied to $\widehat{\gG}$, and \ref{G}(iv), we get the first formula. 
\newline
Using \ref{lemG*}, we get that :
\[[1\underset{A^\ga}{_s\otimes_r}\pi_\ga(1\underset{N}{_b\otimes_\alpha}m)]G^*=
G^*V_{\psi_1}[1\underset{N^o}{_a\otimes_\beta}\hat{J}\hat{R}^c(m^*)\hat{J}]V_{\psi_1}^*\]
and, therefore, using \ref{lemG*}, \ref{notTLG} and \ref{lemW2} :
\begin{eqnarray*}
[1\underset{A^\ga}{_s\otimes_r}\pi_\ga(1\underset{N}{_b\otimes_\alpha}m)]\tilde{G}^*(1\underset{N^o}{_\alpha\otimes_\beta}X)
&=&
G^*V_{\psi_1}[1\underset{N^o}{_a\otimes_\beta}\hat{J}\hat{R}^c(m^*)\hat{J}]V_{\psi_1}^*(X\underset{N}{_b\otimes_\alpha}1)\sigma_{\nu^o}\\
&=&
G^*V_{\psi_1}[1\underset{N^o}{_a\otimes_\beta}\hat{J}\hat{R}^c(m^*)\hat{J}](X\underset{N^o}{_{\hat{\alpha}}\otimes_\beta}1_H)(\sigma W^o\sigma)\sigma\\
&=&
G^*V_{\psi_1}(X\underset{N^o}{_{\hat{\alpha}}\otimes_\beta}1_H)[1\underset{N^o}{_{\hat{\alpha}}\otimes_\beta}\hat{J}\hat{R}^c(m^*)\hat{J}](\sigma W^o\sigma)\sigma\\
&=&
G^*((X\underset{N}{_\beta\otimes_\alpha}1_H)(\sigma W^o\sigma)^*[1\underset{N^o}{_a\otimes_\beta}\hat{J}\hat{R}^c(m^*)\hat{J}](\sigma W^o\sigma)\sigma\\
&=&
\tilde{G}^*(1\underset{N^o}{_\alpha\otimes_\beta}X)W^{o*}(\hat{J}\hat{R}^c(m^*)\hat{J}\underset{N}{_\beta\otimes_{\hat{\alpha}}}1)W^o\\
&=&
\tilde{G}^*(1\underset{N^o}{_\alpha\otimes_\beta}X)\widehat{W}(1\underset{N}{_{\hat{\beta}}\otimes_\alpha}m)\widehat{W}^*
\end{eqnarray*}
from which we get the second formula and finish the proof. \end{proof}

\subsection{Proposition}
\label{coinvo}
{\it With the notations of \ref{notTLG}, for any $X\in\widehat{I}$, we have :}
\[\tilde{G}^*(1\underset{N^o}{_\alpha\otimes_\beta}J_{\psi_1}XJ)\widehat{W}=\varsigma_N[(J_{\psi_1}\underset{A^\ga}{_s\otimes_r}J_{\psi_1})\tilde{G}^*(1\underset{N^o}{_\alpha\otimes_\beta}X)\widehat{W}(J\underset{N}{_{\hat{\beta}}\otimes_\alpha}J)]\]

\begin{proof}
Using \ref{lemG} (i), we get that $(J_{\psi_1}\underset{A^{\ga o}}{_r\otimes_s}J_{\psi_1})\sigma_{\psi_0}\tilde{G}^*(1\underset{N^o}{_\alpha\otimes_\beta}X)\widehat{W}\sigma_{\nu^o}(J\underset{N}{_{\hat{\beta}}\otimes_\alpha}J)$
is equal to :
\begin{multline*}
\tilde{G}^*(\hat{J}\underset{N}{_\beta\otimes_a}J_{\psi_1})\sigma_{\nu^o}V_{\psi_1}^*\sigma_{\nu^o}(1\underset{N^o}{_\alpha\otimes_\beta}X)\widehat{W}\sigma_{\nu^o}(J\underset{N}{_{\hat{\beta}}\otimes_\alpha}J)\\
=
\tilde{G}^*(\hat{J}\underset{N}{_\beta\otimes_a}J_{\psi_1})\sigma_{\nu^o}V_{\psi_1}^*(X\underset{N}{_\beta\otimes_\alpha}1)\sigma_{\nu^o}\widehat{W}\sigma_{\nu^o}(J\underset{N}{_{\hat{\beta}}\otimes_\alpha}J)
\end{multline*}
Using now \ref{notTLG}, we get it is equal to :
\[\tilde{G}^*(\hat{J}\underset{N}{_\beta\otimes_a}J_{\psi_1})\sigma_{\nu^o}(X\underset{N^o}{_{\hat{\alpha}}\otimes_\beta}1)(\sigma W^o\sigma)\sigma_{\nu^o}\widehat{W}\sigma_{\nu^o}(J\underset{N}{_{\hat{\beta}}\otimes_\alpha}J)\]
which, by \ref{lemG}(ii), is equal to :
\[\tilde{G}^*(1\underset{N^o}{_\alpha\otimes_\beta}J_{\psi_1}XJ)\widehat{W}\]
\end{proof}

\subsection{Proposition}
\label{coproduct}
{\it (i) With the notations of \ref{notTLG}, for any element $\left(\begin{matrix}{A}&{X}\cr{Y^*}&{m}\end{matrix}\right)$ in $\widehat{Q}$, let us write :
\[\Gamma_{\widehat{Q}}(\left(\begin{matrix}{A}&{X}\cr{Y^*}&{m}\end{matrix}\right))=
\left(\begin{matrix}{\tilde{G}^*(1\underset{N^o}{_\alpha\otimes_b}A)\tilde{G}}&{\tilde{G}^*(1\underset{N^o}{_\alpha\otimes_\beta}X)\widehat{W}}\cr{\widehat{W}^*(1\underset{N^o}{_\alpha\otimes_b}Y^*)\tilde{G}}&{\widehat{\Gamma}(m)}\end{matrix}\right)\]
Then, we define an application $\Gamma_{\widehat{Q}}$ from $\widehat{Q}$ into $\widehat{Q}\underset{\tilde{N}}{_{\widetilde{\hat{\beta}}}*_{\widetilde{\alpha}}}\widehat{Q}$ which is a coproduct. So, $(\tilde{N}, \widehat{Q}, \widetilde{\alpha}, \widetilde{\hat{\beta}}, \Gamma_{\widetilde{Q}})$ is a Hopf-bimodule. 
\newline
(ii) Let us write :
\[R_{\widehat{Q}}(\left(\begin{matrix}{A}&{X}\cr{Y^*}&{m}\end{matrix}\right))=\left(\begin{matrix}{J_{\psi_1}A^*J_{\psi_1}}&{J_{\psi_1}YJ}\cr{JX^*J_{\psi_1}}&{Jm^*J}\end{matrix}\right)\]
Then, we define an involutive $*$-anti isomorphism $R_{\widehat{Q}}$ of $\widehat{Q}$, which a co-involution for the coproduct $\Gamma_{\widehat{Q}}$. 
\newline
(iii) For any $A\in\widehat{P}$, let us write $\Gamma_{\widehat{P}}(A)=\tilde{G}^*(1\underset{N^o}{_\alpha\otimes_b}A)\tilde{G}$, and $R_{\widehat{P}}(A)=J_{\psi_1}A^*J_{\psi_1}$; then $(A^\ga, \widehat{P}, r, s, \Gamma_{\widehat{P}})$ is a Hopf-bimodule, and $R_{\widehat{P}}$ is a co-involution for $\Gamma_{\widehat{P}}$. }

\begin{proof}
We had got in \ref{lemTLG} that $\tilde{G}^*(1\underset{N}{_\alpha\otimes_\beta}X)\widehat{W}$ belongs to $\widehat{I}\underset{N}{_{\hat{\beta}}*_\alpha}\widehat{I}$; so, for any $\xi$, $\eta$ in $D(_\mu H_{\psi_1}, \widehat{\Phi}')$, we get that $\tilde{G}^*(1\underset{N}{_\alpha\otimes_b}\theta^{\mu,\widehat{\Phi}'}(\xi, \eta))\tilde{G}$ commutes with $\mu(\widehat{M}')\underset{A^\ga}{_s\otimes_r}\mu(\widehat{M}')$, and, therefore, belongs to $\widehat{P}\underset{A^\ga}{_s*_r}\widehat{P}$; by continuity and density, this remains true for any $A$ in $\widehat{P}$. So, we had got that $\Gamma_{\widehat{Q}}$ is an injective $*$-homomorphism from $\widehat{Q}$ into $\widehat{Q}\underset{\tilde{N}}{_{\widetilde{\hat{\beta}}}*_{\widetilde{\alpha}}}\widehat{Q}$. The fact that it is a coassociative coproduct is given by \ref{pent}, which gives (i). 
\newline
We have seen in \ref{notTLG} that $J_{\psi_1}YJ$ belongs to $\widehat{I}$; therefore, for any $\xi$, $\eta$ in $D(_\mu H_{\psi_1}, \widehat{\Phi}')$, we get that $J_{\psi_1}\theta^{\mu, \widehat{\Phi}'}(\xi, \eta)J_{\psi_1}$ belongs to $\widehat{P}$, and, by density, that remains true for any $A$ in $\widehat{P}$. The fact that we obtain a co-involution is given by  \ref{coinvo}, which gives (ii). As $\widehat{P}=\widehat{Q}_{\widetilde{\alpha}(e_1)}$, we easily get (iii). \end{proof}

\subsection{Proposition}
\label{deltahat}
{\it Let $(A, b, \ga, \phi, \psi_0)$ be a Galois system for $\gG$; let $\psi_1=\psi_0\circ T_\ga$; let $\delta_A$ be the modulus introduced in \ref{defGalois}, and $P_A$ be the generator of the one-parameter group of unitaries introduced in \ref{tau}. Then :
\newline
(i) there exists a one-parameter group of unitaries $\widehat{\Delta}_A^{it}=P_A^{it}J_{\psi_1}\delta_A^{it}J_{\psi_1}$; 
\newline
(ii), we have, for all $t\in\mathbb{R}$ and $m\in\widehat{M}'$ :
\[\pi_\ga(1\underset{N}{_b\otimes_\alpha}\sigma_t^{\widehat{\Phi}^c}(m))=\widehat{\Delta}_A^{-it}\pi_A(1\underset{N}{_b\otimes_\alpha}m)\widehat{\Delta}_A^{it}\]
and, in particular, for any $n\in N$ :
\[\widehat{\Delta}_A^{-it}b(n)\widehat{\Delta}_A^{it}=b(\sigma_{-t}^\nu(n))\]
(iii) we have :}
\[(\Delta_{\Phi}^{it}\underset{N}{_\alpha\otimes_b}\widehat{\Delta}_A^{it})\tilde{G}=\tilde{G}(\Delta_{\psi_1}^{it}\underset{A^\ga}{_s\otimes_r}\widehat{\Delta}_A^{it})\]
\[(\Delta_{\widehat{\Phi}}^{it}\underset{N}{_\alpha\otimes_b}P_A^{it})\tilde{G}=\tilde{G}(\widehat{\Delta}_A^{it}\underset{A^\ga}{_s\otimes_r}P_A^{it}\delta_A^{it})\]

\begin{proof}
Using \ref{tau}(vi), we get that $P_A^{it}$ commutes with $J_{\psi_1}$; using \ref{tau}(v), we get that $P_A^{it}$ commutes with $\delta_A^{is}$; so, we get (i). 
\newline
Let now $x\in\gN_{\psi_1}$, and let be $\xi$ in $D_\alpha H, \nu)$ and $\eta$ in $D(_\alpha H, \nu)\cap\mathcal D(\delta^{1/2})$, such that $\delta^{1/2}\eta$ belongs to $D(H_\beta, \nu^o)$. We have then, using \ref{tau}(v), \cite{E5}, 8.4(iii), \ref{thgalois}(i) and (iii), and again \ref{tau}(v) and \cite{E5}, 8.4(iii) :
\begin{eqnarray*}
(id*\omega_{\eta, \xi})(V_{\psi_1})\widehat{\Delta}_A^{-it}\Lambda_{\psi_1}(x)
&=&
(id*\omega_{\eta, \xi})(V_{\psi_1})\Lambda_{\psi_1}(\lambda_A^{-t/2}\tau^A_t(x)\delta_A^{-it})\\
&=&
\Lambda_{\psi_1}[(id\underset{N}{_b*_\alpha}\omega_{\delta^{1/2}\eta, \xi})\ga(\lambda_A^{-t/2}\tau^A_t(x)\delta_A^{-it})]\\
&=&
\Lambda_{\psi_1}[\lambda_A^{-t/2}\tau^A_t(id\underset{N}{_b*_\alpha}\omega_{P^{-it}\delta^{-it}\delta^{1/2}\eta, P^{-it}\xi})\ga(x)\delta_A^{-it}]\\
&=&
\widehat{\Delta_A}^{-it}\Lambda_{\psi_1}((id\underset{N}{_b*_\alpha}\omega_{P^{-it}\delta^{-it}\delta^{1/2}\eta, P^{-it}\xi})\ga(x))\\
&=&
\widehat{\Delta_A}^{-it}(id*\omega_{P^{-it}\delta^{-it}\eta, P^{-it}\xi})(V_{\psi_1})\Lambda_{\psi_1}(x)
\end{eqnarray*}
Therefore, using \ref{thintegrable} and \ref{lemW} applied to $\gG^o$, we get that 
\[\widehat{\Delta}_A^{-it}\pi_\ga(1\underset{N}{_b\otimes_\alpha}(\omega_{\xi, \eta}*id)[(\hat{J}\underset{N}{_\beta\otimes_\alpha}\hat{J})W^*(\hat{J}\underset{N}{_\beta\otimes_{\hat{\alpha}}}\hat{J})])\widehat{\Delta}_A^{it}\]
is equal to $\pi_\ga(1\underset{N}{_b\otimes_\alpha}\sigma_t^{\widehat{\Phi}^c}(\omega_{\xi, \eta}*id)[(\hat{J}\underset{N}{_\beta\otimes_\alpha}\hat{J})W^*(\hat{J}\underset{N}{_\beta\otimes_{\hat{\alpha}}}\hat{J})])$
and, by density and continuity, we get (ii). 
\newline
Using \ref{G} and \ref{tau}(iii) and (ii), we get, where $x$ is in $\gN_{\psi_1}$, $\zeta$ in $D((H_{\psi_1})_b, \nu^o)$ and  $(e_i)_{i\in I}$ is an orthogonal $(b, \nu^o)$-basis of $H_{\psi_1}$ :
\begin{eqnarray*}
\tilde{G}(\Delta_{\psi_1}^{it}\Lambda_{\psi_1}(x)\underset{N}{_\alpha\otimes_b}\widehat{\Delta}_A^{it}\zeta)
&=& 
\sum_i\Lambda_\Phi[(\omega_{\widehat{\Delta}_A^{it}\zeta, e_i}\underset{N}{_b*_\alpha}id)\ga(\sigma_t^{\psi_1}(x))]\underset{\nu^o}{_\alpha\otimes_b}e_i\\
&=&
\sum_i\Lambda_\Phi[(\omega_{\widehat{\Delta}_A^{it}\zeta, e_i}\underset{N}{_b*_\alpha}id)(\tau^A_t\underset{N}{_b*_\alpha}\sigma_t^\Phi)\ga(x)]\underset{\nu^o}{_\alpha\otimes_b}e_i\\
&=&
\sum_i\Delta_\Phi^{it}\Lambda_\Phi[(\omega_{P_A^{it}\widehat{\Delta}_A^{it}\zeta, P_A^{-it}e_i}\underset{N}{_b*_\alpha}\ga(x)]\underset{\nu^o}{_\alpha\otimes_b}e_i\\
&=&
\sum_i\Delta_\Phi^{it}\Lambda_\Phi[(\omega_{J_{\psi_1}\delta_A^{it}J_{\psi_1}\zeta, P_A^{-it}e_i}\underset{N}{_b*_\alpha}\ga(x)]\underset{\nu^o}{_\alpha\otimes_b}e_i\\
&=&
\sum_i\Delta_\Phi^{it}\Lambda_\Phi[(\omega_{\zeta, P_A^{-it}J_{\psi_1}\delta_A^{it}J_{\psi_1}e_i}\underset{N}{_b*_\alpha}\ga(x)]\underset{\nu^o}{_\alpha\otimes_b}e_i
\end{eqnarray*}
and, using the fact that $(P_A^{-it}J_{\psi_1}\delta_A^{it}J_{\psi_1}e_i)_{i\in I}$ is another orthogonal $(b, \nu^o)$-basis of $H_{\psi_1}$, and that the sum doesnot depend on the choice of the basis, we get it is equal to :
\[\sum_i\Delta_\Phi^{it}\Lambda_\Phi[(\omega_{\zeta, e_i}\underset{N}{_b*_\alpha}\ga(x)]\underset{\nu^o}{_\alpha\otimes_b}P_A^{-it}J_{\psi_1}\delta_A^{it}J_{\psi_1}e_i
=(\Delta_{\psi_1}^{it}\underset{A^\ga}{_s\otimes_r}\widehat{\Delta}_A^{it})\tilde{G}(\Lambda_{\psi_1}(x)\underset{N}{_\alpha\otimes_b}\zeta)\]
which gives the first formula of (iii). 
\newline
Finally, we have, using similar arguments :
\begin{eqnarray*}
\tilde{G}(\widehat{\Delta}_A^{it}\underset{A^\ga}{_s\otimes_r}P_A^{it}\delta_A^{it})(\Lambda_{\psi_1}(x)\underset{\psi_0}{_s\otimes_r}\zeta)
&=&
\tilde{G}(P_A^{it}J_{\psi_1}\delta_A^{it}J_{\psi_1}\Lambda_{\psi_1}(x)\underset{\psi_0}{_s\otimes_r}P_A^{it}\delta_A^{it}\zeta)\\
&=&
\tilde{G}(\lambda_A^{-t/2}P_A^{it}\Lambda_{\psi_1}(x\delta_A^{-it})\underset{\psi_0}{_s\otimes_r}P_A^{it}\delta_A^{it}\zeta)\\
&=&
\tilde{G}(\Lambda_{\psi_1}(\tau^A_t(x\delta_A^{-it}))\underset{\psi_0}{_s\otimes_r}P_A^{it}\delta_A^{it}\zeta)\\
&=&
\sum_i \Lambda_\Phi[(\omega_{P_A^{it}\delta_A^{it}\zeta, e_i}\underset{N}{_b*_\alpha}id)\ga(\tau^A_t(x\delta_A^{-it})]\underset{\nu^o}{_\alpha\otimes_b}e_i\\
&=&
\sum_i \Lambda_\Phi[(\omega_{P_A^{it}\delta_A^{it}\zeta, e_i}\underset{N}{_b*_\alpha}id)(\tau^A_t\underset{N}{_b*_\alpha}\tau_t)\ga(x\delta_A^{-it})]\underset{\nu^o}{_\alpha\otimes_b}e_i\\
&=&
\sum_i\lambda_A^{-t/2}P^{it}\Lambda_\Phi[(\omega_{\zeta, P_A^{-it}e_i}\underset{N}{_b*_\alpha}id)\ga(x)\delta^{-it}]\underset{\nu^o}{_\alpha\otimes_b}e_i\\
&=&
\sum_i P^{it}J_\Phi\delta^{it}J_\Phi\Lambda_\Phi[(\omega_{\zeta, P_A^{-it}e_i}\underset{N}{_b*_\alpha}id)\ga(x)]\underset{\nu^o}{_\alpha\otimes_b}e_i\\
&=&
\sum_i P^{it}J_\Phi\delta^{it}J_\Phi\Lambda_\Phi[(\omega_{\zeta, e_i}\underset{N}{_b*_\alpha}id)\ga(x)]\underset{\nu^o}{_\alpha\otimes_b}P_A^{it}e_i\\
&=&
(\Delta_{\widehat{\Phi}}^{it}\underset{N^o}{_\alpha\otimes_b}P_A^{it})\tilde{G}(\Lambda_{\psi_1}(x)\underset{\psi_0}{_s\otimes_r}\zeta)
\end{eqnarray*}
which finishes the proof. \end{proof}

\subsection{Proposition}
\label{propweight}
{\it Let $(A, b, \ga, \phi, \psi_0)$ be a Galois system for $\gG$; let $\psi_1=\psi_0\circ T_\ga$, and let $\widehat{\Delta}_A$ be the operator introduced in \ref{deltahat}. Then, we have :
\newline
(i) there exists a normal semi-finite faithful weight $\Phi_{\hat{P}}$ on $\hat{P}$ such that $\frac{d\Phi_{\hat{P}}}{d\widehat{\Phi}^c}=\widehat{\Delta}_A$; 
\newline
(ii) there exists a normal faithful semi-finite operator valued weight $T^{\hat{P}}_L$ from $\hat{P}$ on $r(A^{\ga})$, such that $\Phi_{\hat{P}}=\psi_0\circ r^{-1}\circ T^{\hat{P}}_L$. }

\begin{proof}
Using \ref{deltahat}(ii) and the definition of the spatial derivative (\cite{T}, IX.3.11), one gets (i). Moreover, we then get that, for all $t\in\mathbb{R}$ and $x\in A^\ga$, we have, using \ref{tau} :
\[\sigma_t^{\Phi_{\hat{P}}}(r(x))=P_A^{it}J_{\psi_1}\delta_A^{it}J_{\psi_1}r(x)J_{\psi_1}\delta_A^{-it}J_{\psi_1}P_A^{-it}=P_A^{it}r(x)P_A^{-it}=r(\sigma_t^{\psi_0}(x))\]
which gives (ii). \end{proof}

\subsection{Notations}
\label{notTLG2}
$(A, b, \ga, \phi, \psi_0)$ be a Galois system for $\gG$; let $\psi_1=\psi_0\circ T_\ga$, and let $\widehat{\Delta}_A$ be the operator introduced in \ref{deltahat}; let $\Phi_{\hat{P}}$ the normal semi-finite faithful weight on $\hat{P}$ introduced in \ref{propweight}(i), and let $T^{\hat{P}}_L$ be the normal faithful semi-finite operator valued weight from $\hat{P}$ on $r(A^{\ga})$, introduced in \ref{propweight}(ii), such that $\Phi_{\hat{P}}=\psi_0\circ r^{-1}\circ T^{\hat{P}}_L$. 
\newline
Let us denote $\Phi_{\hat{Q}}$ the diagonal faithful normal semi-finite weight $\Phi_{\hat{P}}\oplus\widehat{\Phi}$ on the von Neumann algebra introduced in \ref{notTLG}. Let us first remark that we can also define a diagonal normal faithful semi-finite operator-valued weight $T^{\hat{Q}}_L$ from $\hat{Q}$ to $\tilde{\alpha}(\tilde{N})$, defined, for any positive element $\left(\begin{matrix}{A}&{X}\cr{Y^*}&{m}\end{matrix}\right)$ in $\widehat{Q}$ (which implies that $A\in\hat{P}^+$, $m\in\widehat{M}^+$ and $Y=X$), by 
\[T^{\hat{Q}}_L(\left(\begin{matrix}{A}&{X}\cr{X^*}&{m}\end{matrix}\right))=T^{\hat{P}}_L(A)\oplus \hat{T}_L(m)\]
and we get that $\Phi_{\hat{P}}\oplus\widehat{\Phi}=(\psi_0\oplus\nu)\circ T^{\hat{Q}}_L$. 
\newline
It is straightforward to get that $\left(\begin{matrix}{A}&{X}\cr{Y^*}&{m}\end{matrix}\right)$ in $\widehat{Q}$ belongs to $\gN_{\hat{Q}}$ if and only if $A$ belongs to $\gN_{\hat{P}}$, $m$ belongs to $\gN_{\widehat{\Phi}}$, $X$ is such that $\widehat{\Phi}(X^*X)< \infty$, and $Y$ is such that $\Phi_{\hat{P}}(YY^*)\leq\infty$.
\newline
Let us consider the polar decomposition $X=u|X|$; then $u$ belongs to $\widehat{I}$, and $|X|$ belongs to $\gN_{\widehat{\Phi}}$. Writing $\xi=u\Lambda_{\widehat{\Phi}}(|X|)$, we get that, for all $m\in\gN_{\widehat{\Phi}}$, we have :
\[J_{\psi_1}\pi_\ga(1\underset{N}{_b\otimes_\alpha}m^*)J_{\psi_1}\xi=J_{\psi_1}\pi_\ga(1\underset{N}{_b\otimes_\alpha}m^*)J_{\psi_1}u\Lambda_{\widehat{\Phi}}(|X|)=u\hat{J}m^*\hat{J}\Lambda_{\widehat{\Phi}}(|X|)=u|X|\hat{J}\Lambda_{\widehat{\Phi}}(m)\]
which means that $\xi\in D((H_{\psi_1})_{\mu^o}, \widehat{\Phi}^o)$ and $X=R^{\mu^o, \widehat{\Phi}^o}(\xi)$.
\newline
If now we suppose that $\left(\begin{matrix}{A}&{X}\cr{Y^*}&{m}\end{matrix}\right)$ belongs to $\gN_{\hat{Q}}\cap\gN_{\hat{Q}}^*$, we get that there exists $\eta$ in $D((H_{\psi_1})_{\mu^o}, \widehat{\Phi}^o)$, such that $Y=R^{\mu^o, \widehat{\Phi}^o}(\eta)$, and $YY^*=\theta^{\mu^o, \widehat{\Phi}^o}(\eta, \eta)$; by definition of the spatial derivative, the fact that $\Phi_{\hat{P}}(YY^*)<\infty$ implies that $\eta\in\mathcal D(\widehat{\Delta}_A^{1/2})$, and $\Phi_{\hat{P}}(\theta^{\mu^o, \widehat{\Phi}^o}(\eta, \eta))=\|\widehat{\Delta}_A^{1/2}\eta\|^2$; more precisely, there exists an antilinear involutive isometry $\widetilde{J}$ on $H_{\psi_1}$ such that $\tilde{J}\widehat{\Delta}_A^{1/2}=\widehat{\Delta}_A^{-1/2}\tilde{J}$, and we can write :
\[\Lambda_{\Phi_{\hat{Q}}}(\left(\begin{matrix}{A}&{R^{\mu, \widehat{\Phi}'}(\xi)}\cr{R^{\mu, \widehat{\Phi}'}(\eta)^*}&{m}\end{matrix}\right))=\Lambda_{\Phi_{\hat{P}}}(A)\oplus \xi\oplus \widetilde{J}\widehat{\Delta}_A^{1/2}\eta\oplus\Lambda_{\widehat{\Phi}}(m)\]
and we identify this way $H_{\Phi_{\widehat{Q}}}$ with $H_{\Phi_{\widehat{P}}}\oplus H_{\psi_1}\oplus H_{\psi_1}\oplus H$; for simplification, we shall identify $\Lambda_{\Phi_{\hat{Q}}}(\left(\begin{matrix}{A}&{0}\cr{0}&{0}\end{matrix}\right))$ with $\Lambda_{\Phi_{\widehat{P}}}(A)$, $\Lambda_{\Phi_{\hat{Q}}}(\left(\begin{matrix}{0}&{0}\cr{0}&{m}\end{matrix}\right))$ with $\Lambda_{\widehat{\Phi}}(m)$. We shall write $p^{1,2}_{H_{\psi_1}}$ for the projection on the first subspace $H_{\psi_1}$ of $H_{\Phi_{\widehat{Q}}}$, and $p^{2,1}_{H_{\psi_1}}$, for the projection on the second subspace $H_{\psi_1}$. 
\newline
If $X\in\widehat{I}$ is such that $\widehat{\Phi}(X^*X)<\infty$, we shall write $\Lambda^{1,2}(X)=\Lambda_{\Phi_{\widehat{Q}}}
(\left(\begin{matrix}{0}&{X}\cr{0}&{0}\end{matrix}\right))$ (and, therefore $\Lambda^{1,2}(R^{\mu, \widehat{\Phi}'}(\xi))=\xi$ for all $\xi\in D(_\mu H_{\psi_1}, \widehat{\Phi}')$).
\newline
If $Y\in\widehat{I}$ is such that $\Phi_{\widehat{P}}(YY^*)<\infty$, we shall write $\Lambda^{2,1}(Y^*)=
 \Lambda_{\Phi_{\hat{Q}}}(\left(\begin{matrix}{0}&{0}\cr{Y^*}&{0}\end{matrix}\right))$, and, therefore, if $\eta\in D(_\mu (H_{\psi_1}), \widehat{\Phi}')\cap\mathcal D(\widehat{\Delta}_A^{1/2})$, we have $\Lambda^{2,1}(R^{\mu, \widehat{\Phi}'}(\eta)^*)=\widetilde{J}\widehat{\Delta}_A^{1/2}\eta$.
 \newline
The identification of $H_{\Phi_{\widehat{Q}}}$ with $H_{\Phi_{\widehat{P}}}\oplus H_{\psi_1}\oplus H_{\psi_1}\oplus H$ leads also to write :
\[\Delta_{\Phi_{\widehat{Q}}}=\Delta_{\Phi_{\widehat{P}}}\oplus \widehat{\Delta}_A^{1/2}\oplus \widehat{\Delta}_A^{1/2}\oplus\Delta_{\widehat{\Phi}}\] and $J_{\Phi_{\widehat{Q}}}=J_{\Phi_{\widehat{P}}}\oplus (\tilde{J}\oplus\tilde{J})\circ\tau\oplus \hat{J}$, where $\tau (\xi\oplus\eta)=\eta\oplus\xi$, for any $\xi$, $\eta$ in $H_{\psi_1}$. 
\newline
For any $n\in N$, $x\in A^\ga$, we get that :
\[\pi_{\Phi_{\widehat{Q}}}(\widetilde{\alpha}(x\oplus n))\Lambda_{\Phi_{\hat{Q}}}(\left(\begin{matrix}{A}&{R^{\mu, \widehat{\Phi}'}(\xi)}\cr{R^{\mu, \widehat{\Phi}'}(\eta)^*}&{m}\end{matrix}\right))=
\Lambda_{\Phi_{\hat{Q}}}(\left(\begin{matrix}{r(x)A}&{r(x)R^{\mu, \widehat{\Phi}'}(\xi)}\cr{\alpha(n)R^{\mu, \widehat{\Phi}'}(\eta)^*}&{\alpha(n)m}\end{matrix}\right))\]
 Using \ref{spatial}, we get, for any $n\in N$, analytical with respect to $\nu$, that :
 \[R^{\mu, \widehat{\Phi}'}(\eta)\alpha(n)=R^{\mu, \widehat{\Phi}'}(\mu(\sigma_{i/2}^{\widehat{\Phi}'}(\beta(n))\eta)=R^{\mu, \widehat{\Phi}'}(\mu(\beta(\sigma_{i/2}^\nu(n))\eta)=R^{\mu, \widehat{\Phi}'}(b(\sigma_{i/2}^\nu(n))\eta)\]
and, therefore, that :
\[\Lambda^{2,1}(\alpha(n)R^{\mu, \widehat{\Phi}'}(\eta)^*)=\tilde{J}\widehat{\Delta}_A^{1/2}b(\sigma_{i/2}^\nu(n^*))\eta\]
and, using \ref{deltahat}(i), we get that :
\[\Lambda^{2,1}(\alpha(n)R^{\mu, \widehat{\Phi}'}(\eta)^*)=\tilde{J}b(n^*)\widehat{\Delta}_A^{1/2}\eta\]
which, by continuity, remains true for all $n\in N$; from which we obtain :
\[\pi_{\Phi_{\widehat{Q}}}(\widetilde{\alpha}(x\oplus n))=\pi_{\Phi_{\widehat{P}}}(r(x))\oplus r(x)\oplus \tilde{a}(n)\oplus \alpha(n)\]
where we define $\tilde{a}(n)=\tilde{J}b(n^*)\tilde{J}$. 
\newline
With similar arguments, we obtain :
\[\pi_{\Phi_{\widehat{Q}}}(\widetilde{\hat{\beta}}(x\oplus n))=\pi_{\Phi_{\widehat{P}}}(s(x))\oplus s(x)\oplus \tilde{b}(n)\oplus \hat{\beta}(n)\]
where we define $\tilde{b}(n)=\tilde{J}a(n^*)\tilde{J}$. 
Therefore, we get that $\pi_{\Phi_{\widehat{Q}}}(e_1)=p_{H_{\Phi_{\widehat{P}}}}+p^{1,2}_{H_{\psi_1}}$ and $\pi_{\Phi_{\widehat{Q}}}(e_2)=p^{2,1}_{H_{\psi_1}}+p_H$.

\subsection{Proposition}
\label{propTLG}
{\it Let's use the notations of \ref{notTLG} and \ref{notTLG2}. Then, we have :
\newline
(i) for any $\eta\in D((H_{\psi_1})_{\mu^o}, \widehat{\Phi}^o)$, $v\in D(_\alpha H, \nu)\cap D(H_{\hat{\beta}}, \nu^o)$, $\xi\in D((H_{\psi_1})_s, \psi_0^o)$, the element $X=(\omega_{v, \xi}*id)[\tilde{G}^*(1\underset{N}{_\alpha\otimes_\beta}R^{\mu^o, \widehat{\Phi}^o}(\eta))\widehat{W}]$, which belongs to $\widehat{I}$ by \ref{lemTLG}, is such that $\left(\begin{matrix}0&X\cr0&0\end{matrix}\right)$ belongs to $\gN_{\Phi_{\widehat{Q}}}$. 
\newline
(ii) let $(\xi_i)_{i\in I}$ be an orthogonal $(s, \psi_0^o)$-basis of $H_{\psi_1}$; there exists $\eta_i\in D((H_{\psi_1})_{\mu^o}, \widehat{\Phi}^o)$ such that :
\[(\omega_{v, \xi_i}*id)[\tilde{G}^*(1\underset{N}{_\alpha\otimes_\beta}R^{\mu^o, \widehat{\Phi}^o}(\eta))\widehat{W}]=R^{\mu^o, \widehat{\Phi}^o}(\eta_i)\]
Moreover, we have :
\[\|v\underset{\nu^o}{_\alpha\otimes_b}\eta\|^2=\sum_i\|\eta_i\|^2=\|\sum_i\xi_i\underset{\psi_0}{_s\otimes_r}\eta_i\|^2\]
(iii) we have :
\[\tilde{G}^*(v\underset{\nu^o}{_\alpha\otimes_b}\eta)=\sum_i\xi_i\underset{\psi_0}{_s\otimes_r}\eta_i\]
(iv) we have :
\[\Lambda^{1,2}(X)=(\omega_{v, \xi}*id)(\tilde{G}^*)\eta\]
}

\begin{proof}
We have :
\begin{eqnarray*}
X^*X
&=&
((\omega_{v, \xi}*id)[\tilde{G}^*(1\underset{N^o}{_\alpha\otimes_\beta}R^{\mu^o, \widehat{\Phi}^o}(\eta))\widehat{W}])^*(\omega_{v, \xi}*id)[\tilde{G}^*(1\underset{N^o}{_\alpha\otimes_\beta}R^{\mu^o, \widehat{\Phi}^o}(\eta))\widehat{W}]\\
&=&
(\omega_v\underset{N}{_{\hat{\beta}}*_\alpha}id)[\widehat{W}^*(1\underset{N^o}{_\alpha\otimes_b}R^{\mu^o, \widehat{\Phi}^o}(\eta))^*\tilde{G}(\theta^{s,\psi_0^o}(\xi, \xi)\underset{A^\ga}{_s\otimes_r}1)\tilde{G}^*(1\underset{N}{_\alpha\otimes_\beta}R^{\mu^o, \widehat{\Phi}^o}(\eta))\widehat{W}]\\
&\leq&
\|R^{s, \psi_0}(\xi)\|^2(\omega_v\underset{N}{_{\hat{\beta}}*_\alpha}id)[\widehat{W}^*(1\underset{N^o}{_\alpha\otimes_\beta}<\eta, \eta>_{\mu^o, \widehat{\Phi}^o})\widehat{W}]\\
&=&
\|R^{s, \psi_0}(\xi)\|^2(\omega_v\underset{N}{_{\hat{\beta}}*_\alpha}id)(\widehat{\Gamma}(<\eta, \eta>_{\mu^o, \widehat{\Phi}^o}))
\end{eqnarray*}
and, therefore, using the left-invariance of $\widehat{T_L}$, then, using \ref{notTLG2}  :
\begin{eqnarray*}
\widehat{\Phi}(X^*X)
&\leq&
 \|R^{s, \psi_0}(\xi)\|^2 (\widehat{T_L}(<\eta, \eta>_{\mu^o, \widehat{\Phi}^o})v|v)\\
 &=& \|R^{s, \psi_0}(\xi)\|^2 \|v\underset{N^o}{_\alpha\otimes_\beta}\Lambda_{\widehat{\Phi}}(R^{\mu^o, \widehat{\Phi}^o}(\eta))\|^2\\
 &=&\|R^{s, \psi_0}(\xi)\|^2 \|v\underset{N^o}{_\alpha\otimes_b}\eta\|^2
 \end{eqnarray*}
from which we get that $X^*X$ belongs to $\gM^+_{\widehat{\Phi}}$, and, using \ref{notTLG2}, we finish the proof of (i). 
\newline
Making now the same calculation with $X_i=(\omega_{v, \xi_i}*id)[\tilde{G}^*(1\underset{N^o}{_\alpha\otimes_\beta}R^{\mu^o, \widehat{\Phi}^o}(\eta))\widehat{W}]$, we get that :
\[\sum_iX_i^*X_i=(\omega_v\underset{N}{_{\hat{\beta}}*_\alpha}id)\widehat{\Gamma}(<\eta, \eta>_{\mu^o, \widehat{\Phi}^o})\]
and, then :
\[\widehat{\Phi}(\sum_iX_i^*X_i)=\|v\underset{N^o}{_\alpha\otimes_b}\eta\|^2\]
Using again \ref{notTLG2}, we get that there exists $\eta_i\in D((H_{\psi_1})_{\mu^o}, \widehat{\Phi}^o)$ such that $X_i=R^{\mu^o, \widehat{\Phi}^o}(\eta_i)$; from which we get that :
\[\sum_i\|\xi_i\underset{\psi_0^o}{_s\otimes_r}\eta_i\|^2=\sum_i\|\eta_i\|^2=\sum_i\|\Lambda_{\widehat{\Phi}}(R^{\mu^o, \widehat{\Phi}^o}(\eta_i))\|^2=\sum_i\widehat{\Phi}(X_i^*X_i)=\|v\underset{N^o}{_\alpha\otimes_b}\eta\|^2\]
which is (ii). Let now $m\in\gN_{\widehat{\Phi}}$; we have :
\begin{eqnarray*}
(\hat{J}m^*\hat{J}\underset{N^o}{_\alpha\otimes_b}1)\tilde{G}(\xi_i\underset{\psi_0^o}{_s\otimes_r}\eta_i)
&=&
\tilde{G}(\xi_i\underset{\psi_0^o}{_s\otimes_r}\mu(\hat{J}m^*\hat{J})\eta_i)\\
&=&
\tilde{G}(\xi_i\underset{\psi_0^o}{_s\otimes_r}R^{\mu^o, \widehat{\Phi}^o}(\eta_i)\hat{J}\Lambda_{\widehat{\Phi}}(m))\\
&=&
\tilde{G}(\xi_i\underset{\psi_0^o}{_s\otimes_r}X_i\hat{J}\Lambda_{\widehat{\Phi}}(m))
\end{eqnarray*}
and, therefore :
\begin{eqnarray*}
(\hat{J}m^*\hat{J}\underset{N^o}{_\alpha\otimes_b}1)\tilde{G}\sum_i(\xi_i\underset{\psi_0^o}{_s\otimes_r}\eta_i)
&=&
\tilde{G}\sum_i(\xi_i\underset{\psi_0^o}{_s\otimes_r}(\omega_{v, \xi_i}*id)[\tilde{G}^*(1\underset{N}{_\alpha\otimes_\beta}R^{\mu^o, \widehat{\Phi}^o}(\eta))\widehat{W}]\hat{J}\Lambda_{\widehat{\Phi}}(m)\\
&=&
(1\underset{N}{_\alpha\otimes_\beta}R^{\mu^o, \widehat{\Phi}^o}(\eta))\widehat{W}(v\underset{N}{_{\hat{\beta}}\otimes_\alpha}\hat{J}\Lambda_{\widehat{\Phi}}(m))
\end{eqnarray*}
Therefore, taking now $\zeta_1\in D(_\alpha H, \nu)$ and $\zeta_2\in D((H_{\psi_1})_b, \nu^o)$, we get that :
\begin{multline*}
((\hat{J}m^*\hat{J}\underset{N^o}{_\alpha\otimes_b}1)\tilde{G}\sum_i(\xi_i\underset{\psi_0^o}{_s\otimes_r}\eta_i)|\zeta_1\underset{\nu^o}{_\alpha\otimes_b}\zeta_2)\\
=
((1\underset{N}{_\alpha\otimes_\beta}R^{\mu^o, \widehat{\Phi}^o}(\eta))\widehat{W}(v\underset{N}{_{\hat{\beta}}\otimes_\alpha}\hat{J}\Lambda_{\widehat{\Phi}}(m))|\zeta_1\underset{\nu^o}{_\alpha\otimes_b}\zeta_2)\\
=
(R^{\mu^o, \widehat{\Phi}^o}(\eta)(\omega_{v, \zeta_1}*id)(\widehat{W})\hat{J}\Lambda_{\widehat{\Phi}}(m)|\zeta_2)
\end{multline*}
and, using now (\cite{E5}, 3.10(ii) applied to $\widehat{\gG}$, and 3.11(iii)), we get it is equal to :
\[(R^{\mu^o, \widehat{\Phi}^o}(\eta)\hat{J}\Lambda_{\widehat{\Phi}}(\omega_{J\zeta_1, Jv}*id)\widehat{\Gamma}(m))|\zeta_2)=(\mu(\hat{J}(\omega_{J\zeta_1, Jv}*id)\widehat{\Gamma}(m)^*\hat{J})\eta|\zeta_2)\]
Taking the limit when $m$ goes to $1$, we get :
\begin{eqnarray*}
(\tilde{G}\sum_i\xi_i\underset{\psi_0^o}{_s\otimes_r}\eta_i|\zeta_1\underset{\nu^o}{_\alpha\otimes_b}\zeta_2)
&=&
(\beta(<Jv, J\zeta_1>_{{\hat{\beta}}, \nu^o})\eta|\zeta_2)\\
&=&
(\beta(<\zeta_1, v>_{\alpha, \nu})\eta|\zeta_2)\\
&=&(v\underset{\nu^o}{_\alpha\otimes_b}\eta|\zeta_1\underset{\nu^o}{_\alpha\otimes_b}\zeta_2)
\end{eqnarray*}
from which we get that $\tilde{G}(\sum_i\xi_i\underset{\psi_0^o}{_s\otimes_r}\eta_i)=v\underset{\nu^o}{_\alpha\otimes_b}\eta$, which is (iii); this can be written :
\[(\omega_{v, \xi_i}*id)(\tilde{G}^*)\eta=\eta_i=\Lambda^{1,2}(X_i)\]
which, by linearity and continuity, gives (iv). 
\end{proof}

\subsection{Proposition}
\label{propTLG2}
{\it Let us use the notations of \ref{notTLG}, \ref{notTLG2}, \ref{propTLG} and take $\eta\in D((H_{\psi_1})_{\mu^o}, \widehat{\Phi}^o)\cap\mathcal D(\widehat{\Delta}_A^{1/2})$; let us define $\tilde{s}$ the anti-representation of $A^\ga$ on $H_{\psi_1}$ defined by $\tilde{s}(x)=\tilde{J}r(x^*)\tilde{J}$, for all $x\in A^\ga$; let us define $\tilde{a}$ the representation of $N$ on $H_{\psi_1}$ defined by $\tilde{a}(n)=\tilde{J}b(n^*)\tilde{J}$, for all $n\in N$; then, for any $v\in D(H_{\hat{\beta}}, \nu^o)$ and $\xi\in D((H_{\psi_1})_s, \psi_0^o)\cap D(_rH_{\psi_1}, \psi_0)$ the element 
\[X=(\omega_{v, \xi}*id)[\tilde{G}^*(1\underset{N^o}{_\alpha\otimes_\beta}R^{\mu^o, \widehat{\Phi}^o}(\eta))\widehat{W}]\]
is such that $\left(\begin{matrix}{0}&{0}\cr{X^*}&{0}\end{matrix}\right)$ belongs to $\gN_{\Phi_{\widehat{Q}}}$, and we have :
\[\Lambda^{2,1}(X^*)=(\omega_{\xi, v}*id)[(J\underset{N^o}{_\alpha\otimes_b}\tilde{J})\tilde{G}(J_{\psi_1}\underset{A^{\ga o}}{_r\otimes_{\tilde{s}}}\tilde{J})]\Lambda^{2,1}(R^{\mu^o, \widehat{\Phi}^o}(\eta))\]
}

\begin{proof}
Let us first take $\eta$ such that $\Lambda_{\Phi_{\widehat{Q}}}(\left(\begin{matrix}{0}&{R^{\mu^o, \widehat{\Phi}^o}(\eta)}\cr{0}&{0}\end{matrix}\right))$ belongs to the Tomita algebra $\mathcal T_{\Phi_{\widehat{Q}}}$, $x$ in the Tomita algebra $\mathcal T_{\psi_1, T_\ga}$, and $y$, $z$ in $\mathcal T_{\Phi, T_L}$. Then, $\Lambda_{\psi_1}(x)$ belongs to $D((H_{\psi_1})_s, \psi_0^o)$, and $J\Lambda_\Phi(y^*z)$ belongs to $D(_\alpha H, \nu)\cap D(H_{\hat{\beta}}, \nu^o)$. Therefore, we can apply \ref{propTLG}(i) to the element $X=(\omega_{J\Lambda_\Phi(y^*z), \Lambda_{\psi_1}(x)}*id)[\tilde{G}^*(1\underset{N^o}{_\alpha\otimes_\beta}R^{\mu^o, \widehat{\Phi}^o}(\eta))\widehat{W}]$.
\newline
Using \ref{notTLG2} and (\cite{E5}, 3.11 applied to $\gG$), we get that the element  $\sigma_t^{\Phi_{\widehat{Q}}}(\left(\begin{matrix}{0}&{X}\cr{0}&{0}\end{matrix}\right))$ is of the form $\left(\begin{matrix}{0}&{X_t}\cr{0}&{0}\end{matrix}\right)$, with $X_t=\widehat{\Delta}_A^{it}X\Delta_{\widehat{\Phi}}^{-it}$. Using now \ref{deltahat}(iii), we get that :
\[X_t=(\omega_{\Delta_\Phi^{it}J\Lambda_\Phi(y^*z), \Delta_{\psi_1}^{-it}\Lambda_{\psi_1}(x)}*id)[\tilde{G}^*(1\underset{N^o}{_\alpha\otimes_\beta}\widehat{\Delta}_A^{it}R^{\mu^o, \widehat{\Phi}^o}(\eta)\widehat{\Delta}^{-it})\widehat{W}]\]
and the hypothesis on $\eta$, $x$, $y$, $z$ give that the function $t\mapsto X_t$ extends to an analytic function; in particular, we get that $\Lambda^{1,2}(X)$ belongs to $\mathcal D(\widehat{\Delta}_A^{1/2})$, and, using \ref{propTLG}(iv) and \ref{tildeG}(iii), we get :
\begin{eqnarray*}
\widehat{\Delta}_A^{1/2}\Lambda^{1,2}(X)
&=&
\Lambda^{1,2}(X_{-i/2})\\
&=&
(\omega_{\Delta_\Phi^{1/2}J\Lambda_\Phi(y^*z), \Delta_{\psi_1}^{-1/2}\Lambda_{\psi_1}(x)}*id)(\tilde{G}^*)\tilde{J}\widehat{\Delta}_A^{1/2}\eta\\
&=&
(\omega_{\Delta_{\psi_1}^{-1/2}\Lambda_{\psi_1}(x), \Delta_\Phi^{1/2}J\Lambda_\Phi(y^*z)}*id)(\tilde{G})^*\tilde{J}\widehat{\Delta}_A^{1/2}\eta\\
&=&
(\omega_{J_{\psi_1}\Lambda_{\psi_1}(x), \Lambda_\Phi(y^*z)}*id)(\tilde{G})\tilde{J}\widehat{\Delta}_A^{1/2}\eta
\end{eqnarray*}
and, therefore :
\begin{eqnarray*}
\Lambda^{2,1}(X^*)
&=&
\tilde{J}\widehat{\Delta}_A^{1/2}\Lambda^{1,2}(X)\\
&=&
\tilde{J}(\omega_{J_{\psi_1}\Lambda_{\psi_1}(x), \Lambda_\Phi(y^*z)}*id)(\tilde{G})\tilde{J}\widehat{\Delta}_A^{1/2}\eta\\
&=&
(\omega_{\Lambda_{\psi_1}(x),J\Lambda_\Phi(y^*z)}*id)[(J\underset{N^o}{_\alpha\otimes_b}\tilde{J})\tilde{G}(J_{\psi_1}\underset{A^{\ga o}}{_r\otimes_{\tilde{s}}}\tilde{J})]\Lambda^{2,1}(R^{\mu^o, \widehat{\Phi}^o}(\eta)^*)
\end{eqnarray*}
By the closedness of $\Lambda_{\widehat{Q}}$, and, therefore, of $\Lambda^{2,1}$, we get, for any $v\in D(H_{\hat{\beta}}, \nu^o)$ and $\xi\in D((H_{\psi_1})_s, \psi_0^o)\cap D(_rH_{\psi_1}, \psi_0)$, that 
$X=(\omega_{v, \xi}*id)[\tilde{G}^*(1\underset{N^o}{_\alpha\otimes_\beta}R^{\mu^o, \widehat{\Phi}^o}(\eta))\widehat{W}]$ is such that $X^*$ belongs to $\mathcal D (\Lambda^{2,1})$ and that :
\[\Lambda^{1,2}(X^*)=(\omega_{\xi, v}*id)[(J\underset{N^o}{_\alpha\otimes_b}\tilde{J})\tilde{G}(J_{\psi_1}\underset{A^{\ga o}}{_r\otimes_{\tilde{s}}}\tilde{J})]\Lambda^{2,1}(R^{\mu^o, \widehat{\Phi}^o}(\eta)^*)\]
Using again the closednes of $\Lambda^{2,1}$, we get that this result remains true for any $\eta$ such that 
$\left(\begin{matrix}{0}&{R^{\mu^o, \widehat{\Phi}^o}(\eta)}\cr {0}&{0}\end{matrix}\right)$ belongs to 
$\gN_{\Phi_{\widehat{Q}}}\cap \gN_{\Phi_{\widehat{Q}}}^*$ (i.e., using \ref{notTLG2}, if $\eta$ belongs to $\mathcal D(\widehat{\Delta}_A^{1/2})$). 
\end{proof}

\subsection{Theorem}
\label{thP}
{\it The operator-valued weight $T_L^{\widehat{P}}$ is left-invariant. }
\begin{proof}
Let $\eta$ in $D((H_{\psi_1})\mu^o, \widehat{\Phi}^o)\cap \mathcal {\Delta}_A^{1/2})$; let $(v_i)_{i\in I}$ a $(\hat{\beta},\nu^o)$ orthogonal basis of $H$, and $\xi$ in $D((H_{\psi_1})_s, \psi_0^o)\cap D(_rH_{\psi_1}, \psi_0)$; let us write :
\[X_i=(\omega_{v_i, \xi}*id)[\tilde{G}^*(1\underset{N^o}{_\alpha\otimes_\beta}R^{\mu^o, \widehat{\Phi}^o}(\eta))\widehat{W}]\]
We then get :
\[\omega_\xi(id*\Phi_{\widehat{P}})(\Gamma_{\widehat{P}}(\theta^{\mu^o, \widehat{\Phi}^o}(\eta, \eta)))
=
\Phi_{\widehat{P}}(\omega_\xi*id)[\tilde{G}(1\underset{N^o}{_\alpha\otimes_b}\theta^{\mu^o, \widehat{\Phi}^o}(\eta, \eta))\tilde{G}^*]\]
is equal to :
\[\sum_i\Phi_{\widehat{P}}((\omega_\xi*id)(\tilde{G}(1\underset{N}{_\alpha\otimes_\beta}R^{\mu^o, \widehat{\Phi}^o}(\eta))\widehat{W}(\theta^{\hat{\beta}}(v_i, v_i)\underset{N}{_{\hat{\beta}}\otimes_\alpha}1)\widehat{W}^*(1\underset{N^o}{_\alpha\otimes_b}R^{\mu^o, \widehat{\Phi}^o}(\eta)^*\tilde{G}^*]\]
which can be written, using \ref{propTLG2}  :
\begin{eqnarray*}
\sum_i\Phi_{\widehat{P}}(X_iX_i^*)
&=&
\sum_i\|\Lambda^{2,1}(X_i^*)\|^2\\
&=&
\sum_i\|(\omega_{\xi, v_i}*id)[(J\underset{N^o}{_\alpha\otimes_b}\tilde{J})\tilde{G}(J_{\psi_1}\underset{A^{\ga o}}{_r\otimes_{\tilde{s}}}\tilde{J})]\Lambda^{2,1}(R^{\mu^o, \widehat{\Phi}^o}(\eta))\|^2\\
&=&
\sum_i \|v_i\underset{N}{_{\hat{\beta}}\otimes_\alpha}(\omega_{\xi, v_i}*id)[(J\underset{N^o}{_\alpha\otimes_b}\tilde{J})\tilde{G}(J_{\psi_1}\underset{A^{\ga o}}{_r\otimes_{\tilde{s}}}\tilde{J})]\Lambda^{2,1}(R^{\mu^o, \widehat{\Phi}^o}(\eta))\|^2\\
&=&
\|[(J\underset{N^o}{_\alpha\otimes_b}\tilde{J})\tilde{G}(J_{\psi_1}\underset{A^{\ga o}}{_r\otimes_{\tilde{s}}}\tilde{J})](\xi\underset{A^{\ga o}}{_r\otimes_{\tilde{s}}}\Lambda^{2,1}(R^{\mu^o, \widehat{\Phi}^o}(\eta))\|^2\\
&=&
\|\xi\underset{A^{\ga o}}{_r\otimes_{\tilde{s}}}\Lambda^{2,1}(R^{\mu^o, \widehat{\Phi}^o}(\eta))\|^2\\
&=&
(T_L^{\widehat{P}}(\theta^{\mu^o, \widehat{\Phi}^o}(\eta, \eta))\xi|\xi)
\end{eqnarray*}
from which we get that $(id*\Phi_{\widehat{P}})(\Gamma_{\widehat{P}}(\theta^{\mu^o, \widehat{\Phi}^o}(\eta, \eta)))=T_L^{\widehat{P}}(\theta^{\mu^o, \widehat{\Phi}^o}(\eta, \eta))$. As any element in $\gM_{\Phi_{\widehat{P}}}^+$ can be approximated by below by finite sums of operators of the form $\theta^{\mu^o, \widehat{\Phi}^o}(\eta, \eta))$, we get the result. \end{proof}

\subsection{Theorem}
\label{thP2}
{\it With the notations of \ref{notTLG} and \ref{notTLG2}, we have :
\newline
(i) $(A^\ga, \widehat{P}, r, s, \Gamma_{\widehat{P}}, T_L^{\widehat{P}}, R_{\widehat{P}}\circ T_L^{\widehat{P}}\circ R_{\widehat{P}}, \psi_0)$ is a measured quantum groupoid. We shall denote this measured quantum groupoid by $\gG_1(A,b,\ga, \phi, \psi_0)$, or simply by $\gG_1(\ga)$. Following \cite{DC1}, its dual $\widehat{\gG_1(\ga)}$ will be called the reflected measured quantum groupoid of $\gG$ throught the Galois system $(A,b,\ga, \phi, \psi_0)$, or simply, throught $\ga$. 
\newline
(ii) $(\tilde{N}, \widehat{Q}, \widetilde{\alpha}, \widetilde{\hat{\beta}}, \Gamma_{\widehat{Q}}, T_L^{\widehat{Q}}, R_{\widehat{Q}}\circ T_L^{\widehat{Q}}\circ R_{\widehat{Q}}, \psi_0\oplus\nu)$ is a measured quantum groupoid. We shall denote this measured quantum groupoid by $\gG_2(A,b,\ga, \phi, \psi_0)$, or simply by $\gG_2(\ga)$. }

\begin{proof}
By \ref{coproduct}(iii), we know that $(A^\ga, \widehat{P}, r, s, \Gamma_{\widehat{P}})$ is a Hopf-bimodule, and,  by \ref{thP}, that $T_L^{\widehat{P}}$ is left-invariant. Using again \ref{coproduct}(iii), we get that  $R_{\widehat{P}}\circ T_L^{\widehat{P}}\circ R_{\widehat{P}}$ is right-invariant. The only result needed is that the modular automorphism groups $\sigma^{\Phi_{\widehat{P}}}$ and $\sigma^{\Phi_{\widehat{P}}\circ R_{\widehat{P}}}$ commute. By definition, we have, for all $A\in\widehat{P}$, we have, using \ref{deltahat}(i):
\begin{eqnarray*}
\sigma^{\Phi_{\widehat{P}}}_t(A)
&=&
\widehat{\Delta}_A^{it}A\widehat{\Delta}_A^{-it}\\
&=&
P_A^{it}J_{\psi_1}\delta_A^{it}J_{\psi_1}AJ_{\psi_1}\delta_A^{-it}J_{\psi_1}P_A^{-it}
\end{eqnarray*}
and, using \ref{deltahat}(i) and \ref{tau}(v) and (vi) :
\begin{eqnarray*}
\sigma^{\Phi_{\widehat{P}}\circ R_{\widehat{P}}}_s(A)
&=&
R_{\widehat{P}}\circ\sigma^{\Phi_{\widehat{P}}}_{-s}\circ R_{\widehat{P}}(A)\\
&=&
J_{\psi_1}\widehat{\Delta}_A^{-is}J_{\psi_1}AJ_{\psi_1}\widehat{\Delta}_A^{is}J_{\psi_1}\\
&=&
J_{\psi_1}P_A^{-is}J_{\psi_1}\delta_A^{-is}A\delta_A^{is}J_{\psi_1}P_A^{is}J_{\psi_1}\\
&=&
P_A^{-is}\delta_A^{-is}A\delta_A^{is}P_A^{is}
\end{eqnarray*}
and, as $P_A^{it}J_{\psi_1}\delta_A^{it}J_{\psi_1}$ commutes with $P_A^{-is}\delta_A^{-is}$, we obtain the result, and we finish the proof of (i). 
\newline
We had obtained in \ref{coproduct}(i), that $(\tilde{N}, \widehat{Q}, \widetilde{\alpha}, \widetilde{\hat{\beta}}, \Gamma_{\widehat{Q}})$ is a Hopf-bimodule; using \ref{thP} and the definition of $T_L^{\widehat{Q}}$ (\ref{notTLG2}), we get that $T_L^{\widehat{Q}}$ is left-invariant; using \ref{coproduct}(ii), we get that $R_{\widehat{Q}}\circ T_L^{\widehat{Q}}\circ R_{\widehat{Q}}$ is right-invariant. The calculation made in (i) proves as well that the automorphism groups $\sigma^{\Phi_{\widehat{Q}}}$ and $\sigma^{\Phi_{\widehat{Q}}\circ R_{\widehat{Q}}}$ commute, which finish the proof. \end{proof}

\subsection{Theorem}
\label{dual}
{\it Let $\gG$ a measured quantum groupoid, and $(A, b, \ga, \phi, \psi_0)$ a Galois system for $\gG$; let us denote $(\tilde{N}, Q, \widetilde{\alpha}, \widetilde{\beta}, \Gamma_Q, T_L^Q, R_QT_L^QR_Q, \psi_0\oplus\nu)$ the dual measured quantum groupoid $\widehat{\gG_2(\ga)}$. This measured quantum groupoid will be denoted the linking measured quantum groupoid between $\gG$ and the reflected measured quantum groupoid $\widehat{\gG_1(\ga)}$. We shall consider that the von Neumann algebra $Q$ acts on $H_{\Phi_{\widehat{Q}}}=H_{\Phi_{\widehat{P}}}\oplus H_{\psi_1}\oplus H_{\psi_1}\oplus H$. Then :
\newline
(i)  $\widetilde{\alpha}(e_1)$, $\widetilde{\alpha}(e_2)$, $\widetilde{\beta}(e_1)$, $\widetilde{\beta}(e_2)$ belong to $Z(Q)$. 
\newline
(ii) we have $p_{H_{\Phi_{\widehat{P}}}}=\tilde{\alpha}(e_1)\tilde{\beta}(e_1)$; $p^{1,2}_{H_{\psi_1}}= \tilde{\alpha}(e_1)\tilde{\beta}(e_2)$; $p^{2,1}_{H_{\psi_1}}=\tilde{\alpha}(e_2)\tilde{\beta}(e_1)$, and $p_H=\tilde{\alpha}(e_2)\tilde{\beta}(e_2)$; all these projections belong to $Z(Q)$. 
\newline
(iii) we have $Q_{p_{H_{\Phi_{\widehat{P}}}}}=P$, $Q_{p^{1,2}_{H_{\psi_1}}}=A$, $Q_{p^{2,1}_{H_{\psi_1}}}=\tilde{J}A\tilde{J}$, and $Q_{p_H}=M$. Therefore, we have $Q=P\oplus A\oplus A^o\oplus M$. 
\newline
(iv) if $x\in P$, $y\in M$, $z\in A$, we have :
\[\Gamma_P(x)=\Gamma_Q(x)_{\tilde{\alpha}(e_1)\tilde{\beta}(e_1)\underset{\tilde{N}}{_{\tilde{\beta}}\otimes_{\tilde{\alpha}}}\tilde{\alpha}(e_1)\tilde{\beta}(e_1)}\]
\[\Gamma(y)=\Gamma_Q(y)_{\tilde{\alpha}(e_2)\tilde{\beta}(e_2)\underset{\tilde{N}}{_{\tilde{\beta}}\otimes_{\tilde{\alpha}}}\tilde{\alpha}(e_2)\tilde{\beta}(e_2)}\]
\[\ga(z)=\Gamma_Q(z)_{\tilde{\alpha}(e_1)\tilde{\beta}(e_2)\underset{\tilde{N}}{_{\tilde{\beta}}\otimes_{\tilde{\alpha}}}\tilde{\alpha}(e_2)\tilde{\beta}(e_2)}\] 
(v) let $R$ (resp. $R_P$, resp. $R_Q$) be the co-inverse of $\gG$ (resp. of the reflected measured quantum groupoid, resp. of the linking measured quantum groupoid); let $\tau_t$, $\tau_t^P$, $\tau_t^Q$ be the scaling groups of these measured quantum groupoids, $\gamma_t$, $\gamma_t^P$, $\gamma_t^Q$ be the automorphism groups on the basis of these measured quantum groupoids, as defined in \ref{MQG} or \cite{E5}, 3.8(i), (ii) and (v); we have, for any $x\in P$, $y\in M$, $z_1$, $z_2$ in $A$, $n\in N$, $u\in A^\ga$ :
\[R_Q(x\oplus z_1\oplus z_2^o\oplus y)=R_P(x)\oplus z_2\oplus z_1^o\oplus R(y)\]
\[\tau_t^Q(x\oplus z_1\oplus z_2^o\oplus y)=\tau_t^P(x)\oplus \widehat{\Delta}_A^{it}z_1\widehat{\Delta}_A^{-it}\oplus (\widehat{\Delta}_A^{it}z_2\widehat{\Delta}_A^{-it})^o\oplus\tau_t(y)\] 
\[\gamma_t^Q(u\oplus n)=\gamma_t^P(u)\oplus \gamma_t(n)\]}

\begin{proof}
As $\widetilde{\alpha}(e_1)=\widetilde{\hat{\beta}}(e_1)$ (\ref{notTLG}), we get that $\widetilde{\alpha}(e_1)$ belongs to $Z(Q)$; so $\widetilde{\alpha}(e_2)=1-\widetilde{\alpha}(e_1)$ belongs also to $Z(Q)$, and, as $\widetilde{\beta}(e_1)=R_Q(\widetilde{\alpha}(e_1))$ and $\widetilde{\beta}(e_2)=R_Q(\widetilde{\alpha}(e_2))$, we get (i). 
\newline
We have seen in \ref{notTLG2} that $\pi_{\Phi_{\widehat{Q}}}(\widetilde{\alpha}(e_1)=p_{H_{\Phi_{\widehat{P}}}}+p^{1,2}_{H_{\psi_1}}$ (as we shall consider that $Q$ is acting on $H_{\Phi_{\widehat{Q}}}$, we shall now skip the representation $\pi_{\Phi_{\widehat{Q}}}$). Using now 
the formula obtained for $J_{\widehat{Q}}$, we obtain :
\[\tilde{\alpha}(e_1)=p_{H_{\Phi_{\widehat{Q}}}}+p^{1,2}_{H_{\psi_1}}\]
\[\tilde{\alpha}(e_2)=p^{2,1}_{H_{\psi_1}}+p_H\]
\[\tilde{\beta}(e_1)=p_{H_{\Phi_{\widehat{Q}}}}+p^{2,1}_{H_{\psi_1}}\] 
\[\tilde{\beta}(e_2)=p^{1,2}_{H_{\psi_1}}+p_H\]
from which we get (ii). 
\newline
Let $W_{\widehat{Q}}$ be the canonical pseudo-multiplicative unitary associated to $\gG_2(A,b,\ga, \Phi, \psi_0)$; then, $Q$ is the weak closure of the linear set generated by all operators of the form $(\omega_{w,v}*id)(W_{\widehat{Q}}^*)$, for all $v\in D(_{\tilde{\alpha}}H_{\Phi_{\widehat{Q}}}, \nu\oplus \psi_0)\cap D((H_{\Phi_{\widehat{Q}}})_{\tilde{\hat{\beta}}}, \nu\oplus \psi_0)$, and $w\in D((H_{\Phi_{\widehat{Q}}})_{\tilde{\hat{\beta}}}, \nu\oplus \psi_0)$. Using now \cite{E5}, 3.10 (ii), we get that, for $A$ in $\gN_{\Phi_{\widehat{P}}}$, $\xi\in D(_\mu H_{\psi_1}, \widehat{\Phi}')$, $\eta\in D(_\mu H_{\psi_1}, \widehat{\Phi}')\cap\mathcal D(\widehat{\Delta}_A^{1/2})$ and $m\in\gN_{\widehat{\Phi}'}$ :
\[p_{H_{\Phi_{\widehat{P}}}}(\omega_{w,v}*id)(W_{\widehat{Q}}^*)p_{H_{\Phi_{\widehat{P}}}}\Lambda_{\Phi_{\widehat{Q}}}(\left(\begin{matrix}{A}&{R^{\mu, \widehat{\Phi}'}(\xi)}\cr{R^{\mu, \widehat{\Phi}'}(\eta)^*}&{m}\end{matrix}\right))\]
is, using \ref{coproduct}(i) and (iii), equal to :
\begin{eqnarray*}
p_{H_{\Phi_{\widehat{P}}}}\Lambda_{\Phi_{\widehat{Q}}}[(\omega_{w,v}\underset{\tilde{N}}{_{\tilde{\hat{\beta}}}*_{\tilde{\alpha}}}id)\Gamma_{\widehat{Q}}(\left(\begin{matrix}{A}&{0}\cr{0}&{0}\end{matrix}\right))]
&=&
\Lambda_{\Phi_{\widehat{P}}}[(\omega_{p_{H_{\Phi_{\widehat{P}}}}w, p_{H_{\Phi_{\widehat{P}}}}v}\underset{A^\ga}{_s*_r}id)\Gamma_{\widehat{P}}(A)]\\
&=&
(\omega_{p_{H_{\Phi_{\widehat{P}}}}w, p_{H_{\Phi_{\widehat{P}}}}v}*id)(W_{\widehat{P}}^*)\Lambda_{\Phi_{\widehat{P}}}(A)
\end{eqnarray*}
from which we get that $Q_{p_{H_{\Phi_{\widehat{P}}}}}=P$. The proof for $Q_{p_H}$ is similar. 
\newline
The same way, we get that :
\[p^{1,2}_{H_{\psi_1}}(\omega_{w,v}*id)(W_{\widehat{Q}}^*)p^{1,2}_{H_{\psi_1}}\Lambda_{\Phi_{\widehat{Q}}}(\left(\begin{matrix}{A}&{R^{\mu, \widehat{\Phi}'}(\xi)}\cr{R^{\mu, \widehat{\Phi}'}(\eta)^*}&{m}\end{matrix}\right))\]
is equal, using \ref{propTLG}(iv) to $(\omega_{p_H w, p^{1,2}_{H_{\psi_1}}v}*id)(\tilde{G}^*)\xi$, and, therefore, using \ref{tildeG}(iv), we get that $A$ is the the weak closure of the linear set generated by all elements of the form $p^{1,2}_{H_{\psi_1}}(\omega_{w,v}*id)(W_{\widehat{Q}}^*)p^{1,2}_{H_{\psi_1}}$. For $Q_{p^{2,1}_{H_{\psi_1}}}$, the proof is the same, using \ref{propTLG2}, which finishes the proof of (iii). 
\newline
The same calculations prove that the restriction of $(p_{H_{\Phi_{\widehat{P}}}}\underset{\tilde{N}}{_{\tilde{\beta}}\otimes_{\tilde{\alpha}}}p_{H_{\Phi_{\widehat{P}}}})W^*_{\widehat{Q}}((p_{H_{\Phi_{\widehat{P}}}}\underset{\tilde{N}^o}{_{\tilde{\alpha}}\otimes_{\tilde{\hat{\beta}}}}p_{H_{\Phi_{\widehat{P}}}})$ to $H_{\Phi_{\widehat{P}}}\underset{\psi_0}{_r\otimes_s}H_{\Phi_{\widehat{P}}}$ is equal to $W^*_{\widehat{P}}$, that the restriction of $(p_H\underset{\tilde{N}}{_{\tilde{\beta}}\otimes_{\tilde{\alpha}}}p_H)W^*_{\widehat{Q}}(p_H\underset{\tilde{N}^o}{_{\tilde{\alpha}}\otimes_{\tilde{\hat{\beta}}}}p_H)$ to $H\underset{\nu^o}{_\alpha\otimes_{\hat{\beta}}}H$ is equal to $\widehat{W}$, and that the restriction of $(p^{1,2}_{H_{\psi_1}}\underset{\tilde{N}}{_{\tilde{\beta}}\otimes_{\tilde{\alpha}}}p^{1,2}_{H_{\psi_1}})W^*_{\widehat{Q}}(p_H\underset{N^o}{_\alpha\otimes_b}p^{1,2}_{H_{\psi_1}})$ to $H\underset{\nu^o}{_\alpha\otimes_b}H_{\psi_1}$ is equal to $\tilde{G}^*$. Then the result (iv) comes from (\cite{E5}, 3.6(ii)) applied to $\widehat{\gG_2(\ga)}$, $\widehat{\gG_1(\ga)}$ and $\gG$, and \ref{G}(iv). 
\newline
For any $X\in Q$, we have $R_Q(X)=J_{\Phi_{\widehat{Q}}}X^*J_{\Phi_{\widehat{Q}}}$ (\cite{E5}, 3.10(v)), and $\tau_t^Q(X)=\Delta_{\Phi_{\widehat{Q}}}^{it}X\Delta_{\Phi_{\widehat{Q}}}^{-it}$ (\cite{E5}, 3.10 (vii)). So, the result about $R_Q$ (resp. $\tau_t^Q$) is then given by the formula about $J_{\Phi_{\widehat{Q}}}$ (resp. $\Delta_{\Phi_{\widehat{Q}}}$) obtained in 
\ref{notTLG2}. Let's look at the automorphism group $\gamma_t^{\widehat{Q}}$; we have, using \ref{notTLG} and \ref{notTLG2} :
\begin{eqnarray*}
\widetilde{\hat{\beta}}(\gamma_t^{\widehat{Q}}(u\oplus n))
&=&
\sigma_t^{\Phi_{\widehat{Q}}}(\widetilde{\hat{\beta}}(u\oplus n))\\
&=&
\sigma_t^{\Phi_{\widehat{Q}}}(s(u)\oplus \hat{\beta}(n))\\
&=&
s(\gamma_t^P(u))\oplus\hat{\beta}(\hat{\gamma}_t(n))\\
&=&
\widetilde{\hat{\beta}}(\gamma_t^P(u)\oplus\hat{\gamma}_t(n))
\end{eqnarray*}
from which we get $\gamma_t^{\widehat{Q}}(u\oplus n)=\gamma_t^P(u))\oplus\hat{\gamma}_t(n)$, and, using \cite {E5}3.10 (vii), $\gamma_t^Q(u\oplus n)=\gamma_t^P(u)\oplus \gamma_t(n)$. 
\end{proof}
 
\subsection{Proposition}
\label{propGalois}
{\it Let $\gG$ be a measured quantum groupoid, and $(A, b, \ga, \phi, \psi_0)$ be a Galois system for $\gG$; let $A\subset \tilde{A}$ be a unital inclusion of von Neumann algebras, and $(b, \tilde{\ga})$ be an action of $\gG$ on $\tilde{A}$; let us suppose that $\tilde{A}^{\tilde{\ga}}=A^\ga$, and that $\tilde{\ga}_{|A}=\ga$; then $\tilde{A}=A$.} 

\begin{proof}
As the restriction of $T_{\tilde{\ga}}$ to $A$ is equal to $T_\ga$, we get clearly that $\tilde{\ga}$ is integrable. Let now $\psi_0$ be a normal faithful semi-finite weight on $A^\ga$, and $\psi_1=\psi_0\circ T_\ga$, $\tilde{\psi_1}=\psi_0\circ T_{\tilde{\ga}}$; clearly, we get that theses two weights are normal faithful semi-finite, and that $\psi_1$ is equal to the restriction of $\tilde{\psi_1}$ to $A$; from which we get that there exists a normal faithful conditional expectation $E$ from $\tilde{A}$ onto $A$, such that $\tilde{\psi_1}=\psi_1\circ E$, and a projection $p$ in $\mathcal L(H_{\tilde{\psi_1}})$ such that $p\Lambda_{\tilde{\psi_1}}(x)=\Lambda_{\tilde{\psi_1}}(Ex)$, for any $x\in\gN_{\tilde{\psi_1}}$; moreover, as $\tilde{\psi_1}$ is $\delta$-relatively invariant and bears the density property (\ref{thintegrable}), we get, using the implementation $V_{\tilde{\psi_1}}$ of $\tilde{\ga}$ recalled in \ref{action}, that, for any $x\in\gN_{\tilde{\psi_1}}$, $\xi\in D(_\alpha H, \nu)$ and $\eta\in D(_\alpha H, \nu)\cap \mathcal D(\delta^{1/2})$ such that $\delta^{1/2}\eta$ belongs to $D(H_\beta, \nu^o)$, we get :
\[\Lambda_{\tilde{\psi_1}}[(id\underset{N}{_b*_\alpha}\omega_{\eta, \xi})\tilde{\ga}(x)]=(id*\omega_{\delta^{1/2}\eta, \xi})(V_{\tilde{\psi_1}})\Lambda_{\tilde{\psi_1}}(x)\]
from which we get that :
\begin{eqnarray*}
(id*\omega_{\delta^{1/2}\eta, \xi})(V_{\tilde{\psi_1}})p\Lambda_{\tilde{\psi_1}}(x)
&=&
(id*\omega_{\delta^{1/2}\eta, \xi})(V_{\tilde{\psi_1}})\Lambda_{\tilde{\psi_1}}(Ex)\\
&=&
\Lambda_{\tilde{\psi_1}}[(id\underset{N}{_b*_\alpha}\omega_{\eta, \xi})\tilde{\ga}(Ex)]\\
&=&
\Lambda_{\tilde{\psi_1}}[E(id\underset{N}{_b*_\alpha}\omega_{\eta, \xi})\tilde{\ga}(Ex)]\\
&=&
p\Lambda_{\tilde{\psi_1}}[(id\underset{N}{_b*_\alpha}\omega_{\eta, \xi})\tilde{\ga}(Ex)]\\
&=&
p(id*\omega_{\delta^{1/2}\eta, \xi})(V_{\tilde{\psi_1}})p\Lambda_{\tilde{\psi_1}}(x)
\end{eqnarray*}
from which we get $(id*\omega_{\delta^{1/2}\eta, \xi})(V_{\tilde{\psi_1}})p=p(id*\omega_{\delta^{1/2}\eta, \xi})(V_{\tilde{\psi_1}})p$. Using now \ref{pia} and \ref{MQG}, we get that $p$ belongs to $\pi_{\tilde{\ga}}(1\underset{N}{_b\otimes_\alpha}\widehat{M}')'$. Turning back to the same calculation, we then get that :
\begin{eqnarray*}
\Lambda_{\tilde{\psi_1}}[(id\underset{N}{_b*_\alpha}\omega_{\eta, \xi})\tilde{\ga}(Ex)]
&=&
(id*\omega_{\delta^{1/2}\eta, \xi})(V_{\tilde{\psi_1}})p\Lambda_{\tilde{\psi_1}}(x)\\
&=&
p(id*\omega_{\delta^{1/2}\eta, \xi})(V_{\tilde{\psi_1}})\Lambda_{\tilde{\psi_1}}(x)\\
&=&
p\Lambda_{\tilde{\psi_1}}[(id\underset{N}{_b*_\alpha}\omega_{\eta, \xi})\tilde{\ga}(x)]\\
&=&
\Lambda_{\tilde{\psi_1}}[E(id\underset{N}{_b*_\alpha}\omega_{\eta, \xi})\tilde{\ga}(x)]
\end{eqnarray*}
from which we get that $\ga\circ E=(E\underset{N}{_b*_\alpha}id)\tilde{\ga}$. Hence, $E\underset{N}{_b*_\alpha}id$ can be extended to a faithful conditional expectation from $\tilde{A}\rtimes_{\tilde{\ga}}\gG$ onto $A\rtimes_\ga\gG$, and we easily get that, for any $X\in \tilde{A}\rtimes_{\tilde{\ga}}\gG$, $\pi_\ga(E\underset{N}{_b*_\alpha}id)(X)=p\pi_{\tilde{\ga}}(X)p$; as the action $\ga$ is Galois by hypothesis, $\pi_\ga$ is faithful, and we then get that $\pi_{\tilde{\ga}}$ is also faithful, and, therefore, that $\tilde{\ga}$ is also Galois. Let $\tilde{G}_{\tilde{\ga}}$ be its Galois unitary, as defined in \ref{defGalois}. 
\newline
Moreover, if $\gT$ is the normal faithful semi-finite operator-valued weight from $A$ onto $b(N)$ such that $\phi=\nu\circ b^{-1}\circ \gT$, we get that $E\circ\gT$ is a normal faithful semi-finite operator-valued weight from $\tilde{A}$ onto $b(N)$, which satisfies $(E\circ\gT\underset{N}{_b*_\alpha}id)\tilde{\ga}=(E\underset{N}{_b*_\alpha}id)\ga\circ \gT=
\ga\circ E\circ \gT$, which gives that $\phi\circ E$ is invariant by $\ga$. 
For all $t\in\mathbb{R}$, using the notations of \ref{defGalois}, we get that $(D\phi\circ E : D\tilde{\psi_1})_t=(D\phi:D\psi_1)_t=\lambda_A^{it^2/2}\delta_A^{it}$, which proves that the modular automorphism groups of $\phi\circ E$ and $\tilde{\psi_1}$ commute, and, therefore, we have obtained that $(\tilde{A}, b, \tilde{\ga}, \phi\circ E, \psi_0)$ is a Galois system for $\gG$. 
\newline
So, using \ref{thP2}, we get that $p\in \widehat{P}$, where $(A^\ga, \widehat{P}, r,s,\Gamma_{\widehat{P}}, T_L^{\widehat{P}}, R_{\widehat{P}}\circ T_L^{\widehat{P}}\circ R_{\widehat{P}}, \psi_0)$ is the measured quantum groupoid $\gG_1(\tilde{A}, b, \tilde{\ga}, \phi\circ E, \psi_0)$; more precisely, using the definition of $p$, we get that $p\in r(A^\ga)'\cap b(N)'$, and that $J_{\tilde{\psi_1}}pJ_{\tilde{\psi_1}}=p$, which gives (\ref{coproduct}(iii)) that $R_{\widehat{P}}(p)=p$, and, therefore, that $p\in s(A^\ga)'$; using now the definition of $G_{\tilde{\alpha}}$ given in \ref{G}, we get that $G_{\tilde{\ga}}(p\underset{A^\ga}{_s\otimes_r}p)=(p\underset{N}{_b\otimes_\alpha}1)G_{\tilde{\ga}}$, which gives then, using again \ref{coproduct}(iii), that $\Gamma_{\widehat{P}}(p)=p\underset{A^\ga}{_s\otimes_r}p\leq 1\underset{A^\ga}{_s\otimes_r}p$. Applying now the faithful operator valued weight $id\underset{A^\ga}{_s*_r}T_L^{\widehat{P}}$ to the positive operator 
$1\underset{A^\ga}{_s\otimes_r}p-\Gamma_{\widehat{P}}(p)$, we obtain that $\Gamma_{\widehat{P}}(p)=1\underset{A^\ga}{_s\otimes_r}p$, and, therefore, that $p\underset{A^\ga}{_s\otimes_r}p=1\underset{A^\ga}{_s\otimes_r}p$, and $p=1$; from which we infer the result. \end{proof}

 \section{Morita equivalence for measured quantum groupoids}
 \label{Morita}
 In that chapter, we begin (\ref{prop1} and \ref{prop2}) by the converse result of \ref{thP2}; starting from a measured quantum groupoid with a basis of the form $N_1\oplus N_2$, we see under which conditions it is a linking measured quantum groupoid between a measured quantum groupoid $\gG_1$ (with basis $N_1$) and a measured quantum groupoid $\gG_2$ (with basis $N_2$). This leads to some technical additional results about the reflected groupoid of a measured quantum groupoid $\gG$ throught some Galois system (\ref{biGalois} and \ref{thleftQ}). Then, we can define Morita equivalence of measured quantum groupoids (\ref{defMorita}), prove it is indeed an equivalence relation (\ref{thMorita}), and give a complete link between Morita equivalence and Galois systems (\ref{thMorita2}). We finish this chapter by giving some examples of Morita equivalences between locally compact quantum groups and measured quantum groupoids (\ref{exLCQG}). 
 
 \subsection{Proposition}
 \label{prop1}
 {\it Let $\gG_{1,2}$ be a measured quantum groupoid with a basis which is a sum $N_1\oplus N_2$, we shall denote $(N_1\oplus N_2, M, \alpha, \beta, \Gamma, T, RTR, \nu_1\oplus\nu_2)$, with a co-inverse $R$; we shall identify $H_{\nu_1\oplus\nu_2}$ with $H_{\nu_1}\oplus H_{\nu_2}$; let us denore by $e_1$ the unit of $N_1$, considered as a projection in $N_1\oplus N_2$, and $e_2=1-e_1$. Let us suppose that $\alpha(e_1)$ belongs to $Z(M)$; let us denote $\alpha_i$ (resp. $\beta_i$) the restriction of $\alpha$ (resp. $\beta$) to $N_i$ $(i=1,2)$. Let us write $\Phi=(\nu_1\oplus\nu_2)\circ\alpha^{-1}\circ T$, and $H=H_\Phi$. Let us write $M_{i,j}=M_{\alpha(e_i)\beta(e_j)}$. Let $\widehat{M}$ be the underlying von Neumann algebra of the dual measured quantum groupoid $\widehat{\gG_{1,2}}$. Then, we have  :
 \newline
(i) the projections $\alpha(e_2)$, $\beta(e_1)$ and $\beta(e_2)$ belong to $Z(M)$. Moreover,  the projection $\alpha(e_1)$ belongs to $Z(\widehat{M})$ if and only if $\alpha(e_1)\beta(e_2)= 0$.
 \newline
 (ii) if $\eta$ belongs to $D(_\alpha H, \nu_1\oplus\nu_2)$, then $\alpha(e_i)\eta$ belongs to $D(_{\alpha_i}H, \nu_i)$ $(i=1,2)$, and :
 \[R^{\alpha, \nu_1\oplus\nu_2}(\eta)=R^{\alpha_1, \nu_1}(\alpha(e_1)\eta)\oplus R^{\alpha_2, \nu_2}(\alpha(e_2)\eta)\]
 and, for $\eta_1$, $\eta_2$ in $D(_\alpha H, \nu)$ :
 \[<\eta_1, \eta_2>_{\alpha, \nu_1\oplus\nu_2}^o=<\alpha(e_1)\eta_1, \alpha(e_1)\eta_2>^o_{\alpha_1, \nu_1}e_1 +<\alpha(e_2)\eta_1, \alpha(e_2)\eta_2>^o_{\alpha_2, \nu_2}e_2\]
 \newline
 (iii) the application which sends (for $\eta\in D(_\alpha H, \nu_1\oplus\nu_2)$ and $\xi$ in $H$) the vector $\xi\underset{\nu_1\oplus\nu_2}{_\beta\otimes_\alpha}\eta$ on $(\beta(e_1)\xi\underset{\nu_1}{_{\beta_1}\otimes_{\alpha_1}}\alpha(e_1)\eta)\oplus(\beta(e_2)\xi\underset{\nu_2}{_{\beta_2}\otimes_{\alpha_2}}\alpha(e_2)\eta)$ extends to an isometry, 
 which leads to the identification of $H\underset{\nu_1\oplus\nu_2}{_\beta\otimes_\alpha}H$ with :
 \[(\beta(e_1)H\underset{\nu_1}{_{\beta_1}\otimes_{\alpha_1}}\alpha(e_1)H)\oplus 
( \beta(e_2)H\underset{\nu_2}{_{\beta_2}\otimes_{\alpha_2}}\alpha(e_2)H)\]
 and, for all $(i,j)=1,2$, we have :
 \[\Gamma(\alpha(e_i)\beta(e_j))=[\alpha(e_i)\beta(e_1)\underset{N_1}{_{\beta_1}\otimes_{\alpha_1}}\alpha(e_1)\beta(e_j)]\oplus [\alpha(e_i)\beta(e_2)\underset{N_2}{_{\beta_2}\otimes_{\alpha_2}}\alpha(e_2)\beta(e_j)]\]
 and $R(M\alpha(e_i)\beta(e_j))=M\alpha(e_j)\beta(e_i)$. 
 \newline
(iii) for $(i,j)=1,2$, let us write $M_{i,j}=M_{\alpha(e_i)\beta(e_j)}$; we can define $*$-anti-isomorphisms $R_{i,j}$ from $M_{i,j}$ onto $M_{j,i}$ by writing $R_{i,j}(x_{\alpha(e_i)\beta(e_j)})=R(x)_{\alpha(e_j)\beta(e_i)}$. So, 
we get that $M_{j,i}$ is isomorphic to $M_{i,j}^o$; using (i), we get that $M_{1,2}\neq\{0\}$ if and only if $\alpha(e_1)$ does not belong to $Z(\widehat{\gG_{1,2}})$. Moreover, we can define, for all $x\in M$ :
 \[\Gamma_{i,j}^1(x_{\alpha(e_i)\beta(e_j)})=\Gamma(x)_{\alpha(e_i)\beta(e_1)\underset{N_1}{_{\beta_1}\otimes_{\alpha_1}}\alpha(e_1)\beta(e_j)}\]
 \[\Gamma_{i,j}^2(x_{\alpha(e_i)\beta(e_j)})=\Gamma(x)_{\alpha(e_i)\beta(e_2)\underset{N_2}{_{\beta_2}\otimes_{\alpha_2}}\alpha(e_2)\beta(e_j)}\]
which satisfies, for $k=1,2$, for any $n_i\in N_i$, and $n_j\in N_j$ :
\[\Gamma_{i,j}^k(\alpha_i(n_i))=\alpha_i(n_i)\underset{N_k}{_{\beta_k}\otimes_{\alpha_k}}1\]
\[\Gamma_{i,j}^k(\beta_j(n_j))=1\underset{N_k}{_{\beta_k}\otimes_{\alpha_k}}\beta_j(n_j)\]
and $\Gamma_{i,j}^k$ are normal injective $*$-homomorphisms from $M_{i,j}$ into $M_{i,k}\underset{N_k}{_{\beta_k}*_{\alpha_k}}M_{k,j}$. 
These homomorphisms satisfy :
\[(\Gamma_{i,i}^i\underset{N_i}{_{\alpha_i}*_{\beta_i}}id)\Gamma_{i,i}^i=(id\underset{N_i}{_{\alpha_i}*_{\beta_i}}\Gamma_{i,i}^i)\Gamma_{i,i}^i\]
\[(\Gamma_{i,j}^j\underset{N_j}{_{\beta_j}*_{\alpha_j}}id)\Gamma_{i,j}^j=(id\underset{N_j}{_{\beta_j}*_{\alpha_j}}\Gamma_{j,j}^j)\Gamma_{i,j}^j\]
\[(\Gamma_{i,i}^i\underset{N_i}{_{\beta_i}*_{\alpha_i}}id)\Gamma_{i,j}^i=(id\underset{N_i}{_{\beta_i}*_{\alpha_i}}\Gamma_{i,j}^i)\Gamma_{i,j}^i\]
\[(\Gamma_{i,j}^i\underset{N_j}{_{\beta_j}*_{\alpha_j}}id)\Gamma_{i,j}^j=(id\underset{N_i}{_{\beta_i}*_{\alpha_i}}\Gamma_{i,j}^j)\Gamma_{i,j}^i\]
and, therefore, $(N_i, M_{i,i}, \alpha_i, \beta_i, \Gamma_{i,i}^i)$ is a Hopf-bimodule, with $R_{i,i}$ as a co-inverse, and, if $\alpha(e_1)$ does not belong to $Z(\widehat{\gG_{1,2}})$, $(\beta_j, \Gamma_{i,j}^j)$ is an action of $(N_j, M_{j,j}, \alpha_j, \beta_j, \Gamma_{j,j}^j)$ on $M_{i,j}$, and $(\alpha_i, \Gamma_{i,j}^i)$ is a left-action of $(N_i, M_{i,i}, \alpha_i, \beta_i, \Gamma_{i,i}^i)$ on $M_{i,j}$. Moreover, these two actions commute.}

\begin{proof}
As $\alpha(e_2)=1-\alpha(e_1)$, $\beta(e_1)=R(\alpha(e_1))$, $\beta(e_2)=R(\alpha(e_2))$, the beginning of (i) is clear. If $\alpha(e_1)$ belongs to $Z(\widehat{\gG_{1,2}})$, we have $\alpha(e_1)=\beta(e_1)$, and, therefore, $\alpha(e_1)\beta(e_2)=0$. Conversely, if $\alpha(e_1)\beta(e_2)=0$, we have $\alpha(e_1)\leq \beta(e_1)$, and $\alpha(e_1)=\alpha(e_1)\beta(e_1)$. Applying $R$, we get $\beta(e_1)=\alpha(e_1)\beta(e_1)$, and therefore $\alpha(e_1)=\beta(e_1)$, from which we get that $\alpha(e_1)$ belongs to $Z(\widehat{\gG_{1,2}})$, which finishes the proof of (i). Result (ii)  and (iii) are straightforward. 
 \end{proof}

\subsection{Proposition}
\label{prop2}
{\it Let's use the notations of \ref{prop1};  then :
\newline
(i) let us remark that, for any $(i,j)=1,2$, we have $T(M\alpha(e_i)\beta(e_j))=\alpha(N_i)$, and $RTR(M\alpha(e_i)\beta(e_j))=\beta(N_j)$; this leads to define normal semi-finite faithful operator valued weights $T_{i,j}$ from $M_{i,j}$ onto $\alpha_i(N_i)$, and $T'_{i,j}=R_{j,i}T_{j,i}R_{i,j}$ from $M_{i,j}$ onto $\beta_j(N_j)$. Moreover, the left-invariance of $T$ (resp. the right-invariance of $RTR$) gives then the following formulae, for any $x_{i,j}\in M_{i,j}^+$ :
\[(id\underset{N_1}{_{\beta_1}*_{\alpha_1}}T_{1,1})\Gamma_{1,1}^1(x_{1,1})=T_{1,1}(x_{1,1})\]
\[(T'_{1,1}\underset{N_1}{_{\beta_1}*_{\alpha_1}}id))\Gamma_{1,1}^1(x_{1,1})=T'_{1,1}(x_{1,1})\]
\[(id\underset{N_2}{_{\beta_2}*_{\alpha_2}}T_{2,2})\Gamma_{2,2}^2(x_{2,2})=T_{2,2}(x_{2,2})\]
\[(T'_{2,2}\underset{N_2}{_{\beta_2}*_{\alpha_2}}id)\Gamma_{2,2}^2(x_{2,2})=T'_{2,2}(x_{2,2})\]
\[(id\underset{N_2}{_{\beta_2}*_{\alpha_2}}T_{2,2})\Gamma_{1,2}^2(x_{1,2})=T_{1,2}(x_{1,2})\]
\[(T'_{1,2}\underset{N_2}{_{\beta_2}*_{\alpha_2}}id)\Gamma_{1,2}^2(x_{1,2})=T'_{1,2}(x_{1,2})\]
\[(id\underset{N_1}{_{\beta_1}*_{\alpha_1}}T_{1,1})\Gamma_{2,1}^1(x_{2,1})=T_{2,1}(x_{2,1})\]
\[(T'_{2,1}\underset{N_1}{_{\beta_1}*_{\alpha_1}}id)\Gamma_{2,1}^1(x_{2,1})=T'_{2,1}(x_{2,1})\]
from which we get that $T_{1,1}$ (resp. $T'_{1,1}$) is left-invariant (resp. right-invariant) with respect to $\Gamma_{1,1}^1$, that $T_{2,2}$ (resp. $T'_{2,2}$) is left-invariant (resp. right-invariant) with respect to $\Gamma_{2,2}^2$, and that (if $\alpha(e_1)$ does not belong to $Z(\widehat{\gG_{1,2}})$) both actions $\Gamma_{1,2}^2$ and $\Gamma_{2,1}^1$ are integrable and have invariant weights. 
\newline
(ii) let us define :
\[\Phi_1=\nu_1\circ\alpha_1^{-1}\circ T_{1,1}\]
\[\Psi_1=\nu_1\circ\beta_1^{-1}\circ T'_{1,1}\]
\[\Phi_2=\nu_2\circ\alpha_2^{-1}\circ T_{2,2}\]
\[\Psi_2=\nu_2\circ\beta_2^{-1}\circ T'_{2,2}\]
\[\psi_{1,2}=\nu_1\circ\alpha_1^{-1}\circ T_{1,2}\]
\[\phi_{1,2}=\nu_2\circ\beta_2^{-1}\circ T'_{1,2}\]
\[\psi_{2,1}=\nu_2\circ\alpha_2^{-1}\circ T_{2,1}\]
\[\phi_{2,1}=\nu_1\circ\beta_1^{-1}\circ T'_{2,1}\]
The fact that $\nu_1\oplus\nu_2$ is relatively invariant leads to the commutation of $\sigma^{\Phi_1}$ and $\sigma^{\Psi_1}$, of $\sigma^{\Phi_2}$ and $\sigma^{\Psi_2}$, of $\sigma^{\psi_{1,2}}$ and $\sigma^{\phi_{1,2}}$, and of $\sigma^{\psi_{2,1}}$ and $\sigma^{\phi_{2,1}}$. 
\newline
(iii) $\gG_1=(N_1, M_{1,1}, \alpha_1, \beta_1, \Gamma_{1,1}^1, T_{1,1}, T'_{1,1}, \nu_1)$ and $\gG_2=(N_2, M_{2,2}, \alpha_2, \beta_2, \Gamma_{2,2}, T_{2,2}, T'_{2,2}, \nu_2)$ are two measured quantum groupoids. Moreover, $R_{i,i}$ is the co-inverse of $\gG_i$. 
\newline
(iv) if $\alpha(e_1)\in Z(\widehat{\gG_{1,2}})$, then $\gG_{1,2}=\gG_1\oplus\gG_2$; if $\alpha(e_1)$ does not belong to $Z(\widehat{\gG_{1,2}})$, then $(M_{1,2}, \beta_2, \Gamma_{1,2}^2, \phi_{1,2}, \nu_1)$ is a Galois system for $\gG_2$ (and $\gG_1$ is the measured quantum groupoid reflected from $\gG_2$ throught this Galois system), and $(M_{2,1}, \beta_1, \Gamma_{2,1}^1, \phi_{2,1}, \nu_2)$ is a Galois system for $\gG_1$ (and $\gG_2$ is the measured quantum groupoid reflected from $\gG_1$ throught this Galois system). Moreover, the left action $(\alpha_1, \Gamma_{1,2}^1)$ of $\gG_1$ on $M_{1,2}$ leads (\ref{action}) to an action $(\alpha_1, \varsigma_{N_1}\Gamma_{1,2}^1)$ of $\gG_1^o$ on $M_{1,2}^o$, which, by the identification of $M_{2,1}$ with $M_{1,2}^o$ made in \ref{prop1}(iii), is equal to $(\beta_1, \Gamma_{2,1}^1)^o$.}

\begin{proof}
Results (i), (ii), (iii) are straightforward. If $\alpha(e_1)\in Z(\widehat{\gG_{1,2}})$, then $M_{1,2}=M_{2,1}=\{0\}$, and we get that  $\gG_{1,2}=\gG_1\oplus\gG_2$ (in the sense of \ref{exMQG}(v)). Otherwise, we had got in (i) that $(\beta_2, \Gamma_{1,2}^2)$ is an integrable action of $\gG_2$ on $M_{1,2}$, with an invariant normal faithful semi-finite weight $\phi_{1,2}$; moreover, the invariant algebra $M_{1,2}^{\Gamma_{1,2}^2}$ is $\alpha_1(N_1)$, and the modular automorphism group of the lifted weight $\psi_{1,2}$ commutes with the modular automorphism group $\phi_{1,2}$, which gives that $\nu_1$ is $\Gamma_{1,2}^2$-relatively invariant, in the sense of \ref{defGalois}. Therefore, to get that $(M_{1,2}, \beta_2, \Gamma_{1,2}^2, \phi_{1,2}, \nu_1)$ is a Galois system for $\gG_2$, we have only to prove that the Galois homomorphism $\pi_{\Gamma_{1,2}^2}$ is faithful, or, equivalently (\ref{G*}(iv)), that the isometry $G$ constructed in \ref{G} from $\Gamma_{1,2}^2$ is surjective. As $\Gamma_{1,2}^2$ is "part of " $\Gamma$, we get, using \ref{exGalois}(ii) that $\varsigma G$ is the restriction and co-restriction of
\[(\alpha(e_1)\beta(e_2)\underset{\nu}{_{\hat{\beta}}\otimes_\alpha}\alpha(e_2)\beta(e_2))W^*(\alpha(e_1)\beta(e_2)\underset{\nu^o}{_\alpha\otimes_\beta}\alpha(e_2)\beta(e_2))\]
which is a unitary. The proof for $(\beta_1, \Gamma_{2,1}^1)$ is identical. \end{proof}

\subsection{Theorem}
\label{biGalois}
{\it Let $\gG$ be a measured quantum groupoid, $(A,b, \ga, \phi, \psi_0)$ a Galois system for $\gG$, and $\gG_1$ be the measured quantum groupoid reflected from $\gG$ throught $(A,b, \ga, \phi, \psi_0)$; let's use the notations of \ref{dual}; then, for $z\in A$, let us write :
\[\gb(z)=\Gamma_Q(z)_{\tilde{\alpha}(e_1)\tilde{\beta}(e_1)\underset{\tilde{N}}{_{\tilde{\beta}}\otimes_{\tilde{\alpha}}}\tilde{\alpha}(e_1)\tilde{\beta}(e_2)}\]
 Then, $\gb(z)$ belongs to $P\underset{A^\ga}{_{\hat{s}}*_r}A$, with $\hat{s}(x)=J_{\Phi_{\widehat{P}}}r(x)^*J_{\Phi_{\widehat{P}}}$ for all $x\in A^\ga$, and $(r, \gb)$ is a left action of $\gG_1$ on $A$, with $A^\gb=b(N)$; the left action $(r, \gb)$ commutes with $\ga$, and leads to a Galois sytem for $\gG_1$. }

\begin{proof}
Let us denote $(\tilde{N}, Q, \tilde{\alpha}, \tilde{\beta}, \Gamma_Q, T_L^Q, R_QT_L^QR_Q, \psi_0\oplus\nu)$ the linking measured quantum groupoid between $\gG$ and $\gG_1$, as in \ref{dual}. Then, the result comes from \ref{prop2}(iv). \end{proof}

\subsection{Theorem}
\label{thleftQ}
{\it Let $\gG$ be a measured quantum groupoid, $(A,b, \ga, \phi, \psi_0)$ a Galois system for $\gG$, and $\gG_1$ be the measured quantum groupoid reflected from $\gG$ throught $(A,b, \ga, \phi, \psi_0)$, and $\gG_2$ the linking measured groupoid between $\gG$ and $\gG_1$; we have, for $x\in P^+$, $y\in M^+$, $z_1$, $z_2$ in $A^+$ :
\newline
(i) $T_L^Q(x\oplus z_1\oplus z_2^o\oplus y)=(T_L(x)+r\circ T_\ga(z_1))\oplus (T_L(y)+\alpha\circ b^{-1}\gT(z_2))$, where $\gT$ is the normal semi-finite faithful operator-valued weight from $A$ onto $b(N)$ defined by $\phi=\nu\circ b^{-1}\circ\gT$. 
\newline
(ii) $\delta_Q=\delta_P\oplus\delta_A\oplus (\delta_A^{-1})^o\oplus \delta$
\newline
(iii) $\lambda_P=\lambda_A=r(b\circ\beta^{-1}(\lambda))$}

\begin{proof}
Using \ref{prop2}(i) to the measured quantum groupoid $\gG_2$, we see that the application $x\mapsto T_L^Q(x)$ is a left-invariant weight on $(P, \Gamma_P)$; moreover, using \ref{dual}(v), we get that, for all $t\in\mathbb{R}$, we have ${\tau_t^Q}_{|P}=\tau_t^P$ and $(\gamma_t^Q)_{|A^\ga}=\gamma_t^P$; therefore, we can use Lesieur's theorem (\cite{L}, 5.21), and we get that there exists a non-singular positive operator $h$ affiliated to $Z(A^\ga)$ such that, for all $x\in P^+$, we have $T_L^Q(x)=T_L^P(r(h)x)$. And, therefore, we have then $\Phi_Q(x)=\Phi_P(r(h)x)$; but using now the link between $W_{\widehat{Q}}^*$ and $W_{\widehat{P}}^*$ found in \ref{dual}, we get, using \cite{E5}, 3.10(v), that, for an operator of the form $x=(\omega*id)(W_{\widehat{P}}^*)$, with $\omega\in I_{\Phi_{\widehat{P}}}$ (with the notations of \cite{E5}, 3.10 (v)), we have $\Phi_Q(x^*x)=\Phi_P(x^*x)$, from which we infer that $h=1$, and $T_L^Q(x)=T_L^P(x)$, for all $x\in P^+$. 
\newline
The fact that $T_L^Q(y)=T_L(y)$, for all $y\in M^+$, is proved by similar arguments. 
\newline
Using now \ref{dual}(iv) and \ref{prop2}(i), we get that $T_L^Q(z_1)=T_\ga(z_1)$; we have obtained that $\Phi_Q(x)=\Phi_P(x)$, $\Phi_Q(y)=\Phi(y)$, $\Phi_Q(z_1)=\psi_1(z_1)$, and using \ref{dual}(v), that $\Phi_Q\circ R_Q(x)=\Phi_P\circ R_P(x)$, $\Phi_Q\circ R_Q(y)=\Phi\circ R(y)$, $\Phi_Q\circ R_Q(z_2)=\psi_1(z_2)$. 
\newline
Let's look now at the operator $P_Q^{it}$ which is the canonical implementation of $\tau_t^Q$; using the results obtained for $\Phi_Q$ and for $\tau_t^Q$ (\ref{prop2}(v)), we easily get that $(P_Q)_{H_{\Phi_{\widehat{P}}}}=P_P$, $(P_Q)_{H_{\Phi}}=P$, and, using \ref{tau}(v), that $(P_Q)_{H_{\psi_1}^{1,2}}=P_A$. With same arguments, we get that $(\lambda_Q)_{H_{\Phi_{\widehat{P}}}}=\lambda_P$, $(\lambda_Q)_{H_{\Phi}}=\lambda$ and $(\lambda_Q)_{H_{\psi_1}^{1,2}}=\lambda_A$. But using now \cite{E5}, 3.10 (vii), and the result about $\Delta_{\widehat{\Phi}}$ obtained in \ref{notTLG2}, we get that $(\delta_Q)_{H_{\Phi_{\widehat{P}}}}=\delta_P$, $(\delta_Q)_{H_{\Phi}}=\delta$ and, using \ref{deltahat}(i), that $(\delta_Q)_{H_{\psi_1}^{1,2}}=\delta_A$.  
\newline
So, we get, for all $t\in\mathbb{R}$, using \ref{defGalois} :
\[(D(\Phi_Q\circ R_Q)_{|A}:D(\Phi_Q)_{|A})_t=\lambda_A^{it^2/2}\delta_A^{it}=(D\phi:D\psi_1)_t\]
from we we infer that $\Phi_Q\circ R_Q (z_1)=\phi(z_1)$, for all positive $z_1$ in $A$; so, we have $\Phi_Q(z_2^o)=\phi(z_2)$ for all positive $z_2$ in $A$, from which we finish the proof of (i). 
\newline
Now we have :
\[(D(\Phi_Q\circ R_Q)_{|A^o}:D(\Phi_Q)_{|A^o})_t=(D\Phi^o:D\psi_1^o)_t=[(\lambda_A)^o]^{it^2/2}[(\delta_A)^o]^{-it}\]
from which we get (ii). Finally, there is $p\in Z(N)$ such that $\lambda=\alpha(p)=\beta(p)$, and $u\in Z(A^\ga)$ such that $\lambda_P=r(u)=\hat{s}(u)$; on the other hand, there are $q\in Z(N)$ and $v\in Z(A^\ga)$ such that $\lambda_Q=r(v)\oplus\alpha(q)=\hat{s}(v)\oplus \beta(q)$. From all our calculations above, we infer that $q=p$, $v=u$, $\lambda_A=r(v)$ and $\lambda_A^o=\alpha(p)$; from which we get (iii). \end{proof}

\subsection{Definition} 
\label{defMorita}
For $i=1,2$, let $\gG_i=(N_i, M_i, \alpha_i, \beta_i, T_i, T'_i, \nu_i)$ be a measured quantum groupoid. we shall say that $\gG_1$ is Morita equivalent to $\gG_2$ if there exists a von Neumann algebra $A$, a Galois action $(b, \ga)$ of $\gG_1$ on $A$, a Galois left-action $(a, \gb)$ of $\gG_2$ on $A$, such that :
\newline
(i) $A^\ga=a(N_2)$, $A^\gb=b(N_1)$, and the actions $(b, \ga)$ and $(a, \gb)$ commute; 
\newline
(ii) the modular automorphisms groups of the normal semi-finite faithful weights $\nu_1\circ b^{-1}\circ T_\gb$ and $\nu_2\circ a^{-1}\circ T_\ga$ commute. 
\newline
Then $A$ (or, more precisely, $(A, b, \ga, a, \gb)$) will be called the imprimitivity bi-comodule for $\gG_1$ and $\gG_2$. 

\subsection{Remark}
\label{remMorita}
Then, using \ref{biaction}, we get that the system $(A, b, \ga, \nu_1\circ b^{-1}\circ T_\gb, \nu_2\circ a^{-1})$ is Galois for $\gG_1$ and that the system $(A, a, \gb, \nu_2\circ a^{-1}\circ T_\ga, \nu_1\circ b^{-1})$ is left-Galois for $\gG_2$. Therefore, we can construct, following \ref{dual}, the reflected measured quantum groupoid $\tilde{\gG_2}$ of $\gG_1$ throught the Galois system $(A, b, \ga, \nu_1\circ b^{-1}\circ T_\gb, \nu_2\circ a^{-1})$, and the reflected measured quantum groupoid $\tilde{\gG_1}$ of $\gG_2$ throught the left-Galois system $(A, a, \gb, \nu_2\circ a^{-1}\circ T_\ga, \nu_1\circ b^{-1})$, and, using \ref{biGalois}, an action $\tilde{\ga_1}$ of $\tilde{\gG_1}$ on $A$, and a left-action of $\tilde{\gG_2}$ on $A$; let us first remark that the basis of $\tilde{\gG_2}$ is $A^\ga=a(N_2)$ and is therefore isomorphic to $N_2$ which is the basis of $\gG_2$. Similarly, the $\gG_1$ and $\tilde{\gG_1}$ has the same basis. 
\newline
As the action $\tilde{\ga_1}$ is Galois, the homomophism $\pi_{\tilde{\ga_1}}$ is an isomorphism from the crossed product $A\rtimes_{\tilde{\ga}}\tilde{\gG_1}$ onto the algebra $\tilde{A_2}$ constructed by basic construction made from the inclusion $A^{\tilde{\ga}}\subset A$; as $A^{\tilde{\ga}}=a(N_2)=A^\ga$, we get that $\tilde{A_2}$ is equal to the algebra $s(A^\ga)'$ constructed by basic construction made from the inclusion $A^\ga\subset A$, which is isomorphic, via $\pi_\ga^{-1}$, to $A\rtimes_\ga\gG_1$; therefore, there exists an isomorphism $\mathcal I_1$ from $A\rtimes_\ga\gG_1$ onto $A\rtimes_{\tilde{\ga}}\tilde{\gG_1}$ such that $\mathcal I_1\circ\ga=\tilde{\ga}$; similarly, there exists an isomorphism $\mathcal I_2$ from $A\ltimes_\gb\gG_2$ onto $A\ltimes_{\tilde{\gb}}\tilde{\gG_2}$ such that $\mathcal I_2\circ\gb=\tilde{\gb}$; using all these remarks, we easily get that $(A, b, \tilde{\ga}, a, \tilde{\gb})$ is an imprimitivity bi-comodule between $\tilde{\gG_1}$ and $\tilde{\gG_2}$; we can prove also that if $A$ is an imprimitivity bi-comodule for $\tilde{\gG_1}$ and $\tilde{\gG_2}$, it is also an imprimitivity bi-comodule for $\gG_1$ and $\gG_2$.

\subsection{Theorem}
\label{thMorita}
{\it Morita equivalence is indeed an equivalence relation}
\begin{proof}
Using the Galois system $(M, \beta, \Gamma, \Phi\circ R, \nu)$(\ref{exGalois}(ii)), we get the left-Galois system $(M, \alpha, \Gamma, \Phi, \nu)$, and that $\gG$ is Morita equivalent to $\gG$, with $M$ as imprimitivity bi-comodule; so, Morita equivalence is indeed reflexive. 
\newline
If $\gG_1$ is Morita equivalent to $\gG_2$, with $A$ as imprimitivity co-bimodule, we get, using \ref{defGalois}, that $(b^o, \sigma_N\ga^o)$ and $(a^o, (\sigma_N\gb)^o)$ make $\gG_2$ be Morita equivalent to $\gG_1$, with $A^o$ as imprimitivity co-bimodule; so, Morita equivalence is indeed symmetric. 
\newline
Let us suppose now that $\gG_1$, $\gG_2$, $\gG_3$ are three measured quantum groupoids, and that $(A_1, b_1, \ga_1, a_1, \gb_1)$ is an imprimitivity bi-comodule for $\gG_1$ and $\gG_2$, and that $(A_2, b_2, \ga_2, a_2, \gb_2)$ is an imprimitivity bi-comodule for $\gG_2$ and $\gG_3$. Using 
\ref{remMorita}, we know there exists an action $(b_1, \tilde{\ga_1})$ of the reflected measured quantum groupoid $\tilde{\gG_1}$ of $\gG_2$ throught the Galois system $(A_1, a_1, \gb_1, \nu_2\circ a_1^{-1}\circ T_\ga, \nu_1\circ b_1^{-1})$ such that $(A_1, b_1, \tilde{\ga_1}, a_1, \gb_1)$ is an imprimitivity bi-comodule between $\tilde{\gG_1}$ and $\gG_2$; similarly, we shall consider $(A_2, b_2, \ga_2, a_2, \tilde{\gb_2})$ which is an imprimitivity bi-comodule between $\gG_2$ and the reflected measured quantum groupoid $\tilde{\gG_3}$ of $\gG_2$ throught the left Galois system  $(A_2, a_2, \gb_2, \nu_3\circ a_2^{-1}\circ T_\ga, \nu_2\circ b_2^{-1})$. 
\newline
Let $A_3=\{X\in A_2\underset{N_2}{_{b_2}*_{a_2}}A_1; (id\underset{N_2}{_{b_2}*_{a_2}}\gb_1)(X)=(\ga_2\underset{N_2}{_{b_2}*_{a_2}}id)(X)\}$. It is straightforward to check that $a_2(N_3)\underset{N_2}{_{b_2}*_{a_2}}1\subset A_3$ and $1\underset{N_2}{_{b_2}*_{a_2}}b_1(N_1)\subset A_3$, and that :
\newline
(i) $(1\underset{N_2}{_{b_2}*_{a_2}}b_1, (id\underset{N_2}{_{b_2}*_{a_2}}\tilde{\ga_1})_{|A_3}$ is an action of $\tilde{\gG_1}$ on $A_3$, we shall denote $(b_3, \ga_3)$ for simplification.
\newline
(ii) $(a_2\underset{N_2}{_{b_2}*_{a_2}}1, (\tilde{\gb_2}\underset{N_2}{_{b_2}*_{a_2}}id)_{|A_3})$ is a left action of $\tilde{\gG_3}$ on $A_3$, we shall denote $(a_3, \gb_3)$ for simplification.
\newline
(iii) we have $A_3^{\ga_3}=a_3(N_3)$, and $A_3^{\gb_3}=b_3(N_1)$, and the actions $\ga_3$ and $\gb_3$ commute. 
\newline
To prove that we get an imprimitivity system, we shall make a detour. 
\newline
So, let us consider a Galois system for $\gG_2$, with $\tilde{\gG_1}$ as reflected measured quantum groupoid, and  another Galois system for $\gG_2$, with $\tilde{\gG_3}$ as reflected measured quantum groupoid. Let us consider now, as in \ref{notTLG}, the representation $\mu_1$ of $\widehat{M_2}'$ on $H_1$ and the representation $\mu_3$ of $\widehat{M_2}'$ on $H_3$, and the representation $\varpi$ of $\widehat{M_2}'$ on $H_3\oplus H_2\oplus H_1$ given by $\mu_3\oplus id\oplus \mu_1$, and $\widehat{Q}=\varpi(\widehat{M}')'$; using again matrix notations for elements in $\widehat{Q}$, we get that :
\[\widehat{Q}=\left(\begin{matrix}{\widehat{\tilde{M_3}}}&{\widehat{Q}_{2,3}}&{\widehat{Q}_{1,3}}\cr{\widehat{Q}_{2,3}^*}&{\widehat{M_2}}&{\widehat{Q}_{1,2}}\cr{\widehat{Q}_{1,3}^*}&{\widehat{Q}_{1,3}^*}&{\widehat{\tilde{M_1}}}\end{matrix}\right)\]
where, for instance :
\[\widehat{Q}_{1,3}=\{X\in \mathcal L(H_1, H_3), X\mu_1(m)=\mu_3(m)X, \forall m\in\widehat{M_2}'\}\]
We have clearly $\widehat{Q}_{2,3}\widehat{Q}_{1,2}\subset \widehat{Q}_{1,3}$; using again an orthogonal basis as in the proof of \ref{coproduct}(i), we get that the linear set generated by the products in $\widehat{Q}_{2,3}\widehat{Q}_{1,2}$ is weakly dense in $\widehat{Q}_{1,3}$. But, as in \ref{lemTLG}, we can construct a coproduct from $\widehat{Q}_{1,2}$ into $\widehat{Q}_{1,2}\underset{N_1}{_{\hat{\beta_1}}*_{\alpha_1}}\widehat{Q}_{1,2}$ and a coproduct from $\widehat{Q}_{2,3}\underset{N_2}{_{\hat{\beta_2}}*_{\alpha_2}}\widehat{Q}_{2,3}$, and, by product, we obtain therefore a coproduct from $\widehat{Q}_{1,3}$ into $\widehat{Q}_{1,3}\underset{N_1}{_{\hat{\beta_1}}*_{\alpha_1}}\widehat{Q}_{1,3}$, then, as in \ref{coproduct}, a coproduct for $\widehat{Q}$. The proof that $\widehat{Q}$ bears a structure of measured quantum groupoid is completely similar to \ref{propTLG}, \ref{propTLG2} and \ref{thP}. So, as in \ref{dual}, we can look at the dual measured quantum groupoid, which will be on the basis $N_1\oplus N_2\oplus N_3$; let us denote $\alpha_Q$ and $\beta_Q$ the canonical homomorphism and antihomomorphism from $N_1\oplus N_2\oplus N_3$ into $Q$; as in \ref{dual}, we can prove that $\alpha_Q(e_i)\in Z(Q)$ and $\beta_Q(e_i)\in Z(Q)$, where, for $(i=1,2,3)$, $e_i$ is the unit of $N_i$, considered as a projection in $N_1\oplus N_2\oplus N_3$. 
Then, it is easy to get that the reduced algebra on $H_1\oplus H_3$ bears a structure of measured quantum groupoid, over the basis $N_1\oplus N_3$. As $\widehat{Q}_{1,3}\neq \{0\}$, we can use 
\ref{prop2}(iv), and we get the existence of a Galois system for $\tilde{\gG_3}$, with $\tilde{\gG_1}$ as reflected measured quantum groupoid, which means that $\tilde{\gG_3}$ is Morita equivalent to $\tilde{\gG_1}$ (and, by the reflexivity, that $\tilde{\gG_1}$ is Morita equivalent to $\tilde{\gG_3}$); using them arguments analoguous to \ref{remMorita}, we get that $\gG_1$ is Morita equivalent to $\gG_3$, which proves the transitivity. To get the imprimitivity bi-comodule, we must look at the dual $Q=\oplus_{i,j=1}^3Q_{i,j}$, which bears a coproduct $\Gamma_Q$, which can be split into maps $(\Gamma_Q)_{i,j}^k: Q_{i,j}\mapsto Q_{i,k}*Q_{k,j}$. 
\newline
We know that $Q_{1,1}=\tilde{M_1}$, $Q_{2,2}=M_2$, $Q_{3,3}=M_3$, $Q_{2,1}=A_1$, $Q_{3,2}=A_2$, and we are looking for $Q_{3,1}$. We know also that $(\Gamma_Q)_{1,1}^1=\tilde{\Gamma_1}$, $(\Gamma_Q)_{2,2}^2=\Gamma_2$, $(\Gamma_Q)_{3,3}^3=\tilde{\Gamma_3}$, $(\Gamma_Q)_{2,1}^1=\tilde{\ga_1}$, $(\Gamma_Q)_{2,1}^2=\gb_1$, $(\Gamma_Q)_{3,2}^2=\ga_2$, $(\Gamma_Q)_{3,2}^3=\tilde{\gb_2}$. 
\newline
So, $(\Gamma_Q)_{3,1}^2$ sends $Q_{3,1}$ into $A_2*A_1$, and it is easy, with the co-associativity condition of $\Gamma_Q$, to get that $(\Gamma_Q)_{3,1}^2$ sends $Q_{3,1}$ into $A_3$, and that $(\Gamma_Q)_{3,1}^2$ sends the action $(\Gamma_Q)_{3,1}^1$ on $id\underset{N_2}{_{b_2}*_{a_2}}\tilde{\ga_1}$ and the left action $(\Gamma_Q)_{3,1}^3$ on $\tilde{\gb_2}\underset{N_2}{_{b_2}*_{a_2}}id$; using then \ref{propGalois}, we get that $A_3$ is the image of $(\Gamma_Q)_{3,1}^2$, which allow us to identify $Q_{3,1}$ with $A_3$, $(\Gamma_Q)_{3,1}^1$ with $\ga_3$, and $(\Gamma_Q)_{3,1}^3$ with $\gb_3$. By these identifications, we prove that $(A_3, \ga_3, \gb_3)$ is an imprimitivity bi-comodule between $\tilde{\gG_1}$ and $\tilde{\gG_3}$. By similar arguments to \ref{remMorita}, we get an imprimitivity bi-comodule between $\gG_1$ and $\gG_3$. 
\end{proof}

\subsection{Notations}
\label{notcomp}
Let $\gG_1$, $\gG_2$, $\gG_3$ be three measured quantum groupoids; let us suppose that $\gG_1$ is Morita equivalent to $\gG_2$, with $(A_1, \ga_1, \gb_1)$ (or $A_1$ for simplification) as imprimitivity bi-comodule and that $\gG_2$ is Morita equivalent to $\gG_3$ with $(A_2, \ga_2, \gb_2)$ (or simply $A_2$) as imprimitivity co-bimodule; we have proved in \ref{thMorita} that $\gG_1$ is Morita equivalent to $\gG_3$, with $(A_3, \ga_3, \gb_3)$ as imprimitivity co-bimodule, with :
\[A_3=\{X\in A_2\underset{N_2}{_{b_2}*_{a_2}}A_1; (id\underset{N_2}{_{b_2}*_{a_2}}\gb_1)(X)=(\ga_2\underset{N_2}{_{b_2}*_{a_2}}id)(X)\}\]
and $\ga_3=(id\underset{N_2}{_{b_2}*_{a_2}}\ga_1)_{|A_3}$, $\gb_3=(\gb_2\underset{N_2}{_{b_2}*_{a_2}}1)_{|A_3}$. 
\newline
We shall write $(A_3, \ga_3, \gb_3)=(A_2, \ga_2, \gb_2)\circ (A_1, \ga_1, \gb_1)$, or, simply $A_3=A_2\circ A_1$; we can check that this product is associative, and that, if we write $M_1$ for the imprimitivity bi-comodule $(M_1, \Gamma_1, \Gamma_1)$ between $\gG_1$ and itself, we easily get that $A_1\circ M_1=A_1$ and $M_2\circ A_1=A_1$. 

\subsection{Proposition}
\label{propsMorita}
{\it Let $\gG$, $\gG_1$, $\gG_2$ be measured quantum groupoids; let us use the notations of \ref{notcomp}; 
\newline
(i) let us suppose that $\gG_1$ is Morita equivalent to $\gG_2$ with an imprimitivity co-bimodule $A$; then, we have $A^o\circ A=M_1$. 
\newline
(ii) let us suppose that $\gG$ is Morita equivalent to $\gG$ with an imprimitivity co-bimodule $A$; then $A=M$.
\newline
(iii) let us suppose that $\gG_1$ is Morita equivalent to $\gG_2$ with an imprimitivity co-bimodule $A_1$, and with another imprimitivity co-bimodule $A_2$; then $A_1=A_2$. 
\newline
(iv) let us suppose that $\gG_1$ is Morita equivalent to $\gG_2$ with an imprimitivity bi-comodule $(A, \ga, \gb)$; then $\gG_2$ is the reflected measured quantum groupoid of $\gG_1$ throught the Galois system $(A, b, \ga, \nu_1\circ b^{-1}\circ T_\gb, \nu_2\circ a^{-1})$.}

\begin{proof}
 Let us use the Galois system $(A, b, \ga, \nu_1\circ b^{-1}\circ T_\gb, \nu_2\circ a^{-1})$, and apply the constructions and results of \ref{thP2} applied to this Galois system; for any $y\in M_1$, the operator $\Gamma_Q(y)_{\tilde{\alpha}(e_2)\tilde{\beta}(e_1)\underset{\tilde{N}}{_{\tilde{\beta}}\otimes_{\tilde{\alpha}}}\tilde{\alpha}(e_1)\tilde{\beta}(e_2)}$ belongs to $A^o\underset{N_2}{_{a_2^o}*_{a_2}}A$, and, more precisely, using the coassociativity of the coproduct $\Gamma_Q$, we can check it belongs to the subagebra $A^o\circ A$; we define this way an injective morphism from $M_1$ into $A^o\circ A$, which sends $\Gamma_1$ on the action (and on the left-action) canonically defined on $A^o\circ A$; therefore, using \ref{propGalois}, we get (i). 
 \newline
 Let us now use the Galois system $(A, b, \ga, \nu\circ b^{-1}\circ T_\gb, \nu\circ a^{-1})$, and apply the constructions and results of \ref{thP2} to this Galois system; for $x\in A$, the operator $\Gamma_Q(x)_{\tilde{\alpha}(e_2)\tilde{\beta}(e_1)\underset{\tilde{N}}{_{\tilde{\beta}}\otimes_{\tilde{\alpha}}}\tilde{\alpha}(e_1)\tilde{\beta}(e_2)}$ belongs to $A^o\circ A$, and, therefore, using (i), to $M$; we define this way an injective morphism from $A$ into $M$, which sends $\ga$ on $\Gamma$; using again \ref{propGalois}, we get (ii). 
 \newline
 As $A_2^o\circ A_1$ is an imprimitivity bi-comodule for a Morita equivalence between $\gG_1$ and $\gG_1$, we get, using (ii), that $A_2^o\circ A_1=M_1$; therefore, we have, using (i) :
 \[A_1=M_2\circ A_1=A_2\circ A_2^o\circ A_1=A_2\circ M_1=A_2\]
 which is (iii). 
 \newline
 Let $\tilde{\gG_2}$ be the reflected measured quantum groupoid of $\gG_1$ throught the Galois system $(A, b, \ga, \nu_1\circ b^{-1}\circ T_\gb, \nu_2\circ a^{-1})$; there exists a left-action $\tilde{\gb}$ of $\tilde{\gG_2}$ on $A$, and $\tilde{A}=(A, \ga, \tilde{\gb})$ is an imprimitivity bi-comodule which makes $\gG_1$ and $\tilde{\gG_2}$; therefore, using \ref{notcomp}, we get that $ \tilde{A}\circ A^o$ (whose underlying von Neumann algebra is $M_2$ by (i), and that we shall denote by $P$) is an imprimitivity bi-comodule between $\tilde{\gG_2}$ and $\gG_2$; we then get, using again (i),  that $P^o\circ P=\tilde{M_2}$ and $P\circ P^o=M_2$, which leads, using \ref{thP2}, to define injective morphisms $M_2\mapsto \tilde{M_2}$ and $\tilde{M_2}\mapsto M_2$ as parts of the coproduct of the same measured quantum groupoid. Using then the co-associativity of this coproduct, we get that these applications are each other's inverse, which leads to the isomorphism of $M_2$ and $\tilde{M_2}$, which is (iv). 
 \end{proof}
 
 \subsection{Theorem}
 \label{thMorita2}
 {\it Let $\gG_i=(N_i, M_i, \alpha_i, \beta_i, \Gamma_i, T_i, T'_i, \nu_i)$ ($i=1,2$) be two measured quantum groupoids; then, are equivalent :
 \newline
 (i) $\gG_1$ and $\gG_2$ are Morita equivalent, with a imprimitivity bi-comodule $(A, \ga, \gb)$; 
 \newline
 (ii) there exists a Galois system $(A, b, \ga, \phi, \psi_0)$ for $\gG_1$, such that $\gG_2$ is the reflected measured groupoid of $\gG_1$ throught this Galois system; 
 \newline
 (iii) there exists a measured quantum groupoid $\gG_{1,2}=(N_1\oplus N_2, M, \alpha, \beta, \Gamma, T, T', \nu_1\oplus\nu_2)$ such that $\alpha(e_1)$ belongs to $Z(M)$, and does not belong to $Z(\widehat{M})$, where $e_1$ is the unit of $N_1$, considered as a projection in $N_1\oplus N_2$, and $(\gG_{1,2})_{\alpha(e_1)}=\gG_1$,  $(\gG_{1,2})_{\alpha(1-e_1)}=\gG_2$.}
 
 \begin{proof}
 The result (i) implies (ii) by \ref{propsMorita}(iv); the result (ii) implies (i) was obtained in \ref{biGalois}; the result (ii) implies (iii) is given by \ref{thP2}(ii), and \ref{prop1} gives that (iii) implies (ii). 
 \end{proof}
 
 \subsection{Remark}
 \label{remMorita}
 A morphism between an action $(b_1, \ga_1)$ of $\gG$ on a von Neumann algebra $A_1$, and an action $(b_2, \ga_2)$ on a von Neumann algebra $A_2$ will be a$*$-homomorphism $h$ from $A_1$ in $A_2$, such that $h\circ b_1=b_2$, and $(h\underset{N}{_{b_1}*_\alpha}id)\ga_1=\ga_2$; clearly this leads to a category $\mathcal A(\gG)$; it is easy to get that, if $\gG_1$ and $\gG_2$ are two measured quantum groupoids which are Morita equivalent, then these categories $\mathcal A(\gG_1)$ and $\mathcal A(\gG_2)$ are equivalent too. 

 \subsection{Examples of locally compact quantum groups Morita equivalent to measured quantum groupoids}
 \label{exLCQG}
 Here we are looking to examples of locally compact quantum groups which are Morita equivalent to measured quantum groupoids. I am indebted to S. Vaes who put my attention on this question.  We first begin to give two constructions in which any locally compact quantum group is Morita equivalent to a measured quantum groupoid, whose basis is a given factor $N$ (\ref{thampliation}, \ref{dualampliation}). More convincing is K. De Commer's example (\ref{KDC}, \cite{DC4}) : he proves that the compact quantum $SU_q(2)$ is Morita equivalent to some measured quantum groupoid (whose basis is a finite sum of type I factors). 
  
 \subsubsection{Ampliation of a locally compact quantum group.} 
 \label{ampliation}
If $\bf{G}$$=(M, \Gamma, \varphi, \psi)$ is a locally compact quantum group, and $N$ is a von Neumann algebra, we shall call the measured quantum groupoid $\gG(N)\otimes\bf{G}$ the ampliation of $\bf{G}$ by $N$, where $\gG(N)$ is the $N$-measured quantum groupoid defined in \ref{exMQG}(viii) and the tensor product of measured quantum groupoids had been defined in \ref{exMQG}(ix). Morover, the measured quantum groupoid $\widehat{\gG(N)}\otimes\bf{G}$ is, using also \ref{exMQG}(viii) and (ix), another measured quantum groupoid, we shall call the dual ampliation of $\gG$ by $N$. 

\subsubsection{Theorem}
\label{thampliation}
{\it Let $\bf{G}$$=(M, \Gamma, \varphi, \psi)$ be a locally compact quantum group, $N$ a factor, $\gG(N)\otimes\bf{G}$ the ampliation of $\gG$ by $N$, as defined in \ref{ampliation}. Then, the locally compact quantum group $\bf{G}$ and the measured quantum groupoid $\gG(N)\otimes\bf{G}$ are Morita equivalent. }

\begin{proof}
Let us consider the von Neumann algebra $N\otimes M$; then, $(id\otimes\Gamma)$ is an action of $\bf{G}$ on this algebra; we get that the invariant subalgebra is $(N\otimes M)^{(id\otimes\Gamma)}=N\otimes \mathbb{C}$, and that the crossed-product is $N\otimes\mathcal L(H_\varphi)$. Therefore, we get also that this action is Galois, and that $T_{id\otimes\Gamma}=id\otimes\varphi$. Let us chose a normal semi-finite faithful trace $\nu$ on $N$; we get 
$\nu\circ T_{id\otimes\Gamma}=\nu\otimes\varphi$. 
\newline
Taking now on this algebra the restriction of the coproduct of $\gG(N)\otimes\bf{G}$, we obtain a left action $\gb$ of $\gG(N)\otimes\bf{G}$ on $N\otimes M$, and we get that $T_\gb=\nu\otimes\psi$. (Taking for $\tau$ the canonical finite trace on $\mathbb{C}=Z(N)$, we get that the operator-valued weight $T_\nu$ defined in \ref{exMQG}(viii) is $\nu$). So, we then get that $\gb$ is ergodic and Galois. Moreover, as the modular groups of $\varphi$ and $\psi$ commute, we get, by the definition (\ref{defMorita}) that the locally compact quantum group $\bf{G}$ and the measured quantum groupoid $\gG(N)\otimes\bf{G}$ are Morita equivalent, with $N\otimes M$ as imprimitivity bi-comodule. 
\end{proof}

\subsubsection{Proposition}
\label{dualampliation}
{\it Let $\bf{G}$$=(M, \Gamma, \varphi, \psi)$ be a locally compact quantum group, $N$ a factor, $\widehat{\gG(N)}\otimes\bf{G}$ the dual ampliation of $\gG$ by $N$, as defined in \ref{ampliation}. Then, the locally compact quantum group $\bf{G}$ and the measured quantum groupoid $\widehat{\gG(N)}\otimes\bf{G}$ are Morita equivalent. }

\begin{proof}
The proof is very similar to \ref{thampliation}. \end{proof}

\subsubsection{Another Example}
\label{KDC}
In \cite{DC2}, \cite{DC3}, Kenny De Commer had studied Morita equivalences between the compact quantum group $SU_q(2)$ and various quantum groups, and, in \cite{DC4}, with a mesurable quantum groupoid. Indeed, he constructs an integrable Galois action of a $SU_q(2)$, which is not ergodic (the subalgebra of invariants is then a finite sum of type I factors), and, therefore, this construction leads to measured quantum groupoid (whose basis is that finite sum of factors), which is Morita equivalent to the initial compact quantum group. This construction is a particular case of \ref{corsum}.

\section{Application to deformation of a measured quantum groupoid by a 2-cocycle}
\label{cocycles}
In this section, we try to answer the problem of deformation of a measured quantum groupoid by a 2-cocycle. With this deformed coproduct constructed in \ref{propcocycle}, does this new Hopf-bimodule still has a left-invariant (and a right-invariant) Haar operator-valued weight, and therefore remains a measured quantum groupoid ? Following De Commer's strategy, we are able to answer positively to this question for any 2-cocycle only in the case when the basis $N$ is a finite sum of factors (\ref{corG}(xii)). In the general case, we can obtain (\ref{thGaloiscocycle}) sufficient conditions, which leads to positive answers in particular cases (\ref{thcas}, \ref{thpart}). 

\subsection{Definition}
\label{defcocycle}
Let $(N, M, \alpha, \beta, \Gamma)$ be a Hopf-bimodule, in the sense \ref{MQG}; a unitary $\Omega$ in $(M\cap\alpha(N)'\underset{N}{_\beta*_\alpha}(M\cap\beta(N)')$ is called a 2-cocycle for 
 $(N, M, \alpha, \beta, \Gamma)$ if $\Omega$ satisfies the following relation :
 \[(1\underset{N}{_\beta\otimes_\alpha}\Omega)(id\underset{N}{_\beta*_\alpha}\Gamma)(\Omega)=
 (\Omega\underset{N}{_\beta\otimes_\alpha}1)(\Gamma\underset{N}{_\beta*_\alpha}id)(\Omega)\]

If $\mathcal G$ is a measured qroupoid, equipped with a left Haar system and a quasi-invariant mesure on the set of units, and if $\Omega$ is a 2-cocycle for the measured quantum groupoid $\gG(\mathcal G)$(\ref{exMQG}(ii)), then $\Omega$ is just a measurable function from $\mathcal G^{(2)}$ to $\mathbb{T}$, such that, for all $(g_1, g_2)$ and $(g_2, g_3)$ in $\mathcal G^{(2)}$ :
\[\Omega(g_2, g_3)\Omega(g_1, g_2g_3)=\Omega(g_1, g_2)\Omega(g_1g_2, g_3)\]
Let $\gG=(N, M, \alpha, \beta, \Gamma, T, T', \nu)$ be a measured quantum groupoid, and let 
$\Omega$ be a 2-cocycle for $\gG$; let us define, for any $t\in\mathbb{R}$ :
\[\Omega_t=(\tau_t\underset{N}{_\beta*_\alpha}\tau_t)(\Omega)\] 
\[\Omega'_t=(\delta^{it}\underset{N}{_\beta\otimes_\alpha}\delta^{it})\Omega(\delta^{-it}\underset{N}{_\beta\otimes_\alpha}\delta^{-it})=(\sigma_t^{\Phi\circ R}\sigma_{-t}^\Phi\underset{N}{_\beta*_\alpha}\sigma_t^{\Phi\circ R}\sigma_{-t}^\Phi)(\Omega)\]
 One can easily check that $\Omega_t$ and $\Omega'_t$ are also 2-cocycles for $\gG$.  
 \subsection{Proposition}
 \label{propcocycle}
 {\it Let $(N, M, \alpha, \beta, \Gamma)$ be a Hopf-bimodule, and let $\Omega$ be a 2-cocycle for 
 $(N, M, \alpha, \beta, \Gamma)$; let us write, for all $x\in M$ :
 \[\Gamma_\Omega (x)=\Omega\Gamma(x)\Omega^*\]
 Then, $(N, M, \alpha, \beta, \Gamma_\Omega)$ is a Hopf-bimodule, that we shall call the deformation of the initial one by $\Omega$. }
 \begin{proof}
 We have, thanks to the definition of a 2-cocycle, for any $n\in N$ :
 \[\Gamma_\Omega (\alpha(n))=\alpha(n)\underset{N}{_\beta\otimes_\alpha}1\]
 \[\Gamma_\Omega (\beta(n))=1\underset{N}{_\beta\otimes_\alpha}\beta(n)\]
 which allows us to write :
 \begin{eqnarray*}
 (\Gamma_\Omega \underset{N}{_\beta*_\alpha}id)\Gamma_\Omega (x)
 &=&
 (\Gamma_\Omega \underset{N}{_\beta*_\alpha}id)(\Omega\Gamma(x)\Omega^*)\\
 &=&
  (\Gamma_\Omega \underset{N}{_\beta*_\alpha}id)(\Omega)(\Gamma_\Omega \underset{N}{_\beta*_\alpha}id)\Gamma(x) (\Gamma_\Omega \underset{N}{_\beta*_\alpha}id)(\Omega)^*
\end{eqnarray*}
  But we have :
\[(\Gamma_\Omega \underset{N}{_\beta*_\alpha}id)(\Omega)
  =
  (\Omega\underset{N}{_\beta\otimes_\alpha}1)(\Gamma\underset{N}{_\beta*_\alpha}id)(\Omega) (\Omega\underset{N}{_\beta\otimes_\alpha}1)^*\]
  and, therefore :
\[(\Gamma_\Omega \underset{N}{_\beta*_\alpha}id)\Gamma_\Omega (x)=  
(\Omega\underset{N}{_\beta\otimes_\alpha}1)(\Gamma\underset{N}{_\beta*_\alpha}id)(\Omega) 
(\Gamma\underset{N}{_\beta*_\alpha}id)\Gamma(x)(\Gamma\underset{N}{_\beta*_\alpha}id)(\Omega) ^*(\Omega\underset{N}{_\beta\otimes_\alpha}1)^*\]
and, by a similar calculation, we get :
\[(id\underset{N}{_\beta*_\alpha}\Gamma_\Omega)\Gamma_\Omega (x)=  (1\underset{N}{_\beta\otimes_\alpha}\Omega)(id\underset{N}{_\beta*_\alpha}\Gamma)(\Omega)(id\underset{N}{_\beta*_\alpha}\Gamma)\Gamma (x)(id\underset{N}{_\beta*_\alpha}\Gamma)(\Omega)^*(1\underset{N}{_\beta\otimes_\alpha}\Omega)^*\]
which is equal, thanks to the definition of a 2-cocycle, and the coassociativity of $\Gamma$. 
\end{proof}
\subsection{Proposition}
\label{algcocycle}
{\it Let $\gG$ be a measured quantum groupoid, and $\Omega$ be a 2-cocycle for $\gG$; let $W$ be the pseudo-multiplicative unitary associated to $\gG$; let us write $\widetilde{W}=W\Omega^*$, which is a unitary from $H\underset{\nu}{_\beta\otimes_\alpha}H$ onto $H\underset{\nu^o}{_\alpha\otimes_{\hat{\beta}}}H$; then :
\newline
(i) the operator $\widetilde{W}$ satisfies :
\[(1\underset{N^o}{_\alpha\otimes_{\hat{\beta}}}\widetilde{W})(\widetilde{W}\underset{N}{_\beta\otimes_\alpha}1)=
(W\underset{N^o}{_\alpha\otimes_{\hat{\beta}}}1)\sigma^{2,3}_{\alpha, \beta}(\widetilde{W}\underset{N}{_{\hat{\beta}}\otimes_\alpha}1)(1\underset{N}{_\beta\otimes_\alpha}\sigma_{\nu^o})(1\underset{N}{_\beta\otimes_\alpha}\widetilde{W})\]
(with the notations of \ref{MQG}). 
\newline
(ii) for all $\xi$, $\xi'$ in $D(H_\beta, \nu^o)$, $\eta$, $\eta'$ in $D(_\alpha H, \nu)$, we have :
\[(\omega_{\xi, \xi'}\underset{N}{_\beta*_\alpha}id)[\widetilde{W}^*((id*\omega_{\eta, \eta'})(\sigma_{\nu^o} W)\underset{N^o}{_\alpha\otimes_{\hat{\beta}}}1)\widetilde{W}]=
\omega_{(\omega_{\xi, \eta'}*id)(\widetilde{W})^*\xi', \eta}*id)(\widetilde{W})^*\]
(iii) the weakly closed linear space generated by the operators of the form $(\omega_{\xi, \eta}*id)(\widetilde{W})$, for all $\xi\in D(H_\beta, \nu^o)$ and $\eta\in D(_\alpha H, \nu)$ is a non degenerate involutive algebra, therefore a von Neumann algebra on $H$, we shall denote $A_\Omega$; 
\newline
(iv) we have $\alpha(N)\subset A_\Omega$, $\hat{\beta}(N)\subset A_\Omega$, and $A_\Omega\subset \beta(N)'$, $A_\Omega\subset \hat{\alpha}(N)'$. 
\newline
(v) a unitary $v$ on $H$ belongs to $A_\Omega '$ if and only if $v\in \alpha(N)'\cap \hat{\beta}(N)'$ and :
\[\widetilde{W}(1\underset{N}{_\beta\otimes_\alpha}v)=(1\underset{N^o}{_\alpha\otimes_{\hat{\beta}}}v)\widetilde{W}\]
(vi) for any $x\in M$, we have :
\[\Gamma_\Omega(x)=\widetilde{W}^*(1\underset{\nu^o}{_\alpha\otimes_{\hat{\beta}}}x)\widetilde{W}\]
and the weakly closed linera space generated by the operators of the form $(id*\omega_{\zeta_1, \zeta_2})(\widetilde{W})$, for $\zeta_1\in D(_\alpha H, \nu)$ and $\zeta_2\in D(H_{\hat{\beta}}, \nu^o)$ is equal to $M$. }
\begin{proof}
We have, using the definition of a $2$-cocycle :
\begin{eqnarray*}
(1\underset{N^o}{_\alpha\otimes_{\hat{\beta}}}\widetilde{W})(\widetilde{W}\underset{N}{_\beta\otimes_\alpha}1)
&=&
(1\underset{N^o}{_\alpha\otimes_{\hat{\beta}}}W)(1\underset{N^o}{_\alpha\otimes_{\hat{\beta}}}\Omega^*)(W\underset{N}{_\beta\otimes_\alpha}1)(\Omega^*\underset{N}{_\beta\otimes_\alpha}1)\\
&=&
(1\underset{N^o}{_\alpha\otimes_{\hat{\beta}}}W)(W\underset{N}{_\beta\otimes_\alpha}1)(W^*\underset{N}{_\beta\otimes_\alpha}1)(1\underset{N^o}{_\alpha\otimes_{\hat{\beta}}}\Omega^*)(W\underset{N}{_\beta\otimes_\alpha}1)(\Omega^*\underset{N}{_\beta\otimes_\alpha}1)\\
&=&
(1\underset{N^o}{_\alpha\otimes_{\hat{\beta}}}W)(W\underset{N}{_\beta\otimes_\alpha}1)(\Gamma\underset{N}{_\beta*_\alpha}id)(\Omega^*)(\Omega^*\underset{N}{_\beta\otimes_\alpha}1)\\
&=&
(1\underset{N^o}{_\alpha\otimes_{\hat{\beta}}}W)(W\underset{N}{_\beta\otimes_\alpha}1)(id\underset{N}{_\beta*_\alpha}\Gamma)(\Omega^*)(1\underset{N}{_\beta\otimes_\alpha}\Omega^*)\\
\end{eqnarray*}
which is equal to :
\[(1\underset{N^o}{_\alpha\otimes_{\hat{\beta}}}W)(W\underset{N}{_\beta\otimes_\alpha}1)(1\underset{N}{_\beta\otimes_\alpha}W^*)(1\underset{N}{_\beta\otimes_\alpha}\sigma_\nu)(\Omega^*\underset{N}{_{\hat{\beta}}\otimes_\alpha}1)(1\underset{N}{_\beta\otimes_\alpha}\sigma_{\nu^o})(1\underset{N}{_\beta\otimes_\alpha}W)(1\underset{N}{_\beta\otimes_\alpha}\Omega^*)\]
and, using the pentagonal equation for $W$, is equal to :
\[(W\underset{N^o}{_\alpha\otimes_{\hat{\beta}}}1)\sigma^{2,3}_{\alpha, \beta}(W\Omega^*\underset{N}{_{\hat{\beta}}\otimes_\alpha}1)(1\underset{N}{_\beta\otimes_\alpha}\sigma_{\nu^o})(1\underset{N}{_\beta\otimes_\alpha}W\Omega^*)\]
which is (i). 
\newline
For $\zeta$, $\zeta'$ in $H$, we get that :
\[((\omega_{\xi, \xi'}\underset{N}{_\beta*_\alpha}id)[\widetilde{W}^*[(id*\omega_{\eta, \eta'})(\sigma_{\nu^o}W)\underset{N^o}{_\alpha\otimes_{\hat{\beta}}}1)]\widetilde{W}]\zeta|\zeta')\]
is equal to :
\begin{multline*}
(\widetilde{W}^*[(id*\omega_{\eta, \eta'})(\sigma_{\nu^o}W)\underset{N^o}{_\alpha\otimes_{\hat{\beta}}}1)]\widetilde{W}(\xi\underset{\nu}{_\beta\otimes_\alpha}\zeta)|\xi'\underset{\nu}{_\beta\otimes_\alpha}\zeta')=\\
[(\omega_{\eta, \eta'}*id)(W\sigma_\nu)\underset{N^o}{_\alpha\otimes_{\hat{\beta}}}1)]\widetilde{W}(\xi\underset{\nu}{_\beta\otimes_\alpha}\zeta)|\widetilde{W}(\xi'\underset{\nu}{_\beta\otimes_\alpha}\zeta'))=\\
((W\underset{N^o}{_\alpha\otimes_{\hat{\beta}}}1)(\sigma_\nu\underset{N^o}{_\alpha\otimes_{\hat{\beta}}}1)(\eta\underset{\nu^o}{_\alpha\otimes_{\hat{\beta}}}\widetilde{W}(\xi\underset{\nu}{_\beta\otimes_\alpha}\zeta)|\eta'\underset{\nu^o}{_\alpha\otimes_{\hat{\beta}}}\widetilde{W}(\xi'\underset{\nu}{_\beta\otimes_\alpha}\zeta'))=\\
((1\underset{N^o}{_\alpha\otimes_{\hat{\beta}}}\widetilde{W}^*)(W\underset{N^o}{_\alpha\otimes_{\hat{\beta}}}1)(\sigma_\nu\underset{N^o}{_\alpha\otimes_{\hat{\beta}}}1)(1\underset{N^o}{_\alpha\otimes_{\hat{\beta}}}\widetilde{W})(\eta\underset{\nu^o}{_\alpha\otimes_{\hat{\beta}}}\xi\underset{\nu}{_\beta\otimes_\alpha}\zeta)|\eta'\underset{\nu^o}{_\alpha\otimes_{\hat{\beta}}}\xi'\underset{\nu}{_\beta\otimes_\alpha}\zeta')
\end{multline*}
which, using (i), is equal to :
\[((\widetilde{W}\underset{N}{_\beta\otimes_\alpha}1)(1\underset{N}{_\beta\otimes_\alpha}\widetilde{W}^*)(\xi\underset{\nu}{_\beta\otimes_\alpha}(\eta\underset{\nu^o}{_\alpha\otimes_{\hat{\beta}}}\zeta)|\eta'\underset{\nu^o}{_\alpha\otimes_{\hat{\beta}}}\xi'\underset{\nu}{_\beta\otimes_\alpha}\zeta')=(\xi\underset{\nu}{_\beta\otimes_\alpha}\widetilde{W}^*(\eta\underset{\nu^o}{_\alpha\otimes_{\hat{\beta}}}\zeta)|\widetilde{W}^*(\eta'\underset{\nu^o}{_\alpha\otimes_{\hat{\beta}}}\xi')\underset{\nu}{_{\hat{\beta}}\otimes_\alpha}\zeta')\]
Let $(f_i)_{i\in I}$ be an orthogonal $(\beta, \nu^o)$basis of $H$; there exist $\delta_i$ such that :
\[\widetilde{W}^*(\eta'\underset{\nu}{_\alpha\otimes_{\hat{\beta}}}\xi')=\sum_i f_i\underset{\nu}{_\beta\otimes_\alpha}\delta_i\]
and, as in (\cite{E3}3.11), we can prove that $\sum_i\alpha(<f_i, \xi>_{\beta, \nu^o})\delta_i=(\omega_{\xi, \eta'}*id)(\widetilde{W})^*\xi'$, and, therefore, we get that :
\begin{eqnarray*}
(\xi\underset{\nu}{_\beta\otimes_\alpha}\widetilde{W}^*(\eta\underset{\nu^o}{_\alpha\otimes_{\hat{\beta}}}\zeta)|\widetilde{W}^*(\eta'\underset{\nu^o}{_\alpha\otimes_{\hat{\beta}}}\xi')\underset{\nu}{_{\hat{\beta}}\otimes_\alpha}\zeta')
&=&
(\xi\underset{\nu}{_\beta\otimes_\alpha}\widetilde{W}^*(\eta\underset{\nu^o}{_\alpha\otimes_{\hat{\beta}}}\zeta)|\sum_if_i\underset{\nu}{_\beta\otimes_\alpha}\delta_i\underset{\nu}{_{\hat{\beta}}\otimes_\alpha}\zeta')\\
&=&
(\widetilde{W}^*(\eta\underset{\nu^o}{_\alpha\otimes_{\hat{\beta}}}\zeta)|\sum_i f_i\underset{\nu}{_\beta\otimes_\alpha}\delta_i\underset{\nu}{_{\hat{\beta}}\otimes_\alpha}\zeta')\\
&=&
(\widetilde{W}^*(\eta\underset{\nu^o}{_\alpha\otimes_{\hat{\beta}}}\zeta)|(\omega_{\xi, \eta'}*id)(\widetilde{W})^*\xi'\underset{\nu}{_{\hat{\beta}}\otimes_\alpha}\zeta')
\end{eqnarray*}
from which we get (ii). 
\newline
We have :
\begin{eqnarray*}
((\omega_{\xi, \eta}*id)(\widetilde{W})(\omega_{\xi', \eta'}*id)(\widetilde{W})\zeta|\zeta')
&=&
(\widetilde{W}(\xi\underset{\nu}{_\beta\otimes_\alpha} (\omega_{\xi', \eta'}*id)(\widetilde{W})\zeta)|\eta\underset{\nu^o}{_\alpha\otimes_{\hat{\beta}}}\zeta')\\
&=&
(\xi\underset{\nu}{_\beta\otimes_\alpha} (\omega_{\xi', \eta'}*id)(\widetilde{W})\zeta|\widetilde{W}^*(\eta\underset{\nu^o}{_\alpha\otimes_{\hat{\beta}}}\zeta'))
\end{eqnarray*}
which is equal to :
\begin{multline*}
((1\underset{N}{_\beta\otimes_\alpha}\sigma_{\nu^o}\widetilde{W})(\xi\underset{\nu}{_\beta\otimes_\alpha}\xi'\underset{\nu}{_\beta\otimes_\alpha}\zeta)|\widetilde{W}^*(\eta\underset{\nu^o}{_\alpha\otimes_{\hat{\beta}}}\zeta')\underset{\nu}{_{\hat{\beta}}\otimes_\alpha}\eta')\\
=
((\widetilde{W}\underset{N}{_{\hat{\beta}}\otimes_\alpha}1)(1\underset{N}{_\beta\otimes_\alpha}\sigma_{\nu^o}\widetilde{W})(\xi\underset{\nu}{_\beta\otimes_\alpha}\xi'\underset{\nu}{_\beta\otimes_\alpha}\zeta)|(\eta\underset{\nu^o}{_\alpha\otimes_{\hat{\beta}}}\zeta')\underset{N}{_\beta\otimes_\alpha}\eta')\\
=
(\sigma^{2,3}_{\alpha, \beta}(\widetilde{W}\underset{N}{_{\hat{\beta}}\otimes_\alpha}1)(1\underset{N}{_\beta\otimes_\alpha}\sigma_{\nu^o})(1\underset{N}{_\beta\otimes_\alpha}\widetilde{W})(\xi\underset{\nu}{_\beta\otimes_\alpha}\xi'\underset{\nu}{_\beta\otimes_\alpha}\zeta)|(\eta\underset{\nu}{_\beta\otimes_\alpha}\eta')\underset{\nu^o}{_\alpha\otimes_{\hat{\beta}}}\zeta')
\end{multline*}
and, using (i), is equal to :
\[((1\underset{N^o}{_\alpha\otimes_{\hat{\beta}}}\widetilde{W})(\widetilde{W}\underset{N}{_\beta\otimes_\alpha}1)(\xi\underset{\nu}{_\beta\otimes_\alpha}\xi'\underset{\nu}{_\beta\otimes_\alpha}\zeta)|(W(\eta\underset{\nu}{_\beta\otimes_\alpha}\eta')\underset{\nu^o}{_\alpha\otimes_{\hat{\beta}}}\zeta')\]
Let now $(e_i)_{i\in I}$ be an orthogonal $(\alpha, \nu)$-basis. As in (\cite{E3}, 3.4), we can prove that there exist a family $(\eta_i)_{i\in I}$ in $D(_\alpha H, \nu)$, such that :
\[W(\eta\underset{\nu}{_\beta\otimes_\alpha}\eta')=\sum_ie_i\underset{\nu^o}{_\alpha\otimes_{\hat{\beta}}}\eta_i\]
and a family $(\xi_i)_{i\in I}$ in $D(H_\beta, \nu^o)$, such that :
\[\widetilde{W}(\xi\underset{\nu}{_\beta\otimes_\alpha}\xi')=\sum_ie_i\underset{\nu^o}{_\alpha\otimes_{\hat{\beta}}}\xi_i\]
and we get that :
\[((\omega_{\xi, \eta}*id)(\widetilde{W})(\omega_{\xi', \eta'}*id)(\widetilde{W})\zeta|\zeta')
=(\sum_i(\omega_{\xi_i, \eta_i}*id)(\widetilde{W})\zeta|\zeta')\]
from which we get that $(\omega_{\xi, \eta}*id)(\widetilde{W})(\omega_{\xi', \eta'}*id)(\widetilde{W})$ is equal to the weak limit of the finite sums $\sum_1^n(\omega_{\xi_i, \eta_i}*id)(\widetilde{W})$. From which we get that the weakly closed linear set generated by the operators of the form $(\omega_{\xi, \eta}*id)(\widetilde{W})$ is an algebra $A_\Omega$. 
\newline
Let us now use (ii); the weak regularity of the pseudo-multiplicative unitary $W$ (\cite{E5}, 3.8) means that $\alpha(N)'$ is the closed linear space generated by the operators $(id*\omega_{\eta, \eta'})(\sigma_{\nu^o}W)$, for all $\eta$, $\eta'$ in $D(_\alpha H, \nu)$ (\cite{E3}, 4.1); in particular, 
there exists a family in the linear space generated by these operators  which weakly converges to $1$. Using then (ii), we get that, for any $\xi$, $\xi'$ in $D(H_\beta, \nu^o)$, there exists a family in the linear space generated by the operators of the form $(\omega_{(\omega_{\xi, \eta'}*id)(\widetilde{W})^*\xi', \eta}*id)(\widetilde{W})$, for all $\eta$, $\eta'$ in $D(_\alpha H, \nu)$,  which is weakly converging to $\alpha(<\xi, \xi'>_{\beta, \nu^o})$; therefore, we get by density 
that $\alpha(N)$ is included in $A_\Omega$, and, therefore, that $1$ belongs to $A_\Omega$. 
\newline
So, there exists a family of operators of the form (with finite sums) $\sum_i(\omega_{\xi_i, \eta'_i}*id)(\widetilde{W})$ which is weakly converging to $1$. But it is straightforward, using the intertwining properties of $W$ and the definition of $\Omega$, that 
$(\omega_{\xi, \eta'}*id)(\widetilde{W})^*$ commutes with $\beta(N)$, and we get that $R^{\beta, \nu^o}(\sum_i
(\omega_{\xi_i, \eta'_i}*id)(\widetilde{W})^*\xi')=\sum_i(\omega_{\xi_i, \eta'_i}*id)(\widetilde{W})^*R^{\beta, \nu^o}(\xi')$ is converging to $R^{\beta, \nu^o}(\xi')$; finally, we get that $A_\Omega$ is the weakly closed linear set generated by  all operators of the type 
$(\omega_{(\omega_{\xi, \eta'}*id)(\widetilde{W})^*\xi', \eta}*id)(\widetilde{W})$, for all $\xi$, $\xi'$ in $D(H_\beta, \nu^o)$ and $\eta$, $\eta'$ in $D(_\alpha H, \nu)$; using again (ii) and the weak regularity of $W$, we get that $A_\Omega$ is closed under the involution,  which finishes the proof of (iii).  
\newline
For any $n\in N$, we have, using (\cite{E5}, 3.2) :
\begin{eqnarray*}
((\omega_{\xi, \eta}*id)(\widetilde{W})\hat{\beta}(n)\zeta_1|\zeta_2)
&=&
(\widetilde{W}(\xi\underset{\nu}{_\beta\otimes_\alpha}\hat{\beta}(n)\zeta_1)|\eta\underset{\nu^o}{_\alpha\otimes_{\hat{\beta}}}\zeta_2)\\
&=&
(\widetilde{W}(\xi\underset{\nu}{_\beta\otimes_\alpha}\zeta_1)|\beta(n^*)\eta\underset{\nu^o}{_\alpha\otimes_{\hat{\beta}}}\zeta_2)\\
&=&
((\omega_{\xi, \beta(n^*)\eta}*id)(\widetilde{W})\zeta_1|\zeta_2)
\end{eqnarray*}
from which we get that $(\omega_{\xi, \eta}*id)(\widetilde{W})\hat{\beta}(n)=(\omega_{\xi, \beta(n^*)\eta}*id)(\widetilde{W})$, which gives that $\hat{\beta}(n)$ belongs to $A_\Omega$. 
\newline
We had seen that $\beta(n)$ commutes with $(\omega_{\xi, \eta}*id)(\Omega W^*)$; using then (\cite{E5}, 3.11 (iii)), we get also that $\hat{\alpha}(n)$ commutes with $(\omega_{\xi, \eta}*id)(\Omega W^*)$, which finishes the proof of (iv). 
\newline
Then, using (iii), the proof of (v) is clear. 
\newline
It is clear that $(id*\omega_{\zeta_1, \zeta_2})(\widetilde{W})$ belongs to $M$. Let us denote $M_\Omega$ the closed linear set generated by these operators. Using (i), we get, for $\zeta'_1\in D(_\alpha H, \nu)$ and $\zeta'_2\in D(H_{\hat{\beta}}, \nu^o)$, that :
\[(id*\omega_{\zeta'_1, \zeta'_2})(W)(id*\omega_{\zeta_1, \zeta_2})(\widetilde{W})
=
[id*\omega_{\widetilde{W}^*(\zeta'_1\underset{\nu}{_\alpha\otimes_{\hat{\beta}}}\zeta_1), \widetilde{W}^*(\zeta'_2\underset{\nu}{_\alpha\otimes_{\hat{\beta}}}\zeta_2)}](\widetilde{W}\underset{N}{_\beta\otimes_\alpha}1)\]
which belongs to $M_\Omega$. By linearity and weak closure, we get that $M_\Omega$ is a left ideal of $M$. 
\newline
Moreover, the formula $\Gamma_\Omega(x)=\widetilde{W}^*(1\underset{\nu^o}{_\alpha\otimes_{\hat{\beta}}}x)\widetilde{W}$ is clear by the definition of $\Gamma_\Omega$ (\ref{propcocycle}) and $\widetilde{W}$. Using that, we get, for any $(\hat{\beta}, \nu^o)$-orthogonal basis $(e_i)_{i\in I}$ of $H$, and any $\eta\in D(_\alpha H, \nu)$, by taking $x=1$ :
\[\beta(<\eta, \eta>_{\alpha, \nu})=(id\underset{N}{_\beta*_\alpha}\omega_\eta)\Gamma_\Omega(1)=\sum_i(id*\omega_{\eta, e_i})(\widetilde{W})^*(id*\omega_{\eta, e_i})(\widetilde{W})\]
from which we get that $\beta(<\eta, \eta>_{\alpha, \nu})$ belongs to $M_\Omega$; by density, we get that $\beta(N)$ belongs to $M_\Omega$, and therefore, that $1\in M_\Omega$, which finishes the proof. \end{proof}

\subsection{Theorem}
\label{thaction}
{\it Let $\gG$ be a measured quantum groupoid, and $\Omega$ be a 2-cocycle for $\gG$; let $W$ be the pseudo-multiplicative unitary associated to $\gG$, $A_\Omega$ the von Neumann algebra on $H$ defined in \ref{algcocycle}; let us write $\widetilde{W}=W\Omega^*$, and, for any $x\in A_\Omega$, let us write :
\[\ga(x)=W^c(x\underset{N^o}{_{\hat{\alpha}}\otimes_{\hat{\beta}}}1)(W^c)^*\]
Then :
\newline
(i) for any $\xi\in D(H_\beta, \nu^o)$ and $\eta\in D(_\alpha H, \nu)\cap D(H_\beta, \nu^o)$, we have :
\[\ga[(\omega_{\xi, \eta}*id)(\widetilde{W})]=(\omega_{\xi, \eta}*id*id)[(1\underset{N^o}{_\alpha\otimes_{\hat{\beta}}}\sigma_{\nu^o})(W\underset{N^o}{_\alpha\otimes_{\hat{\beta}}}1)\sigma^{2,3}_{\beta, \alpha}(\widetilde{W}\underset{N^o}{_{\hat{\beta}}\otimes_\alpha}1)]\]
(ii) $(\hat{\beta}, \ga)$ is an action of $\widehat{\gG}$ on $A_\Omega$. 
\newline
(iii) this action is integrable and Galois. }
\begin{proof}
Using (\cite{E5}, 5.6), we get that $W^c$ is a corepresentation of $\widehat{\gG}$ on $_{\hat{\alpha}}H_{\hat{\beta}}$. Therefore, the formula, for $y\in\hat{\alpha}(N)'$ :
\[\ga(y)=W^c(y\underset{N^o}{_{\hat{\alpha}}\otimes_{\hat{\beta}}}1)(W^c)^*\]
leads to an action $(\hat{\beta}, \ga)$ of $\widehat{\gG}$ on $\hat{\alpha}(N)'$. Using (\cite{E5}, 3.12), we get, for any $n\in N$ that $\ga(\alpha(n))=\widehat{\Gamma}(\alpha(n))=\alpha(n)\underset{N}{_{\hat{\beta}}\otimes_\alpha}1$, and $\ga(\hat{\beta}(n))=\widehat{\Gamma}(\hat{\beta}(n))=1 \underset{N}{_{\hat{\beta}}\otimes_\alpha}\hat{\beta}(n)$. 
\newline
For any orthogonal $(\beta, \nu^o)$-basis $(e_i)_{i\in I}$ of $H$, we get :
\begin{eqnarray*}
\ga((\omega_{\xi, \eta}*id)(\widetilde{W}))
&=&
W^c ((\omega_{\xi, \eta}*id)(\widetilde{W})\underset{N^o}{_{\hat{\alpha}}\otimes_{\hat{\beta}}}1)W^{c*}\\
&=&
\sum_i W^c[(\omega_{e_i, \eta}\underset{N}{_\beta*_\alpha}id)(W)(\omega_{\xi, e_i}*id)(\Omega^*)\underset{N^o}{_{\hat{\alpha}}\otimes_\beta}1]W^{c*}\\
&=&
\sum_i \widehat{\Gamma}((\omega_{e_i, \eta}\underset{N}{_\beta*_\alpha}id)(W))[(\omega_{\xi, e_i}*id)(\Omega^*)\underset{N}{_{\hat{\beta}}\otimes_\alpha}1]
\end{eqnarray*}
Applying \ref{Gamma} to $\widehat{\gG}$, we get :
\[\widehat{\Gamma}[(id*\omega_{\eta, e_i})(\sigma W^*\sigma)]=(id\underset{N}{_{\hat{\beta}}*_\alpha}id*\omega_{\eta, e_i})[\sigma^{2,3}_{\alpha, \hat{\beta}}(\sigma W^*\sigma\underset{N}{_\beta\otimes_\alpha}1)(1\underset{N}{_{\hat{\beta}}\otimes_\alpha}\sigma_{\nu^o})(1\underset{N}{_{\hat{\beta}}\otimes_\alpha}\sigma W^*\sigma)]\]
from which we get that :
\[\widehat{\Gamma}[(\omega_{\eta, e_i}*id)(W^*)]=(\omega_{\eta, e_i}*id\underset{N}{_{\hat{\beta}}*_\alpha}id)[(W^*\underset{N^o}{_\alpha\otimes_{\hat{\beta}}}1)\sigma^{2,3}_{\hat{\beta}, \alpha}(W^*\underset{N^o}{_\alpha\otimes_{\hat{\beta}}}1)(1\underset{N^o}{_\alpha\otimes_{\hat{\beta}}}\sigma_\nu)]\]
and, for any orthogonal $(\alpha, \nu)$-basis $(f_j)_{j\in J}$ of $H$ :
\[\widehat{\Gamma}[(\omega_{\eta, e_i}*id)(W^*)]=\sum_j(\omega_{f_j, e_i}*id)(W^*)\underset{N}{_\beta\otimes_\alpha}1)(\omega_{\eta, f_j}*id*id)[\sigma^{2,3}_{\hat{\beta}, \alpha}(W^*\underset{N}{_\alpha\otimes_{\hat{\beta}}}1)(1\underset{N^o}{_\alpha\otimes_{\hat{\beta}}}\sigma_\nu)]\]
and, therefore $\ga((\omega_{\xi, \eta}*id)(\widetilde{W}))$ is equal to :
\[\sum_{i,j}[(\omega_{f_j, \eta}*id*id)([(1\underset{N^o}{_\alpha\otimes_{\hat{\beta}}}\sigma_{\nu^o}(W\underset{N^o}{_\alpha\otimes_{\hat{\beta}}}1)\sigma^{2,3}_{\beta, \alpha}])(\omega_{e_i, f_j}*id)(W)\underset{N}{_{\hat{\beta}}\otimes_\alpha}1)(\omega_{\xi, e_i}\underset{N}{_\beta*_\alpha}id)(\Omega^*)\underset{N}{_{\hat{\beta}}\otimes_\alpha}1]\]
which is equal to :
\[\sum_j[(\omega_{f_j, \eta}*id*id)([(1\underset{N^o}{_\alpha\otimes_{\hat{\beta}}}\sigma_{\nu^o}(W\underset{N^o}{_\alpha\otimes_{\hat{\beta}}}1)\sigma^{2,3}_{\beta, \alpha}]))(\omega_{\xi, f_j}*id)(\widetilde{W})\underset{N}{_{\hat{\beta}}\otimes_\alpha}1]\]
from which we get (i). 
\newline
For any $\delta_1$, $\delta_2$ in $D(H_{\hat{\beta}}, \nu^o)$, we get that :
\begin{eqnarray*}
(\omega_{\eta, f_j}*\omega_{\delta_1, \delta_2}*id)[\sigma^{2,3}_{\hat{\beta}, \alpha}(W^*\underset{N}{_\alpha\otimes_{\hat{\beta}}}1)(1\underset{N^o}{_\alpha\otimes_{\hat{\beta}}}\sigma_\nu)]
&=&
(\omega_{\eta, f_j}*id)[(\alpha(<\delta_1, \delta_2>_{\hat{\beta}, \nu^o}\underset{N}{_\beta\otimes_\alpha}1)W^*]\\
&=&
(\omega_{\eta, f_j}*id)(W^*)\alpha(<\delta_1, \delta_2>_{\hat{\beta}, \nu^o})
\end{eqnarray*}
as $\hat{\beta}(N)\subset A_\Omega$, any unitary $u\in A'_\Omega$ commutes with $\hat{\beta}(N)$, and we have :
\[<u\delta_1, u\delta_2>_{\hat{\beta}, \nu^o}=<\delta_1, \delta_2>_{\hat{\beta}, \nu^o}\]
from which we get that $(\omega_{\eta, f_j}*id*id)[\sigma^{2,3}_{\hat{\beta}, \alpha}(W^*\underset{N}{_\alpha\otimes_{\hat{\beta}}}1)(1\underset{N^o}{_\alpha\otimes_{\hat{\beta}}}\sigma_\nu)]$ commutes with $A'_\Omega\underset{N}{_{\hat{\beta}}\otimes_\alpha}1$, and, therefore, belongs to $A_\Omega\underset{N}{_{\hat{\beta}}*_\alpha}\mathcal L(H)$, and, more precisely, to $A_\Omega\underset{N}{_{\hat{\beta}}*_\alpha}\hat{M}$. 
\newline
Then, using \ref{algcocycle}, we easily get that $\ga((\omega_{\xi, \eta}*id)(\widetilde{W}))$ belongs to $A_\Omega\underset{N}{_{\hat{\beta}}*_\alpha}\hat{M}$; using again \ref{algcocycle}, we get $\ga(A_\Omega)\subset A_\Omega\underset{N}{_{\hat{\beta}}*_\alpha}\hat{M}$, which gives (ii). 
\newline
Using (\cite{E5}, 11.2), we know that the von Neumann algebra $A_\Omega\underset{N}{_{\hat{\beta}}*_\alpha}\mathcal L(H)$ is isomorphic to the double crossed-product $(A\rtimes_\ga\widehat{\gG})\rtimes_{\tilde{\ga}}\gG^c$ and that this isomorphism sends the bidual action on the action $\underline{\ga}$ defined, for any $X\in A_\Omega\underset{N}{_{\hat{\beta}}*_\alpha}\mathcal L(H)$ by :
\[\underline{\ga}(X)=(1\underset{N}{_{\hat{\beta}}\otimes_\alpha}W^*)(id\underset{N}{_{\hat{\beta}}*_\alpha}\varsigma_N)(\ga\underset{N}{_{\hat{\beta}}*_\alpha}id)(X)(1\underset{N}{_{\hat{\beta}}\otimes_\alpha}W)\]
Let us define $\mathcal I(X)=\varsigma_N(\widetilde{W}^*)X\varsigma_N( \widetilde{W})$; then $\mathcal I$ is an isomorphism from $A_\Omega\underset{N}{_{\hat{\beta}}*_\alpha}\mathcal L(H)$ onto $A_\Omega\underset{N^o}{_\alpha*_\beta}\mathcal L(H)$, and the above calculations show that :
\[(\mathcal J\underset{N}{_{\hat{\beta}}*_\alpha}id)\underline{\ga}(X)=(id\underset{N}{_{\hat{\beta}}*_\alpha}\varsigma_N)(\ga\underset{N^o}{_\alpha*_\beta}id)(\mathcal I(X))\]
from which we get that $(id\underset{N}{_{\hat{\beta}}*_\alpha}\varsigma_N)(\ga\underset{N^o}{_\alpha*_\beta}id)$ is an integrable and Galois action of $\widehat{\gG}$ on $A_\Omega\underset{N^o}{_\alpha*_\beta}\mathcal L(H)$. As we have $(A_\Omega\underset{N^o}{_\alpha*_\beta}\mathcal L(H))^{(id\underset{N}{_{\hat{\beta}}*_\alpha}\varsigma_N)(\ga\underset{N^o}{_\alpha*_\beta}id)}=A_\Omega^{\ga}\underset{N^o}{_\alpha*_\beta}\mathcal L(H)$, we easily get that $T_{(id\underset{N}{_{\hat{\beta}}*_\alpha}\varsigma_N)(\ga\underset{N^o}{_\alpha*_\beta}id)}=T_\ga\underset{N^o}{_\alpha*_\beta}\mathcal L(H)$; as $T_{(id\underset{N}{_{\hat{\beta}}*_\alpha}\varsigma_N)(\ga\underset{N^o}{_\alpha*_\beta}id)}$ is semi-finite, we get that $T_\ga$ is also semi-finite, and $\ga$ is integrable; moreover, as $(A_\Omega\underset{N^o}{_\alpha*_\beta}\mathcal L(H))\rtimes_{(id\underset{N}{_{\hat{\beta}}*_\alpha}\varsigma_N)(\ga\underset{N^o}{_\alpha*_\beta}id)}\widehat{\gG}=(id\underset{N}{_{\hat{\beta}}*_\alpha}\varsigma_N)[(A\rtimes_\ga\widehat{\gG})\underset{N^o}{_\alpha*_\beta}\mathcal L(H)]$, we get that $\pi_{(id\underset{N}{_{\hat{\beta}}*_\alpha}\varsigma_N)(\ga\underset{N^o}{_\alpha*_\beta}id)}=(\pi_\ga\underset{N^o}{_\alpha*_\beta}id)(id\underset{N}{_{\hat{\beta}}*_\alpha}\varsigma_N)$. As $\pi_{(id\underset{N}{_{\hat{\beta}}*_\alpha}\varsigma_N)(\ga\underset{N^o}{_\alpha*_\beta}id)}$ is injective, we easily get that $\pi_\ga$ is injective also, which is (iii). 
\end{proof}
 
 \subsection{Proposition}
 \label{invariant}
 {\it Let $\gG$ be a measured quantum groupoid, $\Omega$ a $2$-cocycle for $\gG$, let $W$ be the pseudo-multiplicative unitary associated to $\gG$, $A_\Omega$ the von Neumann algebra on $H$ defined in \ref{algcocycle} and $(\hat{\beta}, \ga)$ the action of $\widehat{\gG}$ on $A_\Omega$ defined in \ref{action}; let us write $\widetilde{W}=W\Omega^*$. Then, for $\xi\in D(H_\beta, \nu^o)$ and $\eta\in D(_\alpha H, \nu)$ such that $\omega_{\xi, \eta}$ belongs to $I_\Phi$, in the sense of (\cite{E3}, 4.1), which implies (\cite{E3}, 4.6) that $(\omega_{\xi, \eta}*id)(W)$ belongs to $\gN_{\hat{\Phi}}$, we have :
 \newline
(i) let $\mathcal P_\eta$ be the element of the positive extension of $M'$ defined in (\cite{E3}, 4.1); then, $R^{\beta, \nu^o}(\xi)^*\mathcal P_\eta R^{\beta, \nu^o}(\xi)$ belongs to the positive extension of $N$, and we have :
\[T_\ga[((\omega_{\xi, \eta}*id)(\widetilde{W}))^*((\omega_{\xi, \eta}*id)(\widetilde{W}))]=\alpha(R^{\beta, \nu^o}(\xi)^*\mathcal P_\eta R^{\beta, \nu^o}(\xi))\]
\newline
(ii) $A_\Omega^\ga=\alpha(N)$
\newline
(iii) let us write $\psi_1=\nu\circ \alpha^{-1}\circ T_\ga$; then $\psi_1$ is a normal semi-finite faithful weight on $A_\Omega$, $\hat{\delta}$-invariant with respect to $\ga$, bearing the density condition, and we have :
\[\psi_1[((\omega_{\xi, \eta}*id)(\widetilde{W}))^*((\omega_{\xi, \eta}*id)(\widetilde{W}))]=\hat{\Phi}[(\omega_{\xi, \eta}*id)(W)^*(\omega_{\xi, \eta}*id)(W)]\]
For all $n\in N$ and $t\in\mathbb{R}$, we have $\sigma_t^{\psi_1}(\alpha(n))=\alpha(\sigma_t^\nu(n))$ and $\sigma_t^{\psi_1}(\hat{\beta}(n))=\hat{\beta}(\gamma_{-t}(n))$. 
 (iv) there exists a unitary $u$ from $H$ onto $H_{\psi_1}$ such that :
 \[u\Lambda_{\hat{\Phi}}((\omega_{\xi, \eta}*id)(W))=\Lambda_{\psi_1}((\omega_{\xi, \eta}*id)(\widetilde{W}))\]
 and we have, for all $n\in N$ :
 \[u\alpha(n)=\pi_{\psi_1}(\alpha(n))u\]
 \[u\hat{\beta}(n)=\pi_{\psi_1}(\hat{\beta}(n))u\]
 \[u\beta(n)=J_{\psi_1}\pi_{\psi_1}(\alpha(n^*))J_{\psi_1}u\]
 \[u\hat{\alpha}(n)=J_{\psi_1}\pi_{\psi_1}(\hat{\beta}(n^*))J_{\psi_1}u\]
(v) the normal faithful semi-finite weight $\nu\circ\alpha^{-1}$ on $\alpha(N)=A_\Omega^\ga$ bears the Galois density condition defined in \ref{defdens} for the Galois action $(\hat{\beta}, \ga)$ of $\widehat{\gG}$ on $A_\Omega$. 
\newline
(vi) the operator $(u\underset{N}{_{\hat{\beta}}\otimes_\alpha}1)W^c(u^*\underset{N^o}{_{\hat{\alpha}}\otimes_{\hat{\beta}}}1)$ is the standard implementation $V_{\psi_1}$ of the action $(\hat{\beta}, \ga)$ associated to the weight $\psi_1$ on $A_\Omega$. 
\newline
(vii) for any $x\in A_\Omega$, we have $\pi_{\psi_1}(x)=uxu^*$. 
\newline
(viii) if $N$ is a finite sum of factors, then there exists a normal semi-finite faithful weight $\phi$ on $A_\Omega$ such that $(A_\Omega, \hat{\beta}, \ga, \phi, \nu)$ is a Galois system.}
\begin{proof}
Using the calculations made in \ref{action}, with an orthogonal $(\alpha, \nu)$-basis $(f_j)_{j\in J}$ of $H$, we get that $\ga[((\omega_{\xi, \eta}*id)(\widetilde{W}))^*((\omega_{\xi, \eta}*id)(\widetilde{W}))]$ is equal to :
\begin{multline*}
\sum_{j,j'}((\omega_{\xi, f_{j'}}*id)(\widetilde{W})^*\underset{N}{_{\hat{\beta}}\otimes_\alpha}1)(\omega_{\eta, f_{j'}}*id*id)[\sigma^{2,3}_{\hat{\beta}, \alpha}(W^*\underset{N}{_\alpha\otimes_{\hat{\beta}}}1)(1\underset{N^o}{_\alpha\otimes_{\hat{\beta}}}\sigma_{\nu^o})]...\\
...(\omega_{f_j, \eta}*id*id)([\sigma^{2,3}_{\hat{\beta}, \alpha}(W^*\underset{N}{_\alpha\otimes_{\hat{\beta}}}1)(1\underset{N^o}{_\alpha\otimes_{\hat{\beta}}}\sigma_{\nu^o})]^*)((\omega_{\xi, f_j}*id)(\widetilde{W})\underset{N}{_{\hat{\beta}}\otimes_\alpha}1)
\end{multline*}
and, therefore, $T_\ga[((\omega_{\xi, \eta}*id)(\widetilde{W}))^*((\omega_{\xi, \eta}*id)(\widetilde{W}))]$ is equal to :
\begin{multline*}
\sum_{j,j'}((\omega_{\xi, f_{j'}}*id)(\widetilde{W})^*)...
\\(id\underset{N}{_{\hat{\beta}}*_\alpha}\hat{\Phi})[(\omega_{\eta, f_{j'}}*id*id)[\sigma^{2,3}_{\hat{\beta}, \alpha}(W^*\underset{N}{_\alpha\otimes_{\hat{\beta}}}1)(1\underset{N^o}{_\alpha\otimes_{\hat{\beta}}}\sigma_{\nu^o})]
(\omega_{f_j, \eta}*id*id)([\sigma^{2,3}_{\hat{\beta}, \alpha}(W^*\underset{N}{_\alpha\otimes_{\hat{\beta}}}1)(1\underset{N^o}{_\alpha\otimes_{\hat{\beta}}}\sigma_{\nu^o})]^*)]\\
...((\omega_{\xi, f_j}*id)(\widetilde{W})\underset{N}{_{\hat{\beta}}\otimes_\alpha}1)
\end{multline*}
Let $\delta$ be in $D(H_{\hat{\beta}}, \nu^o)$, and let $(\delta_i)_{i\in I}$ be an orthogonal $(\hat{\beta}, \nu^o)$-basis of $H$, we get that :
\[(\omega_\delta\underset{N}{_{\hat{\beta}}*_\alpha}\hat{\Phi})[(\omega_{\eta, f_{j'}}*id*id)[\sigma^{2,3}_{\hat{\beta}, \alpha}(W^*\underset{N}{_\alpha\otimes_{\hat{\beta}}}1)(1\underset{N^o}{_\alpha\otimes_{\hat{\beta}}}\sigma_{\nu^o})]
(\omega_{f_j, \eta}*id*id)([\sigma^{2,3}_{\hat{\beta}, \alpha}(W^*\underset{N}{_\alpha\otimes_{\hat{\beta}}}1)(1\underset{N^o}{_\alpha\otimes_{\hat{\beta}}}\sigma_{\nu^o})]^*)]\]
is equal to :
\begin{multline*}
\sum_i\hat{\Phi}(\omega_{\eta, f_{j'}}*\omega_{\delta_i, \delta}*id)[\sigma^{2,3}_{\hat{\beta}, \alpha}(W^*\underset{N}{_\alpha\otimes_{\hat{\beta}}}1)(1\underset{N^o}{_\alpha\otimes_{\hat{\beta}}}\sigma_{\nu^o})]...\\
...(\omega_{f_j, \eta}*\omega_{\delta, \delta_i}*id)([\sigma^{2,3}_{\hat{\beta}, \alpha}(W^*\underset{N}{_\alpha\otimes_{\hat{\beta}}}1)(1\underset{N^o}{_\alpha\otimes_{\hat{\beta}}}\sigma_{\nu^o})]^*)
\end{multline*}
which, using the calculations made in \ref{action}, is equal to :
\begin{multline*}
\sum_i\hat{\Phi}[(\omega_{\eta, f_{j'}}*id)(W^*)\alpha(<\delta_i, \delta>_{\hat{\beta}, \nu^o})\alpha(<\delta, \delta_i>_{\hat{\beta}, \nu^o})(\omega_{f_j, \eta}*id)(W)]\\
=\hat{\Phi}[(\omega_{\eta, f_{j'}}*id)(W^*)\alpha(<\delta, \delta>_{\hat{\beta}, \nu^o})(\omega_{f_j, \eta}*id)(W)]
\end{multline*}
Therefore, if $f_j$ and $f_{j'}$ are in $\mathcal D(\pi'(\eta))$, we finally get that it is equal to :
\begin{eqnarray*}
(\alpha(<\delta, \delta>_{\hat{\beta}, \nu^o})\Lambda_{\hat{\Phi}}[(\omega_{f_j, \eta}*id)(W)]|
\Lambda_{\hat{\Phi}}[(\omega_{f_{j'}, \eta}*id)(W)])
&=&
(\alpha(<\delta, \delta>_{\hat{\beta}, \nu^o})\pi'(\eta)^*f_j|\pi'(\eta)f_{j'}^*)\\
&=&
(\delta\underset{\nu}{_{\hat{\beta}}\otimes_\alpha}\pi'(\eta)^*f_j|\delta\underset{\nu}{_{\hat{\beta}}\otimes_\alpha}\pi'(\eta)^*f_{j'})\\
&=&
(\pi'(\eta)^*f_j\underset{\nu^o}{_\alpha\otimes_{\hat{\beta}}}\delta|\pi'(\eta)^*f_{j'}\underset{\nu^o}{_\alpha\otimes_{\hat{\beta}}}\delta)
\end{eqnarray*}
Let now $\delta'_k\in D(H_{\hat{\beta}}, \nu^o)$ be an orthogonal basis of $H$; then :
\[<T_\ga[((\omega_{\xi, \eta}*id)(\widetilde{W}))^*((\omega_{\xi, \eta}*id)(\widetilde{W}))]), \omega_\delta>\]
is equal to :
\[\sum_{k,j,j'}(\tilde{W}^*(f'_j\underset{\nu^o}{_\alpha\otimes_{\hat{\beta}}}\delta'_k))|\xi\underset{\nu}{_\beta\otimes_\alpha}\delta)(\pi'(\eta)^*f_j\underset{\nu^o}{_\alpha\otimes_{\hat{\beta}}}\delta'_k|\pi'(\eta)^*f_{j'}\underset{\nu^o}{_\alpha\otimes_{\hat{\beta}}}\delta'_k)(\tilde{W}(\xi\underset{\nu}{_\beta\otimes_\alpha}\delta)|f_{j'}\underset{\nu^o}{_\alpha\otimes_{\hat{\beta}}}\delta'_k)\]
As the family $(f_j\underset{N^o}{_\alpha\otimes_{\hat{\beta}}\delta'_k)_{(k,j)}}$ is an orthogonal basis of $H\underset{\nu^o}{_\alpha\otimes_{\hat{\beta}}}H$, we can use the Plancherel formula, which gives that :
\begin{eqnarray*}
<T_\ga[((\omega_{\xi, \eta}*id)(\widetilde{W}))^*((\omega_{\xi, \eta}*id)(\widetilde{W}))]), \omega_\delta>
&=&
(\tilde{W}(\xi\underset{\nu}{_\beta\otimes_\alpha}\delta)|\tilde{W}(\mathcal P_\eta\xi\underset{\nu}{_\beta\otimes_\alpha}\delta))\\
&=&(\xi\underset{\nu}{_\beta\otimes_\alpha}\delta|(\mathcal P_\eta\xi\underset{\nu}{_\beta\otimes_\alpha}\delta)\\
&=&\|\pi'(\eta)^*\xi \underset{\nu}{_\beta\otimes_\alpha}\delta\|^2
\end{eqnarray*}
from which one gets (i). 
\newline
We have seen in \ref{action} that $\alpha(N)\subset A_\Omega^\ga$; on the other side, if $x\in\gM_{T_\ga}^+$, using \ref{algcocycle}, one gets that $x$ is the upper limit of an increasing positive sum of elements of the form $(\omega_{\xi, \eta}*id)(\widetilde{W}))^*((\omega_{\xi, \eta}*id)(\widetilde{W}))$; therefore, $T_\ga(x)$ is, by (i), the upper limit of an increasing sequence of elements in $\alpha(N)$, and therefore, we get that $T_\ga(x)\in \alpha(N)$; as $T_\ga(\gM_{T_\ga})$ is dense in $A_\Omega$, we get (ii). 
Thanks to (ii), one can define the lifted weight $\psi_1=\nu\circ\alpha^{-1}\circ T_\ga$, which is $\delta$-invariant with respect to $\ga$ by \ref{thintegrable}. Moreover, using (i), one gets that :
\[\psi_1[(\omega_{\xi, \eta}*id)(\widetilde{W}))^*((\omega_{\xi, \eta}*id)(\widetilde{W}))]=\hat{\Phi}[(\omega_{\xi, \eta}*id)(W))^*(\omega_{\xi, \eta}*id)(W))]\]
which gives the first formula of (iii); the formula $\sigma_t^{\psi_1}(\alpha(n))=\alpha(\sigma_t^\nu(n))$ is clear by definition of $\psi_1$; as $\psi_1$ is $\delta$-invariant, using (\cite{E5}, 8,8) and \ref{action}, we get that $\sigma_t^{\psi_1}(\hat{\beta}(n))=\sigma_t^{\hat{\Phi}}(\hat{\beta}(n))=\hat{\beta}(\gamma_{-t}(n))$, which finishes the proof of (iii). 
\newline
Using (iii), we get the existence of an isometry $u$ from $H$ into $H_{\psi_1}$ such that :
\[u\Lambda_{\hat{\Phi}}((\omega_{\xi, \eta}*id)(W))=\Lambda_{\psi_1}((\omega_{\xi, \eta}*id)(\widetilde{W}))\]
Let us write $P$ for the projection on $Im u$; using \ref{algcocycle}, we get that $\pi_{\psi_1}(A_\Omega)P=P\pi_{\psi_1}(A_\Omega)P$, and, therefore, that $P\in \pi_{\psi_1}(A_\Omega)'=J_{\psi_1}\pi_{\psi_1}(A_\Omega)J_{\psi_1}$. Using Kaplansky's theorem, one can find a family $\omega_n$ in $I_\Phi$, such that $\|(\omega_n*id)(\widetilde{W})\|\leq 1$ and $J_{\psi_1}\pi_{\psi_1}[(\omega_n*id)(\widetilde{W})]J_{\psi_1}$ is weakly converging to $1-P$; then, we get that :
\[\pi_{\psi_1}[(\omega_{\xi, \eta}*id)(\widetilde{W})]J_{\psi_1}\Lambda_{\psi_1}[(\omega_n*id)(\widetilde{W})]=
J_{\psi_1}\pi_{\psi_1}[(\omega_n*id)(\widetilde{W})]J_{\psi_1}\Lambda_{\psi_1}[(\omega_{\xi, \eta}*id)(\widetilde{W})]\]
is converging to $0$, because $J_{\psi_1}\pi_{\psi_1}[(\omega_n*id)(\widetilde{W})]J_{\psi_1}$ is weakly converging to $1-P$; using now the weak density of the linear combinations of elements of the form 
$\pi_{\psi_1}[(\omega_{\xi, \eta}*id)(\widetilde{W})]$ in $\pi_{\psi_1}(A_\Omega)$, we get that 
$\Lambda_{\psi_1}[(\omega_n*id)(\widetilde{W})]$ is converging to $0$; from which one gets that $\psi_1(J_{\psi_1}(1-P)J_{\psi_1})=0$ and that $P=1$, which proves that $u$ is a unitary. 
\newline
Moreover, we have :
\begin{eqnarray*}
u\alpha(n)\Lambda_{\hat{\Phi}} ((\omega_{\xi, \eta}*id)(W))
&=&
u\Lambda_{\hat{\Phi}} [(\omega_{\xi, \eta}*id)((1\underset{N^o}{_\alpha\otimes_{\hat{\beta}}}\alpha(n))W)]\\
&=&
u\Lambda_{\hat{\Phi}} [(\omega_{\xi, \eta}*id)(W(\alpha(n)\underset{N}{_\beta\otimes_\alpha}1))]\\
&=&
u\Lambda_{\hat{\Phi}} ((\omega_{\alpha(n)\xi, \eta}*id)(W))\\
&=&
\Lambda_{\psi_1}((\omega_{\alpha(n)\xi, \eta}*id)(\widetilde{W}))\\
&=&
\Lambda_{\psi_1}[(\omega_{\xi, \eta}*id)(\widetilde{W}(\alpha(n)\underset{N}{_\beta\otimes_\alpha}1))]\\
&=&
\Lambda_{\psi_1}[(\omega_{\xi, \eta}*id)((1\underset{N^o}{_\alpha\otimes_{\hat{\beta}}}\alpha(n))\widetilde{W})]\\
&=&
\pi_{\psi_1}(\alpha(n))\Lambda_{\psi_1}((\omega_{\xi, \eta}*id)(\widetilde{W}))
\end{eqnarray*}
Let us suppose now that $n\in N$ is analytic with respect to $\nu$ and let's use (ii); then, we have :
\begin{eqnarray*}
u\hat{\beta}(n)\Lambda_{\hat{\Phi}} ((\omega_{\xi, \eta}*id)(W))
&=&
u\Lambda_{\hat{\Phi}} [(\omega_{\xi, \eta}*id)((1\underset{N^o}{_\alpha\otimes_{\hat{\beta}}}\hat{\beta}(n))W)]\\
&=&
u\Lambda_{\hat{\Phi}} [(\omega_{\xi, \eta}*id)((\alpha(\sigma^\nu_{i/2}(n))\underset{N^o}{_\alpha\otimes_{\hat{\beta}}}1)W)]\\
&=&
u\Lambda_{\hat{\Phi}}[(\omega_{\xi, \alpha(\sigma_{-i/2}(n^*))\eta}*id)(W)]\\
&=&
\Lambda_{\psi_1}[(\omega_{\xi, \alpha(\sigma_{-i/2}(n^*))\eta}*id)(\widetilde{W})]\\
&=&
u\Lambda_{\hat{\Phi}} [(\omega_{\xi, \eta}*id)((\alpha(\sigma^\nu_{i/2}(n))\underset{N^o}{_\alpha\otimes_{\hat{\beta}}}1)\widetilde{W})]\\
&=&
u\Lambda_{\hat{\Phi}} [(\omega_{\xi, \eta}*id)((1\underset{N^o}{_\alpha\otimes_{\hat{\beta}}}\hat{\beta}(n))\widetilde{W})]\\
&=&
\pi_{\psi_1}(\hat{\beta}(n))u\Lambda_{\hat{\Phi}} ((\omega_{\xi, \eta}*id)(W))
\end{eqnarray*}
and :
\begin{eqnarray*}
u\beta(n)\Lambda_{\hat{\Phi}} ((\omega_{\xi, \eta}*id)(W))
&=&
uJ_{\hat{\Phi}}\alpha(n^*)J_{\hat{\Phi}}\Lambda_{\hat{\Phi}} ((\omega_{\xi, \eta}*id)(W))\\
&=&
u\Lambda_{\hat{\Phi}} [(\omega_{\xi, \eta}*id)(W)(1\underset{N}{_\beta\otimes_\alpha}\alpha(\sigma_{-i/2}(n))]\\
&=&
u\Lambda_{\hat{\Phi}} [(\omega_{\xi, \eta}*id)(W)(\beta(n)\underset{N}{_\beta\otimes_\alpha}1)]\\
&=&
u\Lambda_{\hat{\Phi}}[\omega_{\beta(n)\xi, \eta}*id)(W)]\\
&=&
\Lambda_{\psi_1}[\omega_{\beta(n)\xi, \eta}*id)(\widetilde{W})]\\
&=&
\Lambda_{\psi_1}[(\omega_{\xi, \eta}*id)(\widetilde{W})(\beta(n)\underset{N}{_\beta\otimes_\alpha}1)]\\
&=&
\Lambda_{\psi_1}[(\omega_{\xi, \eta}*id)(\widetilde{W})(1\underset{N}{_\beta\otimes_\alpha}\alpha(\sigma_{-i/2}(n))]\\
&=&
J_{\psi_1}\alpha(n^*)J_{\psi_1}\Lambda_{\psi_1}((\omega_{\xi, \eta}*id)(\widetilde{W}))\\
&=&
J_{\psi_1}\alpha(n^*)J_{\psi_1}u\Lambda_{\hat{\Phi}} ((\omega_{\xi, \eta}*id)(W))
\end{eqnarray*}
If now we suppose that $n$ is analytic with respect to $\gamma$, and use again (ii), we shall get :
\begin{eqnarray*}
u\hat{\alpha}(n)\Lambda_{\hat{\Phi}}((\omega_{\xi, \eta}*id)(W))
&=&
uJ_{\hat{\Phi}}\hat{\beta}(n^*)J_{\hat{\Phi}}\Lambda_{\hat{\Phi}} ((\omega_{\xi, \eta}*id)(W))\\
&=&
u\Lambda_{\hat{\Phi}} [(\omega_{\xi, \eta}*id)(W)(1\underset{N}{_\beta\otimes_\alpha}\hat{\beta}(\gamma_{i/2}(n))]\\
&=&
u\Lambda_{\hat{\Phi}} [(\omega_{\xi, \eta}*id)((\beta(\gamma_{i/2}(n)\underset{N}{_\alpha\otimes_{\hat{\beta}}}1)W))]\\
&=&
u\Lambda_{\hat{\Phi}} [(\omega_{\beta(\gamma{i/2}(n))^*\xi, \eta}*id)(W))]\\
&=&
\Lambda_{\psi_1}[(\omega_{\beta(\gamma{i/2}(n))^*\xi, \eta}*id)(\widetilde{W}))]\\
&=&
\Lambda_{\psi_1} [(\omega_{\xi, \eta}*id)((\beta(\gamma_{i/2}(n)\underset{N}{_\alpha\otimes_{\hat{\beta}}}1)\widetilde{W}))]\\
&=&
\Lambda_{\psi_1}  [(\omega_{\xi, \eta}*id)(\widetilde{W})(1\underset{N}{_\beta\otimes_\alpha}\hat{\beta}(\gamma_{i/2}(n))]\\
&=&
J_{\psi_1}\hat{\beta}(n^*)J_{\psi_1}\Lambda_{\psi_1} [(\omega_{\xi, \eta}*id)(\widetilde{W})]\\
&=&
J_{\psi_1}\hat{\beta}(n^*)J_{\psi_1}u\Lambda_{\hat{\Phi}}((\omega_{\xi, \eta}*id)(W))
\end{eqnarray*}
which, by continuity, finishes the proof of (iv). 
\newline
The weight $\nu\circ\alpha^{-1}$ satisfies the density condition if the subspace :
\[D((H_{\psi_1})_{\pi_{\psi_1}\circ\hat{\beta}}, \nu^o)\cap D(_{\pi_{\psi_1}\circ\alpha}H_{\psi_1}, \nu)\] is dense in $H_{\psi_1}$. Using now (iv), we get that this subspace is the image by $u$ of $D(H_{\hat{\beta}}, \nu^o)\cap D(_\alpha H, \nu)$
which is dense in $H_\Phi$ by (\cite{E4}, 2.3), from which we get (v), using (iv) again. 
\newline
In (\cite{E5}, 8.2), one gets that $\widehat{W^o}=\widehat{W}^c$ is the standard implementation of the action $(\beta, \Gamma)$ of $\gG$ on $M$, associated to the $\delta$-invaraint weight $\Phi$. So, $W^c$ is the standard implementation of the action $(\hat{\beta}, \widehat{\Gamma})$ on $\widehat{M}$, associated to the $\hat{\delta}$-invariant weight $\widehat{\Phi}$. Which means that, for any orthogonal $(\alpha, \nu)$-basis of $H$, any $\zeta$ in $D(_\alpha H, \nu)\cap\mathcal D(\hat{\delta}^{1/2})$ such that $\hat{\delta}^{1/2}$ belongs to $D(H_{\hat{\beta}}, \nu^o)$, any $x$ in $\gN_{\widehat{\Phi}}$, we have :
\[W^c(\Lambda_{\widehat{\Phi}}(x)\underset{\nu^o}{_{\hat{\alpha}}\otimes_{\hat{\beta}}}\hat{\delta}^{1/2}\zeta)=\sum_i\Lambda_{\widehat{\Phi}}[(id\underset{N}{_{\hat{\beta}}*_\alpha}\omega_{\zeta, e_i})\widehat{\Gamma}(x)]\underset{\nu}{_{\hat{\beta}}\otimes_\alpha}e_i\]
and, therefore, in particular :
\begin{multline*}
W^c(\Lambda_{\widehat{\Phi}}(\omega_{\xi, \eta}*id)(W))\underset{\nu^o}{_{\hat{\alpha}}\otimes_{\hat{\beta}}}\hat{\delta}^{1/2}\zeta)=\sum_i\Lambda_{\widehat{\Phi}}[(id\underset{N}{_{\hat{\beta}}*_\alpha}\omega_{\zeta, e_i})\widehat{\Gamma}((\omega_{\xi, \eta}*id)(W))]\underset{\nu}{_{\hat{\beta}}\otimes_\alpha}e_i\\
=\sum_i\Lambda_{\widehat{\Phi}}(\omega_{\xi, \eta}*id*\omega_{\zeta, e_i})[(1\underset{N^o}{_\alpha\otimes_{\hat{\beta}}}\sigma_{\nu^o})(W\underset{N^o}{_\alpha\otimes_{\hat{\beta}}}1)\sigma^{2,3}_{\beta, \alpha}(W\underset{N^o}{_{\hat{\beta}}\otimes_\alpha}1)]\underset{\nu}{_{\hat{\beta}}\otimes_\alpha}e_i
\end{multline*}
Using now (iv), we then get that $(u\underset{N}{_{\hat{\beta}}\otimes_\alpha}1)W^c(\Lambda_{\widehat{\Phi}}(\omega_{\xi, \eta}*id)(W))\underset{\nu^o}{_{\hat{\alpha}}\otimes_{\hat{\beta}}}\hat{\delta}^{1/2}\zeta)$ is equal to :
\[\sum_i\Lambda_{\psi_1}(\omega_{\xi, \eta}*id*\omega_{\zeta, e_i})[(1\underset{N^o}{_\alpha\otimes_{\hat{\beta}}}\sigma_{\nu^o})(W\underset{N^o}{_\alpha\otimes_{\hat{\beta}}}1)\sigma^{2,3}_{\beta, \alpha}(\widetilde{W}\underset{N^o}{_{\hat{\beta}}\otimes_\alpha}1)]\underset{\nu}{_{\hat{\beta}}\otimes_\alpha}e_i\]
which, thanks to \ref{action}(i), is equal to :
\[\sum_i\Lambda_{\psi_1}[(id\underset{N}{_{\hat{\beta}}*_\alpha}\omega_{\zeta, e_i})\ga ((\omega_{\xi, \eta}*id)(\widetilde{W}))\underset{\nu}{_{\hat{\beta}}\otimes_\alpha}e_i\]
and, therefore, we have, where we denote $V_{\psi_1}$ the standard implementation of $(\hat{\beta}, \ga)$ associated to the weight $\psi_1$ :
\begin{multline*}
(u\underset{N}{_{\hat{\beta}}\otimes_\alpha}1)W^c(u^*\underset{N^o}{_{\hat{\alpha}}\otimes_{\hat{\beta}}}1)(\Lambda_{\psi_1}(\omega_{\xi, \eta}*id)(\widetilde{W}))\underset{\nu^o}{_{\hat{\alpha}}\otimes_{\hat{\beta}}}\hat{\delta}^{1/2}\zeta)\\
=
\sum_i\Lambda_{\psi_1}[(id\underset{N}{_{\hat{\beta}}*_\alpha}\omega_{\zeta, e_i})\ga ((\omega_{\xi, \eta}*id)(\widetilde{W}))\underset{\nu}{_{\hat{\beta}}\otimes_\alpha}e_i\\
=
V_{\psi_1}(\Lambda_{\psi_1}(\omega_{\xi, \eta}*id)(\widetilde{W}))\underset{\nu^o}{_{\hat{\alpha}}\otimes_{\hat{\beta}}}\hat{\delta}^{1/2}\zeta)
\end{multline*}
from which we get (vi), by density. 
\newline
Let $\xi'\in D(H_\beta, \nu^o)$, $\eta'\in D(_\alpha H, \nu)$, we get, with an orthogonal $(\beta, \nu^o)$-basis $(e_i)_{i\in I}$ of $H$,  that :
\begin{eqnarray*}
(\omega_{\xi', \eta'}*id)(\widetilde{W})\Lambda_{\widehat{\Phi}}[(\omega_{\xi, \eta}*id)(W)]
&=&
\sum_i (\omega_{e_i, \eta'}*id)(W)(\omega_{\xi', e_i}\underset{N}{_\beta*_\alpha}id)(\Omega^*)\pi'(\eta)^*\xi\\
&=&
\sum_i (\omega_{e_i, \eta'}*id)(W)\pi'(\eta)^*(\omega_{\xi', e_i}\underset{N}{_\beta*_\alpha}id)(\Omega^*)\xi\\
&=&
\sum_i (\omega_{e_i, \eta'}*id)(W)\Lambda_{\widehat{\Phi}}(\omega_{(\omega_{\xi', e_i}\underset{N}{_\beta*_\alpha}id)(\Omega^*)\xi, \eta}*id)(W)\\
&=&
\sum_i\Lambda_{\widehat{\Phi}}[(\omega_{e_i, \eta'}*id)(W)(\omega_{(\omega_{\xi', e_i}\underset{N}{_\beta*_\alpha}id)(\Omega^*)\xi, \eta}*id)(W)]
\end{eqnarray*}
Let now $(f_j)_{j\in J}$ be an orthogonal $(\alpha, \nu)$-basis of $H$; we know that there exists $\xi_{i,j}$, $\eta_j$ in $H$ such that :
\[W(e_i\underset{\nu}{_\beta\otimes_\alpha}(\omega_{\xi', e_i}\underset{N}{_\beta*_\alpha}id)(\Omega^*)\xi)=\sum_j f_j\underset{\nu^o}{_\alpha\otimes_{\hat{\beta}}}\xi_{i,j}\]
\[W(\eta'\underset{\nu}{_\beta\otimes_\alpha}\eta)=\sum_j f_j\underset{\nu^o}{_\alpha\otimes_{\hat{\beta}}}\eta_j\]
and then, we get that 
\[(\omega_{\xi', \eta'}*id)(\widetilde{W})\Lambda_{\widehat{\Phi}}[(\omega_{\xi, \eta}*id)(W)]=
\sum_{i,j}\Lambda_{\widehat{\Phi}}[(\omega_{\xi_{i,j}, \eta_j}*id)(W)]\]
which implies that :
\begin{eqnarray*}
\sum_j f_j\underset{\nu^o}{_\alpha\otimes_{\hat{\beta}}}\sum_i\xi_{i,j}
&=&
W\sum_i(
e_i\underset{\nu}{_\beta\otimes_\alpha}(\omega_{\xi', e_i}\underset{N}{_\beta*_\alpha}id)(\Omega^*)\xi)\\
&=& W\Omega^*(\xi'\underset{\nu}{_\beta\otimes_\alpha}\xi)
\end{eqnarray*}
and, finally :
\[(\omega_{\xi', \eta'}*id)(\widetilde{W})\Lambda_{\widehat{\Phi}}[(\omega_{\xi, \eta}*id)(W)]=
\sum_{j}\Lambda_{\widehat{\Phi}}[(\omega_{\xi'_{j}, \eta_j}*id)(W)]\]
where $\widetilde{W}(\xi'\underset{\nu}{_\beta\otimes_\alpha}\xi)=\sum_j f_j\underset{\nu^o}{_\alpha\otimes_{\hat{\beta}}}\xi'_j$.
But, using again the calculation already made in \ref{algcocycle}(iii), we get that :
\begin{eqnarray*}
u(\omega_{\xi', \eta'}*id)(\widetilde{W})\Lambda_{\widehat{\Phi}}[(\omega_{\xi, \eta}*id)(W)
&=&
\sum_j\Lambda_{\psi_1}[(\omega_{\xi'_{j}, \eta_j}*id)(W)]\\
&=&
\Lambda_{\psi_1}[(\omega_{\xi', \eta'}*id)(\widetilde{W})(\omega_{\xi, \eta}*id)(\widetilde{W})]\\
&=&
\pi_{\psi_1}((\omega_{\xi', \eta'}*id)(\widetilde{W}))\Lambda_{\psi_1}[(\omega_{\xi, \eta}*id)(\widetilde{W})]\\
&=&
\pi_{\psi_1}((\omega_{\xi', \eta'}*id)(\widetilde{W}))u\Lambda_{\widehat{\Phi}}[(\omega_{\xi, \eta}*id)(W)
\end{eqnarray*}
from which, by density, we get (vii). 
\newline
Moreover, (viii) is a direct application of \ref{corsum} to (ii). Which finishes the proof. \end{proof}

\subsection{Theorem}
\label{thG}
 {\it Let $\gG$ be a measured quantum groupoid, $\Omega$ a $2$-cocycle for $\gG$; let $W$ be the pseudo-multiplicative unitary associated to $\gG$, $A_\Omega$ the von Neumann algebra on $H$ defined in \ref{algcocycle} and $(\hat{\beta}, \ga)$ the Galois action of $\widehat{\gG}$ on $A_\Omega$ defined in \ref{action} whose invariant subalgebra $A_\Omega^\ga$ is equal to $\alpha(N)$ (\ref{invariant}(ii)); let us write $\widetilde{W}=W\Omega^*$; moreover, the weight $\nu\circ\alpha^{-1}$ on $\alpha(N)$ bears the Galois density property defined in \ref{defdens}, by \ref{invariant}(vi). Let us write $\psi_1=\nu\circ\alpha^{-1}\circ T_\ga$. Let $u$ be the unitary from $H$ onto $H_{\psi_1}$ introduced in \ref{invariant}(iv). 
 \newline
 The canonical representation $r$ of $A_\Omega^\ga$ on $H_{\psi_1}$ is the restriction of $\pi_{\psi_1}$ to $\alpha(N)$; using \ref{invariant}(iv), we get that the canonical antirepresentation $s$ of $\alpha(N)$ (identified to $N$ for simplification) on $H_{\psi_1}$ is $s(n)=u\beta(n)u^*$ ($n\in N$); for simplification again, we shall write $\alpha$ for $\pi_{\psi_1}\circ\alpha
$ and $\hat{\beta}$ for $\pi_{\psi_1}\circ\hat{\beta}$. Then the Galois unitary $\tilde{G}$ is a unitary from $H_{\psi_1}\underset{\nu}{_s\otimes_\alpha}H_{\psi_1}$ onto $H\underset{\nu^o}{_\alpha\otimes_{\hat{\beta}}}H_{\psi_1}$; then, we have :
\[(1\underset{N^o}{_\alpha\otimes_{\hat{\beta}}}u^*)\widetilde{G}(u\underset{N}{_\beta\otimes_\alpha}u)=\widetilde{W}\]
 }
\begin{proof}
Let $\xi\in D(_\alpha H, \nu)$, $\eta \in D(H_{\hat{\beta}}, \nu^o)$; 
let $\xi'\in D(H_\beta, \nu^o)$ and $\eta'\in D(_\alpha H, \nu)$ such that $\omega_{\xi', \eta'}$ belongs to $I_\Phi$ (in the sense of \cite{E3} 4.1), which implies that $(\omega_{\xi', \eta'}*id)(W)$ belongs to $\gN_{\widehat{\Phi}}$, and, using \ref{invariant}(iii), that $(\omega_{\xi', \eta'}*id)(\widetilde{W})$ belongs to $\gN_{\psi_1}$. We have, then, using \ref{invariant}(iv), and \ref{tildeG}(i) :
\begin{eqnarray*}
(id*\omega_{\xi, \eta}\circ\pi_{\psi_1})(\tilde{G})u\Lambda_{\widehat{\Phi}}[(\omega_{\xi', \eta'}*id)(W)]
&=&
(id*\omega_{\xi, \eta}\circ\pi_{\psi_1})(\tilde{G})\Lambda_{\psi_1}[(\omega_{\xi', \eta'}*id)(\tilde{W})]\\
&=&
\Lambda_{\widehat{\Phi}}[(\omega_{\xi, \eta}\underset{N}{_{\hat{\beta}}*_\alpha}id)\ga ((\omega_{\xi', \eta'}*id)(\widetilde{W}))]
\end{eqnarray*}
Using now \ref{action}(i), we get that 
\[(\omega_{\xi, \eta}\underset{N}{_{\hat{\beta}}*_\alpha}id)\ga ((\omega_{\xi', \eta'}*id)(\widetilde{W}))=(\omega_{\xi', \eta'}*\omega_{\xi, \eta}*id)[(1\underset{N^o}{_\alpha\otimes_{\hat{\beta}}}\sigma_{\nu^o})(W\underset{N^o}{_\alpha\otimes_{\hat{\beta}}}1)\sigma^{2,3}_{\beta, \alpha}(\widetilde{W}\underset{N^o}{_{\hat{\beta}}\otimes_\alpha}1)]\]
Let now $(\xi_i)_{i\in I}$ be an orthogonal $(\alpha, \nu)$-basis of $H$; we get then that this last expression is equal to :
\[\sum_i(\omega_{\xi', \eta'}*\omega_{\xi, \xi_i}*\omega_{\xi_i, \eta}*id)(W_{1,4}W_{1,3}\Omega^*_{1,2})\]
where we use the leg numbering notation, for simplification. But we get then that it is equal to :
\[\sum_i(\omega_{(id\underset{N}{_\beta *_\alpha}\omega_{\xi, \xi_i})(\Omega^*)\xi', \eta'}*\omega_{\xi_i, \eta}*id)[(1\underset{N^o}{_\alpha\otimes_{\hat{\beta}}}\sigma_{\nu^o})(W\underset{N^o}{_\alpha\otimes_{\hat{\beta}}}1)\sigma^{2,3}_{\beta, \alpha}(W\underset{N^o}{_{\hat{\beta}}\otimes_\alpha}1)]\]
which is :
\[\sum_i(\omega_{\xi_i, \eta}\underset{N}{_{\hat{\beta}}*_\alpha}id)\hat{\Gamma}(\omega_{(id\underset{N}{_\beta *_\alpha}\omega_{\xi, \xi_i})(\Omega^*)\xi', \eta'}*id)(W))\]
For any $i\in I$, the operator $(\omega_{\xi_i, \eta}\underset{N}{_{\hat{\beta}}*_\alpha}id)\hat{\Gamma}(\omega_{(id\underset{N}{_\beta *_\alpha}\omega_{\xi, \xi_i})(\Omega^*)\xi', \eta'}*id)(W)$ belongs to $\gN_{\widehat{\Phi}}$, and, by \cite{E3}, 3.10 (ii) applied to $\widehat{\gG}$, then, \cite{E3}4.6 and 4.1, we get that :
\begin{multline*}
\Lambda_{\widehat{\Phi}}[(\omega_{\xi_i, \eta}\underset{N}{_{\hat{\beta}}*_\alpha}id)\hat{\Gamma}(\omega_{(id\underset{N}{_\beta *_\alpha}\omega_{\xi, \xi_i})(\Omega^*)\xi', \eta'}*id)(W)]\\
=
(\omega_{\xi_i, \eta}*id)(\widehat{W}^*)(id\underset{N}{_\beta*_\alpha}\omega_{\xi, \xi_i})(\Omega^*)\Lambda_{\widehat{\Phi}}(\omega_{\xi', \eta'}*id)(W))\\
=(id*\omega_{\xi_i, \eta})(W)(id\underset{N}{_\beta*_\alpha}\omega_{\xi, \xi_i})(\Omega^*)\Lambda_{\widehat{\Phi}}(\omega_{\xi', \eta'}*id)(W)
\end{multline*}
whose sum is weakly converging to $(id*\omega_{\xi, \eta})(\widetilde{W})\Lambda_{\widehat{\Phi}}[(\omega_{\xi', \eta'}*id)(W)]$. As the application $\Lambda_{\widehat{\Phi}}$ is closed, we get that :
\[\Lambda_{\widehat{\Phi}}[(\omega_{\xi, \eta}\underset{N}{_{\hat{\beta}}*_\alpha}id)\ga ((\omega_{\xi', \eta'}*id)(\widetilde{W}))]=(id*\omega_{\xi, \eta})(\widetilde{W})\Lambda_{\widehat{\Phi}}(\omega_{\xi', \eta'}*id)(W)\]
from which we deduce that :
\[(id*\omega_{\xi, \eta}\circ\pi_{\psi_1})(\tilde{G})=(id*\omega_{\xi, \eta})(\widetilde{W})\]
which gives the result, thanks to \ref{invariant}(viii). 
 \end{proof}

 \subsection{Corollaries}
 \label{corG}
 {\it Let $\gG$ be a measured quantum groupoid, $\Omega$ a $2$-cocycle for $\gG$; let $W$ be the pseudo-multiplicative unitary associated to $\gG$, $A_\Omega$ the von Neumann algebra on $H$ defined in \ref{algcocycle} and $(\hat{\beta}, \ga)$ the Galois action of $\widehat{\gG}$ on $A_\Omega$ defined in \ref{action} whose invariant subalgebra $A_\Omega^\ga$ is equal to $\alpha(N)$ (\ref{invariant}(ii)); let us write $\widetilde{W}=W\Omega^*$; moreover, the weight $\nu\circ\alpha^{-1}$ on $\alpha(N)$ bears the Galois density property defined in \ref{defdens}, by \ref{invariant}(vi). Let us write $\psi_1=\nu\circ\alpha^{-1}\circ T_\ga$. Let $u$ be the unitary from $H$ onto $H_{\psi_1}$ introduced in \ref{invariant}(iii). 
  \newline
 The canonical representation $r$ of $A_\Omega^\ga$ on $H_{\psi_1}$ is the restriction of $\pi_{\psi_1}$ to $\alpha(N)$; using \ref{invariant}(iv), we get that the canonical antirepresentation $s$ of $\alpha(N)$ (identified to $N$ for simplification) on $H_{\psi_1}$ is $s(n)=u\beta(n)u^*$ ($n\in N$); let $\rho_t$ be the one-parameter group of automorphisms of $s(N)'$ and $K^{it}$ its standard implementation defined in \ref{propK}; for simplification again, we shall write $\alpha$ for $\pi_{\psi_1}\circ\alpha
$ and $\hat{\beta}$ for $\pi_{\psi_1}\circ\hat{\beta}$. 
Then :
\newline
(i) for any $x\in A_\Omega$, we have :
\[\ga(x)=\sigma_{\nu^o} \widetilde{W}\sigma_{\nu^o}(1\underset{N}{_\alpha\otimes_\beta}x)\sigma_\nu \widetilde{W}^* \sigma_\nu\]
(ii) for any $y\in M'$, we have $\pi_\ga(1\underset{N}{_{\hat{\beta}}\otimes_\alpha}y)=uyu^*$.
\newline
(iii) for all $t\in\mathbb{R}$, we have :
\[K^{it}=(u\underset{N}{_\beta\otimes_\alpha}u)\Omega(\hat{J}\hat{\delta}^{it}\hat{J}\underset{N}{_\beta\otimes_\alpha}\hat{\delta}^{it})\Omega^*(u^*\underset{N}{_\beta\otimes_\alpha}u^*)\]
(iv) for all $t\in\mathbb{R}$, we have $P_{A_\Omega}^{it}=\Delta_{\psi_1}^{it}uJ\delta^{it}Ju^*$.
\newline
(v) we have, for all $t\in\mathbb{R}$ :
\[\widetilde{W}(u^*P_{A_\Omega}^{it}u\underset{N}{_\beta\otimes_\alpha}uP_{A_\Omega}^{it}u^*)=(P^{it}\underset{N^o}{_\alpha\otimes_{\hat{\beta}}}uP_{A_\Omega}^{it}u^*)\widetilde{W}\]
(vi) for any $\xi\in D(H_\beta, \nu^o)$ and $\eta\in D(_\alpha H, \nu)$, we have :
\[\tau_t^{A_\Omega}[(\omega_{\xi, \eta}*id)(\widetilde{W})]=(\omega_{u^*P_{A_\Omega}^{it}u\xi, P^{-it}\eta}*id)(\widetilde{W})=(\omega_{u^*\Delta_{\psi_1}^{it}u\xi, \Delta_{\widehat{\Phi}}^{-it}\eta}*id)(\widetilde{W})\]
(vii) for any $x\in \gN_{\psi_1}\cap\gN_{\psi_1}^*$, $y$, $z$ in $\gN_{\widehat{\Phi}}\cap\gN_{\hat{T}}$, we have :
\[(\omega_{u^*\Lambda_{\psi_1}(x), \hat{J}\Lambda_{\widehat{\Phi}}(y^*z)}*id)(\widetilde{W})^*=
(\omega_{u^*\Lambda_{\psi_1}(x^*), \hat{J}\Lambda_{\widehat{\Phi}}(z^*y)}*id)(\widetilde{W})\]
(viii) for any $\zeta_1$, $\zeta_2$ in $D(_\alpha H, \nu)\cap D(H_{\hat{\beta}}, \nu^o)$, and $\xi\in \mathcal D(\Delta_{\psi_1}^{1/2})$, $\eta\in \mathcal D(\Delta_{\widehat{\Phi}}^{-1/2})$, we have :
\[((id*\omega_{\zeta_1, \zeta_2})(\widetilde{W})u^*\xi|\eta)=
((id*\omega_{\zeta_2, \zeta_1})(\widetilde{W})^*\hat{J}\Delta_{\widehat{\Phi}}^{-1/2}\eta|u^*J_{\psi_1}\Delta_{\psi_1}^{1/2}\xi)\]
(ix) for any $\zeta_1$, $\zeta_2$ in $D(_\alpha H, \nu)\cap D(H_{\hat{\beta}}, \nu^o)$, the operator $\Delta_{\widehat{\Phi}}^{1/2}(id*\omega_{\zeta_2, \zeta_1})(\widetilde{W})u^*\Delta_{\psi_1}^{-1/2}u$ is bounded, and we have :
\[(id*\omega_{\zeta_1, \zeta_2})(\widetilde{W})=\hat{J}\Delta_{\widehat{\Phi}}^{1/2}(id*\omega_{\zeta_2, \zeta_1})(\widetilde{W})u^*\Delta_{\psi_1}^{-1/2}uu^*J_{\psi_1}u\]
(x) for any $\xi\in D(H_\beta, \nu^o)\cap \mathcal D(u\Delta_{\psi_1}^{1/2}u^*)$, and $\eta\in D(_\alpha H, \nu)\cap \mathcal D(\Delta_{\widehat{\Phi}}^{-1/2})$, we have :
\[(\omega_{\xi, \eta}*id)(\widetilde{W})^*=(\omega_{u^*J_{\psi_1}\Delta_{\psi_1}^{1/2}u\xi, \hat{J}\Delta_{\widehat{\Phi}}^{-1/2}\eta}*id)(\widetilde{W})\]
(xi) for all $t\in\mathbb{R}$, we have :
\[\widetilde{W}(u^*\Delta_{\psi_1}^{it}u\underset{N}{_\beta\otimes_\alpha}u^*\Delta_{\psi_1}^{it}u)=
[(\delta\Delta_{\widehat{\Phi}})^{it}\underset{N^o}{_\alpha\otimes_{\hat{\beta}}}u^*\Delta_{\psi_1}^{it}u]\widetilde{W}\]
(xii) for all $t\in\mathbb{R}$, we have :
\[\sigma_t^{\psi_1}[(\omega_{\xi, \eta}*id)(\widetilde{W})]=(\omega_{u^*\Delta_{\psi_1}^{it}u\xi, (\delta\Delta_{\widehat{\Phi}})^{-it}\eta}*id)(\widetilde{W})\]
(xiii) if $N$ is a finite sum of factors, there exists a normal semi-finite faithful operator weight $T_\Omega$ from $M$ to $\alpha(N)$, (resp. $T'_\Omega$ from $M$ to $\beta(N)$) such that \[\gG_\Omega=(N, M, \alpha, \beta, \Gamma_\Omega, T_\Omega, T'_\Omega, \nu)\] is a measured quantum groupoid.}
\begin{proof}
Result (i) is just the application of \ref{thG} to \ref{G}(iv). 
\newline
Let us apply \ref{G}(v) to the action $(\hat{\beta}, \ga)$ of $\widehat{\gG}$ on $A_\Omega$. We get that :
\begin{eqnarray*}
\pi_\ga(1\underset{N}{_{\hat{\beta}}\otimes_\alpha}y)\underset{N}{_s\otimes_\alpha}1
&=&
\tilde{G}^*(y\underset{N^o}{_\alpha\otimes_{\hat{\beta}}}1)\tilde{G}\\
&=&
(u\underset{N}{_\beta\otimes_\alpha}1)\widetilde{W}^*(y\underset{N^o}{_\alpha\otimes_{\hat{\beta}}}1)\widetilde{W}(u^*\underset{N}{_\beta\otimes_\alpha}1)\\
&=&
uyu^*\underset{N}{_s\otimes_\alpha}1
\end{eqnarray*}
from which we get (ii). 
\newline
Applying now \ref{thG} to \ref{propK}, and successively \cite{E5} 3.11(iii) and \cite{E5} 3.8 (vi) applied to $\widehat{\gG}$ , one gets :
\begin{eqnarray*}
K^{it}
&=&
\widetilde{G}^*(\hat{J}\hat{\delta}^{it}\hat{J}\underset{N^o}{_\alpha\otimes_{\hat{\beta}}}1)\widetilde{G}\\
&=&
(u\underset{N}{_\beta\otimes_\alpha}u)\Omega W^*(\hat{J}\hat{\delta}^{it}\hat{J}\underset{N^o}{_\alpha\otimes_{\hat{\beta}}}1)W\Omega^*(u^*\underset{N}{_\beta\otimes_\alpha}u^*)\\
&=&
(u\underset{N}{_\beta\otimes_\alpha}u)\Omega (\hat{J}\underset{N^o}{_\alpha\otimes_{\hat{\beta}}}J)W(\hat{J}\underset{N^o}{_\alpha\otimes_{\hat{\beta}}}J)(\hat{J}\hat{\delta}^{it}\hat{J}\underset{N^o}{_\alpha\otimes_{\hat{\beta}}}1)(\hat{J}\underset{N}{_\beta\otimes_\alpha}J)W^*(\hat{J}\underset{N}{_\beta\otimes_\alpha}J)\Omega^*(u^*\underset{N}{_\beta\otimes_\alpha}u^*)\\
&=&
(u\underset{N}{_\beta\otimes_\alpha}u)\Omega (\hat{J}\underset{N^o}{_\alpha\otimes_{\hat{\beta}}}J)W(\hat{\delta}^{it}\underset{N}{_\beta\otimes_\alpha}1)W^*(\hat{J}\underset{N}{_\beta\otimes_\alpha}J)\Omega^*(u^*\underset{N}{_\beta\otimes_\alpha}u^*)\\
&=&
(u\underset{N}{_\beta\otimes_\alpha}u)\Omega (\hat{J}\underset{N^o}{_\alpha\otimes_{\hat{\beta}}}J)(\hat{\delta}^{it}\underset{N^o}{_\alpha\otimes_{\hat{\beta}}}\hat{\delta}^{it})(\hat{J}\underset{N}{_\beta\otimes_\alpha}J)\Omega^*(u^*\underset{N}{_\beta\otimes_\alpha}u^*)\\
&=&
(u\underset{N}{_\beta\otimes_\alpha}u)\Omega(\hat{J}\hat{\delta}^{it}\hat{J}\underset{N}{_\beta\otimes_\alpha}\hat{\delta}^{it})\Omega^*(u^*\underset{N}{_\beta\otimes_\alpha}u^*)
\end{eqnarray*}
which is (iii). 
\newline
Applying (ii) to \ref{tau}(ii), one gets (iv). Applying again (i) to \ref{tau}(vi), one gets (v). Then, the first equality of (vi) is a direct corollary of (v), and the second equality is a corollary of (iv) and \cite{E6}3.10(vii). Result (vii) is a direct corollary of \ref{tildeG}(iii) applied to \ref{thG}. Then (ix) is an easy corollary from (viii), and (x) from (ix). Result (xi) is given by \ref{G}(vii) and (\cite{E6} 3.11(ii)), applied to \ref{thG}, and (xii) is a direct corollary of (xi). 
\newline
If $N$ is a finite sum of factors, we can apply \ref{action}(iv) and, therefore, we obtain, by \ref{thP2}, a measured quantum groupoid $\gG_1(\ga)$, whose underlying Hopf-bimodule had been defined in \ref{coproduct}(iii). Using now (ii), we get that the von Neumann algebra is (up to $u$) equal to $M$; using \ref{invariant}(ii), we get that the basis is (up to $\alpha$) equal to $N$, and, by \ref{invariant}(iv), that the imbedding of $N$ into $M$ are $\alpha$ and $\beta$. Using now \ref{thG}, we get that the coproduct is $\Gamma_\Omega$, as defined in \ref{propcocycle}, which finishes the proof. 
\end{proof}

\subsection{Proposition}
\label{v}
{\it Let $\gG$ be a measured quantum groupoid, $\Omega$ a $2$-cocycle for $\gG$; let $W$ be the pseudo-multiplicative unitary associated to $\gG$, $A_\Omega$ the von Neumann algebra on $H$ defined in \ref{algcocycle} and $(\hat{\beta}, \ga)$ the action of $\widehat{\gG}$ on $A_\Omega$ defined in \ref{action} whose invariant subalgebra $A_\Omega^\ga$ is equal to $\alpha(N)$ (\ref{invariant}(ii)); let us write $\widetilde{W}=W\Omega^*$; moreover, the weight $\nu\circ\alpha^{-1}$ on $\alpha(N)$ bears the Galois density property defined in \ref{defdens}, by \ref{invariant}(v). Let us write $\psi_1=\nu\circ\alpha^{-1}\circ T_\ga$. Let $u$ be the unitary from $H$ onto $H_{\psi_1}$ introduced in \ref{invariant}(iv), and let us write, for all $t\in\mathbb{R}$:
\[v_t^\Omega=u^*\Delta_{\psi_1}^{it}u\Delta_{\widehat{\Phi}}^{-it}\]
For all $t\in\mathbb{R}$, let us consider the 2-cocycle $\Omega_t$ introduced in \ref{defcocycle}, the algebra $A_{\Omega_t}$ associated, the action $(\hat{\beta}, \ga_t)$ of $\widehat{\gG}$ on $A_{\Omega_t}$, whose invariant is also equal to $\alpha(N)$. Let us denote $\psi_{1,t}$ the weight $\nu\circ\alpha^{-1}\circ T_{\ga_t}$, and $u_t$ the canonical unitary from $H$ to $H_{\psi_{1,t}}$, which, by \ref{invariant}(vii) applied to $\Omega_t$, implements $\pi_{\psi_{1,t}}$. Let us write $\widetilde{W_t}=W\Omega_t^*$. 
Then :
\newline
(i) $v_t^\Omega$ is a unitary in $M\cap\alpha(N)'\cap\beta(N)'$; moreover, the application $t\mapsto v_t^\Omega$ is a $\tau_t$-cocycle. 
\newline
(ii) we have :
\[\widetilde{W}(v_t^\Omega\underset{N}{_\beta\otimes_\alpha}v_t^\Omega)=(1\underset{N}{_\alpha\otimes_{\hat{\beta}}}v_t^\Omega)\widetilde{W_t}\]
\[\Gamma(v_t^\Omega)\Omega_t^*=\Omega^*(v_t^\Omega\underset{N}{_\beta\otimes_\alpha}v_t^\Omega)\]
(iii) the application $\mathcal I_t : x\mapsto v_t^\Omega x(v_t^\Omega)^*$ is an isomorphism from $A_{\Omega_t}$ to $A_\Omega$, and we have, for all $\xi\in D(H_\beta, \nu^o)$ and $\eta\in D(_\alpha H, \nu)$ :
\[\mathcal I_t[(\omega_{\xi, \eta}*id)(\widetilde{W_t})]=(\omega_{v_t^\Omega\xi, \eta}*id)(\widetilde{W})\]
(iv) moreover, we have $\psi_1\circ \mathcal I_t=\psi_{1,t}$; then $u v_t^\Omega u_t^*$ is the standard implementation of $\mathcal I_t$. 
\newline
(v) we have $v_s^{\Omega_t}=\tau_t(v_s^\Omega)$.
\newline
(vi) if, for all $t\in\mathbb{R}$, we have $\Omega=\Omega_t$, then there exists a positive non singular operator $k_\Omega$ affiliated to $M$, such that $\tau_t(k_\Omega)=k_\Omega$ and $v_t^\Omega=k_\Omega^{it}$.}

\begin{proof}
By definition of $\psi_1$, we get that, for all $t\in\mathbb{R}$ and $n\in N$, we have $\sigma_t^{\psi_1}(\alpha(n))=\alpha(\sigma_t^\nu(n))=\sigma_t^{\widehat{\Phi}}(\alpha(n))$; therefore, we get that $v_t^\Omega\in \alpha(N)'$. 
\newline
We have, using first \cite{E5}3.10 (iv) and 3.8(i), then \ref{invariant}(vii), \cite{E5}3.8(ii) and 3.10(vii) :
\begin{eqnarray*}
u^*\Delta_{\psi_1}^{it}u\beta(n)u^*\Delta_{\psi_1}^{-it}u
&=&
u^*\Delta_{\psi_1}^{it}u\hat{J}\alpha(n^*)\hat{J}u^*\Delta_{\psi_1}^{-it}u\\
&=&
u^*\Delta_{\psi_1}^{it}uu^*J_{\psi_1}uX_\Omega\alpha(n^*)X_\Omega^*u^*J_{\psi_1}uu^*\Delta_{\psi_1}^{-it}u\\
&=&
u^*J_{\psi_1}uu^*\Delta_{\psi_1}^{it}\alpha(n^*)u^*\Delta_{\psi_1}^{-it}uu^*J_{\psi_1}u\\
&=&
X_\Omega \hat{J}\alpha(\sigma_t^\nu(n^*))\hat{J}X_\Omega^*\\
&=&
X_\Omega \beta(\sigma_t^\nu(n))X_\Omega^*\\
&=&
\beta(\sigma_t^\nu(n))\\
&=&
\tau_t(\beta(n))\\
&=&
\Delta_{\widehat{\Phi}}^{it}\beta(n)\Delta_{\widehat{\Phi}}^{-it}
\end{eqnarray*}
from which we get that $v_t^\Omega\in\beta(N)'$. 
\newline
Using \ref{invariant}(vi), and \cite{E5}8.8(ii), we get that :
 \[W^c(u^*\Delta_{\psi_1}^{it}u\underset{N^o}{_\alpha\otimes_{\hat{\beta}}}\hat{\delta}^{-it}P^{-it})=
(u^*\Delta_{\psi_1}^{it}u\underset{N}{_{\hat{\beta}}\otimes_\alpha}\hat{\delta}^{-it}P^{-it})W^c\]
and, as we have also, applying \cite{E5}8.8(ii) to the weight $\widehat{\Phi}$ :
\[W^c(\Delta_{\widehat{\Phi}}^{it}\underset{N^o}{_\alpha\otimes_{\hat{\beta}}}\hat{\delta}^{-it}P^{-it})=
(\Delta_{\widehat{\Phi}}^{it}\underset{N^o}{_\alpha\otimes_{\hat{\beta}}}\hat{\delta}^{-it}P^{-it})W^c\]
we get that :
\[W^c(v_t^\Omega\underset{N^o}{_\alpha\otimes_{\hat{\beta}}}1)=(v_t^\Omega\underset{N^o}{_\alpha\otimes_{\hat{\beta}}}1)W^c\]
from which one gets that $v_t^\Omega$ belongs to $M$, using \cite{E5}3.10(ii) applied to $\widehat{\gG^c}$. 
\newline
We have, for $s$, $t$ in $\mathbb{R}$ :
\[v_{s+t}^\Omega=u^*\Delta_{\psi_1}^{i(s+t)}u\Delta_{\widehat{\Phi}}^{-i(s+t)}=u\Delta_{\psi_1}^{is}u\Delta_{\widehat{\Phi}}^{-is}\Delta_{\widehat{\Phi}}^{is}u\Delta_{\psi_1}^{it}\Delta_{\widehat{\Phi}}^{-it}\Delta_{\widehat{\Phi}}^{-is}=v_s^\Omega\tau_s(v_t^\Omega)\]
which finishes the proof of (i). 
\newline
Using now \ref{corG}(iii) and \cite{E5}3.10(vii), we get that $u^*P_{A_\Omega}^{it}u=v_t^\Omega P^{it}$; therefore, \ref{corG}(iv) can be written :
\[\widetilde{W}(v_t^\Omega\underset{N}{_\beta\otimes_\alpha}v_t^\Omega)(P^{it}\underset{N}{_\beta\otimes_\alpha}P^{it})=(P^{it}\underset{N^o}{_\alpha\otimes_{\hat{\beta}}}v_t^\Omega P^{it})\widetilde{W}\]
or, using \cite{E5}3.8(vii) :
\[\widetilde{W}(v_t^\Omega\underset{N}{_\beta\otimes_\alpha}v_t^\Omega)=(1\underset{N^o}{_\alpha\otimes_{\hat{\beta}}}v_t^\Omega)W(\tau_t\underset{N}{_\beta*_\alpha}\tau_t)(\Omega^*)\]
from which we get the first formula of (ii). The second formula of (ii) is just a straightforward corollary of the first formula. We then get that :
\[(\omega_{v_t^\Omega\xi, \eta}*id)(\widetilde{W})=v_t^\Omega(\omega_{\xi, \eta}*id)(\widetilde{W_t})(v_t^\Omega)¬*\]
from which we get (iii). Using now the definitions of $\ga$ and $\ga_t$, we get that $(\mathcal I_t\underset{N}{_{\hat{\beta}}*_\alpha}id)\ga_t=\ga\circ\mathcal I_t$, $T_{\ga_t}=T_\ga\circ \mathcal I_t$ and $\psi_{1,t}=\psi_1\circ\mathcal I_t$. 
\newline
If we suppose now that $\omega_{\xi, \eta}$ belongs to $I_\Phi$, we get :
\begin{eqnarray*}
\Lambda_{\psi_1}(\mathcal I_t[\omega_{\xi, \eta}*id)(\widetilde{W_t})])
&=&
\Lambda_{\psi_1}[(\omega_{v_t^\Omega\xi, \eta}*id)(\widetilde{W})]\\
&=&
u\Lambda_{\widehat{\Phi}}(\omega_{v_t^\Omega\xi, \eta}*id)(W)]\\
&=&
u\pi'(\eta)^*v_t^\Omega\xi\\
&=&
uv_t^\Omega\pi'(\eta)^*\xi\\
&=&
uv_t^\Omega\Lambda_{\widehat{\Phi}}[(\omega_{\xi, \eta}*id)(W)]\\
&=&
uv_t^\Omega u_t^*\Lambda_{\psi_{1,t}}[(\omega_{\xi, \eta}*id)(\widetilde{W_t})]
\end{eqnarray*}
which finishes the proof of (iv). 
\newline
Using (iv), we get $uv_t^\Omega u_t^*\Delta_{\psi_{1,t}}^{is}=\Delta_{\psi_1}^{is}uv_t^\Omega u_t^*$, from which we infer :
\begin{eqnarray*}
v_t^\Omega v_s^{\Omega_t}
&=&
v_t^\Omega u_t^*\Delta_{\psi_{1,t}}^{is}u_t\Delta_{\widehat{\Phi}}^{-is}\\
&=&
u^*\Delta_{\psi_1}^{is}u\Delta_{\widehat{\Phi}}^{-is}\Delta_{\widehat{\Phi}}^{is}v_t^\Omega\Delta_{\widehat{\Phi}}^{-is}\\
&=&
v_s^\Omega\tau_s(v_t^\Omega)\\
&=&
v_{s+t}^\Omega\\
&=&
v_t^\Omega\tau_t(v_s^\Omega)
\end{eqnarray*}
from which we get (v). 
\newline
Using (v), we get that, if $\Omega=\Omega_t$, $v_t^\Omega$ is invariant under $\tau_s$, and is a one parameter group of unitaries, which, with (vi), finishes the proof. 
\end{proof}

\subsection{Theorem}
\label{thGaloiscocycle}
{\it Let $\gG$ be a measured quantum groupoid, $\Omega$ a $2$-cocycle for $\gG$; let $W$ be the pseudo-multiplicative unitary associated to $\gG$, $A_\Omega$ the von Neumann algebra on $H$ defined in \ref{algcocycle} and $(\hat{\beta}, \ga)$ the action of $\widehat{\gG}$ on $A_\Omega$ defined in \ref{action} whose invariant subalgebra $A_\Omega^\ga$ is equal to $\alpha(N)$ (\ref{invariant}(ii)); let us write $\widetilde{W}=W\Omega^*$; moreover, the weight $\nu\circ\alpha^{-1}$ on $\alpha(N)$ bears the Galois density property defined in \ref{defdens}, by \ref{invariant}(vi). Let us write $\psi_1=\nu\circ\alpha^{-1}\circ T_\ga$. Let $u$ be the unitary from $H$ onto $H_{\psi_1}$ introduced in \ref{invariant}(iii). 
\newline
Then, are equivalent :
\newline
(i) there exists a normal semi-finite faithful weight $\phi$ on $A_\Omega$ such that $(A_\Omega, \hat{\beta}, \ga, \phi, \nu)$ is a Galois system. 
\newline
(ii) there exists a one-parameter group of unitaries $\delta_\Omega^{it}$ on $H$, such that it is possible to define a one parameter group of unitaries $uJ_{\psi_1}u^* \delta_\Omega^{it}uJ_{\psi_1}u^*\underset{N}{_\beta\otimes_\alpha}\delta_\Omega^{it}$, with natural values on elementary tensors, and such that :
\[uJ_{\psi_1}u^* \delta_\Omega^{it}uJ_{\psi_1}u^*\underset{N}{_\beta\otimes_\alpha}\delta_\Omega^{it}=\Omega(\hat{J}\hat{\delta}^{it}\hat{J}\underset{N}{_\beta\otimes_\alpha}\hat{\delta}^{it})\Omega^*\]
(iii) there exists a $\tau_{-s}\sigma_{-s}^{\Phi\circ R}$-cocycle $t\mapsto u_t^\Omega$ in $M\cap \beta(N)'$, such that :
\[\Gamma(u_t^\Omega)=\Omega^*(u_t^\Omega\underset{N}{_\beta\otimes_\alpha}1)(\tau_{-t}\sigma_{-t}^{\Phi\circ R}\underset{N}{_\beta*_\alpha}id)(\Omega)\]
and $u_t^\Omega$ is linked with the $\tau_s$-cocycle $v_t^\Omega$ introduced in \ref{v} by the formula, for all $s$, $t$ in $\mathbb{R}$ :
\[u_t^\Omega\tau_{-t}\sigma_{-t}^{\Phi\circ R}(v_s^\Omega)=v_s^\Omega\tau_s(u_t^\Omega)\]
In that situation, $u^*\delta_\Omega u$ is the modulus of the action $(\hat{\beta}, \ga)$, and we have $\delta_\Omega^{it}=u_t^\Omega\hat{\delta}^{it}$. Moreover, there exists then a normal semi-finite faithful operator weight $T_\Omega$ from $M$ to $\alpha(N)$, (resp. $T'_\Omega$ from $M$ to $\beta(N)$) such that \[\gG_\Omega=(N, M, \alpha, \beta, \Gamma_\Omega, T_\Omega, T'_\Omega, \nu)\] is a measured quantum groupoid.
\newline
Moreover, if $N$ is a finite sum of factors, then any 2-cocycle satisfies these equivalent conditions.  }

\begin{proof}
The equivalence between (i) and (ii) is just an application of \ref{corgalois}, thanks to \ref{corG}(iii). We then get that $u^*\delta_\Omega u$ is the modulus of the action  $(\hat{\beta}, \ga)$ of $\widehat{\gG}$ on $A_\Omega$, and, using \ref{invariant}(vii), we get that $\delta_\Omega$ is affiliated to $A_\Omega$, and that, there exists a one-parameter group of unitaries $\delta_\Omega^{it}\underset{N}{_{\hat{\beta}}\otimes_\alpha}\delta^{it}$ such that, for all $t\in\mathbb{R}$, we have $\ga(\delta_\Omega^{it})=\delta_\Omega^{it}\underset{N}{_{\hat{\beta}}\otimes_\alpha}\hat{\delta}^{it}$. Let us write $u_t^\Omega=\delta_\Omega^{it}\hat{\delta}^{-it}$. Using \ref{algcocycle}(iv), we get that $u_t\in\beta(N)'$; moreover, using \ref{action}, and \cite{E5}3.12(v) and 3.8(vi) applied to $\widehat{\gG}$, we gat that $W^c(u_t^\Omega\underset{N^o}{_{\hat{\alpha}}\otimes_{\hat{\beta}}}1)(W^c)^*=u_t^\Omega\underset{N}{_{\hat{\beta}}\otimes_\alpha}1$, which gives that $u_t^\Omega$ belongs to $M$, thanks to \cite{E5}3.10(ii) applied to $\widehat{\gG}^c$. As, for any $x\in M$, and $t\in\mathbb{R}$, we have, using \cite{E5}3.11(ii), $\hat{\delta}^{it}x\hat{\delta}^{-it}=\tau_{-t}\sigma_{-t}^{\Phi\circ R}(x)$, we get that $t\mapsto u_t^\Omega$ is indeed a $\tau_{-s}\sigma_{-s}^{\Phi\circ R}$-cocycle. 
\newline
Using now \ref{corG}(i), we get that $\hat{\delta}^{it}\underset{N^o}{_\alpha\otimes_{\hat{\beta}}}\delta_\Omega^{it}=\widetilde{W}(\delta_\Omega^{it}\underset{N}{_\beta\otimes_\alpha}1)\widetilde{W}^*$. And, therefore :
\begin{eqnarray*}
1\underset{N^o}{_\alpha\otimes_\beta}u_t^\Omega
&=&
(\hat{\delta}^{it}\underset{N^o}{_\alpha\otimes_{\hat{\beta}}}\delta_\Omega^{it})(\hat{\delta}^{-it}\underset{N^o}{_\alpha\otimes_{\hat{\beta}}}\hat{\delta}^{-it})\\
&=&
\widetilde{W}(\delta_\Omega^{it}\underset{N}{_\beta\otimes_\alpha}1)\widetilde{W}^*(\hat{\delta}^{-it}\underset{N^o}{_\alpha\otimes_{\hat{\beta}}}\hat{\delta}^{-it})\\
&=&
\widetilde{W}(\delta_\Omega^{it}\underset{N}{_\beta\otimes_\alpha}1)\Omega W^*(\hat{\delta}^{-it}\underset{N^o}{_\alpha\otimes_{\hat{\beta}}}\hat{\delta}^{-it})\\
&=&
\widetilde{W}(\delta_\Omega^{it}\underset{N}{_\beta\otimes_\alpha}1)\Omega(\hat{\delta}^{-it}\underset{N}{_\beta\otimes_\alpha}1)W^*\\
&=&
\widetilde{W}(u_t^\Omega\underset{N}{_\beta\otimes_\alpha}1)(\tau_{-t}\sigma_{-t}^{\Phi\circ R}\underset{N}{_\beta*_\alpha}id)(\Omega)W^*
\end{eqnarray*}
and, therefore 
\[\Gamma(u_t^\Omega)=W^*(1\underset{N^o}{_\alpha\otimes_\beta}u_t^\Omega)W=\Omega^*(u_t^\Omega\underset{N}{_\beta\otimes_\alpha}1)(\tau_{-t}\sigma_{-t}^{\Phi\circ R}\underset{N}{_\beta*_\alpha}id)(\Omega)\]
which gives the first formula of (iii). 
\newline
Moreover, using \ref{thgalois}(iii), we get that $\sigma_t^{\psi_1}(\delta_\Omega^{it})=\lambda^{ist}\delta_\Omega^{it}$. Using \ref{v}, and \cite{E5}3.8(vi) applied to $\widehat{\gG}$, we have :
\begin{eqnarray*}
\sigma_t^{\psi_1}(\delta_\Omega^{it})
&=&
v_s^\omega\widehat{\Delta}^{is}u_t^\Omega\hat{\delta}^{it}\widehat{\Delta}^{-is}(v_s^\Omega)^*\\
&=&
v_s^\Omega\tau_s(u_t^\Omega)\sigma_s^{\widehat{\Phi}}(\hat{\delta}^{it})(v_s^\Omega)^*\\
&=&
v_s^\Omega\tau_s(u_t^\Omega)\lambda^{ist}\hat{\delta}^{it}(v_s^\Omega)^*\\
&=&
\lambda^{ist}v_s^\Omega\tau_s(u_t^\Omega)\hat{\delta}^{it}(v_s^\Omega)^*
\end{eqnarray*}
from which we get $u_t^\Omega\hat{\delta}^{it}=v_s^\Omega\tau_s(u_t^\Omega)\hat{\delta}^{it}(v_s^\Omega)^*$, which gives the the second formula of (iii). 
\newline
Conversely, if we have (iii), we can define a one parameter group of unitaries $\delta_\Omega^{it}$ by writing $\delta_\Omega^{it}=u_t^\Omega\hat{\delta}^{it}$. Now, from the first formula of (iii), taking the same calculation upside down, we get that $\hat{\delta}^{it}\underset{N^o}{_\alpha\otimes_{\hat{\beta}}}\delta_\Omega^{it}=\widetilde{W}(\delta_\Omega^{it}\underset{N}{_\beta\otimes_\alpha}1)\widetilde{W}^*$, which gives, by \ref{algcocycle}(iii), that $\delta_\Omega$ is affiliated to $A_\Omega$; so we had obtained that $\ga(\delta_\Omega^{it})=\delta_\Omega^{it}\underset{N}{_{\hat{\beta}}\otimes_\alpha}\hat{\delta}^{it}$. 
\newline
From the second formula of (iii), using again the same calculation upside down, we get that $\sigma_t^{\psi_1}(\delta_\Omega^{it})=\lambda^{ist}\delta_\Omega^{it}$, which proves that $\lambda$ is affiliated to $A_\Omega$; by the definition of $\delta_\Omega$, we see that the operators $\delta_\Omega$ and $\lambda$ strongly commute. Therefore, by \cite{V1}, 5.1, there exists a normal semi-finite faithful weight $\phi$ on $A_\Omega$ such that $(D\phi:D\psi_1)_t=\lambda^{it^2/2}\delta_\Omega^{it}$. 
\newline
Using now \ref{invariant}(iv) and \cite{E5}8.1, we get, for all $x\in\gN_\phi$ such that $x\delta_\Omega^{1/2}$ is bounded (its closure, denoted $\overline{x\delta_\Omega^{1/2}}$ belongs then to $\gN_{\psi_1}$, and we identify $\Lambda_\phi(x)$ with $\Lambda_{\psi_1}(\overline{x\delta_\Omega^{1/2}})$), for all $\eta$ in $D(H_{\hat{\beta}}, \nu^o)\cap \mathcal D(\delta_\Omega^{-1/2})$, such that $\delta_\Omega^{-1/2}\eta$ belongs to $D(_\alpha H, \nu)$ :
\begin{eqnarray*}
\|\Lambda_\phi(x)\underset{N^o}{_{\hat{\alpha}}\otimes_{\hat{\beta}}}\eta\|^2
&=&
\|\Lambda_{\psi_1}(\overline{x\delta_\Omega^{1/2}})\underset{N^o}{_{\hat{\alpha}}\otimes_{\hat{\beta}}}\eta\|^2\\
&=&
(\psi_1\underset{N}{_{\hat{\beta}}*_\alpha}\omega_{\delta_\Omega^{-1/2}\eta})\ga((\overline{x\delta_\Omega^{1/2}})^*\overline{x\delta_\Omega^{1/2}})\\
&=&
(\phi\underset{N}{_{\hat{\beta}}*_\alpha}\omega_\eta)\ga(x^*x)
\end{eqnarray*}
which, by continuity, remains true for all $\eta\in D(_\alpha H, \nu)\cap D(H_{\hat{\beta}}, \nu^o)$ and all $x\in\gN_\phi$, which proves that $\phi$ is invariant by $\ga$. But now, we are in the situation of \ref{thgalois}, which gives that $\lambda$ is affiliated to the center of $A_\Omega$; we then have (i). 

\end{proof}

\subsection{Corollaries}
\label{corGaloiscocycle}
{\it Let $\gG$ be a measured quantum groupoid, $\Omega$ a $2$-cocycle for $\gG$; let $W$ be the pseudo-multiplicative unitary associated to $\gG$; then, are equivalent :
\newline
(i) $\Omega$ satisfies the equivalent conditions of \ref{thGaloiscocycle}.
\newline
(ii) for all $t\in\mathbb{R}$, $\Omega_t$ (resp. $\Omega'_t$) satisfies the equivalent conditions of \ref{thGaloiscocycle}.
\newline
(iii) there exists $t\in\mathbb{R}$ such that $\Omega_t$ (resp. $\Omega'_t$) satisfies the equivalent conditions of \ref{thGaloiscocycle}.}
\begin{proof}
We can easily check that we can write $\tau_s(u_t^\Omega)=u_t^{\Omega_s}$, and $\delta^{is}u_t^\Omega\delta^{-is}=u_t^{\Omega'_s}$, then \ref{thGaloiscocycle} gives the result. 
\end{proof}

\subsection{Theorem}
\label{thcas}
{\it Let $\gG$ be a measured quantum groupoid, and $\Omega$ a 2-cocycle for $\gG$; let us suppose that, for any $t\in\mathbb{R}$, we have $(\tau_t\sigma_{-t}^\Phi\underset{N}{_\beta*_\alpha}\tau_t\sigma_t^{\Phi\circ R})(\Omega)=\Omega$. Then, the cocycle $\Omega$ satisfies the equivalent conditions of \ref{thGaloiscocycle}. In particular, there exists a normal semi-finite faithful operator weight $T_\Omega$ from $M$ to $\alpha(N)$, (resp. $T'_\Omega$ from $M$ to $\beta(N)$) such that \[\gG_\Omega=(N, M, \alpha, \beta, \Gamma_\Omega, T_\Omega, T'_\Omega, \nu)\] is a measured quantum groupoid. Moreover, we get, for all $t\in\mathbb{R}$, that $\tau_{-t}\sigma_{-t}^{\Phi\circ R}(v_s^\Omega)=v_s^\Omega$ and $(\tau_{-t}\sigma_{-t}^{\Phi\circ R} \underset{N}{_\beta*_\alpha}id)(\Omega)=\Omega$.}

\begin{proof}
Using \ref{corG}(iii), we get that 
\[\Omega(\hat{J}\hat{\delta}^{it}\hat{J}\underset{N}{_\beta\otimes_\alpha}\hat{\delta}^{it})\Omega^*=\hat{J}\hat{\delta}^{it}\hat{J}\underset{N}{_\beta\otimes_\alpha}\hat{\delta}^{it}\]
from which, using \ref{invariant}(vii), we get that $\hat{\delta}^{it}$ belongs to $A_\Omega$, and, by \ref{propK}(v), that $\ga(\hat{\delta}^{it})=\hat{\delta}^{it}\underset{N}{_{\hat{\beta}}\otimes_\alpha}\hat{\delta}^{it}$. Using now \cite{E5} 8.8(iii), one gets that, for any $s$, $t$ in $\mathbb{R}$, we have :
\[\ga(\sigma_s^{\psi_1}(\hat{\delta}^{it}))=\sigma_s^{\psi_1}(\hat{\delta}^{it})\underset{N}{_{\hat{\beta}}\otimes_\alpha}\hat{\delta}^{it}\]
from which one gets that $\sigma_s^{\psi_1}(\hat{\delta}^{it})\hat{\delta}^{-it}$ belongs to $A_\Omega^\ga=\alpha(N)$ by \ref{invariant}. 
\newline
More precisely, if $n\in N$, we get that :
\begin{multline*}
\sigma_s^{\psi_1}(\hat{\delta}^{it})\hat{\delta}^{-it}\alpha(n)=\sigma_s^{\psi_1}(\hat{\delta}^{it})\alpha(\sigma_t^\nu\gamma_t(n))\hat{\delta}^{-it}=\sigma_s^{\psi_1}(\hat{\delta}^{it}\alpha(\sigma_{t-s}^\nu\gamma_t(n)))\hat{\delta}^{-it}\\=\sigma_t^{\psi_1}(\alpha(\sigma_{-s}^\nu(n))\hat{\delta}^{it})\hat{\delta}^{-it}=\alpha(n)\sigma_s^{\psi_1}(\hat{\delta}^{it})\hat{\delta}^{-it}
\end{multline*}
and, therefore, we get that $\sigma_s^{\psi_1}(\hat{\delta}^{it})\hat{\delta}^{-it}$ belongs to $\alpha(Z(N))$. 
\newline
But, on the other hand, using \ref{v}, we get that $\sigma_s^{\psi_1}(\hat{\delta}^{it})\hat{\delta}^{-it}=v_s^\Omega\sigma^{\widehat{\Phi}}_s(\hat{\delta}^{it})(v_s^\Omega)^*\hat{\delta}^{-it}=v_s^\Omega\lambda^{ist}\hat{\delta}^{it}(v_s^\Omega)^*\hat{\delta}^{-it}=\lambda^{ist}v_s^\Omega \hat{\delta}^{it}(v_s^\Omega)^*\hat{\delta}^{-it}$, from which, using \cite{E5}8.11 (ii), we get :
\[\sigma_s^{\psi_1}(\hat{\delta}^{it})\hat{\delta}^{-it}=\lambda^{ist}v_s^\Omega \tau_{-t}\sigma_{-t}^{\Phi\circ R}(v_s^\Omega)^*\]
and, for all $s$, $t$ in $\mathbb{R}$, $\tau_{-t}\sigma_{-t}^{\Phi\circ R}(v_s^\Omega)(v_s^\Omega)^*$ belongs to $\alpha(Z(N))$. Therefore, there exists a one-parameter group of unitaries $t\mapsto \mu_s^{it}$ in $Z(N)$ such that $\tau_{-t}\sigma_{-t}^{\Phi\circ R}(v_s^\Omega)=\alpha(\mu_s^{it})v_s^\Omega$; and, therefore, $\sigma_s^{\psi_1}(\hat{\delta}^{it})=\lambda^{ist}\alpha(\mu_s^{-it})\hat{\delta}^{it}$. So, there exists a positive non singular operator $\mu$ affiliated to $Z(N)$ such that $\sigma_s^{\psi_1}(\hat{\delta}^{it})=\lambda^{ist}\alpha(\mu^{-ist})\hat{\delta}^{it}$ and $\tau_t\sigma_t^{\Phi\circ R}(v_s^\Omega)=\alpha(\mu^{ist})v_s^\Omega$. But now, as, for all $u\in\mathbb{R}$, we have $\tau_u\sigma_u^{\Phi\circ R}(\alpha(\mu^{ist}))=\alpha(\gamma_u(\mu^{ist}))$, we get that $\gamma_u(\mu)=\mu$, and, therefore, that $\hat{\delta}$ and $\lambda\alpha(\mu)$ strongly commute. Therefore, by \cite{V1}, 5.1, there exists a normal semi-finite faithful weight $\phi$ on $A_\Omega$ such that $(D\phi:D\psi_1)_t=(\lambda\alpha(\mu))^{it^2/2}\hat{\delta}^{it}$. 
\newline
Using now \ref{invariant}(iv) and \cite{E5}8.1, as in \ref{thGaloiscocycle} that $\phi$ is invariant by $\ga$. But now, we are in the situation of \ref{thgalois}, which gives that $\mu=1$, and proves that we are in the situation of \ref{thGaloiscocycle}, with, moreover, $u_t^\Omega=1$; we then infer from \ref{thGaloiscocycle} that $\tau_{-t}\sigma_{-t}^{\Phi\circ R}(v_s^\Omega)=v_s(\Omega)$ and that $(\tau_{-t}\sigma_{-t}^{\Phi\circ R} \underset{N}{_\beta*_\alpha}id)(\Omega)=\Omega$. \end{proof}

\subsection{Theorem}
\label{thpart}
{\it Let $\gG$ be a measured quantum groupoid, and $\Omega$a 2-cocycle for $\gG$; let us suppose that, for any $t\in\mathbb{R}$, we have $(\tau_{-t}\sigma_{-t}^{\Phi\circ R} \underset{N}{_\beta*_\alpha}id)(\Omega)=\Omega$. Then, the cocycle $\Omega$ satisfies the equivalent conditions of \ref{thGaloiscocycle}. In particular, there exists a normal semi-finite faithful operator weight $T_\Omega$ from $M$ to $\alpha(N)$, (resp. $T'_\Omega$ from $M$ to $\beta(N)$) such that \[\gG_\Omega=(N, M, \alpha, \beta, \Gamma_\Omega, T_\Omega, T'_\Omega, \nu)\] is a measured quantum groupoid. Moreover, we get, for all $t\in\mathbb{R}$, that $\tau_{-t}\sigma_{-t}^{\Phi\circ R}(v_s^\Omega)=v_s^\Omega$ and $(\tau_t\sigma_{-t}^\Phi\underset{N}{_\beta*_\alpha}\tau_t\sigma_t^{\Phi\circ R})(\Omega)=\Omega$. }

\begin{proof}
The proof is similar to \ref{thcas}.
\end{proof}

\section{Examples, at least}
\label{atleast}
In this last chapter, we construct  a general situation in which the deformations of a measured quantum groupoid by some 2-cocycles are still measured quantum groupoids. 

\subsection{Measured quantum groupoids associated to a matched pair of groupoids}
\label{G1G2}
In \cite{Val6} was decribed a procedure for constructing measured quantum groupoids : 
\newline
Let $\mathcal G$ be a locally compact groupoid, with $\mathcal G^{(0)}$ as set of units, and $r: \mathcal G\mapsto \mathcal G^{(0)}$ (resp. $s:\mathcal G\mapsto \mathcal G^{(0)}$) as range (resp. source) application, equipped with a Haar system $(\lambda^u)_{u\in\mathcal G^{(0)}}$ and a quasi-invariant measure $\nu$ on $\mathcal G^{(0)}$. Let us write $\mu=\int_{\mathcal G^{(0)}}\lambda^u d\nu(u)$. 
\newline
Let $\mathcal G_1$, $\mathcal G_2$ two closed subgroupoids of $\mathcal G$, (with $r_1=r_{|\mathcal G_1}$, etc) equipped with their Haar systems $(\lambda_1^u)_{u\in\mathcal G^{(0)}}$, $(\lambda_2^u)_{u\in\mathcal G^{(0)}}$. 
\newline
Then $(\mathcal G_1, \mathcal G_2)$ is called a matched pair of groupoids if :
\newline
(i) $\mathcal G_1\cap\mathcal G_2=\mathcal G^{(0)}$
\newline
(ii) the set $\mathcal G_1\mathcal G_2=\{g_1g_2, g_1\in\mathcal G_1, g_2\in\mathcal G_2^{s(g_1)}\}$ is $\mu$-conegligeable in $\mathcal G$.
\newline
(iii) there exists a measure $\nu$ on $\mathcal G^{(0)}$ with is quasi-invariant for the three Haar systems. 
\newline
Then, Vallin has constructed an action $(s_2, \ga)$ of $\gG(\mathcal G_1)$ on $L^\infty (\mathcal G_2, \mu_2)$, and put a measured quantum groupoid structure on the crossed product $L^\infty(\mathcal G_2, \mu_2)\rtimes_\ga\gG(\mathcal G_1)$. 
\newline
Let us denote $\gG(\mathcal G_1, \mathcal G_2)=(L^\infty(\mathcal G^{(0)},\nu), L^\infty(\mathcal G_2, \mu_2)\rtimes_\pi\gG(\mathcal G_1), m, s, \Gamma, T_L, T_R, \nu)$  this measured quantum groupoid. 
\newline
Moreover, there exists a right action $(r_1, \hat{\ga})$ of $\gG(\mathcal G_2)$ on $L^\infty (\mathcal G_1, \mu_1)$, which leads to a measured quantum groupoid structure on $L^\infty(\mathcal G_1, \mu_1)\ltimes_{\hat{\ga}}\gG(\mathcal G_2)$, we shall write $\gG(\mathcal G_2, \mathcal G_1)$; we have $\gG(\mathcal G_2, \mathcal G_1)=\widehat{\gG(\mathcal G_1, \mathcal G_2)}$
\newline
This measured quantum groupoid $\gG(\mathcal G_1, \mathcal G_2)$ bears some properties :
\newline
(i) the scaling operator $\lambda$ is equal to $1$. \newline
(ii)  for any $f\in L^\infty(\mathcal G_2, \mu_2)$, $\ga(f)$ is invariant under $\sigma_t^\Phi$ (\cite{Val6}, 4.3.5). 
\newline
(iii) for any $f\in L^\infty(\mathcal G_2, \mu_2)$, we have $R(\ga(f))=\ga(\check{f})$, where $R$ is the co-inverse of $\gG(\mathcal G_1, \mathcal G_2)$, and $\check{f}(g_2)=f(g_2^{-1})$, for any $g_2\in\mathcal G_2$. Therefore, using (i), we get that $\ga(f)$ is also invariant under $\sigma_t^{\Phi\circ R}$.
\newline
(iv) using \cite{Val6}4.1.1 and \cite{E5}3.8(ii), one can easily check that, for all $t\in\mathbb{R}$, and $f\in L^\infty(\mathcal G_2, \mu_2)$, we have $\tau_t(\ga(f))=\ga(f)$. Namely we have, using (ii) :
\begin{eqnarray*}
(\ga\underset{L^\infty(\mathcal G^{(0)}, \nu)}{_{s_2}*_{r_2}}\ga)\Gamma_{\mathcal G_2}(f)
&=&
\Gamma(\ga(f))\\
&=&
\Gamma(\sigma_t^\Phi(\ga(f)))\\
&=&
(\tau_t\underset{L^\infty(\mathcal G^{(0)}, \nu)}{_s*_m}\sigma_t^\Phi)\Gamma(\ga(f))\\
&=&
(\tau_t\circ\ga\underset{L^\infty(\mathcal G^{(0)}, \nu)}{_{s_2}*_{r_2}}\sigma_t^\Phi\circ\ga)\Gamma_{\mathcal G_2}(f)\\
&=&
(\tau_t\circ\ga\underset{L^\infty(\mathcal G^{(0)}, \nu)}{_{s_2}*_{r_2}}\ga)\Gamma_{\mathcal G_2}(f)
\end{eqnarray*}
from which we get the result. 
\newline
We refer to \cite{Val6} for all details.

\subsection{Theorem}
\label{thG1G2}
{\it Let $\gG(\mathcal G_1, \mathcal G_2)$ be the measured quantum groupoid constructed from a matched pair $(\mathcal G_1, \mathcal G_2)$ of groupoids. Let us use all notations of \ref{G1G2}. Let $\Omega$ be a 2-cocycle for $\gG(\mathcal G_2)$, as defined in \ref{defcocycle}. Then :
\newline
(i) $(\ga\underset{L^\infty(\mathcal G^{(0)}, \nu)}{_{s_2}*_{r_2}}\ga)(\Omega)$ is a 2-cocycle for $\gG(\mathcal G_1, \mathcal G_2)$, we shall write $\Omega_\ga$ for simplification. 
\newline
(ii) There exists a left-invariant operator-valued weight $T_\Omega$ and a right-invariant operator-valued weight $T'_\Omega$ such that :
\[\gG(\mathcal G_1, \mathcal G_2)_{\Omega_\ga}=(L^\infty(\mathcal G^{(0)},\nu), L^\infty(\mathcal G_2, \mu_2)\rtimes_\ga\gG(\mathcal G_1), m, s, \Gamma_{\Omega_\ga}, T_\Omega, T'_\Omega, \nu)\]
is a measured quantum groupoid. }

\begin{proof}
Using \cite{Val6} 4.1.1, one gets (i). As, for all $t\in\mathbb{R}$, $\tau_t\sigma_{-t}^\Phi\circ\ga=\ga$, and $ \tau_t\sigma_t^{\Phi\circ R}\circ\ga=\ga$, we get that this cocycle $\Omega_\ga$ satisfies the conditions of \ref{thcas} or \ref{thpart}. So, we get (ii). 
\end{proof}

\subsection{Matched pair of groups acting on a space}
\label{G1G2X}
As a particular case of \ref{G1G2}, we can study, following (\cite{Val6}5.1) the case where $G$ is a locally compact group acting (on the right) on a locally compact space $X$, and $G_1$, $G_2$ a matched pair of closed subgroups of $G$, in the sense of \cite{BSV}. Then, we can define almost everywhere Borel functions $p^G_1$ from $G$ to $G_1$ and $p^G_2$ from $G$ to $G_2$, such that :
\[g=p^G_1(g)p^G_2(g)\]
Following \cite{VV}, we can construct an action $a_1$ of $G_1$ on $L^\infty(G_2)$, and put on the crossed product $L^\infty(G_2)\rtimes_{a_1}G_1$ a structure of a locally compact quantum group we shall denote by $\bf{G}$$(G_1, G_2)$; Let us denote $\widetilde{\Gamma}$ the coproduct of this locally compact quantum group. 
\newline
Let us denote now $\mathcal{G}$ (resp. $\mathcal{G}_1$, resp. $\mathcal{G}_2$) the locally compact groupoid given by the action of $G$ (resp. $G_1$, resp. $G_2$) on $X$. Then, it is easy to get that $\mathcal{G}_1$ and $\mathcal{G}_2$ are two closed subgroupoids of $\mathcal{G}$, which are a matched pair of groupoids as defined in \ref{G1G2}. So, there is an action $\ga$ of the measured quantum groupoid $\gG(\mathcal G_1)$ on $L^\infty(\mathcal G_2)$, and a measured quantum groupoid structure $\gG(\mathcal G_1, \mathcal G_2)$ on the crossed product $L^\infty(\mathcal G_2)\rtimes_\ga\gG(\mathcal G_1)$. The action $\ga$ can be identified with an action $\tilde{\ga}$ of $G_1$ on $L^\infty (X\times G_2)$ (\cite{Val6}, 5.1.2) and the crossed product $L^\infty(\mathcal G_2)\rtimes_\ga\gG(\mathcal G_1)$ can be identified with the crossed-product $L^\infty (X\times G_2)\rtimes_{\tilde{\ga}}G_1$, which will be considered as bounded operators on $L^2(X\times G\times G)$ (\cite{Val6}, 5.1.1). 
\newline
We can identify $L^2(X\times G_2)\underset{L^\infty(X)}{_{s_2}\otimes_{r_1}}L^\infty(X\times G_1)$ with $L^2(X\times G_2)\otimes L^2(G_1)$ (\cite{Val6}, 5.1.1); using these identifications,  are given in (\cite{Val6} 5.1.2) the formulae of the coproduct $\Gamma$ we can put on this crossed-product. For any $f\in L^\infty(X\times G_2)$, $h\in L^\infty(X)$, $k\in L^\infty(G_1)$, we have :
\[\Gamma(\ga(f))(x, g, g')=f(x.p^G_1(g), p^G_2(g)p^G_2(g'))\]
\[\Gamma(1\underset{L^\infty(X)}{_{s_2}\otimes_{r_1}}\rho(h\otimes k))=M(h)(1\otimes\widetilde{\Gamma}(1\otimes \rho_1(k)))\]
where $M(h)$ is the function $M(h)(x,g,g')=h(x. gp^G_2(g'))$. 
\newline
Let's see now how this coproduct can be deformed by a 2-cocycle $\Omega_2$ for $\gG(\mathcal G_2)$ to a new coproduct $\Gamma_\Omega$. For simplification, we shall restrict to a 2-cocycle $\Omega_2$ for $G_2$, which can be easily considered as a 2-cocycle for $\gG(\mathcal G_2)$. Using then \cite{VV}, we can put on the crossed-product $L^\infty(G_2)\rtimes_{a_1}G_1$ another structure of locally compact quantum group we shall denote by $\bf{G}$$(G_1, G_2)_{\Omega_2}$, with a deformed coproduct we shall denote $\widetilde{\Gamma}_{\Omega_2}$. 
\newline
By construction, we have :
\[\Gamma_{\Omega_2}(\ga(f))(x, g, g')=\Gamma(\ga(f))(x, g, g')=f(x.p^G_1(g), p^G_2(g)p^G_2(g'))\]
and :
\[\Gamma_{\Omega_2}(1\underset{L^\infty(X)}{_{s_2}\otimes_{r_1}}\rho(h\otimes k))=M(h)(1\otimes\widetilde{\Gamma}_{\Omega_2}(1\otimes \rho_1(k)))\]

\subsection{Looking back to Kac-Paljutkin's examples}
\label{KP}
Following (\cite{VV}, 5.1.1), let's look at the particular case of \ref{G1G2X} where $G_2$ is a normal subgroup of $G$; then $G_1$ acts on $G_2$ by (inner in $G$) automorphisms, the action of $G_2$ on $G_1$ is trivial, the application $p^G_1$ is an homomorphism and $G$ is the semi-direct product $G_2\rtimes G_1$. Then we well know that the old Kac-Paljutkin's examples can be obtained as locally compact quantum groups of the form $\bf{G}$$(G_1, G_2)_{\Omega_2}$. 
\newline
(i) taking $G_1=\mathbb{Z}/2\mathbb{Z}$ acting on $G_2=(\mathbb{Z}/2\mathbb{Z})^2$ by permutations, the cocycle $\Omega$ had been computed in (\cite{BS}, 8.26.1), in order to get that $\bf{G}$$(G_1, G_2)_{\Omega_2}$ is then the dimension 8 example constructed in \cite{KP1}. Taking now an action of the semi-direct product $G=G_2\rtimes G_1$ on a locally compact space $X$, we obtain, by \ref{G1G2X} applied to this particular case, a measured quantum groupoid given by dimension 8 Kac-Paljutkin's example and a right action of $(\mathbb{Z}/2\mathbb{Z})^2\rtimes \mathbb{Z}/2\mathbb{Z}$ on a space $X$. 
\newline
(ii) taking $G_1=\mathbb{R}$ acting on $G_2=\mathbb{R}^2$ by $a_g(x)=exp(gK)(x)$ ($x\in\mathbb{R}^2$, $K$ is a real $2\times 2$ matrix). Then the cocycle had been computed in (\cite{VV},8.26.2) and leads to the infinite dimension Kac-Paljutkin's example (\cite{KP2}). So, starting from this example, and some right action of the Heisenberg group $H_3(\mathbb{R})=\mathbb{R}^2\rtimes_a\mathbb{R}$ on $X$, we get, by \ref{G1G2X}, another example of a measured quantum groupoid.


\end{document}